\documentstyle[11pt,leqno]{article}
\begin{document}
 \textheight 193mm
 \textwidth  130mm
\baselineskip 16pt
\parindent=22pt

\pagestyle{empty}

\vspace*{15mm}

\begin{center}

\vspace{5mm}

 {\bf\huge English Translation of Chapter 9 of the book}

\vspace{8mm}

{\bf\huge Adams spectral sequence and stable homotopy groups of spheres}

\vspace{8mm}

{\bf\huge (In Chinese)}

\vspace{8mm}

{\bf\huge (Sciences Press, Beijing 2007)}

\vspace{8mm}

{\bf\huge By Jinkun Lin}

\vspace{60mm}

{\bf\huge  September, 2007}

\end{center}

\newpage

\pagestyle{plain}

\begin{center}

{\bf\huge Chapter 9}

\vspace{3mm}

{\bf\huge A  sequence  of new families in the stable homotopy groups of spheres}

\vspace{3mm}

\end{center}

In this Chapter, we will state and prove the existence of a sequence of new families in the stable homotopy groups of
spheres which is in base on several papers of the author (especially [7][8][9][24]).As a preminilaries, in \S 1 we introduce
some spectra which is closely related to Moore spectrum and Smith-Toda spectrum $V(1)$ and state some of their properties.
In \S 2 we state and prove a general result on the convergence of $h_0\sigma$ and $h_0\sigma'$ for a pair of $a_0$-related
elements $\sigma$ and $\sigma'$. ( the generalization of [8] Theorem A). \S 3 is devoted to state and prove a general result on
the convergence of $(i'i)_*(h_0\sigma)$ induces the convergence of $(i'i)_*(g_0\sigma)$ in the stable homotopy groups of
Smith-Toda spectrum $V(1)$ (the generalization of [7] Theorem II).  In \S 4 we prove a pull back Theorem in the Adam
spectral sequence and as a corollary of the main results in \S 2 and \S 4 , in \S 5
we obtain the convergence of a sequence of $h_0h_n,
h_0b_n, h_0h_nh_m, h_0(h_nb_{m-1}-h_mb_{n-1})$ new families in the stable homotopy groups of spheres.  \S 6 concerns with
the convergence of a sequence of $h_0\sigma\widetilde{\gamma}_s, g_0\sigma\widetilde{\gamma}_s$-elements. In \S 7, we first
prove $h_n$ Theorem and then obtain the third periodicity $\gamma_{p^n/s}$-families ([24] Theorem I and Theorem II)).
At last , in \S 8, the second periodicity $\beta_{tp^n/j,i+1}$-families in the stable homotopy groups of spheres are detected.

\quad

\begin{center}

{\bf\large \S 1.\quad Some spectra closely related to the Moore spectrum and Smith-Toda spectrum $V(1)$}

\vspace {2mm}

\end{center}

Let $M$  be the Moore spectrum given by the cofibration\\
{\bf (9.1.1)} \qquad\quad $S\stackrel{p}{\longrightarrow} S\stackrel{i}{\longrightarrow}
M\stackrel{j}{\longrightarrow}\Sigma S$\\
Let $\alpha : \Sigma^qM\rightarrow M$ be Adams map and $K$ be its cofibre given by the cofibration\\
{\bf (9.1.2)}\qquad\quad $\Sigma^qM\stackrel{\alpha}{\longrightarrow} M\stackrel{i'}
{\longrightarrow}K\stackrel{j'}{\longrightarrow}\Sigma^{q+1}M$\\
The above spectrum $K$ which we briefly write as $K$ actually is the Smith-Toda spectrum $V(1)$ in Chapter 6 \S 2.

Now we introduce some spectra closely related to $S, M$ or $K$ .  Let $L$ be the cofibre of $\alpha_{1} =
j\alpha i :\Sigma^{q-1}S\rightarrow S$ given by the cofibration\\
{\bf (9.1.3)} $\qquad\quad \Sigma^{q-1}S\stackrel{\alpha_{1}}{\longrightarrow} S\stackrel{i''}
{\longrightarrow}L\stackrel{j''}{\longrightarrow} \Sigma^{q}S$.\\
Let $Y$ be the cofibre of $i'i : S\rightarrow K$ given by the cofibration\\
{\bf (9.1.4)} $\quad\qquad S\stackrel{i'i}{\longrightarrow} K\stackrel{\overline{r}}
{\longrightarrow}Y\stackrel{\epsilon}{\longrightarrow} \Sigma S$.\\
 $Y$ actually is the Toda spectrum $V(1\frac{1}{2})$, and it also  is the cofibre of $j\alpha :
\Sigma^{q}M \rightarrow \Sigma S$  given by the cofibration\\
{\bf (9.1.5)} $\quad\qquad \Sigma^{q}M\stackrel{j\alpha}{\longrightarrow} \Sigma S\stackrel{\overline{w}}
{\longrightarrow} Y\stackrel{\overline{u}}{\longrightarrow}\Sigma^{q+1}M$,\\
 This can be seen by the following homotopy commutative diagram of $3\times 3$-Lemma in the stable homotopy category
 (cf. Chapter 3 \S 7)

$\qquad\qquad S\quad\stackrel{i'i}{\longrightarrow}\quad K\quad\stackrel{j'}
{\longrightarrow}\quad\Sigma^{q+1}M$

$\qquad\qquad\quad\searrow i\quad\nearrow i'\quad\searrow\overline{r}\quad\nearrow
\overline{u}$\newline
{\bf (9.1.6)} $\qquad\qquad\qquad M\quad\qquad\qquad Y$

$\qquad\qquad\quad\nearrow\alpha\quad\searrow j\quad \nearrow \overline{w}\quad
\searrow \epsilon $

$\qquad\qquad \Sigma^{q}M\quad\stackrel{j\alpha}{\longrightarrow}\quad\Sigma S
\quad\stackrel{p}{\longrightarrow}\quad \Sigma S$\\

Note that $\alpha_1\cdot p = p\cdot \alpha_{1} = 0$, then there exist $\pi\in [\Sigma^qS, L]$ and $\xi\in [L,S]$ such that
 $ p = j''\pi$
and $p = \xi i''$. since $\pi_qS$ = 0, then $\pi_qL\cong Z_{(p)}\{\pi\}$. Moreover, $i''\xi i'' =
i''\cdot p = (p\wedge 1_L)i''$, then $p\wedge 1_L = i''\xi + \lambda \pi j''$ for some $\lambda\in Z_{(p)}$.  By composing
$j''$ on the above equation we have $p\cdot j'' = j''(p\wedge 1_L)
= \lambda j''\pi\cdot j'' = \lambda p\cdot j''$ so that $\lambda = 1$ and we have\\
(9.1.7)\qquad $p\wedge 1_L \quad = \quad i'' \xi + \pi j''$.\\
By the following homotopy commutative diagram of $3\times 3$-Lemma

$\qquad\qquad \Sigma^{q}S\quad\stackrel{p}{\longrightarrow}\quad \Sigma^{q}S\quad\stackrel
{\alpha_{1}}{\longrightarrow}\quad \Sigma S$

$\qquad\qquad\quad \searrow \pi\quad\nearrow j''\quad\searrow i
\quad\nearrow j\alpha\quad\searrow i''$\newline
{\bf (9.1.8)} $\qquad\qquad\qquad L\qquad\qquad\quad \Sigma^q M\qquad\qquad\Sigma^{q+1}L$

$\qquad\qquad\quad\nearrow i''\quad\searrow\overline{h}\quad\nearrow
\overline{u}\quad\searrow j\quad\nearrow \pi$

$\qquad\qquad S\quad\stackrel{\overline{w}}{\longrightarrow}\quad \Sigma^{-1}Y\quad
\stackrel{j\overline{u}}{\longrightarrow}\quad \Sigma^{q+1}S$\\
we obtain the following cofibration\\
{\bf (9.1.9)} $\quad \Sigma^{q}S\stackrel{\pi}{\longrightarrow}L\stackrel{\overline{h}}
{\longrightarrow}\Sigma^{-1}Y\stackrel{j\overline{u}}{\longrightarrow}
\Sigma^{q+1}S$\\
and there are equations\quad $\overline{u}\overline{h} = i\cdot j'',
\quad \bar h i'' = \overline{w},\quad \pi\cdot j = i'' j\alpha$
By $2\alpha ij\alpha = ij\alpha^{2} + \alpha^{2}ij$ (cf. (6.5.3)),
then we have  $\alpha_{1}\alpha_{1}$ = 0 and so there are $\phi\in [\Sigma^{2q-1}S, L]$
and $(\alpha_{1})_{L}\in [\Sigma^{q-1}L, S]$ such that\\
{\bf (9.1.10)} $\quad j''\phi\quad = \quad\alpha_{1}\quad = \quad (\alpha_{1})_{L} \cdot i''$.

Let $W$ be the cofibre of $\phi : \Sigma^{2q-1}S\rightarrow L$, then $W$ also is the cofibre of
 $(\alpha_{1})_{L} : \Sigma^{q-1}L\rightarrow S$
, This can be seen by the following homotopy commutative diagram of $3\times 3$-Lemma

$\qquad\quad \Sigma^{2q-1}S\qquad\stackrel{\alpha_{1}}{\longrightarrow}\quad \Sigma^{q}S
\quad\stackrel{\alpha_{1}}{\longrightarrow}\qquad\Sigma S$

$\qquad\qquad\qquad\searrow\phi\quad\nearrow j''\quad\searrow i''\quad\nearrow
(\alpha_{1})_{L}$\newline
{\bf (9.1.11)} $\qquad\qquad\qquad L\qquad\qquad\qquad\Sigma^{q}L$

$\qquad\qquad\qquad\nearrow i''\quad\searrow w\quad\nearrow u\quad\searrow j''$

$\qquad\qquad S\qquad\stackrel{wi''}{\longrightarrow} \qquad W\qquad\stackrel{j''u}
{\longrightarrow} \qquad\Sigma^{2q}S$.\\
that is, we have two cofibrations\newline
{\bf (9.1.12)} $\quad \Sigma^{2q-1}S\stackrel{\phi}{\longrightarrow} L\stackrel{w}{\longrightarrow}
W\stackrel{j''u}{\longrightarrow}\Sigma^{2q}S$\newline
{\bf (9.1.13)} $\quad \Sigma^{q-1}L\stackrel{(\alpha_{1})_{L}}{\longrightarrow} S
\stackrel{wi''}{\longrightarrow} W\stackrel{u}{\longrightarrow}\Sigma^{q}L$.

We write the Toda spectrum $V(\frac{1}{2})$ as $K'$, it is the cofibre of $jj' : \Sigma^{-1}K\rightarrow
\Sigma^{q+1}S$ given by the cofibration\\
{\bf (9.1.14)}$\quad \Sigma^{-1}K\stackrel{jj'}{\longrightarrow}\Sigma^{q+1}S\stackrel
{z}{\longrightarrow}K'\stackrel{x}{\longrightarrow}K$\\
 $K'$ also is the cobibre of $\alpha i : \Sigma^q S\rightarrow M$ given by the cofibration\\
{\bf (9.1.15)}$\quad \Sigma^qS\stackrel{\alpha i}{\longrightarrow}M\stackrel{v}{\longrightarrow}
K'\stackrel{y}{\longrightarrow}\Sigma^{q+1}S$\\
This can be seen by the following homotopy commutative diagram of $3\times 3$-Lemma

$\qquad\qquad\Sigma^qS\quad\stackrel{\alpha i}{\longrightarrow}\quad M\quad
\stackrel{i'}{\longrightarrow}\quad K$

$\qquad\qquad\quad\searrow i\quad\nearrow \alpha\quad\searrow v\quad\nearrow x$\\
{\bf (9.1.16)}\quad\qquad\qquad $\Sigma^qM\qquad\qquad K'$

$\qquad\qquad\quad\nearrow j'\quad\searrow j\quad\nearrow z\quad\searrow y$

$\qquad\quad\Sigma^{-1}K\stackrel{jj'}{\longrightarrow}\quad\Sigma^{q+1}S\quad
\stackrel{p}{\longrightarrow}\Sigma^{q+1}S$

By $\alpha_1\wedge 1_M = ij\alpha - \alpha ij$, let $\alpha' = \alpha_1\wedge 1_K
\in [\Sigma^{q-1}K,K]$, then $j'\alpha' = - (\alpha_1\wedge 1_M)j' = \alpha ijj'\in [\Sigma^{-2}K,M]$.
By (9.1.16) we have  $y\cdot z = p$, then $y\cdot z\cdot y = p\cdot y = y(1_{K'}\wedge p)$,
so that $z\cdot y = 1_{K'}\wedge p$. This is because $[K',M] = 0$ which can be seen by the following exact sequence induced
by (9.1.15)\\
\centerline{$0 = [\Sigma^{q+1}S,M]\stackrel{y^*}{\longrightarrow} [K',M]\stackrel
{v^*}{\longrightarrow}[M,M]\stackrel{(\alpha i)^*}{\longrightarrow}$}\\
where $[M,M]\cong Z_p\{1_M\}$ so that the above $(\alpha i)^*$ is monic. Moreover, by the following homotopy commutative
diagram of $3\times 3$-Lemma we know that the cofibre of $\alpha ijj' : \Sigma^{-1}K\rightarrow\Sigma M$
is $K'\wedge M$ given by the following cofibration\\
{\bf (9.1.17)}$\quad\Sigma^{-1}K\stackrel{\alpha ijj'}{\longrightarrow}\Sigma M\stackrel
{\psi}{\longrightarrow} K'\wedge M\stackrel{\rho}{\longrightarrow} K$

\newpage

$\qquad\qquad\Sigma^{-1}K\quad\stackrel{\alpha ijj'}{\longrightarrow}\Sigma M\quad
\stackrel{v}{\longrightarrow}\quad \Sigma K'$

$\qquad\qquad\quad\searrow jj'\quad\nearrow\alpha i\quad\searrow \psi\quad\nearrow _{1_{K'}\wedge j}$\\
{\bf (9.1.18)}$\qquad\qquad\quad \Sigma^{q+1}S\qquad\qquad K'\wedge M$

$\qquad\qquad\quad\nearrow y\quad\searrow z\quad\nearrow _{1_{K'}\wedge i}\quad\searrow \rho$

$\qquad\qquad K'\quad\stackrel{1_{K'}\wedge p}{\longrightarrow}\quad K'\qquad\stackrel
{x}{\longrightarrow}\qquad K$

From $(1_{K'}\wedge j)(v\wedge 1_{M})\overline{m}_{M} = v(1_{M}\wedge j)
\overline{m}_{M} = v = (1_{K'}\wedge j)\psi$ we have $(v\wedge 1_{M})\overline{m}_{M}
= \psi$ and  $d(\psi)  \in [\Sigma^{2}M, K'\wedge M]$ = 0. Similarly, from $m_{K}(x\wedge 1_{M})(1_{K'}\wedge i) = m_{K}(1_{K}\wedge i)x = x =
\rho (1_{K'}\wedge i)$ we have $\rho = m_{K}(x \wedge 1_{M})$ and
$d(\rho) \in [\Sigma K'\wedge M , K]$ = 0. Concludingly, up to sign we have \\
{\bf (9.1.19)} \quad $\rho = m_{K}(x\wedge 1_{M}), \quad\psi =
(v\wedge 1_{M})\overline{m}_M$,\quad  $d(\rho ) = 0 , d(\psi ) =
0$.

Let $\alpha ' = \alpha_{1}\wedge 1_{K} \in [\Sigma^{q-1}K,K]$, where $\alpha_{1} = j \alpha i
\in \pi_{q-1}S$, then  $j'\alpha '\alpha' = 0$ and so by (9.1.17), there exists $\alpha_{K'\wedge M}'
\in [\Sigma^{q-1}K,K'\wedge M]$ such that $\rho \alpha_{K'\wedge M}' = \alpha '$.
and $d(\alpha_{K'\wedge M}') \in [\Sigma^{q}K,K'\wedge M]$ = 0. Hence
 $\rho \alpha_{K'\wedge M}' i' = \alpha 'i'
= i'(\alpha_{1}\wedge 1_{M}) = \rho (vi\wedge 1_{M})(\alpha_{1}\wedge 1_{M})$
and so we have $\alpha_{K'\wedge M}'i' = (vi\wedge 1_{M})(\alpha_{1}\wedge 1_{M}) + \lambda \psi (ij\alpha ij)$
with $\lambda \in Z_{p}$ ,this is because $[\Sigma^{q-2}M,M] \cong Z_{p}\{ij \alpha ij\}$.
 Since the derivation $d(\alpha_{K'\wedge M}') = 0 , d(i') = 0, d(vi\wedge 1_{M}) = 0, d(\alpha_{1}
\wedge 1_{M}) = 0,  d(\psi ) = 0$ and $ d(ij\alpha ij) = - \alpha_{1}\wedge 1_{M}$,
then by applying $d$ to the above eqaution we have $\lambda\psi(\alpha_{1}\wedge 1_{M})$
= 0 so that $\lambda$ = 0.  Concludingly we have\\
{\bf (9.1.20)}$\quad \rho\alpha_{K'\wedge M}' = \alpha',\quad  \alpha_{K'\wedge M}' i' =
(vi\wedge 1_{M})(\alpha_{1}\wedge 1_{M})$,\quad   $d(\alpha_{K'\wedge M}') = 0$,

$\qquad \rho (1_{K'}\wedge ij)\alpha'_{K'\wedge M} = - \alpha ''\in
[\Sigma^{q-2}K,K]$ \\where we use  $d(\alpha'' ) = - \alpha'$ (cf. (6.5.5) ).

\vspace{2mm}

{\bf Proposition 9.1.21}\quad Let $p\geq 5$, $V$ be any spectrum and  $f : \Sigma^{t}K'\rightarrow V\wedge K$
be any map, then $f\cdot z = 0 \in [\Sigma^{t+q+1}S,V\wedge K]$.

{\bf Proof}: By Theorem 6.5.16 and Theorem 6.5.19, there is a commutative multiplication $\mu : K\wedge K
\rightarrow K$ such that $\mu (i'i\wedge 1_{K}) = 1_{K} = \mu (1_{K}\wedge i'i)$
and there is an injection $\nu : \Sigma^{q+2}K \rightarrow K\wedge K$ such that
$(jj'\wedge 1_{K})\nu = 1_{K}$.  Then by (9.1.14) we have
$  z\wedge 1_{K} = (z\wedge 1_{K})(jj'\wedge 1_{K})\nu = 0$
  and $f\cdot z = (1_V\wedge\mu )(1_{V\wedge K}\wedge i'i)f\cdot z = (1_V\wedge \mu )
(f\cdot z\wedge 1_{K}) i'i$ = 0. Q.E.D.

\quad

By (9.1.6) we have $\epsilon\cdot\overline{w} = p $(up to sign),
Then it is easy to proof that $\overline{w}\cdot\epsilon = (1_Y\wedge p)$.
By the following homotopy commutative diagram of $3\times 3$-Lemma in the stable
homotopy category

\newpage

$\qquad\quad \Sigma^qM\quad\stackrel{\alpha'i'}{\longrightarrow}\qquad\Sigma K\quad
\stackrel{r}{\longrightarrow}\qquad\Sigma Y$

$\qquad\qquad\quad\searrow j\alpha\quad\nearrow i'i\quad\searrow ^{(r\wedge 1_M)\overline{m}_K}
\nearrow _{1_Y\wedge j}$\\
{\bf (9.1.22)}$\quad\qquad\qquad \Sigma S\quad\qquad\qquad Y\wedge M$

$\qquad\qquad\quad\nearrow\epsilon\qquad\searrow\overline{w}\quad\nearrow
_{1_Y\wedge i}\quad\searrow ^{m_M(\overline{u}\wedge 1_M)}$

$\qquad\qquad Y\qquad\stackrel{1_Y\wedge p}{\longrightarrow}\quad Y\qquad
\stackrel{\overline{u}}{\longrightarrow}\quad\Sigma^{q+1}M$\newline
we know that the cofibre of $\alpha' i' : \Sigma^qM\rightarrow \Sigma K$ is $Y\wedge M$
given by the following cofibration\\
{\bf (9.1.23)} $\quad \Sigma^qM\stackrel{\alpha'i'}{\longrightarrow}\Sigma K\stackrel
{(r\wedge 1_M)\overline{m}_K}{\longrightarrow} Y\wedge M\stackrel{m_M(\overline{u}
\wedge 1_M)}{\longrightarrow}\Sigma^{q+1}M$.\\
By (9.1.10), $\alpha_1 = j''\cdot \phi$, where $\phi \in \pi_{2q-1}L$. Then we have\\
{\bf (9.1.24)} $\quad m_M(\overline{u}\bar h\wedge 1_M)(\phi\wedge 1_M) = m_M(ij''\wedge 1_M)(\phi\wedge 1_M)
= \alpha_1\wedge 1_M$ .\\
Since $\alpha'i' \cdot\alpha = 0$, then by the cofibration
(9.1.23), there is $\alpha_{Y\wedge M} \in [\Sigma^{2q+1}M,
\\Y\wedge M]$ such that $\alpha = m_M(\overline{u}\wedge
1_M)\alpha_{Y\wedge M}$. In addition, $m_M(\overline{u}\wedge
1_M)\alpha_{Y\wedge M}m_M\\(\overline{u}\wedge 1_M) = \alpha
m_M(\overline{u}\wedge 1_M) = m_M(\overline{u}\wedge
1_M)(1_Y\wedge \alpha )$ so that by (9.1.23) we have
$\alpha_{Y\wedge M}m_M(\overline{u}\wedge 1_M) = 1_Y\wedge\alpha$
modulo  $(r\wedge 1_M)\overline{m}_K*[\Sigma^{q-1}Y \wedge M , K]$
= 0, this is because $[\Sigma^qK,K] = 0 = [\Sigma^{2q}M,K]$ (cf.
Theorem 6.5.9 and Theorem 6.2.11). Note that $d(\alpha_{Y\wedge
M})\in [\Sigma^{2q+2}M,Y\wedge M]$ = 0, this is because
$[\Sigma^{q+1}M,M] = 0$ and $[\Sigma^{2q+1}M,K]$ = 0 (cf. Theorem
6.2.11 ), then $(j\overline{u}\wedge 1_M)\alpha_{Y\wedge M}\in
[\Sigma^{q-1}M,M] \cap (ker d)\cong Z_p\{\alpha_1\wedge 1_M\}$ and
so $(j\overline{u}\wedge 1_M)\alpha_{Y\wedge M} = \alpha_1\wedge
1_M$ (up to scalar). In addition, $m_M(\overline{u}\wedge
1_M)\alpha_{Y\wedge M}\cdot (\alpha_1\wedge 1_M) = \alpha
(\alpha_1\wedge 1_M) = (\alpha_1\wedge 1_M)\alpha =
m_M(\overline{u} \wedge 1_M)(\bar h\phi\wedge 1_M)\alpha$, Then by
(9.1.23) we have $\alpha_{Y\wedge M}(\alpha_1\wedge 1_M) = (\bar
h\phi \wedge 1_M)\alpha$, this is because $[\Sigma^{3q-1}M,K]$ =
0. Moreover we have, $(1_Y\wedge\alpha ) (\bar h\phi\wedge 1_M) =
\alpha_{Y\wedge M}m_M(\overline{u}\wedge 1_M)(\bar h\phi\wedge
1_M)
= \alpha_{Y\wedge M}(\alpha_1\wedge 1_M)$. Hence, we have the following relations\\
{\bf (9.1.25)} \quad $ m_M(\overline{u}\wedge 1_M)\alpha_{Y\wedge M} = \alpha,\quad$
$ \alpha_{Y\wedge M}m_M(\overline{u}\wedge 1_M) = 1_Y\wedge\alpha,\quad$

$\qquad (j\overline{u}\wedge 1_M)\alpha_{Y\wedge M} = \alpha_1\wedge 1_M$\qquad
(up to scalar)

$\qquad\alpha_{Y\wedge M}(\alpha_1\wedge 1_M) = (1_Y\wedge\alpha)(\bar h\phi\wedge 1_M)
= (\bar h\phi\wedge 1_M)\alpha$\\
where $\alpha_{Y\wedge M}\in [\Sigma^{2q+1}M, Y\wedge M]\cap (ker d)$ and $\phi\in
\pi_{2q-1}L$.

Now recall the ring spectrum properties of the spetrum $K$.  By Theorem 6.5.16 and (6.5.17), there is
a homotopy equivalence $K\wedge K = K\vee\Sigma L\wedge K\vee\Sigma^{q+2}K$  and there are
projections and injections\\
{\bf (9.1.26)}\quad $ \mu : K\wedge K\rightarrow K, \quad \mu_2 : K\wedge K\rightarrow \Sigma L\wedge K,
\quad jj'\wedge 1_K : K\wedge K\rightarrow \Sigma^{q+2}K$

$\quad i'i\wedge 1_K : K\rightarrow K\wedge K, \quad \nu_2 : \Sigma L\wedge K \rightarrow K
\wedge K, \quad \nu : \Sigma^{q+2}K\rightarrow K\wedge K$\\
such that (cf. Theorem 6.5.16, Theorem 6.5.19 )

$\quad\mu(i'i\wedge 1_K) = 1_K = \mu(1_K\wedge i'i)$, \quad $(jj'\wedge 1_K)\nu = 1_K =
(1_K\wedge jj')\nu$ ,

$\quad (i'i\wedge 1_K)\mu + \nu_2\mu_2 + (jj'\wedge 1_K)\nu = 1_{K
\wedge K}$ ,\quad $\mu_2(i'i\wedge 1_K)$ = 0.\\
Hence ,by (9.1.4), there exists $\overline{\mu}_2\in
[Y\wedge K, \Sigma L\wedge K]$ such that  $ \overline{\mu}_2(r\wedge
1_K) = \mu_2$  and $d(\overline{\mu}_2)= 0\in [Y\wedge K, L\wedge K]$,
this can be obtained from $d(\overline{\mu}_2(r\wedge 1_K) ) = d(\mu_2) = 0$
(cf. Theorem 6.5.19(H)).
By the first equation of (9.1.25) , (9.1.23)(9.1.3) and the following homotopy
commutative diagram (9.1.28) of $3\times 3$-Lemma
we know that the cofibre of $\alpha_{Y\wedge M} : \Sigma^{2q+1}M\rightarrow Y\wedge M$
is  $\Sigma L\wedge K$ given by the following cofibration\\
{\bf (9.1.27)} $\quad \Sigma^{2q+1}M\stackrel{\alpha_{Y\wedge M}}{\longrightarrow} Y
\wedge M\stackrel{\overline{\mu}_2(1_Y\wedge i')}{\longrightarrow} \Sigma
L\wedge K\stackrel{j'(j''\wedge 1_K)}{\longrightarrow}\Sigma^{2q+2}M$\\

$\qquad\qquad \Sigma^qM\quad\stackrel{\alpha'i'}{\longrightarrow}\quad\Sigma K\quad
\stackrel{i''\wedge 1_K}{\longrightarrow}\Sigma L\wedge K$

$\qquad\qquad\quad\searrow i'\quad\nearrow \alpha'\quad\searrow ^{(r\wedge 1_M)\overline{m}_K}
\nearrow\overline{\mu}_2(1_Y\wedge i')$\\
{\bf (9.1.28)}$\qquad\qquad\quad\Sigma^qK\qquad\qquad\quad Y\wedge M$

$\qquad\qquad\quad\nearrow _{(j''\wedge 1_K)}\searrow j'\quad\nearrow _{\alpha_{Y\wedge M}}\searrow
^{m_M(\overline{u}\wedge 1_M)}$

$\qquad\quad L\wedge K\stackrel{j'(j''\wedge 1_K)}{\longrightarrow}
\Sigma^{2q+1}M\quad\stackrel{\alpha}{\longrightarrow}\quad\Sigma^{q+1}M$

Since $\epsilon\wedge 1_K = \mu (i'i\wedge 1_K)(\epsilon\wedge 1_K) = 0$, then
the cofibration (9.1.4) induces a split cofibration $K\stackrel{i'i\wedge 1_K}
{\longrightarrow}K\wedge K\stackrel{r\wedge 1_K}{\longrightarrow}Y\wedge K$.
that is, there is a homotopy equivalence $K\wedge K = K\vee Y\wedge K$ so that
$Y\wedge K = \Sigma L\wedge K \vee\Sigma^{q+2}K$and there are projections
$\overline{\mu}_2 : Y\wedge K\rightarrow\Sigma L\wedge K$ , $j\overline{u}
\wedge 1_K : Y\wedge K\rightarrow\Sigma^{q+2}K$ and injections $\nu_Y : \Sigma^{q+2}K
\rightarrow Y\wedge K$, $\overline{\nu}_2 : \Sigma L\wedge K\rightarrow Y\wedge K$
such that $\nu_Y = (r\wedge 1_K)\nu$ and\\
{\bf (9.1.29)} $ \quad (j\overline{u}\wedge 1_K)\nu_Y = 1_K,\quad \overline{\mu}_2\overline
{\nu}_2  = 1_{L\wedge K},\quad \nu_Y(j\overline{u}\wedge 1_K) + \overline{\nu}_2
\overline{\mu}_2 = 1_{Y\wedge K}$.

By (9.1.1)(9.1.15)(9.1.3) and homotopy commutative diagram of $3\times 3$-Lemma we can easily know that
the cofibre of $vi : S\rightarrow K'$ is $\Sigma L$ given by the following cofibration\\
{\bf (9.1.30)} $\quad S\stackrel{vi}{\longrightarrow} K'\stackrel{k}{\longrightarrow}\Sigma L
\stackrel{\xi}{\longrightarrow}\Sigma S$\\
with relations $\xi\cdot i'' = p$ so that $\xi i''\wedge 1_M =
p\wedge 1_M = 0$ and so $\xi\wedge 1_M = \alpha (j''\wedge 1_M)$.
In addition, $\xi i''\wedge 1_K = p\wedge 1_K = 0$ so that $\xi
\wedge 1_K\in (j''\wedge 1_K)^* [\Sigma^qK,K]$ = 0. Then , the
cofibration (9.1.30) induces a split cofibration
$K\stackrel{vi\wedge 1_K} {\longrightarrow}K'\wedge
K\stackrel{k\wedge 1_K}{\longrightarrow}\Sigma L\wedge K$. That is
to say, $K'\wedge K$ splits into $K\vee \Sigma L\wedge K$ so that
there is $\nu'_2 : \Sigma L\wedge K\rightarrow K'\wedge K$ such
that $(k\wedge 1_K) \nu'_2 = 1_{L\wedge K}$ and $\mu (x\wedge 1_K)
(vi\wedge 1_K) = 1_K$, $(vi\wedge 1_K) \mu(x\wedge 1_K) + \nu'_2
(k\wedge 1_K) = 1_{K'\wedge K}$. Moreover, $x(1_{K'}\wedge
\epsilon)  = (1_K\wedge\epsilon)(x\wedge 1_Y) = 0 \in
[\Sigma^{-1}K'\wedge Y , K]$, Hence , by (9.1.14) we have,
$1_{K'}\wedge\epsilon = z\cdot\omega$,  $\omega \in [K'\wedge Y,
\Sigma^{q+2}S]$. We claim that $K'\wedge Y$ splits into
$\Sigma^{q+2}S\vee \Sigma L\wedge K$, this can be seen by the
following homotopy commutative diagram of $3\times 3$-Lemma in the
stable homotopy category

$\qquad\qquad K'\wedge Y\quad\stackrel{1_{K'}\wedge\epsilon}{\longrightarrow}\quad
\Sigma K'\quad\stackrel{x}{\longrightarrow}\quad \Sigma K$

$\qquad\qquad\qquad\searrow\widetilde{\nu}\qquad\nearrow z\quad\searrow ^{1_{K'}\wedge i'i}
\nearrow ^{\mu(x\wedge 1_K)}$\\
{\bf (9.1.31)}$\quad\qquad\qquad\Sigma^{q+2}S\qquad\qquad\Sigma K'\wedge K$

$\qquad\qquad\qquad\nearrow jj'\quad\searrow 0\quad\nearrow \nu'_2\qquad\searrow
^{1_{K'}\wedge r}$

$\quad\qquad\quad K\qquad\stackrel{0}{\longrightarrow}\quad\Sigma^2 L\wedge K\quad\stackrel{\tilde{\nu}_2}
{\longrightarrow}\quad\Sigma K'\wedge Y$. \\
That is, we have a split cofibration $\Sigma L\wedge K\stackrel{\tilde{\nu}_2}{\longrightarrow}
K'\wedge Y\stackrel{\tilde{\nu}}{\longrightarrow}\Sigma^{q+2}S$ so that there are $\widetilde{\tau} :
\Sigma^{q+2}S\rightarrow K'\wedge Y, \widetilde{\mu}_2 : K'\wedge Y\to\Sigma L\wedge K$
such that\\
{\bf (9.1.32)}\qquad $\widetilde{\nu}\cdot \widetilde{\tau} = 1_S$,\quad $\widetilde{\mu}_2\widetilde{\nu}_2 = 1_{L\wedge K}$,
\quad $\widetilde{\tau}\widetilde{\nu} + \widetilde{\nu}_2\widetilde{\mu}_2 = 1_{K'\wedge Y}$.

\vspace{2mm}

{\bf Proposition 9.1.33}\quad Let $V$ be any spectrum, then there is a direct sum decomposition

$\qquad\qquad [\Sigma^*M, V\wedge K] = (ker d) i' \oplus (ker d) i'ij$\\
where $ker d = [\Sigma^*K, V\wedge K]\cap (ker d).$

{\bf Proof} : For any $f\in [\Sigma^*M, V\wedge K]$ we have $(1_V\wedge\mu)(fi\wedge 1_K)i'i
= (1_V\wedge\mu(1_K\wedge i'i))fi = fi$, where $\mu : K\wedge K\rightarrow K$ is
the multiplication of the ring spectrum $K$ such that $\mu(i'i\wedge 1_K) = 1_K = \mu(1_K\wedge i'i)$
(cf. (9.1.26)). Then $f = (1_V\wedge\mu)(fi\wedge 1_K)i' + f_2\cdot j$
for some $f_2\in [\Sigma^{*+1}S,V\wedge K]$. It follows that $f = (1_V\wedge\mu)
(fi\wedge 1_K)i' + (1_V\wedge\mu)(f_2\wedge 1_K)i'ij$ which proves the result,
where $d(fi\wedge 1_K) = fi\wedge d(1_K) = 0$, $d(1_V\wedge\mu) = 1_V\wedge d(\mu)$ = 0
(cf. Theorem 6.5.19(G)). Q.E.D.

\quad

\begin{center}

{\bf \large \S 2.\quad A general result on convergence of $a_0$-related elements}

\vspace{2mm}

\end{center}

From [12] p. 11 Theorem 1.2.14, there is a nontrivial secondary differential in the Adams spectral
sequence $d_2(h_n) = a_0b_{n-1}, n\geq 1$,
where $d_2 : Ext_A^{1,p^nq}(Z_p,Z_p)\rightarrow Ext_A^{3,p^nq+1}(Z_p,Z_p)$ is the secondary diffenrential
of the Adams spectral sequence.  We call $h_n\in Ext_A^{1,p^nq}(Z_p,Z_p)$ and $b_{n-1}\in Ext_A^{2,p^nq}(Z_p,Z_p)$
is a pair of $a_0$-related elements. In this section, we prove a general result on convergence of
$a_0$-related elements in the Adams spectral sequence of sphere spectrum and Moore spectrum.

{\bf Definition 9.2.1}\quad Let $p\geq 7, s\leq 4$, and there is a nontrivial
secondary differential of the Adams spectral sequence
$d_2(\sigma) = a_0\sigma'$, we call $\sigma\in Ext_A^{s,tq}(Z_p,Z_p)$ and $\sigma'\in Ext_A^{s+1,tq}(Z_p,Z_p)$
is a pair of $a_0$-related elements. We have the following general result.

{\bf The main Theorem A} (the generalization of [8] Theorem A)\quad Let $p\geq 7, s\leq 4$, $\sigma$ be
the unique generator of $Ext_A^{s,tq}(Z_p,Z_p)$
and there is a nontrivial secondary differential $d_2(\sigma) = a_0\sigma'$ in the
ASS, where $\sigma'$ is the unique generator (or the linear combination of the two generators) of
$Ext_A^{s+1,tq}(Z_p,Z_p)$. Moreover, suppose that
\\(I)\quad $Ext_A^{s,tq+rq-u}(Z_p,Z_p) = 0 (r = 2,3,4, u = 1,2).$

     $Ext_A^{s+1,tq+q}(Z_p,Z_p)\cong Z_p\{h_0\sigma\},\quad Ext_A^{s+1,tq+1}(Z_p,Z_p)\cong Z_p\{a_0\sigma\}$,

     $Ext_A^{s+1,tq-q}(Z_p,Z_p) = 0$

     $Ext_A^{s+1,tq+kq+r-1}(Z_p,Z_p) = 0 (k = 2,3,4, r = 0,1)$,

     $Ext_A^{s+1,tq+kq+r-2}(Z_p,Z_p) = 0 (k = 1,2,3, r = 0,1)$.\\
(II)\quad $Ext_A^{s+2,tq+rq+u}(Z_p,Z_p)$ = 0,$r = 2,3,4, u =-1,0$ or $r = 3,4, u = 1$,

    $Ext_A^{s+2,tq}(Z_p,Z_p)$ = 0 or has unique generator $\iota$ such that
    $a_0^2\iota\neq 0$,

    $Ext_A^{s+2,tq+q}(Z_p,Z_p)\cong Z_p\{h_0\sigma'\}$ or $Z_p\{h_0\sigma'_1, h_o\sigma'_2\}$.
\\(III) $Ext_A^{s+3,tq+rq+1}(Z_p,Z_p) = 0 ( r = 1,3,4)$,

     $Ext_A^{s+3,tq+rq}(Z_p,Z_p) = 0 (r = 2,3)$

     $Ext_A^{s+3,tq+2q+1}(Z_p,Z_p)\cong Z_p\{\widetilde{\alpha}_2\sigma'\}$ or
     $Z_p\{\widetilde{\alpha}_2\sigma'_1, \widetilde{\alpha}_2\sigma'_2\}$

     $Ext_A^{s+3,tq+2}(Z_p,Z_p)\cong Z_p\{a_0^2\sigma'\}$ or $Z_p\{a_0^2\sigma'_1, a_0^2\sigma'_2\}$

     $Ext_A^{s+3,tq+1}(Z_p,Z_p)\cong Z_p\{a_0\iota\}$ or 0,\\
Then $h_0\sigma'\in Ext_A^{s+2,tq+q}(Z_p,Z_p)$ and $i_*(h_0\sigma)\in Ext_A^{s+1,tq+q}(H^*M,Z_p)$
are permanent cycles in the ASS.

\vspace{2mm}

To prove the main Theorem A, we need some preminilaries as follows.

For $(\alpha_1)_L\in [\Sigma^{q-1}L,S]$ in (9.1.10) we have $\alpha_1\cdot
(\alpha_1)_L\in [\Sigma^{2q-2}L,S]$ = 0 which is obtained from $\pi_{rq-2}S$ = 0 $(r = 2,3)$.
Then there is $\bar\phi \in [\Sigma^{2q-1}L,L]$ such that $j'' \bar\phi = (\alpha_1)_L
\in [\Sigma^{q-1}L,S]$ and $\bar\phi \cdot i''\in \pi_{2q-1}L$. Since $\pi_{rq-1}S$
has unique generator  $\alpha_1 = j\alpha i , \alpha_2 = j\alpha^2 i$ for $r = 1,2$
respectively and $j''\phi \cdot p = \alpha_1\cdot p = 0$, then $\phi\cdot p =
i''\alpha_2$ (up to scalar). That is, $i''_*\pi_{2q-1}S$ also is generated by
$\phi$ , so that $\pi_{2q-1}L\cong Z_{p^s}\{\phi\}$, for some $s\geq 1$.
Hence, $\overline{\phi} i'' = \lambda \phi$ with $\lambda\in Z_{(p)}$ and we have
$\lambda\alpha_1 = \lambda j''\phi = j''\overline{\phi} i'' = (\alpha_1)_L i''
= \alpha_1$ so that $\lambda = 1$ (mod $p$).
Moreover, $(\alpha_1)_L\overline{\phi}\in [\Sigma^{3q-2}L,S] = 0$, this is because
$\pi_{rq-2}S$ = 0  $(r = 3,4)$,  then, by (9.1.13), there is $\bar\phi_W\in [\Sigma^{3q-1}L,W]$
 such that $u\overline{\phi}_W = \overline{\phi}$. Concludingly, we have elements
$\overline{\phi}\in [\Sigma^{2q-1}L,L], \overline{\phi}_W\in [\Sigma^{3q-1}L,
W]$ such that\\
{\bf (9.2.2)}\quad $j''\bar\phi = (\alpha_1)_L , \quad \overline{\phi} i''
= \lambda\phi$, $\lambda = 1$ (mod $p$), \quad
$u\overline{\phi}_W = \overline{\phi}$.

{\bf Proposition 9.2.3} \quad Let $p\geq 7$, then\\
(1) Up to nonzero scalar we have $\phi\cdot p = i''\alpha_2 = \pi\cdot\alpha_1\neq 0$ ,
$(\alpha_1)_L\cdot \pi = \alpha_2$ ,

$ p\cdot (\alpha_1)_L = \alpha_2\cdot j'' = (\alpha_1)_L\pi j''\neq 0$,

$ [\Sigma^{2q-1}L,L]$ has unique generator $\overline{\phi}$ modulo some elements of
filtration $\geq 2$. \\
(2)  $\bar h\overline{\phi}(p\wedge 1_L)\neq 0\in [\Sigma^{2q}L,Y]$\\
(3)  $\bar h\widetilde{\phi}(\pi\wedge 1_L)(p\wedge 1_L)\neq 0\in
[\Sigma^{3q}L, Y]$ , $j''\widetilde{\phi}(\pi\wedge 1_L)\pi = j\alpha^3 i\in \pi_{3q-1}S$
(up to mod $p$ nonzero scalar) , and $\bar h\widetilde{\phi}(\pi\wedge 1_L)\pi\neq 0\in \pi_{4q}Y$
, where $\widetilde{\phi}\in [\Sigma^{2q-1}L\wedge L, L]$]
such that $\widetilde{\phi}(1_L\wedge i'') = \bar\phi$. \\
(4)  $\pi_{4q}Y$ has unique generator $\bar h\widetilde{\phi}(\pi\wedge 1_L)\pi$
such that $\bar h\widetilde{\phi}(\pi\wedge 1_L)\pi\cdot p $ = 0.

{\bf Proof}:\quad  (1) Since $j''\phi\cdot p = \alpha_1\cdot p = 0 = j''\pi
\cdot\alpha_1$, and $\pi_{2q-1}S\cong Z_{p}\{\alpha_2\}$, then $\phi\cdot p
= i''\alpha_2 = \pi\cdot\alpha_1$ (up to scalar).
We claim that $\phi\cdot p\neq 0$, this can be proved as follows.
Consider the following exact sequence

$\quad Z_p\{j\alpha^2\}\cong [\Sigma^{2q-1}M,S]\stackrel{i''_*}{\rightarrow}
[\Sigma^{2q-1}M,L]\stackrel{j''_*}{\rightarrow}[\Sigma^{q-1}M,S]\stackrel{(\alpha_1)_* }
{\rightarrow}$\\
induced by (9.1.3). The right group has unique generator $j\alpha$ satisfying
$(\alpha_1)_*j\alpha = j\alpha ij\alpha = \frac{1}{2}j\alpha\alpha ij\neq 0$, then
the above  $(\alpha_1)_*$ is monic, im$j''_*$ = 0 so that $[\Sigma^{2q-1}M,L]
\cong Z_p\{i''j\alpha^2\}$. Suppose in contrast that $\phi\cdot p = 0$, then
$\phi\in i^*[\Sigma^{2q-1}M,L]$ and so $\phi = i''j\alpha^2 i$ , $\alpha_1 = j''\phi
= j'' i''\alpha_2$ = 0 which is a contradiction.
This shows that $\phi\cdot p\neq 0$ so that the above scalar is nonzero (mod $p$).

The proof of the second result is similar. To prove the last result, let $x$ be any element of
$[\Sigma^{2q-1}L,L]$, then $j'' x\in [\Sigma^{q-1}L,S]\cong Z_{p^s}
\{(\alpha_1)_L\}$ for some $s\geq 2$. Hence, $j'' x = \lambda j''\overline{\phi}$
with $\lambda\in Z_{p^s}$ so that $x = \lambda\overline{\phi} + i'' x'$
,where $x'\in [\Sigma^{2q-1}L,S]$. Since $x' i''\in \pi_{2q-1}S\cong Z_p\{j\alpha^2i\}$
and $\pi_{3q-1}S\cong Z_p\{j\alpha^3i\}$ , then $x'$ is an element of
 filtration $\geq 2$ which shows the result.

(2) Suppose incontrast that  $\bar h\bar\phi(p\wedge 1_L) = 0$, then by (9.1.9) we have
$\bar \phi(p\wedge 1_L) = \lambda'\pi \cdot (\alpha_1)_L$ ,where $\lambda' \in Z_{(p)}$.
Since $j''\pi\wedge 1_M = p\wedge 1_M = 0$, then $\pi\wedge 1_M = (i''\wedge 1_M)\alpha$,
and so $\lambda' (\pi\wedge 1_M)i \cdot (\alpha_1)_L$ = $\lambda' (1_L\wedge i)
\pi(\alpha_1)_L$ = 0. Moreover we have $\lambda'
(i''\wedge 1_M)\alpha i(\alpha_1)_L = \lambda'(\pi\wedge 1_M)i
(\alpha_1)_L$ = 0 , then $\lambda'\alpha i(\alpha_1)_L \in (\alpha_1\wedge 1_M)_*
[\Sigma^qL,M]$ and so $\lambda'\alpha i\alpha_1 \in (\alpha_1\wedge 1_M) (i'')^*[\Sigma^{q}L,M]$ = 0
which can be obtained by the following exact sequence

$\quad [\Sigma^{2q}S,M]\stackrel{(j'')^*}{\rightarrow}[\Sigma^{q}L,M]\stackrel
{(i'')^*}{\rightarrow}[\Sigma^qS,M]\stackrel{(\alpha_1)^*}{\rightarrow}$\\
induced by (9.1.3),where the right group has unique generaator $\alpha i$ satisfying $(\alpha_1)^*\alpha i
= \alpha ij\alpha i\neq 0$ so that  $(i'')^*[\Sigma^qL,M]$ = 0.  The above equation implies that
$\lambda'$ = 0 so that we have $\bar\phi(p\wedge 1_L)$ = 0, this contradicts with the
result in (1) on $j''\overline{\phi}
(p\wedge 1_L) = p\cdot (\alpha_1)_L\neq 0$.  This shows that $\bar h\bar\phi(p\wedge 1_L)\neq 0$.

(3) Since $\pi_{rq-2}S$ = 0( $r = 2,3,4$), then $\bar\phi
(1_L\wedge\alpha_1) \in [\Sigma^{3q-2}L,L]$ = 0, and so there is
$\widetilde{\phi}\in [\Sigma^{2q-1}L\wedge L,L]$ such that
$\widetilde{\phi}(1_L\wedge i'') = \bar \phi$. We first prove
$\widetilde{\phi}(\pi\wedge 1_L)(p\wedge 1_L)\neq 0$. For
otherwise, if it is zero, then $\overline{\phi}\pi\cdot p =
\widetilde{\phi}(\pi\wedge 1_L)(p\wedge 1_L)i''$ = 0 so that
$\overline{\phi}\pi\in i^*[\Sigma^{3q-1}M,L]$. However,
$(j'')_*[\Sigma^{3q-1}M,L] \subset [\Sigma^{2q-1}M,S]$ the last of
which has unique generator $j\alpha^2$ satisfying
$(\alpha_1)_*(j\alpha^2) = j\alpha ij\alpha^2\neq 0$, then
$(j'')_*[\Sigma^{3q-1}M,L]$ = 0 and so we have $(\alpha_1)_L\pi =
j''\overline{\phi}\pi\in i^*(j'')_*[\Sigma^{3q-1}M,\\ L]$ = 0,
this contradicts with the result in (1).

Now suppose in contrast that $\bar h\widetilde{\phi}(\pi\wedge 1_L)(p\wedge 1_L) = 0$,
then by (9.1.9) we have, $\widetilde{\phi}(\pi\wedge 1_L)(p\wedge 1_L) = \pi\cdot \omega$
, where $\omega\in [\Sigma^{2q-1}L,S]$ satisfying $\omega i'' = \lambda_1
\alpha_2$ for some $\lambda_1\in Z_{p}$. It follows that $(i''\wedge 1_M)\alpha i\omega =
(1_L\wedge i)\pi\cdot\omega $ = 0, then $\alpha i\omega \in (\alpha_1\wedge 1_M)_*[\Sigma^{2q}L,M]$
and so $\lambda_1\alpha i\alpha_2 = \alpha i\omega i'' \in (\alpha_1\wedge 1_M)_*
(i'')^*[\Sigma^{2q}L,M] = (\alpha_1)^*(i'')^*[\Sigma^{2q}L.M] = 0$. This shows that
$\lambda_1$ = 0 (since $\alpha i\alpha_2 = \alpha ij\alpha^2 i\neq 0$). Hence,
$\omega = \lambda_2j\alpha^3i\cdot j''$ and $\widetilde{\phi}(\pi\wedge 1_L)(p\wedge 1_L) =
\lambda_2\pi\cdot j\alpha^3 i\cdot j''$ for some $\lambda_2\in Z_{(p)}$.
It follows that $\overline{\phi}\pi\cdot p = \widetilde{\phi}(\pi\wedge 1_L)(p\wedge 1_L) i''$ = 0
, then $\overline{\phi}\pi\in i^*[\Sigma^{3q-1}M,L]$ and so $(\alpha_1)_L\pi =
j''\overline{\phi}\pi\in i^*(j'')_*[\Sigma^{3q-1}M,L]$ = 0. This contradicts with
the result in (1) on $(\alpha_1)_L\pi\neq 0$.

For the second result, by (9.1.9) we have $\pi\cdot j = i'' j\alpha$ ,  then $j''\widetilde{\phi}(\pi\wedge 1_L)\pi\cdot j
= j''\widetilde{\phi}(\pi\wedge 1_L)i''j\alpha = j''\overline{\phi}\pi j\alpha  =
(\alpha_1)_L\pi j\alpha = \alpha_2 j\alpha = j\alpha^3ij$  (up to mod $p$
nonzero scalar) .  Consequently we have $j''\widetilde{\phi}(\pi\wedge 1_L)\pi
= j\alpha^3i$ (up to nonzero scalar),
This is because $\pi_{3q-1}S\cong Z_{p}\{\alpha_3\}$ so that $p^* \pi_{3q-1}S$ = 0.

For the last result, we first prove  $\widetilde{\phi}(\pi\wedge 1_L)\pi\neq 0$.
For otherwise , if it is zero, then $0 = \widetilde{\phi}(\pi\wedge 1_L)\pi\cdot j
= \widetilde{\phi}(\pi\wedge 1_L)i''j\alpha = \overline{\phi}\pi j\alpha$
and so $\alpha_2j\alpha = (\alpha_1)_L\pi j\alpha =
= j''\overline{\phi}\pi j\alpha$ = 0 , this is a contradiction (since
$\alpha_2j\alpha = j\alpha^2ij\alpha\neq 0\in [\Sigma^{3q-2}M,S]$).
Now suppose incontrast that $\bar h\widetilde{\phi}(\pi\wedge 1_L)\pi = 0$,
Then , by (9.1.9)  and $\pi_{3q-1}S\cong Z_p\{\alpha_3\}$ we have
$\widetilde{\phi}(\pi\wedge 1_L)\pi = \lambda \pi\cdot j\alpha^3 i =
\lambda i''j\alpha^4i$ for some $\lambda\in Z_{p}$ and so $j''\widetilde{\phi}(\pi\wedge 1_L)\pi$ = 0,
this contradicts with the second result.

(4) Since $(\overline{u})_*\pi_{4q}Y\subset \pi_{3q-1}M$ and the last of which has
unique generator
$ij\alpha^3i = ij''\widetilde{\phi}(\pi\wedge 1_L)\pi = \overline{u}\bar h
\widetilde{\phi}(\pi\wedge 1_L)\pi$(up to nonzero scalar) and $\pi_{4q-1}S\cong Z_p\{j\alpha^4i\}$
such that $(\overline{w})_*\pi_{4q-1}S$ = 0,
then $\pi_{4q}Y$ has unique generator $\bar h\widetilde{\phi}(\pi\wedge 1_L)\pi$.
Moreover by (9.1.7) we have, $\bar h\widetilde{\phi}(\pi\wedge 1_L)\pi\cdot p = \bar h(p\wedge 1_L)
\widetilde{\phi}(\pi\wedge 1_L)\pi = \bar hi''\xi\widetilde{\phi}(\pi\wedge 1_L)\pi
= \overline{w} j\alpha^4i$ = 0. This shows the result. Q.E.D.

\vspace{2mm}

{\bf Proposition 9.2.4}\quad Let $p\geq 7$, then under the supposition of the main
Theorem A we have

$Ext_A^{s+1,tq+q}(H^*L,Z_p)$ = 0, $Ext_A^{s+1,tq}(H^*L,H^*L)\cong
Z_{p}\{(\sigma')_L\}$ or\\
$Z_p\{(\sigma'_1)_L, (\sigma'_2)_L\} $, where $L$ is the spectrum in (9.1.3) and there are
relations $(i'')^*(\sigma')_L =  (i'')_*(\sigma')$  or $(i'')^*(\sigma'_1)_L = (i'')_*(\sigma'_1),
(i'')^*(\sigma'_2)_L = (i'')_*(\sigma'_2)$.

{\bf Proof}: Consider the following exact sequence

$\quad Ext_A^{s+1,tq+q}(Z_p,Z_p)\stackrel{i''_*}{\longrightarrow} Ext_A^{s+1,tq+q}(H^*L,Z_p)$

$\qquad\quad\stackrel{j''_*}{\longrightarrow}Ext_A^{s+1,tq}(Z_p,Z_p)\stackrel{(\alpha_1)_*}{\longrightarrow}$\\
induced by (9.1.3). The right group has unique generator $\sigma'$ or has two generators
$\sigma'_1, \sigma'_2$
satisfying $(\alpha_1)_*(\sigma') = h_0\sigma' \neq 0$ or $(\alpha_1)_*(\sigma'_1)
= h_0\sigma'_1\neq 0, (\alpha_1)_*(\sigma'_2) = h_0\sigma'_2\neq 0\in Ext_A^{s+2,tq+q}(Z_p,Z_p)$
(cf. the supposition II),  then the above $(\alpha_1)_*$
is monic so that im $j''_*$ = 0. Moreover, the left group has unique generator
$h_0\sigma = (\alpha_1)_*(\sigma)$ , then im $i''_*$ = 0 so that
 $Ext_A^{s+1,tq+q}(H^*L,Z_p)$ = 0. Look at the following exact sequence

$\quad 0 = Ext_A^{s+1,tq+q}(H^*L,Z_p)\stackrel{(j'')^*}{\longrightarrow}
Ext_A^{s+1,tq}(H^*L,H^*L)$

$\qquad\quad\stackrel{(i'')^*}{\longrightarrow}Ext_A^{s+1,tq}(H^*L,Z_p)\stackrel{(\alpha_1)^*}
{\longrightarrow}$\\
induced by (9.1.3). Since $Ext_A^{s+1,tq}(Z_p,Z_p)\cong Z_p\{\sigma'\}$ or $Z_{p}\{\sigma'_1, \sigma'_2\}$
and\\ $Ext_A^{s+1,tq-q}(Z_p,Z_p)$ = 0, then the right group has unique generator $(i'')_*(\sigma')$
or has two generators $(i'')_*(\sigma'_1)$, $ (i'')_*(\sigma'_2)$, the image of which
under $(\alpha_1)^*$ is zero. Then, the result on the middle group is proved. Q.E.D.

\vspace{2mm}

{\bf Proposition 9.2.5}\quad Let $p\geq 7$, then under the supposition of the main Theorem A we have\\
(1) $Ext_A^{s+3,tq+3q+1}(H^*L,Z_p)\cong Z_{p}\{\bar\phi_*\pi_*
(\sigma'_1), \bar\phi_*\pi_*(\sigma'_2)\}$  or has unique generator $\overline{\phi}_*\pi_*\sigma'$ \\
(2) $Ext_A^{s+3,tq+3q+2}(H^*Y,H^*L)\cong Z_{p}\{\bar h_*\widetilde{\phi}_*(\pi\wedge 1_L)_*(\sigma'_1)_L,
\bar h_*\widetilde{\phi}_*(\pi\wedge 1_L)_*(\sigma'_2)_L\}$ or has unique generator $\bar h_*\widetilde{\phi}_*
(\pi\wedge 1_L)_*(\sigma')_L$, where $\widetilde{\phi}\in
[\Sigma^{2q-1}L\wedge L, L]$ such that $\widetilde{\phi}(1_L\wedge i'') = \bar \phi
\in [\Sigma^{2q-1}L,L]$ (cf. Prop. 9.2.3(3)).

{\bf Proof}:  (1) Consider the following exact sequence

$\quad Ext_A^{s+3,tq+3q+1}(Z_p,Z_p)\stackrel{i''_*}{\longrightarrow}Ext_A^{s+3,tq+3q+1}
(H^*L,Z_p)$

$\qquad\quad\stackrel{j''_*}{\longrightarrow}Ext_A^{s+3,tq+2q+1}(Z_p,Z_p)\stackrel
{(\alpha_1)_*}{\longrightarrow}$\\
induced by (9.1.3). The left group is zero and the right group has
unique generator $\widetilde{\alpha}_2\sigma'$ or has two generators
$\widetilde{\alpha}_2\sigma'_1, \widetilde{\alpha}_2\sigma'_2$(cf. the supposition III).
Note that $j\alpha\alpha i = (\alpha_1)_L\cdot \pi = j''\bar\phi\cdot \pi\in \pi_{2q-1}S$,
(cf. Prop. 9.2.3), then $\widetilde{\alpha}_2\sigma'_1 = j_*\alpha_*\alpha_*i_*
(\sigma'_1) = j''_*\bar\phi_*\pi_*(\sigma'_1)$ and $\widetilde{\alpha}_2(\sigma'_2)
= j''_*\bar\phi_*\pi_*(\sigma'_2)$ so that the result on the middle group follows.

(2) Consider the following exact sequence

$\quad 0 = Ext_A^{s+3,tq+4q+1}(H^*L,Z_p)\stackrel{(j'')^*}{\longrightarrow}
Ext_A^{s+3,tq+3q+1}(H^*L,H^*L)$

$\qquad\qquad\stackrel{(i'')^*}{\longrightarrow}Ext_A^{s+3,tq+3q+1}(H^*L,Z_p)\stackrel
{(\alpha_1)^*}{\longrightarrow}$\\
induced by (9.1.3). By the supposition III,
$Ext_A^{s+3,tq+rq+1}(Z_p,Z_p)$ = 0 ( $r = 3,4$). By (1) and
$\bar\phi = \widetilde{\phi}(1_L\wedge i'')$, the right group has
unique generator $\bar\phi_*\pi_*(\sigma') =
(i'')^*(\widetilde{\phi}_*(\pi\wedge 1_L)_*(\sigma')_L$ or has two
generators $\bar\phi_*\pi_*(\sigma'_1) =
(i'')^*\widetilde{\phi}_*(\pi\wedge 1_L)_*(\sigma'_1)_L,
\bar\phi_*\pi_*(\sigma'_2) = (i'')^*\widetilde{\phi}_*(\pi\wedge
1_L)_*(\sigma'_2)_L$ the image of which under $(\alpha_1)^*$ is
zero,  then the middle group has unique generator
$\widetilde{\phi}_*(\pi\wedge 1_L)_*(\sigma')_L$ or has two
generators $\widetilde{\phi}_*(\pi\wedge 1_L)_*(\sigma'_1)_L ,
\widetilde{\phi}_*(\pi\wedge 1_L)_*(\sigma'_2)_L$. Moreover, by
$Ext_A^{s+3,tq+rq}(Z_p,Z_p)$ = 0( $r = 2,3$)  we know that
$Ext_A^{s+3,tq+2q}(Z_p,\\ H^*L)$ = 0,
 then by (9.1.9),  $Ext_A^{s+3,tq+3q+2}(H^*Y,H^*L) =  \bar
 h_*Ext_A^{s+3,tq+3q+1}\\
(H^*L,H^*L)$ and the result follows as desired. Q.E.D.

\vspace{2mm}

{\bf Proposition 9.2.6}\quad Let $p\geq 7$, then under the supposition of the main Theorem A we have

(1) $Ext_A^{s+2,tq+3q+1}(H^*Y,H^*L) = 0$,\quad $Ext_A^{s+2,tq+4q+2}
(H^*Y,Z_p)$ = 0.

(2) $Ext_A^{s+1,tq+3q+r}(H^*Y, H^*L) = 0$, $r = 0,1$.

{\bf Proof}:  (1) Consider the following exact sequence

$ Ext_A^{s+2,tq+3q}(H^*L,H^*L)\stackrel{(\bar h)_*}{\longrightarrow} Ext_A
^{s+2,tq+3q+1}(H^*Y,H^*L)$

$\qquad\quad\stackrel{(j\overline{u})_*}{\longrightarrow} Ext_A^{s+2,
tq+2q-1}(Z_p,H^*L)\stackrel{(\pi)_*}{\longrightarrow}$\\
induced by (9.1.9). \quad By the supposition II on $Ext_A^{s+2,tq+rq}(Z_p,Z_p)$ = 0\\ ( $r = 2,3,4$) we know that
the left group is zero.
By the supposition II on  $Ext_A^{s+2,tq+rq-1}(Z_p,Z_p)$ = 0($r = 2,3$) also know that
the right group is zero.  Then the middle group is zero as desired.

For the second result, consider the following exact sequence

$ Ext_A^{s+2,tq+4q+1}(H^*L,Z_p)\stackrel{(\bar h)_*}{\longrightarrow} Ext_A
^{s+2,tq+4q+2}(H^*Y,Z_p)$

$\qquad\quad\stackrel{(j\overline{u})_*}{\longrightarrow}Ext_A^{s+2,
tq+3q}(Z_p,Z_p)$\\
induced by (9.1.9).  By the supposition II on  $Ext_A^{s+2,tq+rq+1}(Z_p,Z_p)$ = 0
( $r = 3,4$) we know that the left group is zero.  Similarly, the right group also is zero.
Then the middle group is zero as desired.

(2) Consider the following exact sequence ($r = 0,1$)

$Ext_A^{s+1,tq+3q+r-1}(H^*L,H^*L)\stackrel{(\bar h)_*}{\longrightarrow}Ext_A^{s+1,tq+3q+r}
(H^*Y,H^*L)$

$\qquad\quad\stackrel{(j\overline{u})_*}{\longrightarrow} Ext_A^{s+1,tq+2q+r-2}(Z_p,H^*L)$\\
induced by (9.1.9). By the supposition I on $Ext_A^{s+1,tq+kq+r-1}(Z_p,Z_p)$ = 0
($k = 2,3,4, r = 0,1$) we know that the left group is zero.  By the supposition
 II on $Ext_A^{s+1,tq+kq+r-2}(Z_p,Z_p)$ = 0 ($k = 2,3, r = 0,1$) also know that the
right group is zero, then the middle group is zero as desired. Q.E.D.

\vspace{2mm}

{\bf Proposision 9.2.7}\quad Let $p\geq 7$, then under the supposition of the main Theorem A we have

(1) $Ext_A^{s+1,tq+3q}(H^*W,H^*L)\cong Z_p\{(\overline{\phi}_W)_*(\sigma)_L\}$,
 where $\overline{\phi}_W\in [\Sigma^{3q-1}L,\\W]$ satisfying $u\overline{\phi}_W =
\overline{\phi}\in [\Sigma^{2q-1}L,L]$ and $(\sigma)_L\in Ext_A^{s,tq}(H^*L,H^*L)$
such that $ (i'')^*(\sigma)_L = (i'')_*(\sigma)\in Ext_A^{s,tq}(H^*L,Z_p)$.

(2) $Ext_A^{s,tq+3q}(H^*Y,H^*L) = 0 ,\quad Ext_A^{s,tq+q-1}(H^*M,H^*L)$ = 0

{\bf Proof}: (1) Consider the following exact sequence

$Ext_A^{s+1,tq+3q}(H^*L,H^*L)\stackrel{w_*}{\longrightarrow}Ext_A^{s+1,tq+3q}(H^*W,H^*L)$

$\qquad\quad\stackrel{(j''u)_*}{\longrightarrow}Ext_A^{s+1,tq+q}(Z_p,H^*L)\stackrel{\phi_*}
{\longrightarrow}$\\
induced by (9.1.12). By the supposition I on  $Ext_A^{s+1,tq+rq}(Z_p,Z_p)$ = 0 $(r = 2,3,4$) we know that the left group is zero.
By $(i'')^*\cdot Ext_A^{s+1,tq+q}(Z_p,H^*L)\subset
Ext_A^{s+1,tq+q}(Z_p Z_p)$ and the last of which has unique generator $h_0\sigma = (\alpha_1)^*\cdot
(\sigma) = (i'')^*((\alpha_1)_L)^*(\sigma)$ and $Ext_A^{s+1,tq+2q}(Z_p,Z_p)$ = 0 ,
then the right group has unique generator $((\alpha_1)_L))^*(\sigma) = ((\alpha_1)_L)_*
(\sigma)_L = (j''u)_*(\overline{\phi}_W)_*(\sigma)_L$,
where $(\sigma)_L\in Ext_A^{s,tq}(H^*L,H^*L)$ satisfying $(i'')^*(\sigma)_L =
(i'')_*(\sigma)\in Ext_A^{s,tq}(H^*L,Z_p)$. Moreover, $\phi_*((\alpha_1)_L)_*
(\sigma)_L = 0\in Ext_A^{s+2,tq+3q}(H^*L,H^*L)$, then the middle group has unique generator
$(\overline{\phi}_W)_*(\sigma)_L$.

(2) Consider the following exact sequences

$Ext_A^{s,tq+3q-1}(H^*L,H^*L)\stackrel{\bar h_*}{\longrightarrow}Ext_A^{s,tq+3q}(H^*Y,H^*L)$

$\qquad\quad\stackrel{(j\overline{u})_*}{\longrightarrow}Ext_A^{s,tq+2q-2}(Z_p,H^*L)$

$Ext_A^{s,tq+q-1}(Z_p,H^*L)\stackrel{i_*}{\longrightarrow}Ext_A^{s,tq+q-1}(H^*M,H^*L)$

$\qquad\quad\stackrel{j_*}{\longrightarrow}Ext_A^{s,tq+q-2}(Z_p,H^*L)$\\
induced by (9.1.9) and (9.1.1) respectively.\quad  By the
supposition I on\\ $Ext_A^{s,tq+rq-1}(Z_p,Z_p)$ = 0( $r = 2,3,4$)
we know that the upper left group is zero. By the supposition I on
$Ext_A^{s,tq+rq-2}(Z_p,Z_p)$ = 0 ( $r = 2,3$ ), the upper right
group is zero. Then the upper middle group is zero as desired.
Similarly, the lower middle also is zero. Q.E.D.

\vspace{2mm}

{\bf Proposition 9.2.8}\quad Let $p\geq 7$, then under the supposition (I)(III) of the main Theorem A we have

 $Ext_A^{s+3,tq+2}(H^*M, Z_p)$ = 0,

 $Ext_A^{s+1,tq+q+1}(H^*M\wedge L, Z_p)\cong Z_p\{(i\wedge 1_L)_*\pi_*(\sigma)\}$.

{\bf Proof}: Consider the following exact sequence

$Ext_A^{s+3,tq+2}(Z_p,Z_p)\stackrel{i_*}{\longrightarrow}Ext_A^{s+3,tq+2}(H^*M,Z_p)$

$\qquad\quad\stackrel{j_*}{\longrightarrow}Ext_A^{s+3,tq+1}(Z_p,Z_p)\stackrel{p_*}{\longrightarrow}$\\
induced by (9.1.1). By the supposition III, the right group is zero or has unique generator  $a_0\iota$
which satisfies $p_*(a_0\iota) = a_0^2\iota\neq 0\in Ext_A^{s+4,tq+2}(Z_p,Z_p)$,
then  im $j_*$ = 0. By the supposition III, the left group has unique generator $a_0^2\sigma'$or
has two generators  $a_0^2\sigma'_1
= p_*(a_0\sigma'_1), a_0^2\sigma'_2 = p_*(a_0\sigma'_2$ so that we have im $i_*$ = 0.
Then, the middle group is zero as desired.

For the second result, consider the following exact sequence

$Ext_A^{s+1,tq+q+1}(H^*L,Z_p)\stackrel{(i\wedge 1_L)_*}{\longrightarrow}Ext_A
^{s+1,tq+q+1}(H^*M\wedge L, Z_p)$

$\qquad\quad\stackrel{(j\wedge 1_L)_*}{\longrightarrow}
Ext_A^{s+1,tq+q}(H^*L,Z_p)$\\
induced by (9.1.1). By Prop. 9.2.4, the right group is zero.
Since $(j'')_*\\Ext_A^{s+1,tq+q+1}(H^*L,Z_p)\subset Ext_A^{s+1,tq+1}(Z_p,Z_p)\cong
Z_p\{a_0\sigma = (j'')_*\pi_*(\sigma)\}$ and $Ext_A^{s+1,tq+q+1}(Z_p,Z_p)$ = 0
then the left group has unique generator $\pi_*(\sigma)$ and the result follows. Q.E.D.

\vspace{2mm}

Now we begin to prove the main Theorem A. The proof will be done by processing
an argument processing in the Adams resolution of some spectra related to $S$.
Let\\
{\bf (9.2.9)} $\quad\qquad\cdots \stackrel{\bar a_{2}}{\longrightarrow}\quad\Sigma^{-2}E_{2}\quad\stackrel
{\bar a_{1}}{\longrightarrow}\quad\Sigma^{-1}E_{1}\quad\stackrel{\bar a_{0}}{\longrightarrow} E_{0} = S$

$\qquad\qquad\qquad\qquad\qquad\qquad\big\downarrow\bar b_{2}\qquad\qquad\quad\big\downarrow\bar b_{1}
\qquad\qquad\big\downarrow\bar b_{0}$

$\qquad\qquad\qquad\qquad\qquad\Sigma^{-2}KG_{2}\qquad\quad\Sigma^{-1}KG_{1}\qquad\quad KG_{0}$\\
be the minimal Adams resolution of the sphere spectrum $S$ which satisfies\\
(1) $E_{s}\stackrel{\bar b_{s}}{\rightarrow} KG_{s}\stackrel{\bar c_{s}}{\rightarrow} E_{s+1}
\stackrel{\bar a_{s}}{\rightarrow}\Sigma E_{s}$ are cofibrations for all $s\geq 0$,
which induce short exact sequences in $Z_p$-cohomology $0\rightarrow H^{*}E_{s+1}\stackrel
{\bar c_{s}^{*}}{\rightarrow} H^{*}KG_{s}\stackrel{\bar b_{s}^{*}}{\rightarrow}
H^{*}E_{s}\rightarrow 0$.\\
(2) $KG_{s}$ is a graded wedge sum of Eilenberg-Maclane spectra $KZ_p$ of type $Z_p$.\\
(3) $\pi_{t}KG_{s}$ are the $E_{1}^{s,t}$-terms of the Adams spectral sequence,
$(\bar b_{s}\bar c_{s-1})_{*} : \pi_{t}KG_{s-1}\rightarrow \pi_{t}KG_{s}$ is the
$d_1^{s-1,t}$-differentials of the Adams spectral sequence ,
and $\pi_{t}KG_{s}\cong Ext_{A}^{s,t}(Z_{p},Z_{p})$ (cf. [3] p.180).\\
Then, an Adams resolution of an arbitrary spectrum $V$ can be obtained by smashing
$V$ to (9.2.9). We first prove some Lemmas.

\vspace{2mm}

{\bf Lemma 9.2.10}\quad Let $p\geq 7$, then under the supposition of the main Theorem A we have\\
\centerline{$\bar c_{s+1}\cdot h_0\sigma = (1_{E_{s+2}}\wedge\alpha_1)\kappa$\qquad (up to scalar)}\\
where $\kappa\in \pi_{tq+1}E_{s+2}$ such that $\bar c_{s+1}\cdot \sigma = \bar a_{s+1}\cdot
\kappa$ and $\bar b_{s+2}\cdot \kappa = a_0\sigma' \in \pi_{tq+1}
KG_{s+2}\cong Ext_A^{s+2,tq+1}(Z_p,Z_p)$.

{\bf Proof}: The $d_1$-cycle $(1_{KG_{s+1}}\wedge i'')h_0\sigma\in \pi_{tq+q}KG_{s+1}\wedge L$
represents an element in $Ext_A^{s+2,tq+q}(H^*L,Z_p)$  and this group is zero
by Prop. 9.2.4 ,
then it is a $d_1$-boundary and so $(\bar c_{s+1}\wedge 1_L)(1_{KG_{s+1}}
\wedge i'')h_0\sigma$ = 0 , $\bar c_{s+1}\cdot h_0\sigma
= (1_{E_{s+2}}\wedge\alpha_1)f''$ for some $f''\in \pi_{tq+1}E_{s+2}$.
It follows that $\bar a_{s+1}\cdot (1_{E_{s+2}}\wedge\alpha_1)f''$ = 0
, then $\bar a_{s+1}\cdot f'' = (1_{E_{s+1}}\wedge j'')f''_2$ with $f''_2\in \pi_{tq+q}(E_{s+1}\wedge L)$.
The $d_1$-cycle $(\bar b_{s+1}\wedge 1_L)f''_2\in \pi_{tq+q}KG_{s+1}\wedge L$ represents
an element in $Ext_A^{s+1,tq+q}(H^*L, Z_p)$ and this group is zero,
then $(\bar b_{s+1}\wedge 1_L)f''_2 = (\bar b_{s+1}\bar c_{s}\wedge 1_L)g''$ with $g''\in\pi_{tq+q}(KG_s\wedge L)$.
Hence, $f''_2 = (\bar c_s\wedge 1_L)g''+ (\bar a_{s+1}\wedge 1_L)f''_3$, for some $f''_3\in \pi_{tq+q+1}E_{s+2}\wedge L$ and we have
$\bar a_{s+1}\cdot f'' = \bar a_{s+1} (1_{E_{s+2}}\wedge j'')f''_3 + \bar c_s (1_{KG_s}
\wedge j'')g'' = \bar a_{s+1} (1_{E_{s+2}}\wedge j'')f''_3 + \lambda \bar c_s\cdot \sigma$
$ = \bar a_{s+1} (1_{E_{s+2}}\wedge j'')f''_3 + \lambda \bar a_{s+1}\cdot \kappa$
for some $\lambda \in Z_{p}$, this is because $(1_{KG_s}\wedge j'')g''\in \pi_{tq}KG_s
\cong Ext_A^{s,tq}(Z_p,Z_p)\cong Z_p\{\sigma\}$. Then, $f'' = (1_{E_{s+2}}\wedge j'')f''_3
+ \lambda\kappa + \bar c_{s+1}\cdot g''_2$ for some $g''_2\in \pi_{tq+1}KG_{s+1}$
and so $\bar c_{s+1}\cdot h_0\sigma =
(1_{E_{s+1}}\wedge\alpha_1)\kappa$ (up to scalar). Q.E.D.

\vspace{2mm}

Since $\bar h\phi\cdot p = \bar h i''j\alpha^2i$ = 0 (cf. Prop. 9.2.3(1) and (9.1.9)(9.1.5)),
then $\bar h\phi = (1_Y\wedge j)\alpha_{Y\wedge M}i$ , where $\alpha_{Y\wedge M}
\in [\Sigma^{2q+1}M, Y\wedge M]$. Let $\Sigma U$ be the cofibre of $\bar h\phi =
(1_Y\wedge j)\alpha_{Y\wedge M}i : \Sigma^{2q}S\rightarrow  Y$
given by the cofibration\\
{\bf (9.2.11)}\qquad $\Sigma^{2q}S\stackrel{\bar h\phi}{\longrightarrow}Y\stackrel{w_2}
{\longrightarrow}\Sigma U\stackrel{u_2}{\longrightarrow}\Sigma^{2q+1}S$\\
Moreover,  $w_2(1_Y\wedge j)\alpha_{Y\wedge M} = \widetilde{w}\cdot j$, where
$\widetilde{w} : \Sigma^{2q}S\rightarrow U$ whose cofibre is $X$ given by the cofibration
$\Sigma^{2q}S\stackrel{\tilde{w}}{\longrightarrow}
U\stackrel{\tilde{u}}{\longrightarrow}X\stackrel{j\tilde{\psi}}{\longrightarrow}
\Sigma^{2q+1}S$. Then, $\Sigma X$ also is the cofibre of $\omega = (1_Y\wedge j)
\alpha_{Y\wedge M} : \Sigma^{2q}M\rightarrow Y$ given by the cofibration\\
{\bf (9.2.12)}\qquad $\Sigma^{2q}M\stackrel{(1_Y\wedge j)\alpha_{Y\wedge M}}{\longrightarrow}
 Y\stackrel{\tilde{u}w_2}
{\longrightarrow}\Sigma X\stackrel{\tilde{\psi}}{\longrightarrow}\Sigma^{2q+1}M$\\
This can be seen by the following homotopy commutative diagram of $3\times 3$-Lemma

\quad

\quad

$\qquad\quad\Sigma^{2q}S\qquad\stackrel{\bar h\phi}{\longrightarrow}\qquad Y\qquad
\stackrel{\tilde{u}w_2}{\longrightarrow}\qquad\Sigma X$

$\qquad\qquad\quad\searrow i\quad\nearrow _{(1_Y\wedge j)\alpha_{Y\wedge M}}
\searrow w_2\quad\nearrow\widetilde{u}\quad\searrow \widetilde{\psi}$\\
{\bf (9.2.13)}\qquad\qquad\quad $\Sigma^{2q}M\qquad\qquad\qquad\Sigma U\qquad\qquad\Sigma^{2q+1}M$

$\qquad\qquad\quad\nearrow\widetilde{\psi}\qquad\searrow j\qquad\nearrow\widetilde{w}\qquad
\searrow u_2\quad\nearrow i$

$\qquad\qquad X\quad\stackrel{j\tilde{\psi}}{\longrightarrow}\qquad\Sigma^{2q+1}S\quad
\stackrel{p}{\longrightarrow}\qquad\Sigma^{2q+1}S$\\

Since $j\overline{u}(\bar h\phi)$ = 0, then by (9.2.11) we have, $j\overline{u} = u_3w_2$,
for some $u_3\in [U, \Sigma^{q+1}S]$. Hence, the spectrum $U$ in (9.2.11) also is the cofire of
$w\pi : \Sigma^qS\rightarrow W$ given by the cofibration\\
{\bf (9.2.14)}\qquad $\Sigma^qS\stackrel{w\pi}{\longrightarrow}W\stackrel{w_3}{\longrightarrow}
U\stackrel{u_3}{\longrightarrow} \Sigma^{q+1}S$\\
This can be seen by the following homotopy commutative diagram of $3\times 3$- Lemma in the
stable homotopy category

$\qquad\quad\Sigma^{-1}Y\quad\stackrel{j\overline{u}}{\longrightarrow}\quad
\Sigma^{q+1}S\quad\stackrel{w\pi}{\longrightarrow}\quad \Sigma W$

$\qquad\qquad\quad\searrow w_2\quad\nearrow u_3\quad\searrow \pi\quad\nearrow w\quad\searrow w_3$\\
{\bf (9.2.15)} $\quad\qquad\qquad U\qquad\qquad\qquad \Sigma L\qquad\qquad\Sigma U$

$\qquad\qquad\quad\nearrow w_3\quad\searrow u_2\quad\nearrow \phi\quad\searrow \bar h\quad\nearrow w_2$

$\qquad\quad W\qquad\stackrel{j''u}{\longrightarrow}\quad\Sigma^{2q}S\qquad\stackrel
{\bar h\phi}{\longrightarrow}\qquad Y$\\
Moreover, by $u_3\widetilde{w} = \alpha_1$, the cofibre of $\widetilde{u}w_3 : W\rightarrow X$
is $\Sigma^{q+1}L$ given by the cofibration\\
{\bf (9.2.16)} \qquad $W\stackrel{\tilde{u}w_3}{\longrightarrow}X\stackrel{u''}{\longrightarrow}
\Sigma^{q+1}L\stackrel{w'(\pi\wedge 1_L)}{\longrightarrow} \Sigma W$\\
where $w'\in [L\wedge L, W]$ such that $w'(1_L\wedge i'') = w$.
This can be seen by the following homotopy commutative diagram of
$3\times 3$-Lemma in the stable homotopy category

$\qquad\qquad W\quad\stackrel{\tilde{u}w_3}{\longrightarrow}\quad X\quad
\stackrel{j\tilde{\psi}}{\longrightarrow}\quad \Sigma^{2q+1}S$

$\qquad\qquad\quad\searrow w_3\quad\nearrow\widetilde{u}\quad\searrow u''\quad\nearrow j''$\newline
{\bf (9.2.17)}\quad\qquad\qquad $U\qquad\qquad\quad\Sigma^{q+1}L$

$\qquad\qquad\quad\nearrow \widetilde{w}\quad\searrow u_3\quad\nearrow i''
\quad \searrow ^{w'(\pi\wedge 1_L)}$

$\qquad\qquad \Sigma^{2q}S\quad\stackrel{\alpha_1}{\longrightarrow}\quad\Sigma^{q+1}S\quad
\stackrel{w\pi}{\longrightarrow}\quad\Sigma W$

\vspace{2mm}

{\bf Lemma 9.2.18} \quad (1) Let $\overline{\phi}_W\in [\Sigma^{3q-1}L,W]$
be the map in Prop. 9.2.7 such that $u\overline{\phi}_W = \overline{\phi}\in [\Sigma^{2q-1}L,L]$,
then

(1) $\widetilde{u}w_3\overline{\phi}_W(p\wedge 1_L)\neq 0\in [\Sigma^{3q-1}L,X]$.

(2) $Ext_A^{s,tq+3q-1}(H^*X,H^*L)$ = 0,

\qquad $ Ext_A^{s+1,tq+3q}(H^*X,H^*L) = (\widetilde{u}w_3)_*Ext_A^{s+1,tq+3q}(H^*W,H^*L)$.

{\bf Proof}: (1) Suppose in contrast that $\widetilde{u}w_3\overline{\phi}_W(p\wedge 1_L)$
= 0, then by (9.2.16) and the result of Prop. 9.2.3(1) on $[\Sigma^{2q-1}L,L]$ we have\\
(9.2.19) $\quad \overline{\phi}_W(p\wedge 1_L) = \lambda w'(\pi\wedge 1_L)\overline{\phi}$
\quad modulo $F_3[\Sigma^{3q-1}L,W]$\\
for some $\lambda\in Z_{(p)}$, where $F_3[\Sigma^{3q-1}L,W]$
denotes the subgroup of $[\Sigma^{3q-1}L,W]$ consisting by all
elements of filtration $\geq 3$. Moreover, note that
$uw'(\pi\wedge 1_L)\in [L,L]$ and this group has two generators
$(p\wedge 1_L), \pi j''$ which has filtration one, then
$uw'(\pi\wedge 1_L) = \lambda_1(p\wedge 1_L) + \lambda_2 \pi j''$
with $\lambda_1,\lambda_2\in Z_{(p)}$. By (9.1.13) we have
 $\lambda_1 p\cdot (\alpha_1)_L + \lambda_2(\alpha_1)_L \pi j''$ = 0 so that
 $\lambda_2 = \lambda_0\lambda_1$, here we use the equation
$(\alpha_1)_L\pi j'' = - \lambda_0 p\cdot (\alpha_1)_L$,  $\lambda_0\neq 0\in Z_p$.
Then , by composing $u$ on (9.2.19) we have$\overline{\phi}
(p\wedge 1_L) = u\overline{\phi}_W(p\wedge 1_L) = \lambda uw'(\pi\wedge 1_L)
\overline{\phi} = \lambda\lambda_1\overline{\phi}(p\wedge 1_L) + \lambda\lambda_0\lambda_1 \pi j''\overline{\phi}$
(mod $F_3[\Sigma^{2q-1}L,L]$) and so by (9.1.9) $\bar h\overline{\phi}
(p\wedge 1_L) = \lambda\lambda_1\bar h\overline{\phi}(p\wedge 1_L)$ (mod $F_3[\Sigma^{2q}L,Y]$).
This implies that $\lambda\lambda_1 = 1$ (mod $p$) (cf. the following Remark 9.2.20).

Hence we have $\lambda\lambda_1\lambda_0\pi j''\overline{\phi}$ = 0
(mod $F_3[\Sigma^{2q-1}L,L]$) and by the same reason as shown in the following Remark 9.2.20,
this implies that $\lambda\lambda_1\lambda_0$ = 0 (mod $p$) which yields a contradiction.

(2) Consider the following exact sequence

$ Ext_A^{s,tq+3q}(H^*Y,H^*L)\stackrel{(\tilde{u}w_2)_*}{\longrightarrow}
Ext_A^{s,tq+3q-1}(H^*X,H^*L)$

$\qquad\quad\stackrel{(\tilde{\psi})_*}{\longrightarrow}
Ext_A^{s,tq+q-1}(H^*M,H^*L)$ \\
induced by (9.2.12). By Prop. 9.2.7(2), both sides of groups are zero ,so tha the middle
group is zero as desired.  Look at the following exact sequence

$ Ext_A^{s+1,tq+3q}(H^*W,H^*L)\stackrel{(\tilde{u}w_3)_*}{\longrightarrow}
Ext_A^{s+1,tq+3q}(H^*X,H^*L)$

$\qquad\quad\stackrel{(u'')_*}{\longrightarrow}Ext_A^{s+1,tq+2q-1}(H^*L,H^*L)$\\
induced by (9.2.16). By the supposition on $Ext_A^{s+1,tq+rq-1}(Z_p,Z_p)$ = 0
( $r = 1,2,3$) we know that the right group is zero and so the result follows. Q.E.D.

\vspace{2mm}

{\bf Remark 9.2.20} \quad Here we give an explanation on the reason why the coefficient
in the equation $(1 - \lambda\lambda_1)\bar h\overline{\phi}(p\wedge 1_L) = 0 $
(mod $F_3[\Sigma^{2q}L,Y]$) must be zero (mod $p$). For otherwise ,
if $1 - \lambda\lambda_1 \neq 0$ (mod $p$), then $(1 - \lambda\lambda_1)
\bar h\overline{\phi}(p\wedge 1_L)$ must be represented by some nonzero element
 $x\in Ext_A^{2,2q+2}(H^*Y,H^*L)$ in the ASS. However , it also equals
 to an element of filtration $\geq 3$, then
$x$ must be a $d_2$-boundary, that is, $x = d_2(x')\in d_2Ext_A^{0,2q+1}(H^*Y,H^*L)$
= 0, this is because $Ext_A^{0,2q+1}(H^*Y,H^*L) = Hom_A^{2q+1}(H^*Y,H^*L)$ = 0 which is
obtained by $H^rL \neq 0$ only for $r = 0,q$.
This is a contradiction so that we have $1 - \lambda\lambda_1 = 0$ (mod $p$).

\vspace{2mm}

{\bf Lemma 9.2.21}\quad For the element $\kappa\in \pi_{tq+1}E_{s+2}$ in Lemma 9.2.19,
it is known that $\bar a_{s+1}\cdot\kappa = \bar c_s\cdot \sigma$ and $\bar b_{s+2}\cdot\kappa
= a_0\sigma'\in \pi_{tq+1}KG_{s+2}\cong Ext_A^{s+2,tq+1}(Z_p,\\Z_p)$, then there exists
$f\in \pi_{tq+3}E_{s+4}\wedge M$ and $g\in \pi_{tq+1}(KG_{s+1}\wedge M)$ such that\\
(A)\qquad $ (1_{E_{s+2}}\wedge i)\kappa = (\bar c_{s+1}\wedge 1_M)g + (\bar a_{s+2}\bar a_{s+3}\wedge 1_M)f$\newline
and\\
(B) \qquad $(1_{E_{s+4}}\wedge (1_Y\wedge j)\alpha_{Y\wedge M})f\cdot
(\alpha_1)_L$ = 0 $\in [\Sigma^{tq+4q+2}L, E_{s+4}\wedge Y]$,\\
where $\alpha_{Y\wedge M}\in [\Sigma^{2q+1}M,Y\wedge M]$ satisfying $(1_Y
\wedge j)\alpha_{Y\wedge M}i = \bar h\phi\in \pi_{2q}Y$.

{\bf Proof}: Note that the $d_1$-cycle
$(\bar b_{s+2}\wedge 1_M)(1_{E_{s+2}}\wedge i)\kappa\in \pi_{tq+1}KG_{s+2}\wedge M$
represents an element $i_*(a_0\sigma') = i_*p_*(\sigma') = 0 \in Ext_A^{s+2,tq+1}
(H^*M,Z_p)$ so that it is a $d_1$-boundary. That is, $(\bar b_{s+2}\wedge 1_M)(1_{E_{s+2}}
\wedge i)\kappa = (\bar b_{s+2}\bar c_{s+1}\wedge 1_M)g$ for some $g\in \pi_{tq+1}KG_{s+1}\wedge M$
. Then, by $Ext_A^{s+3,tq+2}(H^*M,Z_p)$ = 0 (cf. Prop. 9.2.8) we have
$(1_{E_{s+2}}\wedge i)\kappa = (\bar c_{s+1}\wedge 1_M)g + (\bar a_{s+2}\bar a_{s+3}\wedge 1_M)f$
for some $f\in \pi_{tq+3}E_{s+4}\wedge M$. This shows (A).
For (B), by Prop. 9.2.3(1) we have $\phi\cdot p
= i''j\alpha^2i$ (up to nonzero scalar), then $\bar h\phi\cdot p = \bar h i''j\alpha^2 i = 0$
and so $\bar h \phi = (1_Y\wedge j)\alpha_{Y\wedge M} i$ for some $\alpha_{Y\wedge M}
\in [\Sigma^{2q+1}M,Y\wedge M]$. Then, by composing  $1_{E_{s+2}}\wedge (1_Y\wedge j)\alpha_{Y\wedge M}$
on the equation (A) we have\\
{\bf (9.2.22)}\qquad $(1_{E_{s+2}}\wedge \bar h\phi)\kappa = (1_{E_{s+2}}\wedge (1_Y\wedge j)\alpha_{Y\wedge M}i)
\kappa $

$ \quad = (\bar a_{s+2}\bar a_{s+3}\wedge 1_Y)(1_{E_{s+4}}\wedge (1_Y\wedge j)
\alpha_{Y\wedge M})f$\\
where $(1_Y\wedge j)\alpha_{Y\wedge M}$  induces zero homomorphism in $Z_p$-cohomology so that
 $(\bar c_{s+1}\wedge 1_Y)(1_{KG_{s+1}}\wedge (1_Y\wedge j)\alpha_{Y\wedge M})g$ = 0.

By composing $(\alpha_1)_L$ on (9.2.22) we have  $(\bar a_{s+2}\bar a_{s+3}\wedge 1_Y)
(1_{E_{s+4}}\wedge (1_Y\wedge j)\alpha_{Y\wedge M})f\cdot (\alpha_1)_L = (1_{E_{s+2}}
\wedge\bar h)(\kappa\wedge 1_L)\phi\cdot (\alpha_1)_L$ = 0, this is because $\phi\cdot (\alpha_1)_L
\in [\Sigma^{3q-2}L,L]$ = 0 which is obtained by $\pi_{rq-2}S$ = 0( $r = 2,3,4$). Then we have\\
\centerline{$(\bar a_{s+3}\wedge 1_Y)(1_{E_{s+4}}\wedge (1_Y\wedge j)\alpha_{Y\wedge M})f\cdot (\alpha_1)_L
= (\bar c_{s+2}\wedge 1_Y)g_1$ = 0}\\ where the  $d_1$-cycle $g_1\in [\Sigma^{tq+3q+1}L,
KG_{s+2}\wedge Y]$ represents an element in $Ext_A^{s+2,tq+3q+1}(H^*Y,H^*L)$ and this group is zero (cf. Prop.
9.2.6(1)) so that it is a $d_1$-boundary and we have $(\bar c_{s+2}\wedge 1_Y)g_1$ = 0.
Briefly write $(1_Y\wedge j)\alpha_{Y\wedge M} = \omega$  and let  $V$ be the cofibre of
$(1_Y\wedge (\alpha_1)_L)(\omega\wedge 1_L) = \omega\cdot
(1_M\wedge (\alpha_1)_L ): \Sigma^{3q-1}M\wedge L\rightarrow Y$ given by the cofibration\\
{\bf (9.2.23)}\qquad $\Sigma^{3q-1}M\wedge L\stackrel{(1_Y\wedge (\alpha_1)_L)(\omega\wedge 1_L)}
{\longrightarrow}Y\stackrel{w_4}{\longrightarrow}V\stackrel{u_4}{\longrightarrow}
\Sigma^{3q}M\wedge L$\\
It follows that $(\bar a_{s+3}\wedge 1_Y)(1_{E_{s+4}}\wedge 1_Y\wedge (\alpha_1)_L)
(\omega\wedge 1_L)(f\wedge 1_L) = (\bar a_{s+3}\wedge 1_Y)(1_{E_{s+4}}\wedge (1_Y\wedge j)\alpha_{Y\wedge M})f\cdot
(\alpha_1)_L$ = 0, then by (9.2.23) we have $(\bar a_{s+3}\wedge 1_{M\wedge L})
(f\wedge 1_L) = (1_{E_{s+3}}\wedge u_4)f_2$ for some $f_2 \in [\Sigma^{tq+3q+2}L, E_{s+3}\wedge V]$.
Consequently, $(\bar b_{s+3}\wedge 1_{M\wedge L})(1_{E_{s+3}}\wedge u_4)f_2$ = 0  so that we have\\
{\bf (9.2.24)}\qquad $(\bar b_{s+3}\wedge 1_V)f_2 = (1_{KG_{s+3}}\wedge w_4)g_2$ \\
with $g_2\in [\Sigma^{tq+3q+2}L, KG_{s+3}\wedge Y]$. Then, $(\bar b_{s+4}
\bar c_{s+3}\wedge 1_V)(1_{KG_{s+3}}\wedge w_4)g_2$ = 0
and so $(\bar b_{s+4}\bar c_{s+3}\wedge 1_Y)g_2 \in (1_{KG_{s+4}}\wedge (1_Y\wedge
(\alpha_1)_L(\omega\wedge 1_L))_*[\Sigma^*L, KG_{S+4}\wedge M\wedge L]$ = 0. That is,
$g_2$ is a  $d_1$-cycle which represents an element $[g_2] \in Ext_A^{s+3,tq+3q+2}
(H^*Y\\,H^*L)$ and this group has two generators as shown in Prop. 9.2.5(2)), then we have\\
{\bf (9.2.25)} \qquad $[g_2] = \bar h_*\widetilde{\phi}_*(\pi\wedge 1_L)_*(\lambda_1[\sigma'_1\wedge 1_L]
+ \lambda_2[\sigma'_2\wedge 1_L])$ \\ with $\lambda_1, \lambda_2\in Z_p$.
By (9.2.24) we know that \\
\centerline{$(w_4)_*[g_2]\in E_2^{s+3,tq+3q+2}(V) = Ext_A^{s+3,tq+3q+2}
(H^*V,H^*L)$}\\ is a permanent cycle in the ASS. However, $(1_Y\wedge (\alpha_1)_L)(\omega\wedge 1_L)$
is a map of filtration 2, then the cofibration (9.2.23)  induces an exact sequence
in $Z_p$-cohomology which is split as $A$-module.  That is , it induces a split exact sequence in the $E_1$-term of the ASS
: $E_1^{s+3,*}(Y)\stackrel{(w_4)_*}
{\longrightarrow}E_1^{s+3,*}(V)\stackrel{(u_4)_*}{\longrightarrow}E_1^{s+3,*-3q}(M\wedge L)$.
It follows that  it induces a split exact sequence in the $E_r$-term of the ASS for all $(r\geq 2$)\\
{\bf (9.2.26)} \qquad\qquad\quad $E_r^{s+3,*}(Y)\stackrel{(w_4)_*}{\longrightarrow}
E_r^{s+3,*}(V)\stackrel{(u_4)_*}{\longrightarrow}E_r^{s+3,*-3q}(M\wedge L)$\\
 Then , $d_r((w_4)_*[g_2]) = 0$ implies that $d_r
([g_2]) = 0$( $r\geq 2$). That is , (9.2.24) implies that $[g_2]$ also is a permanent cycle
in the ASS. Since the secondary differential
$ d_2[g_2]$ = 0 and $d_2(\sigma) = a_0\sigma'$ in which $\sigma'$ is the linear combination of
$\sigma'_1,\sigma'_2$, then
 $\lambda_1 ,\lambda_2$ linearly dependent. That is, (9.2.25)  becomes\\
\qquad\quad $[g_2] = \lambda_1\bar h_*\widetilde{\phi}_*(\pi\wedge 1_L)_*[\sigma'\wedge 1_L]$.
\\Now we consider the case $\lambda_1$ is nonzero or zero respectively.

If $\lambda_1\neq 0$ , (9.2.24) imp;ies $[g_2] $ and so $\bar h_*\widetilde{\phi}_*
(\pi\wedge 1_L)_*[\sigma'\wedge 1_L]\in E_2^{s+3,tq+3q+2}\\(Y) = Ext_A^{s+3,tq+3q+2}(H^*Y,H^*L)$
is a permanent cycle in the ASS. Moreover, by
$(\bar a_{s+3}\wedge 1_Y)(1_{E_{s+4}}\wedge (1_Y\wedge j)\alpha_{Y\wedge M})f
\cdot (\alpha_1)_L$ = 0 we have

\quad\qquad $(1_{E_{s+4}}\wedge (1_Y\wedge j)\alpha_{Y\wedge M})f\cdot
(\alpha_1)_L = (\bar c_{s+3}\wedge 1_Y)g_3$
\\ for some $d_1$-cycle $g_3\in [\Sigma^{tq+3q+2}L, KG_{s+3}\wedge Y]$
and it represents an element  $[g_3]\in Ext_A^{s+3,tq+3q+2}(H^*Y,H^*L)$ so that we have
$[g_3] = \bar h_*\widetilde{\phi}_*(\pi\wedge 1_L)_*(\lambda_3[\sigma'_1\wedge 1_L] +
\lambda_4[\sigma'_2\wedge 1_L])$ for some $\lambda_3, \lambda_4\in Z_p$.
By the above equation and  $(1_Y\wedge (\alpha_1)_L)(\omega\wedge 1_L)$
has filtration 2 we know that the secondary differential $d_2([g_3])$ = 0 so that by the similar reason as above
, $\lambda_3 , \lambda_4 $ is linearly dependent . That is,
$[g_3] = \lambda_3\bar h_*\widetilde{\phi}_*(\pi\wedge 1_L)_*
[\sigma'\wedge 1_L]$ so that we have $(1_{E_{s+4}}\wedge (1_Y\wedge (\alpha_1)_L)
(\omega\wedge 1_L))(f\wedge 1_L) = (\bar c_{s+3}\wedge 1_Y)g_3 = 0$
and the result follows.

If $\lambda_1$ = 0, then $g_2 = (\bar b_{s+3}\bar c_{s+2}\wedge 1_Y)g_4$ for some
$g_4\in [\Sigma^{tq+3q+2}L, KG_{s+2}\wedge Y]$ and (9.2.24) becomes $(\bar b_{s+3}\wedge 1_V)f_2
= (\bar b_{s+3}\bar c_{s+2}\wedge 1_V)(1_{KG_{s+2}}\wedge w_4)g_4$. Consequently we have
$f_2 = (\bar c_{s+2}\wedge 1_V)(1_{KG_{s+2}}\wedge w_4)g_4 + (\bar a_{s+3}\wedge 1_V)f_3$
with $f_3\in [\Sigma^{tq+3q+3}L, E_{s+4}\wedge V]$ and so
$(\bar a_{s+3}\wedge 1_{M\wedge L})(f\wedge 1_L) = (1_{E_{s+3}}\wedge u_4)f_2 =
(\bar a_{s+3}\wedge 1_{M\wedge L})(1_{E_{s+4}}\wedge u_4)f_3$.
Hence, $(f\wedge 1_L) = (1_{E_{s+4}}\wedge u_4)f_3 + (\bar c_{s+3}\wedge 1_{M\wedge L})g_5$
for some $g_5\in [\Sigma^{tq+3q+3}L, KG_{s+3}\wedge M\wedge L]$ and so by (9.2.23) we have
$(1_{E_{s+4}}\wedge (1_Y\wedge (\alpha_1)_L)(\omega\wedge 1_L))(f\wedge 1_L) =
(\bar c_{s+3}\wedge 1_Y)(1_{KG_{s+3}}\wedge (1_Y\wedge (\alpha_1)_L)(\omega\wedge 1_L))g_5$ = 0
(this is because $(\alpha_1)_L$ induces zero homomorphsm is  $Z_p$-cohomology). Q.E.D.

\vspace{2mm}

{\bf Proof of the main Theorem A:}\quad We will continue the argument in Lemma 9.2.21.
Note that the spectrum $V$ in (9.2.23) also is the cofibre of $(1_M\wedge wi'')\widetilde{\psi}
: X\rightarrow\Sigma^{2q}M\wedge W$ given by the cofibration\\
{\bf (9.2.27)}\qquad $X\stackrel{(1_M\wedge wi'')\tilde{\psi}}
{\longrightarrow}\Sigma^{2q}M\wedge W\stackrel{w_5}{\longrightarrow}V\stackrel
{u_5}{\longrightarrow}\Sigma X$\\
this can be seen by the following homotopy commutative diagram of $3\times 3$-Lemma

$\qquad\Sigma^{3q-1}M\wedge L\quad\longrightarrow \qquad Y\qquad\stackrel
{\tilde{u}w_2}{\longrightarrow}\quad \Sigma X$

$\qquad\qquad\quad\searrow ^{1_M\wedge (\alpha_1)_L}\quad\nearrow \omega\quad
\searrow w_4\quad\nearrow u_5\quad\searrow \widetilde{\psi}$\newline
{\bf (9.2.28)}\quad\qquad\qquad$\Sigma^{2q}M\quad\qquad\qquad\qquad V\qquad\qquad\Sigma^{2q+1}M$

$\qquad\qquad\qquad\nearrow \widetilde{\psi}\quad\searrow ^{1_M\wedge wi''}\quad
\nearrow w_5\quad\searrow u_4\quad\nearrow _{1_M\wedge (\alpha_1)_L}$

$\quad\qquad\quad X\quad\longrightarrow\quad\Sigma^{2q}M\wedge W
\quad\stackrel{1_M\wedge u}{\longrightarrow}\quad\Sigma^{3q}M\wedge L$\\
By Lemma 9.2.21(B) and (9.2.23), $f\wedge 1_L = (1_{E_{s+4}}\wedge u_4)f_5$
for some $f_5\in [\Sigma^{tq+3q+3}L\\, E_{s+4}\wedge V]$  and so by Lemma 9.2.21(A)  we have \\
{\bf (9.2.29)}\quad $(\bar a_{s+2}\bar a_{s+3}\wedge 1_{M\wedge L})
(1_{E_{s+4}}\wedge u_4)f_5\quad  = \quad (\bar a_{s+2}\bar a_{s+3}\wedge 1_{M\wedge L})(f\wedge 1_L)$

$\qquad = \quad (1_{E_{s+2}}\wedge i\wedge 1_L)(\kappa\wedge 1_L) - (\bar c_{s+1}\wedge
1_{M\wedge L})(g\wedge 1_L)$.\\
It follows that  $(\bar a_s\bar a_{s+1}\bar a_{s+2}\bar
a_{s+3}\wedge 1_{M\wedge L}) (1_{E_{s+4}}\wedge u_4)f_5 = 0$ and
so $(\bar a_s\bar a_{s+1}\bar a_{s+2}\\\bar a_{s+3}\wedge 1_V)f_5
= (1_{E_s}\wedge w_4)f_6$ with $f_6\in [\Sigma^{tq+3q-1}L,
E_s\wedge Y]$. Clearly we have $(\bar b_s\wedge 1_V)(1_{E_s}\wedge
w_4)f_6$ = 0, then $(\bar b_s\wedge 1_Y)f_6$ = 0 and by
$Ext_A^{s+1+r, tq+3q+r}(H^*Y,\\H^*L)$ = 0( $r = 0,1$,  cf. Prop.
9.2.6) we have $(\bar a_s\bar a_{s+1}\bar a_{s+2}\bar
a_{s+3}\wedge 1_V)f_5 = (\bar a_s\bar a_{s+1}\\ \bar a_{s+2}\wedge
1_V)(1_{E_{s+3}}\wedge w_4)f_7$, with $f_7\in [\Sigma^{tq+3q+2}L,
E_{s+3}\wedge Y]$. Consequently we have\newline {\bf
(9.2.30)}\qquad $(\bar a_{s+1}\bar a_{s+2}\bar a_{s+3}\wedge
1_V)f_5$

\qquad\quad $ = (\bar a_{s+1}\bar a_{s+2}\wedge 1_V)(1_{E_{s+3}}\wedge w_4)f_7
+ (\bar c_s\wedge 1_V)g_6$\\
with $d_1$-cycle $g_6\in [\Sigma^{tq+3q}L, KG_s\wedge V]$ which represents an element
$[g_6]\in Ext_A^{s,tq+3q}(H^*V,H^*L)$. Note that the $d_1$-cycle $(\bar
b_{s+3}\wedge 1_Y)f_7\in [\Sigma^{tq+3q+2}L,\\ KG_{s+3}\wedge Y]$ represents an element
$[(\bar b_{s+3}\wedge 1_Y)f_7]\in Ext_A^{s+3,tq+3q+2}(H^*Y,H^*L)$ which has two generators
(cf. Prop. 9.2.5)(2)),  then $[(\bar b_{s+3}\wedge 1_Y)f_7] = \lambda'\bar h_*\widetilde{\phi}_*(\pi\wedge 1_L)_*
[\sigma'_1\wedge 1_L] + \lambda''\bar h_*\widetilde{\phi}_*(\pi\wedge 1_L)_*[\sigma'_2\wedge 1_L]$
for some $\lambda', \lambda''\in Z_p$. By the vanishes of the secondary differential :
$ 0 = d_2[(\bar b_{s+3}\wedge 1_Y)f_7]$ we know that $\lambda' , \lambda''$ is linearly dependent
.  Then we have\\
{\bf (9.2.31)}\qquad $[(\bar b_{s+3}\wedge 1_Y)f_7] = \lambda'\bar
h_*\widetilde{\phi}_*(\pi\wedge 1_L)_*[\sigma'\wedge 1_L]$

$\qquad\qquad\qquad\in Ext_A^{s+3,tq+3q+2}(H^*Y,H^*L)$\\
We claim that  the scalar $\lambda'$ in (9.2.31) is zero. This can be proved as follows.

The equation (9.2.30) means that the secondary differential of the ASS $d_2[g_6]
= 0\in E_2^{s+2,tq+3q+1}(L,V) = Ext_A^{s+2,tq+3q+1}(H^*V,H^*L)$, then $[g_6]\in
E_3^{s,tq+3q}(L,V)$  and

 The third differential\quad $d_3[g_6] = (w_4)_*[(\bar b_{s+3}\wedge 1_Y)f_7]
\in E_3^{s+3,tq+3q+2}(L,V)$\\
Note that $(\omega\wedge 1_L)(1_M\wedge (\alpha_1)_L)(i\wedge 1_L)\pi
= (1_Y\wedge j)\alpha_{Y\wedge M}i(\alpha_1)_L\pi = \bar h\phi(\alpha_1)_L\pi$ = 0
, this is because $\phi(\alpha_1)_L\in [\Sigma^{3q-2}L,L]$ = 0 which is obtained by $\pi_{rq-2}S$ = 0
($r = 2,3,4$). Hence , by (9.2.23), $(i\wedge 1_L)\pi = u_4\tau$  with
$\tau \in [\Sigma^{4q}S, V]$ which has  filtration 1.
Moreover, $u_4\tau\cdot p = (i\wedge 1_L)\pi\cdot p$ = 0 , then , by Prop. 9.2.3(4),
$\tau \cdot p = \widetilde{\lambda}w_4\bar h\widetilde{\phi}(\pi\wedge 1_L)
\pi$ for some $\widetilde{\lambda}\in Z_{(p)}$. This scalar $\widetilde{\lambda}$  must be zero
(mod $p$),this is because the left hand side of the equation has filtration 2 and the right hand side has
filtration 3 (cf. Remark 9.2.20  and $Ext_A^{0,4q+1}(H^*V,Z_p) = 0$ which is obtained by
$Ext_A^{0,4q+1}(H^*Y,Z_p) = 0 = Ext_A^{0,q+1}(H^*M\wedge L,Z_p)$).
Consequently, by Prop. 9.2.3(4) we have $\tau\cdot p$ = 0  and so $\tau = \overline{\tau}i$ with $\overline{\tau}\in
[\Sigma^{4q}M,V]$. Since $(u_4)_*(\pi)^*[g_6]\in Ext_A^{s+1,tq+q+1}(H^*M\wedge L, Z_p)\cong Z_p\{
(i\wedge 1_L)_*(\pi)_*(\sigma)\}$ ( cf. Prop. 9.2.8),  then $(u_4)_*\pi^*[g_6] = \lambda_0(i\wedge 1_L)_*\pi_*(\sigma)
= \lambda_0(u_4)_*(\overline{\tau} i)_*(\sigma)$ for some $\lambda_0\in Z_p$ and so by (9.2.23)
 we have $\pi^*[g_6] = \lambda_0\overline{\tau}_* i_*(\sigma)
\in Ext_A^{s+1,tq+3q+1}(H^*V, Z_p)$, this is because $Ext_A^{s+1,tq+3q+1}\\(H^*Y,H^*L)$ = 0 (cf. Prop.
9.2.6).  By the supposition on $d_2(\sigma) = a_0\sigma' = p_*(\sigma')\in Ext_A^{s+2,tq+1}(Z_p,Z_p)$,
we have $d_2i_*( \sigma)$ = 0 so that $i_*(\sigma)\in E_3^{s+2,tq+1}(S,\\M)$.
Moreover, $E_2^{s+3,tq+2}(S,M) = Ext_A^{s+3,tq+2}(H^*M,Z_p)$ = 0 (cf. Prop. 9.2.8)
then $E_3^{s+3,tq+2}(S,M)$ = 0  so that the third differential
 $d_3i_*(\sigma) \in E_3^{s+3,tq+2}(S,M)$ = 0. Since $\pi^*[g_6]
= \lambda_0(\overline{\tau})_*i_*(\sigma)\in
E_2^{s+1,tq+4q+1}(S,V)$, then $\pi^*[g_6] =
\lambda_0\overline{\tau}_*(i_*(\sigma))\in E_3^{s+1,tq+4q+1}
(S,V)$ and so

$ d_3\pi^*[g_6] = \lambda_0d_3(\overline{\tau})_*(i_*(\sigma)) =
\lambda_0(\overline{\tau})_*d_3(i_*(\sigma)) = 0\in
E_3^{s+4,tq+4q+3}(S,V)$\\
Hence, $(w_4)_*\pi^*[(\bar b_{s+3}\wedge 1_Y)f_7]
= d_3\pi^*[g_6]$ = 0 $\in E_3^{s+4,tq+4q+2}(S,V)$. In addition, by the split exact sequence
(9.2.26) we have $\pi^*[(\bar b_{s+3}\wedge 1_Y)f_7] = 0\in
E_3^{s+4,tq+4q+3}(S,Y)$.  Then, in the $E_2$-term, $\pi^*[(\bar b_5
\wedge 1_Y)f_7]$ must be a $d_2$-boundary, that is \\
\centerline{$\pi^*[(\bar b_{s+3}\wedge 1_Y)f_7] \in d_2E_2^{s+2,tq+4q+2}(S,Y) = d_2
Ext_A^{s+2,tq+4q+2}(H^*Y,Z_p)$ = 0}\\
(cf. Prop.  9.2.6(1)). Hence, by (9.2.31), $\lambda'\bar h_*\widetilde{\phi}_*
(\pi\wedge 1_L)_*\pi_*(\sigma')$ = 0. This implies that the scalar
 $\lambda'$ is zero  (cf. Prop. 9.2.9(3)) which shows the above claim.

So,(9.2.30) becomes $(\bar a_{s+1}\bar a_{s+2}\bar a_{s+3}\wedge
1_V)f_5 = (\bar a_{s+1}\bar a_{s+2}\bar a_{s+3}\wedge
1_V)(1_{E_{s+4}}\wedge w_4)f_8 + (\bar c_s\wedge 1_V)g_6$ for some
$f_8\in [\Sigma^{tq+3q+3}L, E_{s+4}\wedge Y]$. By composing
$1_{E_{s+1}}\wedge u_5$ on the above equation we have $(\bar
a_{s+1}\bar a_{s+2}\bar a_{s+3}\wedge 1_{X}) (1_{E_{s+4}}\wedge
u_5)f_5 = (\bar a_{s+1}\bar a_{s+2}\bar a_{s+3}\wedge
1_{X})(1_{E_{s+4}}\wedge\widetilde{u}w_2)f_8$ ( cf. (9.2.28)),
this is because $(1_{KG_s}\wedge u_5)g_6\in [\Sigma^{tq+3q-1}L,
KG_s\wedge X]$ represents an element in $Ext_A^{s,tq+3q-1}(H^*X,\\
H^*L)$ = 0 (cf. Lemma 9.2.18(2)) so that it is a $d_1$-boundary
and $(\bar c_s\wedge 1_{X})(1_{KG_s}\wedge u_5)g_6$ = 0.
Consequently we have \\
{\bf (9.2.32)}\qquad $(\bar a_{s+2}\bar a_{s+3}\wedge 1_{X})(1_{E_{s+4}}\wedge u_5)f_5$

\quad $ = (\bar a_{s+2}\bar a_{s+3}\wedge 1_{X})(1_{E_{s+4}}\wedge\widetilde{u}w_2)f_8 +
(\bar c_{s+1}\wedge 1_{X})g_7$ \\ for some  $d_1$-cycle $g_7\in [\Sigma^{tq+3q+1}L,
KG_{s+1}\wedge X]$ such that $[g_7]\in Ext_A^{s+1,tq+3q}\\(H^*X, H^*L)$.

Now we prove $(\bar c_{s+1}\wedge 1_X)g_7$ = 0 as follows.  By Lemma 9.2.18(2) and Prop. 9.2.7(1),
$[g_7] = \lambda_3(\widetilde{u}w_3)_*(\overline{\phi}_W)_*[\sigma\wedge 1_L]$
and the equation (9.2.32) means that the secondary differential $d_2[g_7] = 0$.
Since $d_2(\sigma) = a_0\sigma' = p_*(\sigma')\in Ext_A^{s+2,tq+1}(Z_p,Z_p)$
, then $\lambda_3(\widetilde{u}w_3)_*
(\overline{\phi}_W)_*(p\wedge 1_L)_*[\sigma'\wedge 1_L] = d_2[g_7] =
 0 \in Ext_A^{s+3,tq+3q+1}(H^*X,H^*L)$. By Lemma 9.2.18(1), this implies that
$\lambda_3 = 0 $ so that $g_7$ is a  $d_1$-boundary and
$(\bar c_{s+1}\wedge 1_X)g_7$ = 0.

Hence, (9.2.32) becomes  $(\bar a_{s+2}\bar a_{s+3}\wedge
1_{Y\wedge W}) (1_{E_{s+4}}\wedge u_5)f_5 = (\bar a_{s+2}\bar
a_{s+3}\wedge 1_{X})(1_{E_{s+4}}\wedge \widetilde{u}w_2)f_8$ and
so by (9.2.28)(9.2.12) we have, $(\bar a_{s+2}\bar a_{s+3}\wedge
1_{M})(1_{E_{s+4}}\\\wedge (1_M\wedge (\alpha_1)_L)u_4)f_5 = (\bar
a_{s+2}\bar a_{s+3}\wedge
1_M)(1_{E_{s+4}}\wedge\widetilde{\psi}u_5)f_5 = 0$. On the other
hand, by composing $(1_{E_{s+2}}\wedge 1_M\wedge (\alpha_1)_L)$ on
the equation (9.2.29) we have $(1_{E_{s+2}}\wedge i)\kappa\cdot
(\alpha_1)_L = (1_{E_{s+2}}\wedge 1_M\wedge
(\alpha_1)_L)(1_{E_{s+2}}\wedge i\wedge 1_L)(\kappa\wedge 1_L) =
(\bar a_{s+2}\bar a_{s+3}\wedge 1_{M})(1_{E_{s+4}}\wedge
(1_M\wedge (\alpha_1)_L)u_4)f_5 = 0$.

It follows that\\
{\bf (9.2.33)}\qquad\qquad $\kappa\cdot (\alpha_1)_L = (1_{E_{s+2}}\wedge p)f_9$\\
for some $f_9\in [\Sigma^{tq+q}L, E_{s+2}]$.
Since $\bar b_{s+2}\cdot\kappa = a_0\sigma' = p_*(\sigma')\in Ext_A^{s+2,tq+1}(Z_p,\\Z_p)$,
then $\kappa\cdot (\alpha_1)_L$ lifts to a map  $\tilde{f}\in [\Sigma^{tq+q+1}L,
E_{s+3}]$ such that $\bar b_{s+3}\cdot \tilde{f} $ represents $p_*((\alpha_1)_L)_*[\sigma'\wedge 1_L]\neq 0
\in Ext_A^{s+3,tq+q+1}(Z_p, H^*L)$ ( cf. Prop. 9.2.3(1)). Then , by (9.2.33),
$p_*[\bar b_{s+2}\cdot f_9] = p_*((\alpha_1)_L)_*[\sigma'\wedge 1_L]$ so that
$[\bar b_{s+2}\cdot f_9]\in Ext_A^{s+2,tq+q}(Z_p,H^*L)$  must equal to
$((\alpha_1)_L)_*[\sigma'\wedge 1_L]$
, this is because the group has two generators $((\alpha_1)_L)_*[\sigma'_1\wedge 1_L]$,
$((\alpha_1)_L)_*[\sigma'_2\wedge 1_L]$.  Write $\xi_{n,s+2} = f_9 i''$,
then \\
{\bf (9.2.34)} \qquad\quad $\kappa\cdot\alpha_1 = (1_{E_s+2}\wedge p)\xi_{n,s+2}$\\
such that $\bar b_{s+2}\cdot \xi_{n,s+2} = h_0\sigma'\in Ext_A^{s+2,tq+q}(Z_p,Z_p)$
and by Lemma 9.2.10 we have $(\bar c_{s+1}\wedge 1_M)(1_{KG_{s+1}}\wedge i)h_0\sigma =
(1_{E_{s+2}}\wedge i)\kappa\cdot\alpha_1$ = 0. This shows the second result of the
main Theorem.

In addition, by (9.2.34) and Lemma 9.2.10(2), $\bar a_0\bar a_1\cdots\bar a_{s+1}(1_{E_{s+2}}
\wedge p)\xi_{n,s+2}\\ = 0$, this shows that $\xi_n = \bar a_0\bar a_1\cdots\bar a_{s+1}
\cdot\xi_{n,s+2}\in \pi_{tq+q-s-2}S$  is an element of order $p$ and it is represented
by $h_0\sigma'\in Ext_A^{s+2,tq+q}(Z_p,Z_p)$ in the ASS. Q.E.D.

{\bf Remark 9.2.35}\quad  In the proof of the main Theorem A, we
obtain a stronger result. By (9.2.33), $\kappa\cdot (\alpha_1)_L =
(1_{E_{s+2}}\wedge p)f_9$, then $(1_{E_{s+2}}\wedge i)\kappa\cdot
(\alpha_1)_L$ = 0, and so $(1_{E_{s+2}}\wedge 1_L\wedge
i)(\kappa\wedge 1_L)\phi = (1_{E_{s+2}}\wedge 1_L\wedge i)
(\kappa\wedge 1_L)((\alpha_1)_L\wedge 1_L)\widetilde{i''}$ = 0,
where $\widetilde{i}''\in \pi_qL\wedge L$ such that
$((\alpha_1)_L\wedge 1_L)\widetilde{i}'' = \phi$. It can be easily
proved that $(\kappa\wedge 1_L)\phi = (\bar c_{s+1}\wedge
1_L)\sigma\phi$, where $\sigma\phi\in \pi_{tq+2q}(KG_{s+1}\wedge
L)$ is a $d_1$-cycle which represents $(\phi)_*(\sigma)\in
Ext_A^{s+1,tq+2q}(H^*L,Z_p)$.  Then we obtain that $(\bar
c_{s+1}\wedge 1_{L\wedge M})(1_{KG_{s+1}}\wedge i)\sigma\phi$ = 0.
That is to say, $(1_L\wedge i)_*(\phi)_*(\sigma)\in
Ext_A^{s+1,tq+2q}(H^*L\wedge M,Z_p)$ is a permanent cycle in the
ASS . Moreover, by (9.2.34) we have $\xi_{n,s+2} = f_9i''$, then
$(1_{KG_{s+2}}\wedge\alpha_1)\xi_{n,s+2} = (1_{KG_{s+2}}
\wedge\alpha_1)f_9i'' = f_9i''\cdot\alpha_1 = 0$ and so
$\xi_{n,s+2} = (1_{E_{s+2}}\wedge j'')\tilde{f}_9$ with
$\tilde{f}_9\in\pi_{tq+2q}E_{s+2}\wedge L$. Since $\xi_{n,s+2}$ is
represented by $h_0\sigma' = (j'')_*\phi_*(\sigma')$ in the ASS ,
then $\tilde{f}_9$ is represented by $(\phi)_*(\sigma')\in
Ext_A^{s+2,tq+2q}(H^*L,Z_p)$ in the ASS.
 That is to say $(\phi)_*(\sigma')\in Ext_A^{s+2,tq+2q}(H^*L,Z_p)$ is a permanent cycle in the ASS.
This is a stronger result obtained in the main Theorem A.

\quad

\begin{center}

{\bf \large \S 3.\quad A general result on convergence in the spectrum $V(1)$}

\vspace{2mm}

\end{center}

In this section we will prove , under some suppositions, a general result on
 the convergence of $i'_*i'*(h_0\sigma)\in Ext_A^{s+1,tq+q}
(H^*V(1),Z_p)$ to the homotopy groups of the spectrum $V(1)$ can implies the
convergence of $i'_*i_*(g_0\sigma)\in\\ Ext_A^{s+2,tq+pq+2q}(H^*V(1),Z_p)$ in the ASS.
We have the following main Theorem.

\vspace{2mm}

{\bf The main Thoerem B} (generalization of [7] Theorem II)\quad Let $p\geq 5,s\leq 4$,
$Ext_A^{s,tq}(Z_p,Z_p)\cong Z_p\{\sigma\}$,
$Ext_A^{s+1,tq+q}(Z_p,Z_p)\cong Z_p\{h_0\sigma\}$,\\
$Ext_A^{s+2,tq+2q+1}(Z_p,Z_p)\cong Z_p\{\widetilde{\alpha}_2\sigma\}$ and suppose that

(I) $Ext_A^{s+1, tq+rq+u}(Z_p,Z_p)$ = 0, for $r =1, u = -1,1,2,3$ or $r = 2, u = $

$\qquad\quad -1,0,1,2,3$.

    $Ext_A^{s+1,tq}(Z_p,Z_p)$ is zero or has (one or two) generator $\sigma'$ satisfying
    $a_0\sigma'\neq 0$.

    $Ext_A^{s+1,tq+r}(Z_p,Z_p)$ = 0 for $r = -2,-1,2,3$ and has unique generator $a_0\sigma$
    for $r = 1$ satisfying $a_0^2\sigma\neq 0$.

    $Ext_A^{s,tq+q}(Z_p,Z_p)$ = 0 or $= Z_p\{h_0\tau'\}, Ext_A^{s,tq+1}(Z_p,Z_p)$ = 0 or $  Z_p\{a_0\tau'\}$.

    $Ext_A^{s,tq+rq+u}(Z_p,Z_p)$ = 0, $r = 1,u = 1,2$, $r = -1, u = -1,0$.

(II) $i'_*i_*(h_0\sigma)\in Ext_A^{s+1,tq+q}(H^*K,Z_p)$ is a permanent cycle in the ASS,\\
then $i'_*i_*(g_0\sigma)\in Ext_A^{s+2,tq+pq+2q}(H^*K,Z_p)$ also is a permanent cycle in the ASS
and it conveges to a nontrivial element in $\pi_{tq+pq+2q-s-2}K$.

\vspace{2mm}

To prove the main Theorem B, we need some knowledge on derivation of maps between $M$-
module spectra and some lower dimensional Ext groups. These preminilaries will be
used in the proof of the main Theorem B and especially in the proof of Theorem 9.3.9 below.

\vspace{2mm}

{\bf Prop. 9.3.0}\quad Let $p\geq 5, s\leq 4$, then under the supposition of the main Theorem B we have

\quad (1) $Ext_A^{s+1,tq+r}(H^*M,H^*M)$ = 0 for $ r = 1,2$.

(2) $Ext_A^{s+1,tq+r}(Z_p,H^*M)$ = 0  for $r = 0,1$,

\qquad $Ext_A^{s,tq+r}(H^*M,Z_p)$ = 0 for $r = 1,2$,

{\bf Proof}\quad (1) Consider the following exact sequence ( $r =
1,2,3$) \\$Ext_A^{s+1,tq+r}(Z_p,Z_p)\stackrel{i_*}{\rightarrow}
Ext_A^{s+1,tq+r}(H^*M,Z_p)
\stackrel{j_*}{\rightarrow} Ext_A^{s+1,tq+r-1}(Z_p,Z_p)\stackrel{p_*}{\rightarrow}$\\
induced by (9.1.1).  By the supposition, the left group is zero for $r = 2,3$ and has unique generator $a_0\sigma = p_*(\sigma)$ for $r = 1$ so that
im $i_*$ = 0.  By the supposition, the right group is zero for $r = 3$ and has unique generator $a_0\sigma$ for $r = 2$ which satisfies
$p_*(a_0\sigma) = a_0^2\sigma\neq 0$.  By the supposition, the right group is zero for $r = 1$ and has (one or two)
generator $\sigma'$ for $r = 2$ (both) satisfying $p_*(\sigma') = a_0\sigma'\neq 0$. Then , the above $p_*$ is monic so that
im $j_*$ = 0. This shows that the middle group is zero which shows the first result.
The second result can be obtained immediately by the first result.

(2) Consider the following exact sequence  $(r = 0,1)$

$Ext_A^{s+1,tq+r+1}(Z_p,Z_p)\stackrel{j^*}{\longrightarrow} Ext_A^{s+1,tq+r}(Z_p,H^*M)$

$\qquad\quad\stackrel{i^*}{\longrightarrow} Ext_A^{s+1,tq+r}(Z_p,Z_p)\stackrel{p^*}{\longrightarrow}$\\
induced by (9.1.1). By the supposition, the left group is zero for $r = 1$ and has unique generator $a_0\sigma = p^*(\sigma)$ for $r = 0$ so that
im $j^*$ = 0. The right group has unique generator $a_0\sigma$ for $r = 1$ which satisfies $p^*(a_0\sigma) = a_0^2\sigma\neq 0$.
The right group is zero for $r = 0$ or has (one or two) generator $\sigma'$ satisfying
$p^*(\sigma') = a_0\sigma'\neq 0$. Then
im $i^*$ = 0 so that the middle group is zero as desired. The proof of the second result is similar. Q.E.D.

\vspace{2mm}

{\bf Proposition 9.3.1}\quad Let $p\geq 5, s\leq 4$, then under the supposition of the main Theorem B we have

(1) $Ext_{A}^{s,tq}(H^{*}M,H^{*}M)\cong Z_p\{\tilde{\sigma}\} $
satisfying $i^*(\tilde{\sigma}) = i_*(\sigma)\in
Ext_A^{s,tq}\\(H^*M,Z_p)$, $j_*(\tilde{\sigma}) = j^*(\sigma)\in
Ext_A^{s,tq-1}(Z_p,H^*M)$.

(2) $Ext_{A}^{s+1,tq+q}(H^{*}M,H^{*}M) \cong Z_{p} \{(ij)_{*}\alpha_{*}(\tilde{\sigma}) , $
$\alpha_{*}(ij)^{*}(\tilde{\sigma})\}$ ,

(3) $Ext_A^{s+1,tq+q+1}(H^*M,H^*M)\cong Z_p\{\alpha_*(\tilde{\sigma}) = \alpha^*(\tilde{\sigma})\}$,

(4) $Ext_{A}^{s+1,tq+q}(H^{*}K,H^{*}M)\cong Z_{p}
\{i'_{*}(ij)_{*}\alpha_{*}(\tilde{\sigma}) = i'_{*}(\alpha_{1}\wedge 1_{M})_{*}
(\tilde{\sigma})\}$,\\ where $\alpha_{1} = j\alpha i : \Sigma^{q-1} S\rightarrow
S $ and $\alpha_* : Ext_A^{s,tq}(H^*M,H^*M)\to Ext_A^{s+1,tq+q+1}\\(H^*M,H^*M)$ is the connecting (or boundary)
homomorphism induced by $\alpha : \Sigma^qM\to M$.

{\bf Proof}: \quad (1) Consider the following exact sequence\\
$\quad 0 = Ext_A^{s,tq+1}(H^*M,Z_p)\stackrel{j^*}{\rightarrow} Ext_A^{s,tq}(H^*M,H^*M)
\stackrel{i^*}{\rightarrow} Ext_A^{s,tq}(H^*M,Z_p)\stackrel{p^*}{\rightarrow}$\\
induced by (9.1.1). The right group has unique generator  $i_*(\sigma)$ , this is because
 $Ext_A^{s,tq-r}(Z_p,Z_p)$ = 0 for $r = 1$) and has unique generator $\sigma$ for $r = 0$.
Moreover, $p^*i_*(\sigma) = i_*p^*(\sigma) = i_*(a_0\sigma) = i_* p_*(\sigma)$
= 0, then the middle has unique generator $\tilde{\sigma}$ such that $i^*(\tilde{\sigma}) = i_*(\sigma)$.
This shows the result and the second relation can be similarly proved.

(2) By the supposition, $Ext_{A}^{s+1,tq+q}(Z_{p},Z_{p})$ has unique generator
$h_{0}\sigma = j_{*}\alpha_{*}i_{*}(\sigma)$ = $j_{*}\alpha_{*}i^{*}(\tilde{\sigma})$,
Then the result follows by the following exact sequence

$ \stackrel{p^{*}}{\rightarrow}Ext_{A}^{s+1,tq+q+1}(H^{*}M,Z_{p})
\stackrel{j^{*}}{\rightarrow}Ext_{A}^{s+1,tq+q}(H^{*}M,H^{*}M)$

$\qquad\quad\stackrel{i^{*}}{\rightarrow}Ext_{A}^{s+1,tq+q}(H^{*}M,Z_{p})\stackrel{p^{*}}
{\longrightarrow} $\\
induced by (9.1.1), where the right group has unique generator
 $i^{*}(ij)_{*}\alpha_{*}(\tilde{\sigma}) = (ij)_{*}\alpha_{*}i_{*}(\sigma)$ satisftying $p^{*}(ij)_{*}
\alpha_{*}i_{*}(\sigma) = (ij)_{*}\alpha_{*}i_{*}p_{*}(\sigma)$ = 0 and
the left group has unique generator $\alpha_{*}i_{*}(\sigma) = i^{*}\alpha_{*}(\tilde{\sigma})$.

(3) Consider the following exact sequence

$ Ext_A^{s+1,tq+q+2}(H^*M,Z_p)\stackrel{j^*}{\longrightarrow} Ext_A^{s+1,tq+q+1}
(H^*M,H^*M)$

$\qquad\quad \stackrel{i^*}{\longrightarrow} Ext_A^{s+1,tq+q+1}(H^*M,Z_p)\stackrel
{p^*}{\longrightarrow}$\\
induced by (9.1.1). The left group is zero, this is because by the supposition, $Ext_A^{s+1,tq+q+r}(Z_p,Z_p)$ = 0 for $r = 1,2,3$. The right group
 has unique generator $(\alpha i)_*(\sigma) = i^*\alpha_*(\tilde{\sigma})$, this is because $Ext_A^{s+1,tq+q+r}
 (Z_p,Z_p)$ is zero for $r = 1$ and has unique generator $h_0\sigma = j_*(\alpha i)_*(\sigma)$ for $r = 0$.
Since $p^*(\alpha i)_*(\sigma) = (\alpha i)_*p_*(\sigma)$ = 0, then the middle
group  has unique generator $\alpha_*(\tilde{\sigma})$ as desired. Moreover we have $\alpha_*(\tilde{\sigma})
= \alpha^*(\tilde{\sigma})$, this is because $i^*j_*\alpha_*(\tilde{\sigma}) = j_*\alpha_*i_*(\sigma)
= h_0\sigma = (j\alpha i)^*(\sigma) = i^*j_*\alpha^*(\tilde{\sigma})$.

(4) Consider the following exact sequence

$  Ext_{A}^{s+1,tq+q}(H^{*}M,H^{*}M)
\stackrel{i'_{*}}{\rightarrow}Ext_{A}^{s+1,tq+q}(H^{*}K,H^{*}M)$

$\qquad\quad\stackrel{j'_{*}}
{\rightarrow}Ext_{A}^{s+1,tq-1}(H^{*}M,H^{*}M)\stackrel{\alpha_{*}}{\longrightarrow} $\\
induced by (9.1.2). By the supposition,
$Ext_{A}^{s+1,tq-r}(Z_{p},Z_{p})$ = 0 for $r = 1,2$ and has unique
generator  $\sigma'$ for $r = 0$ , then the right group has unique
generator $(ij)^{*}(\tilde{\sigma'})$ satisfying
$\alpha_{*}(ij)^{*}(\tilde{\sigma'}) =
j^{*}\alpha_{*}i_{*}(\sigma') \neq 0\in
Ext_{A}^{s+2,tq+q}(H^{*}M,\\H^{*}M)$ . Hence,
$Ext_{A}^{s+1,tq+q}(H^{*}K,H^{*}M) = i'_{*}Ext_{A}
^{s+1,tq+q}\\(H^{*}M,H^{*}M)$ has unique generator $
(i')_{*}(ij)_{*}\alpha_{*}\tilde{\sigma}$ =
$i'_{*}(\alpha_{1}\wedge 1_{M})_{*}(\tilde{\sigma})$, this is
because $(\alpha_{1}\wedge 1_{M})_{*}(\tilde{\sigma})$ =
$(ij)_{*}\alpha_{*}(\tilde{\sigma}) - \alpha_{*}(ij)_{*}
(\tilde{\sigma})$ which is obtained by $\alpha_{1}\wedge 1_{M} =
ij\alpha - \alpha ij$.  Q.E.D.

\vspace{2mm}

{\bf Proposition 9.3.2}\quad Let $p\geq 5, s\leq 4$ , then under the supposition of the main Theorem B we have

(1) $Ext_{A}^{s+1,tq+2q+r}(H^{*}K,H^{*}M) $ = 0, $r = 0,1,2$,

\qquad $Ext_{A}^{s+1,tq+2q+1}(H^{*}K,Z_{p})$ = 0.

(2) $Ext_{A}^{s+1,tq+q+r}(H^{*}K,Z_{p})$ = 0 , $r = 1,2,3$,

\qquad $Ext_{A}^{s+1,tq+q+r}(H^{*}K,H^{*}M)$ = 0 ,$ r = 1,2$.

(3) $Ext_{A}^{s+1,tq+q}(H^{*}K,H^{*}K)$  $\cong Z_{p}\{(h_{0}\sigma)'\}$
with $(i')^{*}(h_{0}\sigma)' = (i'ij\alpha )_{*}(\tilde{\sigma})$.

{\bf Proof}: (1) Consider the following exact sequence

$\qquad Ext_{A}^{s+1,tq+2q+r}(H^{*}M,H^{*}M)\stackrel{i'_{*}}{\rightarrow}
Ext_{A}^{s+1,tq+2q+r}(H^{*}K,H^{*}M)$

$\qquad\qquad \stackrel{j'_{*}}{\rightarrow}Ext_{A}
^{s+1,tq+q+r-1}(H^{*}M,H^{*}M)\stackrel{\alpha_{*}}{\rightarrow}$\\
induced by (9.1.2). The left group is zero by the supposition on
$Ext_{A}^{s+1,tq+2q+u}\\(Z_{p}, Z_{p})$ = 0 for $u = -1, 0, 1, 2$
. The right group has unique generator
$(ij)^{*}(ij)_{*}\alpha_{*}(\tilde{\sigma})$ for $r = 0$ and is
generated by two generators $(ij)_{*}\alpha_{*}(\tilde{\sigma})$ ,
$(ij)^{*}\alpha_{*}(\tilde{\sigma})$ for $r = 1$. Moreovre, the
right group has unique generator $\alpha_{*}(\tilde{\sigma})$ for
$r = 2$ (cf. Prop. 9.3.1(3)). We claim that (i)
$\alpha_{*}(ij)^{*}(ij)_{*}\alpha_{*}(\tilde{\sigma}) \neq 0$ .
(ii) $\alpha_{*}[\lambda_{1}(ij)_{*}\alpha_{*}(\tilde{\sigma}) +
\lambda_{2} \alpha_{*}(ij)^{*}(\tilde{\sigma})] \neq 0$. (iii)
$\alpha_{*}\alpha_{*} (\tilde{\sigma}) \neq 0$.  Then the above
$\alpha_{*}$ is monic and so $im j'_{*}$ = 0.  This shows
$Ext_{A}^{s+1,tq+2q+r}(H^{*}K,H^{*}M)$ = 0  with $r = 0, 1, 2$ and
consequently we have $Ext_{A}^{s+1,tq+2q+1}$ $(H^{*}K,Z_{p})$ = 0.

To prove the claim,  recall from the supposition that
$\widetilde{\alpha}_{2}\sigma =
j_{*}\alpha_{*}\alpha_{*}i_{*}(\sigma)\\ \neq 0 \in
Ext_{A}^{s+2,tq+2q+1}(Z_{p},Z_{p})$, then
$i_{*}(\widetilde{\alpha}_{2}\sigma) \neq 0 \in
Ext_{A}^{s+2,tq+2q+1}(H^{*}M,Z_{p})$ , this is because
$Ext_{A}^{s+1,tq+2q}(Z_{p},Z_{p})$ = 0 from the supposition. In
addition, we also have $j^{*}i_{*}(\alpha_{2}\sigma)\neq 0 \in
Ext_{A}^{s+2,tq+2q}(H^{*}M, H^{*}M)$ ,this is because
$Ext_{A}^{s+1,tq+2q}$ $(H^{*}M,Z_{p})$ = 0. Hence, by $2 \alpha
ij\alpha =
ij \alpha^{2} + \alpha^{2}ij$ we have\\
{\bf (9.3.3)} $\qquad \alpha_{*}(ij)^{*}(ij)_{*}\alpha_{*}(\tilde{\sigma}) =
j^{*}\alpha_{*}(ij)_{*}\alpha_{*}i_{*}(\sigma)$

$\qquad\quad =  \frac{1}{2} j^{*}(ij)_{*}\alpha_{*}\alpha_{*}i_{*}(\sigma) =
\frac{1}{2}j^{*}i_{*}(\alpha_{2}\sigma) \neq 0$\\
This shows the claim (i).  For the claim  (ii),

$\alpha_{*}[\lambda_{1}(ij)_{*}\alpha_{*}(\tilde{\sigma}) + \lambda_{2}\alpha_{*}
(ij)^{*}(\tilde{\sigma})]$

$\quad = \frac{1}{2}\lambda_{1}(ij)_{*}\alpha_{*}\alpha_{*}(\tilde{\sigma})
+ (\frac{1}{2}\lambda_{1}+\lambda_{2})\alpha_{*}\alpha_{*}(ij)^{*}(\tilde{\sigma})\neq 0$\newline
, this is because this two terms is linearly independent which can be obtained from
$(ij)_{*}\alpha_{*}\alpha_{*}(ij)^{*}(\tilde{\sigma})\neq 0$ ( cf. (9.3.3)).  The claim (iii) is immediate
, this is because $i^{*}j_{*}\alpha_{*}\alpha_{*}(\tilde{\sigma}) = j_{*}\alpha_{*}\alpha_{*}
i_{*}(\sigma) = \widetilde{\alpha}_{2}\sigma \neq 0$.

(2) Consider the following exact sequence ($r = 1,2,3)$

$ Ext_{A}^{s+1,tq+q+r}(H^{*}M,Z_{p})\stackrel{i'_{*}}{\longrightarrow}
Ext_{A}^{s+1,tq+q+r}(H^{*}K,Z_{p})$

$\qquad\quad\stackrel{j'_{*}}{\longrightarrow}
Ext_{A}^{s+1,tq+r-1}(H^{*}M,Z_{p})\stackrel{\alpha_{*}}{\longrightarrow}$\\
induced by (9.1.2). The left group is zero for  $r = 2,3$ which can be obtained from
the supposition of $Ext_{A}^{s+1,tq+q+u}(Z_{p},Z_{p})$ = 0( $u = 1,2,3$).
The left group has unique generator $\alpha_{*}i_{*}(\sigma)$ for $r = 1$,
then in any case we have  im  $i'_{*}$ = 0.  The right group is zero for $r = 2,3$
(cf. Prop. 9.3.0) and has unique generator $i_{*}(\sigma')$ for $r = 1$ which satisfies
$\alpha_{*}i_{*}(\sigma') \neq 0 \in Ext_{A}^{s+2,tq+q+1}(H^{*}M,Z_{p})$
so that $j'_{*} $ = 0 and the result follows.

(3) Consider the following exact sequence

$Ext_{A}^{s+1,tq+2q+1}(H^{*}K,H^{*}M)\stackrel{(j')^{*}}{\longrightarrow}
Ext_{A}^{s+1,tq+q}(H^{*}K,H^{*}K)$

$\qquad\quad\stackrel{(i')^{*}}{\longrightarrow}Ext_{A}
^{s+1,tq+q}(H^{*}K,H^{*}M)$\\
induced by (9.1.2). The left group is zero by (1) and the right group has unique generator
 $i'_{*}(ij)_{*}\alpha_{*}(\tilde{\sigma})$ (cf. Prop. 9.3.1(4)) which satisfies
$\alpha^{*}i'_{*}(ij)_{*}\\\alpha_{*}(\tilde{\sigma}) =
i'_{*}(ij)_{*} \alpha_{*}\alpha^{*}(\tilde{\sigma}) =
i'_{*}(ij)_{*}\alpha_{*}\alpha_{*}(\tilde{\sigma})$ = 0, this is
because  $i'ij\alpha^{2}\\ = 2 i'\alpha ij\alpha - i'\alpha^{2} ij
= 0 \in [\Sigma^{2q-1}M,K]$. Then the result follows. Q.E.D.

\vspace{2mm}

{\bf Proposition 9.3.4}\quad Let $p\geq 5, s\leq 4$, then under the supposition of the main Theorem B we have

$\qquad\quad Ext_{A}^{s+1,tq+q-1}(H^{*}K,H^{*}K)$ $\cong Z_{p}\{(h_{0}\sigma)''\}$
\\ satisfying $(i')^{*}(h_{0}\sigma)'' = i'_{*}(ij)_{*}(\alpha_{1}\wedge 1_{M})_{*}(\tilde{\sigma})$.

{\bf Proof}:  Consider the following exact sequence

$Ext_{A}^{s+1,tq+2q}(H^{*}K,H^{*}M)\stackrel{(j')^{*}}{\longrightarrow}
Ext_{A}^{s+1,tq+q-1}(H^{*}K,H^{*}K)$

$\qquad\quad\stackrel{(i')^{*}}{\longrightarrow}Ext_{A}^{s+1,tq+q-1}(H^{*}K,
H^{*}M) $\\
induced by (9.1.2). The left group is zero by Prop. 9.3.2(1) and similar to that in
Prop. 9.3.1, the right group has unique generator $(ij)^{*}i'_{*}
 (ij)_{*}\alpha_{*}(\tilde{\sigma})
= i'_{*}(ij)_{*}(\alpha_{1}\wedge 1_{M})_{*}(\tilde{\sigma})$ which satisfies
$\alpha^{*} i'_{*}(ij)_{*}(\alpha_{1}\wedge 1_{M})_{*}(\tilde{\sigma})$ =
$i'_{*}(ij)_{*}(\alpha_{1}\wedge 1_{M})_{*}\alpha_{*}(\tilde{\sigma})$ = 0
$ \in Ext_{A}^{s+2,tq+2q}(H^{*}K,$ $H^{*}M)$, this is because $i'ij(\alpha_{1}\wedge 1_{M})
\alpha $ = 0 $\in [\Sigma^{2q-2}M,K]$. Then the result follows. Q.E.D.

\vspace{2mm}

{\bf Proposition 9.3.5} \quad Let $p\geq 5, s\leq 4,$ then under the supposition of the main Theorem B we have

 $Ext_{A}^{s+1,tq+q+1}(H^{*}K'\wedge M , H^{*}M)\cong Z_{p}\{\psi_{*}(ij)_{*}\alpha_{*}
(\tilde{\sigma}) , \psi_{*}(ij)^{*}\alpha_{*}(\tilde{\sigma})\}.$
\\where $\psi : \Sigma M\rightarrow K'\wedge M$ is the map in (9.1.17).

{\bf Proof}: Consider the following exact sequence

$ Ext_{A}^{s+1,tq+q}(H^{*}M,H^{*}M)\stackrel{\psi_{*}}{\longrightarrow}
Ext_{A}^{s+1,tq+q+1}(H^{*}K'\wedge M , H^{*}M)$

$\qquad\quad \stackrel{\rho_{*}}{\longrightarrow}
Ext_{A}^{s+1,tq+q+1}(H^{*}K,H^{*}M) = 0$\\
induced by (9.1.17). The result follows immediately form Prop. 9.3.1(2) and Prop.9.3.2.(Note: By the supposition,
similar to that given in Prop. 9.3.2(2),  we can prove that $Ext_A^{s,tq+q+1}(H^*K,H^*M)$ = 0
so that the above $\psi_*$ is monic). Q.E.D.

\vspace{2mm}

{\bf Proposition 9.3.6}\quad  Let $p\geq 5, s\leq 4$, then under the supposition of the main Theorem B we have

 $Ext_{A}^{s,tq}(H^{*}K,H^{*}K)$ $\cong Z_{p}\{(\sigma)'\}$
 satisfying $(i')^{*}(\sigma)' = (i')_{*}(\tilde{\sigma})$.

{\bf Proof}:  Consider the following exact sequence

$Ext_{A}^{s,tq+q+1}(H^{*}K,H^{*}M)\stackrel{(j')^{*}}{\longrightarrow}
Ext_{A}^{s,tq}(H^{*}K,H^{*}K)$

$\qquad\quad\stackrel{(i')^{*}}{\longrightarrow} Ext_{A}^{s,tq}(H^{*}K,H^{*}M)$\\
induced by (9.1.2). Since
$j'_{*}Ext_{A}^{s,tq+q+1}(H^{*}K,H^{*}M)$ $\subset
Ext_{A}^{s,tq}(H^{*}M,H^{*}M)\\\cong Z_{p}\{\tilde{\sigma}\}$ and
$\alpha_{*}(\tilde{\sigma}) \neq 0 \in
Ext_{A}^{s+1,tq+q+1}(H^{*}M,H^{*}M)$, then im $(j')_*$ = 0 and so
$Ext_{A}^{s,tq+q+1}(H^{*}K,H^{*}M) = i'_{*}
Ext_{A}^{s,tq+q+1}(H^{*}M,H^{*}M) $.
 Moreover, by the supposition on $Ext_{A}^{s,tq+q+r}(Z_{p},Z_{p})$ = 0 , $ r = 1,2$,
and $Ext_A^{s,tq+q}(Z_p,Z_p)$ is zero or $\cong Z_p\{h_0\sigma''\}$ we have $Ext_A^{s,tq+q+1}(H^*M,H^*M)
\cong Z_p\{\alpha_*(\tilde{\sigma''})\}$, then $Ext_A^{s,tq+q+1}(H^*K,H^*M) =
(i')_*Ext_A^{s,tq+q+1}(H^*M,H^*M)$ = 0.
On the other hand, it is easily seen  that $Ext_{A}^{s,tq}(H^{*}K,H^{*}M)$ has unique generator
 $i'_{*}(\tilde{\sigma})$ which satisfies $\alpha^{*}i'_{*}(\tilde{\sigma})
= i'_{*}\alpha_{*}(\tilde{\sigma})$ = 0.  Then the result follows. Q.E.D.

\vspace{2mm}

By (6.5.5), there is $\alpha'' \in [\Sigma^{q-2}K,K]$ such that $\alpha'' i'
= i'ij\alpha ij$. Let $X$  be the cofibre of  $\alpha'' : \Sigma^{q-2}K \rightarrow
K$ given by the cofibration\\
{\bf (9.3.7)} $\qquad\quad \Sigma^{q-2}K\stackrel{\alpha''}{\longrightarrow} K\stackrel{w}
{\longrightarrow}X \stackrel{u}{\longrightarrow} \Sigma^{q-1}K$,\\
Then, $\alpha''$ induces a boundary homomorphism (or connecting homomotphism) $(\alpha'')^{*}$ :
$Ext_{A}^{s,tq}(H^{*}K,H^{*}K)) \rightarrow Ext_{A}^{s+1,tq+q-1}$
$(H^{*}K,H^{*}K)$.  Since $\alpha'' i' = i'
ij\alpha ij = i' ij(\alpha_{1}\wedge 1_{M})$ , then
$(i')^{*}(\alpha'')^{*}(\sigma)' = (\alpha'' i')^{*}(\sigma)' =
(i'ij(\alpha_{1}\wedge 1_{M}))^{*}(\sigma)' = (\alpha_{1}\wedge
1_{M})^{*}(ij)^{*}(i')^{*}(\sigma)'= (i'ij)_{*}(\alpha_{1}\wedge
1_{M})_{*}(\tilde{\sigma}) = (i')^{*}(h_{0}\sigma)''$ (cf. Prop. 9.3.4)
. Then we have\\
{\bf (9.3.8)} $ \quad (h_{0}\sigma)''  = (\alpha'')^{*}(\sigma)' \in
Ext_{A}^{s+1,tq+q-1}(H^{*}K,H^{*}K)$ \\
 this is because the above  $(i')^{*}$ is monic which can be obtained by
 $Ext_{A}^{s+1,tq+2q}\\(H^{*}K,H^{*}M)$ = 0 (cf. Prop. 9.3.2).

After finishing the above preminilaries , we now turn to prove the following Theorem 9.3.9. It is proved by some argument
processing in the Adams resolution (9.2.9) of some spectra related to the sphere spectrum $S$.

\vspace{2mm}

{\bf  Theorem 9.3.9}\quad Let $p\geq 5, s\leq 4$ , then under the supposition of the main Theorem B we have
$(\bar c_{s+1}\wedge 1_K)(h_0\sigma)''$ = 0, where $(h_{0}\sigma)''\in [\Sigma^{tq+q-1}K,\\
KG_{s+1}\wedge K]$ is a  $d_{1}$-cycle which represents the unique generator $(h_0\sigma)''$ of
$Ext_A^{s+1,tq+q-1}(H^*K,H^*K)$ (cf. Prop. 9.3.4).

Before proving Theorem 9.3.9, we first prove the following Lemma.

\vspace{2mm}

{\bf Lemma 9.3.10}\quad Let $p\geq 5, s\leq 4$, then under the supposition of the main Theorem B we have

(1) $(\bar c_{s+1}\wedge 1_{K})(h_{0}\sigma)'' $ = $(1_{E_{s+2}}
\wedge \alpha'')(\kappa\wedge 1_{K})$,

(2) $(\bar c_{s+1}\wedge 1_K)(h_0\sigma\wedge 1_K) = (1_{E_{s+2}}\wedge\alpha_1\wedge 1_K)(\kappa\wedge 1_K)$\\
where $\kappa\in\pi_{tq+1}E_{s+2}$ such that $\bar a_{s+1}\kappa =
\bar c_{s} \sigma$  with $\sigma\in \pi_{tq}KG_{s}\\ \cong
Ext_{A}^{s,tq}(Z_{p},Z_{p})$.

{\bf Proof}:  Recall that $X$ is the cofibre of $\alpha'' :
\Sigma^{q-2}K \rightarrow K$ given by the cofibration (9.3.7).
Since $(h_{0}\sigma)'' \in [\Sigma^{tq+q-1}K,KG_{s+1}\wedge K]$
represents $(h_{0}\sigma)'' = (\alpha'')^{*}(\sigma)' \in
Ext_{A}^{s+1,tq+q-1} (H^{*}K,$ $H^{*}K)$,  then $(h_{0}\sigma)''u
\in [\Sigma^{tq}X, \\KG_{s+1}\wedge K]$ is a $d_{1}$-boundary so
that $(\bar c_{s+1}\wedge 1_{K})(h_{0}\sigma)''u$
 = 0 and $(\bar c_{s+1}\wedge 1_{K})(h_{0}\sigma)'' = f' \alpha''$ for some $f' \in
\Sigma^{tq+1}K,E_{s+2}\wedge K].$  It follows that $(\bar a_{s+1}\wedge 1_{K})
f' \alpha'' = 0$  and so $(\bar a_{s+1}\wedge 1_{K})f' = f'_{2}w$ with $f'_{2} \in
[\Sigma^{tq}X, E_{s+1}\wedge K]$.  Then, $(\bar b_{s+1}\wedge 1_{K})f'_{2}w$ = 0
and $(\bar b_{s+1}\wedge 1_{K})f'_{2} = g'\cdot u$ for some $g' \in [\Sigma^{tq+q-1}K,
KG_{s+1}\wedge K]$.  $g'$ is a $d_{1}$-cycle, this is because $(\bar b_{s+2}\bar c_{s+1}
\wedge 1_{K})g' = g'_{2}\alpha''$ (with $g'_{2} \in [\Sigma^{tq+1}K,KG_{s+2}\wedge K])$
= 0 since $\alpha''$ induces zero homomorphism in $Z_{p}$-cohomology.
Then,  by Prop. 9.3.4 and (9.3.8), $g' $ represents $(h_{0}\sigma)'' = (\alpha'')^{*}(\sigma)'
\in Ext_{A}^{s+1,tq+q-1}(H^{*}K,H^{*}K)$ and so $g'\cdot u$ is a  $d_{1}$-boundary
, that is $g'\cdot u = (\bar b_{s+1}\bar c_{s}\wedge 1_{K})g'_{3}$ with $g'_{3} \in
[\Sigma^{tq}X,KG_{s}\wedge K]$. Then $(\bar b_{s+1}\wedge 1_{K})f'_{2}
= (\bar b_{s+1}\bar c_{s}\wedge 1_{K})g'_{3}$ and so $f'_{2} = (\bar c_{s}\wedge 1_{K})g'_{3}
+ (\bar a_{s+1}\wedge 1_{K})f'_{3}$ for some $f'_{3} \in [\Sigma^{tq+1}X,E_{s+2}\wedge K]$
and we have $(\bar a_{s+1}\wedge 1_{K})f' = f'_{2} w = (\bar c_{s}\wedge 1_{K})g'_{3}w
+ (\bar a_{s+1}\wedge 1_{K})f'_{3}w$. Clearly, $g'_{3}w \in [\Sigma^{tq}K,KG_{s}\wedge K]$
is a  $d_{1}$-cycle which represents $Ext_{A}^{s,tq}(H^{*}K,H^{*}K)
\cong Z_{p}\{(\sigma)'\}$ (cf. Prop. 9.3.6).  Then $g'_{3}w = \sigma\wedge 1_{K}$
(up to scalar and modulo $d_1$-boundary), where $\sigma \in \pi_{tq}KG_{s} \cong Ext_{A}^{s,tq}(Z_{p},Z_{p})$.
So we have $(\bar a_{s+1}\wedge 1_{K})f' = (\bar c_{s}\wedge 1_{K})(\sigma\wedge 1_{K}) +
(\bar a_{s+1}\wedge 1_{K})f'_{3}w$ = $(\bar a_{s+1}\wedge 1_{K})(\kappa\wedge 1_{K}) +
(\bar a_{s+1}\wedge 1_{K})f'_{3}w$ , where $\kappa \in \pi_{tq+1}E_{s+2}$
satisfying $\bar a_{s+1}\kappa = \bar c_{s}\sigma$.  It follows that $ f' = \kappa
\wedge 1_{K} +  f'_{3}w + (\bar c_{s+1}\wedge 1_{K})g'_{4}$ for some $g'_{4} \in
[\Sigma^{tq+1}K, KG_{s+1}\wedge K]$  and we have $(\bar c_{s+1}\wedge 1_{K})
(h_{0}\sigma)'' = f'\alpha'' = (\kappa\wedge 1_{K})\alpha'' = (1_{E_{s+2}}\wedge \alpha'')
(\kappa\wedge 1_{K}).$  This shows (1). The proof of (2) is similar. Q.E.D.

\vspace{2mm}

{\bf Proof of Theorem 9.3.9}\quad At first, by the supposition of the main Theorem B on $i'_*i_*(h_0\sigma)\in
Ext_A^{s+1,tq+q}(H^*K,Z_p)$ in a permanent cycle in the ASS we have $(\bar c_{s+1}\wedge 1_K)(h_0\sigma\wedge 1_K)$ = 0.
The there exists $\eta'_{n,s+1}\in [\Sigma^{tq+q}K, E_{s+1}\wedge K]$ such that $(\bar b_{s+1}\wedge 1_K)\eta'_{n,s+1}
= (h_0\sigma\wedge 1_K)$.

By Lemma 9.3.10, it suffices to prove $(1_{E_{s+2}}\wedge \alpha'')(\kappa
\wedge 1_{K})$ = 0.  Note that, by $\bar a_{s+1}\kappa = \bar c_{s} \sigma$ we have
 $\bar a_{s+1}(1_{E_{s+2}}\wedge \alpha_{1})\kappa =\bar c_{s} (1_{KG_{s}}
\wedge\alpha_{1})\sigma$ = 0 and so  $(1_{E_{s+2}}\wedge \alpha_{1})\kappa =
\bar c_{s+1}(h_{0}\sigma)$ (up to scalar),  this is because $\pi_{tq+q}KG_{s+1} \cong Ext_{A}^{s+1,tq+q}
(Z_{p},Z_{p}) \cong Z_{p}\{h_{0}\sigma\}$.   Then, by Lemma 9.3.10
we have\\
{\bf (9.3.11)}\qquad $(1_{E_{s+2}}\wedge\alpha_{1}\wedge 1_{K})
(\kappa\wedge 1_{K}) = (\bar c_{s+1}\wedge 1_K)(h_0\sigma\wedge 1_K)$ = 0.\\
Moreover, by (9.1.20) we have
$\quad (1_{E_{s+2}}\wedge \rho\alpha'_{K'\wedge M})(\kappa\wedge
1_{K})i'$ $\quad = \quad (1_{E_{s+2}}\wedge \alpha')(\kappa\wedge 1_{K})i'$ = 0,
Then,  by (9.1.17), $(1_{E_{s+2}}\wedge\alpha'_{K'\wedge M})(\kappa\wedge 1_{K})i' =
(1_{E_{s+2}}\wedge (v\wedge 1_{M})\overline{m}_{M})f$
for some $f\in[\Sigma^{tq+q-1}M, E_{s+2}\wedge K'\wedge M]$ and so
$(1_{E_{s+2}}\wedge i')f = (1_{E_{s+2}}\wedge\rho(1_{K'}\wedge ij)\alpha'_{K'\wedge
M})(\kappa\wedge 1_{K})i' = (1_{E_{s+2}}\wedge \alpha'')(\kappa\wedge
1_{K})i' = (1_{E_{s+2}}\wedge \alpha')(\kappa\wedge 1_{K})i'ij$ = 0.
Hence $f = (1_{E_{s+2}}\wedge\alpha )f_{2}$ for some $f_{2}\in
[\Sigma^{tq-1}M, E_{s+2}\wedge M]$ and we have $(1_{E_{s+2}}\wedge
(x\wedge 1_{M})\alpha'_{K'\wedge M})(\kappa\wedge 1_{K})i' =
(1_{E_{s+2}}\wedge (i'\wedge 1_{M})\overline{m}_{M}\alpha )f_{2}$ = 0
 and $(1_{E_{s+2}}\wedge (x\wedge 1_{M})\alpha'_{K'\wedge
M})(\kappa\wedge 1_{K})\rho(v\wedge 1_{M}) = (1_{E_{s+2}}\wedge (x\wedge
1_{M})\alpha'_{K'\wedge M})(\kappa\wedge 1_{K})\rho(vi\wedge
1_{M})m_{M} +
 + (1_{E_{s+2}}\wedge (x\wedge 1_{M})\alpha'_{K'\wedge M})(\kappa\wedge 1_{K})
\rho (v\wedge 1_{M})\overline{m}_{M}(j\wedge 1_{M})$ = 0, this is because $\rho (v\wedge 1_{M})
\overline{m}_{M}$ = 0 , $\rho (vi\wedge 1_{M}) = i'$.  Then we have\\
{\bf (9.3.12)}   $\qquad\quad (1_{E_{s+2}}\wedge (x\wedge 1_{M})\alpha'_{K'\wedge
M})(\kappa\wedge 1_{K})\rho = f_{3}(y\wedge 1_{M})$ \\
with $f_{3}\in [\Sigma^{tq+2q+1}M, E_{s+2}\wedge K\wedge M]\cap
(ker d)$ (cf. (9.1.15) and Cor. 6.4.15).  It follows  that\\
{\bf (9.3.13)} $\qquad\quad (\bar a_{s+1}\wedge 1_{K\wedge M})f_{3} = f_{4}(\alpha i\wedge
1_{M})$\\
with $f_{4}\in [\Sigma^{tq+q}M\wedge M , E_{s+1}\wedge K\wedge
M]\cap (ker d)$(cf. (9.1.15) and Cor. 6.4.15).  Note that the $d_{1}$-cycle $(\bar b_{s+1}\wedge
1_{M})(1_{E_{s+1}}\wedge jj'\wedge 1_{M})f_4 \in [\Sigma^{tq-2}M\wedge M
, KG_{s+1}\wedge M]\cong Z_{p}\{ (\sigma'\wedge 1_{M})ij(j\wedge
1_{M})\}$, then $(\bar b_{s+1}\wedge 1_{M})(1_{E_{s+1}}\wedge jj'\wedge
1_{M})f_{4} = \lambda\cdot (\sigma'\wedge 1_{M})ij(j\wedge 1_{M})$ and
by applying the derivation $d$ we have $\lambda\cdot (\sigma'\wedge 1_{M})(j\wedge
1_{M})$ = 0  and this implies that $\lambda$ = 0.  That is to say $(\bar b_{s+1}\wedge
1_{M})(1_{E_{s+1}}\wedge jj'\wedge 1_{M})f_{4}$ = 0 , then $(\bar
b_{s+1}\wedge 1_{K\wedge M})f_{4} = (1_{KG_{s+1}}\wedge x\wedge 1_{M})g$
with $d_{1}$-cycle $g \in [\Sigma^{tq+q}M\wedge M, E_{s+1}\wedge K'\wedge M]\cap
(ker d)$( cf. Cor. 6.4.15).

By Theorem 6.4.3, $g = g (i\wedge 1_{M})m_{M} +
 g \overline{m}_{M}(j\wedge 1_{M})$. Now we claim that
$g(i\wedge 1_{M}) = \lambda_{1}(1_{KG_{s+1}}\wedge vi\wedge 1_{M})(h_{0}\sigma\wedge 1_{M})$
and $g \overline{m}_{M} = \lambda_{2}
(1_{KG_{s+1}}\wedge (v\wedge 1_{M})\overline{m}_{M})(h_{0}\sigma\wedge 1_{M})$
(mod $d_{1}$-boundary), where $\lambda_{1}, \lambda_{2} \in Z_{p}$.

To prove the claim, note that the $d_{1}$-cycle $g(i\wedge 1_{M})$
represents an element $[g(i\wedge 1_{M})] \in
Ext_{A}^{s+1,tq+q}(H^{*}K'\wedge M , H^{*}M)$ and $[(1_{KG_{s+1}}
\wedge \rho )g(i\wedge 1_M)]\in Ext_A^{s+1,tq+q}(H^*K,H^*M)$
$\cong Z_{p}\{[(1_{KG_{s+1}}\wedge i')(h_{0}\sigma\wedge 1_{M})]\}$ (cf. Prop. 9.3.2(4)).
Then $(1_{KG_{s+1}}\wedge \rho )g(i\wedge 1_{M}) =
\lambda_{1}(1_{KG_{s+1}}\wedge \rho (vi\wedge 1_{M}))(h_{0}\sigma\wedge 1_{M})
+ (\bar b_{s+1}\bar c_{s}\wedge 1_{K})g_{2}$ for some $ g_{2} \in [\Sigma^{tq+q}M,
KG_{s}\wedge K]$.  Since $(1_{KG_{s}}\wedge j'\alpha')g_{2}$ = 0, then
$g_{2} =
(1_{KG_{s}}\wedge \rho )g_{3}$ with $g_{3} \in [\Sigma^{tq+q}M, KG_{s}\wedge K'
\wedge M]$.  Then $g(i\wedge 1_{M})$
$= \lambda_{1}(1_{KG_{s+1}}\wedge vi\wedge 1_{M})(h_{0}\sigma\wedge 1_{M}) +
(\bar b_{s+1}\bar c_{s}\wedge 1_{K'\wedge M})g_{3} + (1_{KG_{s+1}}\wedge
\psi )g_{4}$
with $g_{4} \in [\Sigma^{tq+q-1}M, KG_{s+1}\wedge M] \cong Z_{p}\{(h_{0}\sigma
\wedge 1_{M})ij \}$ and so $g_{4} = \lambda'(h_{0}\sigma\wedge 1_{M})ij$ for some
 $\lambda' \in Z_{p}$.  However, $d(i\wedge 1_{M}) = 0 $ and $d(g)$ = 0
this implies that $d(g(i\wedge 1_{M}))$ = 0, then, by applying the derivation $d$ to the above equation
we have  $(1_{KG_{s+1}}\wedge \psi )d(g_{4}) + (\bar b_{s+1}\bar c_{s}\wedge
1_{K'\wedge M})d(g_{3}) = 0$, that is $\lambda'(1_{KG_{s+1}}\wedge\psi )(h_{0}\sigma
\wedge 1_{M}) = (\bar b_{s+1}\bar c_{s}\wedge 1_{K'\wedge M})d(g_{3})$ and this means that
the scalar $\lambda'$ = 0, this is because $\psi_{*}[h_{0}\sigma\wedge 1_{M}]\neq 0 \in
Ext_{A}^{s+1,tq+q+1}(H^{*}K'\wedge M , H^{*}M)$(cf. Prop. 9.3.5).
This shows that $g(i\wedge 1_{M}) =
\lambda_{1}(1_{KG_{s+1}}\wedge vi\wedge 1_{M})(h_{0}\sigma\wedge 1_{M})$ (mod
$d_{1}$-boundary). In addition, by $d(\overline{m}_{M})$  $\in [\Sigma^{2}M,
M\wedge M]\cong [\Sigma^{2}M,M] + [\Sigma M, M]$ = 0, then similarly we have
$g\overline{m}_{M} = \lambda_{2}(1_{KG_{s+1}}\wedge \psi )
(h_{0}\sigma\wedge 1_{M})$ (mod $d_{1}$-boundary).  This proves the above claim.

Then,  modulo $d_1$-boundary we have\\
{\bf (9.3.14)} $\qquad g  = g(i\wedge 1_{M})m_{M} +  g \overline{m}_{M}(j\wedge 1_{M})$

$ = \lambda_{1}(1_{KG_{s+1}}\wedge vi\wedge 1_{M})(h_{0}\sigma\wedge 1_{M})m_{M}
+ \lambda_{2}(1_{KG_{s+1}}\wedge \psi )(h_{0}\sigma\wedge 1_{M})(j\wedge 1_{M})$

$ = \lambda_{1}(h_{0}\sigma\wedge 1_{K'\wedge M})(vi\wedge 1_{M})m_{M}
+ \lambda_{2}(h_0\sigma\wedge 1_{K'\wedge M})(v\wedge 1_{M})\overline{m}_{M}(j\wedge 1_{M})$\newline
We claim that\\
{\bf (9.3.15)} \qquad The scalar in (9.3.14) $\lambda_{1} = \lambda_{2}$.\\
This will be proved in the last. Then , $g = \lambda_{1} (1_{KG_{s+1}}\wedge v\wedge 1_{M})(h_{0}\sigma\wedge
1_{M}\wedge 1_{M})$ and so we have $(\bar b_{s+1}\wedge 1_{K\wedge M})f_{4} =
(1_{KG_{s+1}}\wedge x\wedge 1_{M})g = \lambda_1 (1_{KG_{s+1}}\wedge i'\wedge
1_{M})(h_0\sigma\wedge 1_M\wedge 1_M) = \lambda_1 (h_0\sigma\wedge 1_K\wedge
1_M)(i'\wedge 1_M) + (\bar b_{s+1}\bar c_{s}\wedge 1_{K\wedge M})g_{5} =
\lambda_1 (\bar b_{s+1}\wedge 1_{K\wedge M})(\eta'_{n,s+1}\wedge 1_M)(i'\wedge 1_M) + (\bar b_2\bar
c_1\wedge 1_{K\wedge M})g_5$  and $f_4 = \lambda_1
(\eta'_{n,s+1}i'\wedge 1_M) + (\bar c_s\wedge 1_{K\wedge M})g_5 + (\bar
a_{s+1}\wedge 1_{K\wedge M})f_5 $ with $f_5 \in
[\Sigma^{tq+q+1}M\wedge M, E_{s+2}\wedge K\wedge M]$. It follows that
$(\bar a_{s+1}\wedge 1_{K\wedge M})f_3 = f_4(\alpha i\wedge 1_M) = (\bar
a_{s+1}\wedge 1_{K\wedge M})f_5(\alpha i\wedge 1_M)$ and so $f_3 = f_5
(\alpha i\wedge 1_M) + (\bar c_{s+1}\wedge 1_{K\wedge M})g_6$
for some $g_6\in [\Sigma^{tq+2q+1}M, E_{s+2}\wedge K\wedge M]$. So
$(1_{E_{s+2}}\wedge \alpha'')(\kappa\wedge 1_{K})\rho =
(1_{E_{s+2}}\wedge(1_K\wedge j)(x\wedge 1_M)\alpha'_{K'\wedge
M})(\kappa\wedge 1_{K})\rho = (1_{E_{s+2}}\wedge 1_K\wedge j)f_3 (y\wedge 1_M) =
(\bar c_{s+1}\wedge 1_K)(1_{KG_{s+1}}\wedge 1_K\wedge j)g_6 (y\wedge 1_M)$ =
0, this is because  the $d_1$-cycle $(1_{KG_{s+1}}\wedge 1_K\wedge j)g_6 \in
[\Sigma^{tq+2q}M, KG_{s+1}\wedge K]$ represents an element in
$Ext_{A}^{s+1,tq+2q}(H^{*}K, H^{*}M)$ = 0 (cf. Prop. 9.3.2(1)).

It follows from (9.1.11) that $(1_{E_{s+2}}\wedge \alpha'')(\kappa\wedge 1_{K}) = f_{6}\alpha ijj'$
for some $f_{6} \in [\Sigma^{tq+q+1}M, E_{s+2}\wedge K]$ and $(\bar a_{s+1}\wedge 1_{K})
f_{6}\alpha ijj' = (\bar a_{s+1}\wedge 1_{K})(1_{E_{s+2}}\wedge \alpha'')(\kappa\wedge 1_{K}) =
(\bar c_{1}\wedge 1_{K})(1_{KG_{s}}\wedge\alpha'')(\sigma\wedge 1_{K})$ = 0.
Then, by (9.1.14) we have $(\bar a_{s+1}\wedge 1_{K})f_{6}\alpha i =
f_{7} z $ with
$f_{7} \in [\Sigma^{tq+q-1}K', E_{s+1}\wedge K]$.  Moreover, by Prop. 9.1.21,
 $f_{7} z $ = 0, then $f_6\alpha i = (\bar
c_{s+1}\wedge 1_K)g_7$ for some $g_7 \in \pi_{tq+2q+1}(KG_{s+1}\wedge
K)$ and so $(1_{E_{s+2}}\wedge\alpha'')(\kappa\wedge 1_K) = f_6\alpha ijj'
= (\bar c_{s+1}\wedge 1_K)g_7jj'$ = 0, this is because the $d_1$-cycle $g_7\in
\pi_{tq+2q+1}(KG_{s+1}\wedge K)$ reprresents an element in
$Ext_{A}^{s+1,p^{n}q+2q+1}(H^{*}K , Z_p)$ = 0.
This shows the result of the Theorem and the remaining work is to prove the claim (9.3.15).

To prove (9.3.15), Note that by Theorem 6.4.3 and (9.1.15) we have $(v\wedge 1_M)\overline{m}_M(\alpha_1\wedge 1_M) =
(v\wedge 1_M)\overline{m}_M(j\wedge 1_M)(\alpha i\wedge 1_M) = -
(v\wedge 1_M)(i\wedge 1_M)m_M(\alpha i\wedge 1_M) = - (vi\wedge 1_M)\alpha$.
Similarly we have $\alpha (j\overline{u}\wedge 1_M) = - (\alpha_1\wedge
1_M)m_{M}(\overline{u}\wedge 1_M)$, where $\overline{u} : Y\rightarrow
\Sigma^{q+1}M$  and $v : \Sigma M\rightarrow K'$  are the map (9.1.5)(9.1.15).

Then, modulo  $d_1$-boundary we have

$\quad (1_{KG_{s+1}}\wedge vi\wedge 1_{M})(\widetilde{h_0\sigma}) = -
(1_{KG_{s+1}}\wedge v\wedge 1_M)\overline{m}_M(h_0\sigma\wedge 1_M)$

$\quad (\widetilde{h_0\sigma})(j\overline{u}\wedge 1_M) = - (h_0\sigma\wedge
1_M)m_M(\overline{u}\wedge 1_M)$\\
where $\widetilde{h_0\sigma}\in [\Sigma^{tq+q+1}M, KG_{s+1}\wedge M]$ is a $d_1$-cycle which represents
$\alpha_*(\tilde{\sigma})\in Ext_A^{s+1,tq+q+1}(H^*M,H ^*M)$.
So, by (9.3.14), modulo $d_1$-boundary we have

$\quad g (\overline{u}\wedge 1_{M})$
$ = \lambda_1(1_{KG_{s+1}}\wedge
v\wedge 1_M)(h_0\sigma\wedge 1_M\wedge 1_M)(\overline{u}\wedge 1_M) +
(\lambda_2 - \lambda_1)$

$\quad\qquad\qquad\qquad (1_{KG_{s+1}}\wedge v\wedge 1_{M})\overline{m}_{M}(h_0\sigma\wedge 1_M)(j\overline{u}\wedge 1_M)$

$\qquad\quad\qquad = (\lambda_1 - \lambda_2)(1_{KG_{s+1}}\wedge vi\wedge
1_M)(\widetilde{h_0\sigma})(j\overline{u}\wedge 1_M)$\\
this is because $(1_{KG_{s+1}}\wedge v)(h_0\sigma\wedge 1_M)\overline{u} =
(1_{KG_{s+1}}\wedge v)[\widetilde{h_0\sigma})ij + (1_{KG_{s+1}}\wedge
ij)(\widetilde{h_0\sigma})\overline{u}$ = 0 (mod $d_1$-boundary).
On the other hand, modulo $d_1$-boundary we have

$\quad g (\overline{u}\wedge 1_M) = \lambda_2 (1_{KG_{s+1}}\wedge v\wedge
1_M)(h_0\sigma\wedge 1_M\wedge 1_M)(\overline{u}\wedge 1_M) +
(\lambda_1 - $

$\qquad\qquad\qquad\quad\lambda_2)(1_{KG_{s+1}}\wedge vi\wedge 1_M)(h_0\sigma\wedge
1_M)m_M(\overline{u}\wedge 1_M)$

$\quad\qquad\qquad = (\lambda_2 - \lambda_1)(1_{KG_{s+1}}\wedge vi\wedge
1_M)\widetilde{h_0\sigma}(j\overline{u}\wedge 1_M)$.\\
Moreover, $(1_{KG_{s+1}}\wedge vi\wedge
1_M)\widetilde{h_0\sigma}(j\overline{u}\wedge 1_M)$ represents an nonzero element in the Exr group,
this is because $(1_{KG_{s+1}}\wedge (1_{K}\wedge i)(x\wedge 1_M))(1_{KG_{s+1}}\wedge vi\wedge
1_M)\widetilde{h_0\sigma}(j\overline{u}\wedge 1_M)(\overline{r}\wedge
1_M)(1_K\wedge i) = (1_{KG_{s+1}}\wedge i'ij)\widetilde{h_0\sigma}ijj' =
(1_{KG_{s+1}}\wedge i')(h_0\sigma\wedge 1_M)ijj'$ represents a nonzero element in the Ext group.
Then, by comparison to the above two equations we have $\lambda_1 - \lambda_2 =
\lambda_2 - \lambda_1$ so that  $\lambda_1 = \lambda_2$. This shows the claim (9.3.15). Q.E.D.

{\bf Remark}\quad In the last of section 4, we will also give another proof of Theorem 9.3.9.

\vspace{2mm}

{\bf Proof of the main Theorem B}\quad  By Theorem 9.3.9, there
exists \\$(\eta_{n,s+1})'' \in [\Sigma^{tq+q-1}K, E_{s+1}\wedge
K]$ such that $(\bar b_{s+1}\wedge 1_K)(\eta_{n,s+1})'' =
(h_0\sigma)''\in [\Sigma^{tq+q-1}K, KG_{s+1}\wedge K]$. Let
$(\eta_n)'' = (\bar a_0\cdots \bar a_s\wedge 1_K)
(\eta_{n,s+1})''\in [\Sigma^{tq+q-s-2}K,\\K]$ and  consider the map\\
\centerline{$(\eta_n)''\beta i'i\in\pi_{tq+pq+2q-s-2}K$}\\
 where $\beta\in [\Sigma^{(p+1)q}K,K]$ is the known second periodicity element which has filtration 1.
  Since $(\eta_n)''$ is represented by $(h_0\sigma)''
\in Ext_A^{s+1,tq+q-1}(H^*K,H^*K)$ in the ASS,  then
$(\eta_n)''\beta i'i\in \pi_{tq+pq+2q-s-2}K$ is represented by
$(\beta i'i)^*(h_0\sigma)''\\ = (\beta i'i)^*
(\alpha'')^*(\sigma)' = \alpha''_*\beta_*(i'i)_*(\sigma)\in
Ext_A^{s+2,tq+pq+2q}(H^*K,Z_p)$.  By [14] Theorem 3.2 and [15]
Theorem 5.2 we know that $\alpha''\beta i'i\in \pi_{pq+2q-2}K$ is
represented by $\alpha''_*\beta_*(i'i)_*(1) = (i'i)_*(g_0)\in
Ext_A^{2,pq+2q}(H^*K,Z_p)$(up to nonzero scalar)  in the ASS so
that $(\eta_n)''\beta i'i$ is represented by $\alpha''_*
\beta_*(i'i)_*(\sigma) = (i'i)_*(g_0\sigma)\in
Ext_A^{s+2,tq+pq+2q}(H^*K,Z_p)$. Q.E.D.

\vspace{2mm}

Using the stronger result of the main Theorem A which is stated in the Remark 9.2.35, the result of the main Theorem B
also can be obtained by the following main Theorem B'.

{\bf The main Theorem B'}\quad  Let $\sigma\in
Ext_A^{s,tq}(Z_p,Z_p), \sigma'\in Ext_A^{s+1,tq}(Z_p,\\Z_p)$ be a
pair of $a_0$-related elements, that is, there is a secondary
differential $d_2(\sigma) = a_0\sigma'$. Suppose that all the
supposition of the main Theorem A hold, then
$(i'i)_*(g_0\sigma)\in Ext_A^{s+2,tq+pq+2q}(H^*K,Z_p),
(i'i)_*(g_0\sigma')\in Ext_A^{s+3,tq+pq+2q}\\(H^*K,Z_p)$ are
permanent cycles in the ASS.

{\bf Proof}\quad  Let $\phi\sigma\in\pi_{tq+2q}KG_{s+1}, \phi\sigma'\in\pi_{tq+2q}KG_{s+2}$ be $d_1$-cycles which represent
$\phi_*(\sigma)\in Ext_A^{s+1,tq+2q}(H^*L,Z_p), \phi_*(\sigma')\in Ext_A^{s+2,tq+2q}(H^*L,Z_p)$ respectively. By the stronger result of the
main Theorem A (cf. Remark 9.2.35) we have $(\bar c_{s+2}\wedge 1_L)\phi\sigma' = 0, (\bar c_{s+1}
\wedge 1_{L\wedge M})(1_{KG_{s+1}}\wedge 1_L\wedge i)\phi\sigma = 0$.  Then $(\bar c_{s+2}\wedge 1_{L\wedge K})
(\phi\sigma'\wedge 1_K)$ = 0 and by using the multiplication of the ring spectrum $K$ we have
$(\bar c_{s+1}\wedge 1_{L\wedge K})(\phi\sigma)\wedge 1_K) = 0$.  In addition,
$(h_0\sigma)'' = (1_{KG_{s+1}}\wedge\overline{\Delta})(\phi\sigma\wedge 1_K), (h_0\sigma')'' = (1_{KG_{s+2}}\wedge
\overline{\Delta})(\phi\sigma'\wedge 1_K)$, this is because $\overline{\Delta}(\phi\wedge 1_K) = \alpha''\in [\Sigma^{q-2}K,K]$.
Then we have $(\bar c_{s+1}\wedge 1_K)(h_0\sigma)'' = 0,(\bar c_{s+2}\wedge 1_K)(h_0\sigma')'' = 0$. The remaining steps is similiar
to that given in the proof of the main Theorem B. Q.E.D.

\quad

\begin{center}

{\bf\large \S 4.\quad A general result on pull back convergence of $h_0\sigma$}

\end{center}

\vspace{2mm}

In this section, we will prove that, under some suppositions, the convergence of the element
$(1_L\wedge i)_*\phi_*(\sigma)\in Ext_A^{s+1,tq+2q}(H^*L
\wedge M,Z_p)$ can be pull backed to obtain the convergence of $h_0\sigma\in Ext_A^{s+1,tq+q}(Z_p,Z_p)$
in the stable homotopy groups of spheres. We have the following main Theorem.

\vspace {2mm}

{\bf  The main Theorem C} ( generalization of [24] Theorem A)\quad Let $p\geq 5, s\leq 4$ and suppose that

(I)(a) $Ext_A^{s,tq}(Z_p,Z_p)\cong Z_p\{\sigma\}$, $Ext_A^{s+1,tq+q}(Z_p,Z_p)\cong Z_p\{h_0\sigma\}$

\qquad\quad $Ext_A^{s+2,tq+2q+1}(Z_p,Z_p)\cong
Z_p\{\widetilde{\alpha}_2\sigma\}$ satisfying $a_0^2\sigma\neq 0$.

\quad (b) $Ext_A^{s+1,tq+u}(Z_p,Z_p)\cong Z_p\{a_0\sigma\}$  for $u = 1$ and is zero for $u = 2,3$.

\qquad\quad  $Ext_A^{s+1,tq}(Z_p,Z_p)$ is zero or has (one or two) generator $\sigma'$ such that
(both) satisfies

\qquad\quad $h_0\sigma'\neq 0, a_0\sigma'\neq 0$,

\qquad\quad $Ext_A^{s+1,tq+rq+u}(Z_p,Z_p)$ = 0 for $r = -1,2,3, u = -2,-1,0,1,2,3$ or

\qquad\quad for $r = 1, u = -2,-1,1,2,3$

\quad (c) $Ext_A^{s,tq+u}(Z_p,Z_p)$ = 0 for $u = -1,1,2,3$

\qquad $Ext_A^{s,tq+rq+u}(Z_p,Z_p)$ = 0 for $r = -2,-1,1,2, u = -2,-1,0,1,2,3$

(II)   $(1_L\wedge i)_*(\phi)_*(\sigma)\in Ext_A^{s+1,tq+2q}(H^*L\wedge M,Z_p)$
is a permanent cycle in the ASS, then $(\alpha i)_*(\sigma)\in Ext_A^{s+1,tq+q+1}(H^*M,Z_p)$
also is a permanent cycle in tha ASS so that $h_0\sigma = j_*(\alpha i)_*(\sigma)\in Ext_A^{s+1,tq+q}
(Z_p,Z_p)$ converges to an element  in $\pi_{tq+q-s-1}S$ of order $p$.

Note that the supposition (I) of the main Thoerem C contains
the supposition I of the main Thoerem B, then some results on Ext groups in $\S 3$ also hold
under the supposition of the main Theorem C. Before proving the main Theorem C,
we first recall the properties of some spectra related to $K$ and $M$ and prove some
results on low dimensional Ext groups.

By (9.1.27), $((1_Y\wedge j)\alpha_{Y\wedge M}\wedge 1_M)\overline{m}_M = \alpha_{Y\wedge M}$
, (9.2.12) and the following homotopy commutative diagram of  $3\times 3$-Lemma

$\qquad\quad X\wedge M\stackrel{m_M(\tilde{\psi}\wedge 1_M)}{\longrightarrow}
\Sigma^{2q}M\qquad\stackrel{0}{\longrightarrow}\qquad\Sigma^{2q+2}M$

$\qquad\qquad\quad\searrow ^{\tilde{\psi}\wedge 1_M}\quad\nearrow m_M\qquad\searrow
^{(\phi\wedge 1_K) i'}\nearrow _{j'(j''\wedge 1_K)}\quad\searrow\overline{m}_M$\\
{\bf (9.4.1)}\qquad\qquad\qquad $\Sigma^{2q}M\quad\qquad\qquad\qquad
\Sigma L\wedge K\quad\qquad\qquad\Sigma^{2q+2}M\wedge M$

$\qquad\qquad\quad\nearrow \overline{m}_M\quad\searrow ^{(1_Y\wedge j)\alpha_{Y\wedge M}\wedge 1_M}
\nearrow \quad\qquad\searrow u'\quad\nearrow \widetilde{\psi}\wedge 1_M$

$\qquad\quad\Sigma^{2q+1}M\quad \stackrel{\alpha_{Y\wedge M}}{\longrightarrow}\quad
Y\wedge M\quad\stackrel{\tilde{u}w_2\wedge 1_M}{\longrightarrow}\quad\Sigma X\wedge M$\\
 we know that the cofibre of  $m_M(\widetilde{\psi}\wedge 1_M) :
X \wedge M\rightarrow \Sigma^{2q}M$ is $\Sigma L\wedge K$ given by the following cofibration\\
{\bf (9.4.2)}\qquad $X\wedge M\stackrel{m_M(\tilde{\psi}\wedge 1_M)}{\longrightarrow}
\Sigma^{2q}M\stackrel{(\phi\wedge 1_K)i'}{\longrightarrow}\Sigma L\wedge K
\stackrel{u'}{\longrightarrow}\Sigma X\wedge M$\\

Since $(1_L\wedge i')(\phi\wedge 1_M)m_M(\widetilde{\psi}\wedge 1_M)$ = 0,
then by $[\Sigma^{-q-1}X\wedge M, L\wedge M]\cap (ker d)\cong Z_p\{u''\wedge 1_M\}$ and (9.1.2) we have
$(\phi\wedge 1_M)m_M(\widetilde{\psi}\wedge 1_M)
= (1_L\wedge\alpha)(u''\wedge 1_M)$ (up to nonzero scalar).
Since $(\phi\wedge 1_K)i' \alpha $ = 0, then by (9.4.2), there exists $\alpha_{X\wedge M}
\in [\Sigma^{3q}M, X\wedge M]$ such that $m_M(\widetilde{\psi}\wedge 1_M)\alpha_{X\wedge M} = \alpha$.
In addition, $m_M(\widetilde{\psi}\wedge 1_M)\alpha_{X\wedge M}m_M(\widetilde{\psi}\wedge 1_M)
= \alpha m_M(\widetilde{\psi}\wedge 1_M) = m_M(\widetilde{\psi}\wedge 1_M)(1_X\wedge\alpha )$
so that by (9.4.2) we have $\alpha_{X\wedge M}m_M(\widetilde{\psi}\wedge 1_M) = 1_X\wedge\alpha$
modulo $(u')_*[\Sigma^qX\wedge M, L\wedge K]$ = 0 , this is because $[\Sigma^qL\wedge K, L\wedge K]$ = 0  and $[\Sigma^{3q}M, L\wedge K]$ = 0.  Concludingly we have\\
{\bf (9.4.3)}\qquad $(\phi\wedge 1_M)m_M(\widetilde{\psi}\wedge 1_M) =
(1_L\wedge\alpha)(u''\wedge 1_M),$

$\qquad\quad\alpha_{X\wedge M}m_M(\widetilde{\psi}\wedge 1_M) = 1_X\wedge\alpha $

The cofibre of the map $\alpha_{X\wedge M} : \Sigma^{3q}M\rightarrow X\wedge M$ is $W\wedge K$ given by the cofibration\\
{\bf (9.4.4)}\qquad $\Sigma^{3q}M\stackrel{\alpha_{X\wedge M}}{\longrightarrow} X\wedge M
\stackrel{\mu_{X\wedge M}}{\longrightarrow} W\wedge K\stackrel{j'(j''u\wedge 1_K)}
{\longrightarrow}\Sigma^{3q+1}M$\\
This can be seen by the following homotopy commutative diagram of $3\times 3$-Lemma

$\qquad\qquad\Sigma^{3q}M\quad\stackrel{\alpha}{\longrightarrow}\quad \Sigma^{2q}M
\quad\stackrel{(\phi\wedge 1_K)i'}{\longrightarrow}\quad\Sigma L\wedge K$

$\qquad\qquad\quad\searrow \alpha_{X\wedge M}\nearrow _{m_M(\tilde{\psi}\wedge 1_M)}\searrow i'\quad
\nearrow (\phi\wedge 1_K)$\\
{\bf (9.4.5)}$\qquad\qquad\quad X\wedge M\qquad\qquad\quad \Sigma^{2q}K$

$\qquad\qquad\quad\nearrow u'\quad\searrow \mu_{X\wedge M}\quad\nearrow  j''u\wedge 1_K
\searrow j'$

$\qquad\qquad L\wedge K\quad\stackrel{w\wedge 1_K}{\longrightarrow}\quad W\wedge K\stackrel{j'(j''u\wedge 1_K)}
{\longrightarrow}\quad\Sigma^{3q+1}M$

By (9.2.13), $ijm_M(\widetilde{\psi}\cdot\widetilde{u}\wedge 1_M) = ij(u_2\wedge 1_M)
= (u_2\wedge 1_M)(1_U\wedge ij) = m_M(\widetilde{\psi}\cdot\widetilde{u}\wedge 1_M)(1_U\wedge ij)
= m_M(\widetilde{\psi}\wedge 1_M)(1_X\wedge ij)(\widetilde{u}\wedge 1_M)$,
then we have $ijm_M(\widetilde{\psi}\wedge 1_M) = m_M(\widetilde{\psi}\wedge 1_M)
(1_X\wedge ij) + \lambda (j\widetilde{\psi}\wedge 1_M)$ for some $\lambda\in Z_p$. It follows that
$\lambda j(j\widetilde{\psi}\wedge 1_M) = - jm_M(\widetilde{\psi}\wedge 1_M)(1_X\wedge ij)
= - j\widetilde{\psi}(1_X\wedge j) = j(j\widetilde{\psi}\wedge 1_M)$ and so $\lambda = 1$.
In addition, $i'(\alpha_1\wedge 1_M)m_M(\widetilde{\psi}\wedge 1_M) = (j''\wedge 1_K)(1_L\wedge i')(\phi\wedge 1_M)m_M(\widetilde{\psi}\wedge 1_M)$ = 0,
then by (9.1.23) we have  $m_M(\widetilde{\psi}\wedge 1_M) = m_M(\overline{u}\wedge 1_M)\psi_{X\wedge M}$
, where $\psi_{X\wedge M}\in [\Sigma^{-q+1}X\wedge M, Y\wedge M]$. In addition,
$[\Sigma^{-q+1}X\wedge M,Y\wedge M]\cong
Z_p\{\psi_{X\wedge M}\}$ , this can be obtained from $[\Sigma^{-2q}X\wedge M, M]\cong Z_p\{m_M(\widetilde{\psi}\wedge 1_M)\}$
, (9.1.23) and $[\Sigma^{-q}X\wedge M, K]$ = 0.  Then, by $j'(j''u\wedge 1_K)\cdot \mu_{X\wedge M}$ = 0
and (9.1.27) we have $(u\wedge 1_K)\mu_{X\wedge M} = \overline{\mu}_2(1_Y\wedge i')\psi_{X\wedge M}$
(up to nonzero scalar).  Concludingly we have\\
{\bf (9.4.6)}\quad $j\widetilde{\psi}\wedge 1_M = ijm_M(\widetilde{\psi}\wedge 1_M) - m_M(\widetilde{\psi}\wedge 1_M)(1_X\wedge ij)$,

$\qquad (u\wedge 1_K)\mu_{X\wedge M} = \overline{\mu}_2(1_Y\wedge i')\psi_{X\wedge M}$\quad
(up to nonzero scalar)

$\qquad [\Sigma^{-q+1}X\wedge M,Y\wedge M]\cong Z_p\{\psi_{X\wedge M}\},$

$\qquad\quad m_M(\overline{u}\wedge 1_M)\psi_{X\wedge M} = m_M(\widetilde{\psi}\wedge 1_M)$,\newline

By the following homotopy commutative diagram of  $3\times 3$-Lemma

\quad

$\qquad\qquad L\wedge K\quad\stackrel{(1_X\wedge j)u'}{\longrightarrow}\quad\Sigma X\qquad
\stackrel{1_X\wedge p}{\longrightarrow}\qquad\Sigma X$

$\qquad\qquad\quad\searrow u'\qquad\nearrow 1_X\wedge j\qquad\searrow\omega\qquad\nearrow\widetilde{u}w_2$

\quad\qquad\qquad\qquad $X\wedge M\qquad\qquad\qquad\qquad Y$

$\qquad\qquad\quad\nearrow 1_X\wedge i\quad\searrow ^{m_M(\tilde{\psi}\wedge 1_M)}
\nearrow _{(1_Y\wedge j)\alpha_{Y\wedge M}}\searrow ^{\overline{\mu}_2(1_Y\wedge i'i)}$

$\qquad\qquad X\qquad\stackrel{\tilde{\psi}}{\longrightarrow}\qquad\Sigma^{2q}M\qquad
\stackrel{(\phi\wedge 1_K)i'}{\longrightarrow}\qquad\Sigma L\wedge K$\\
we know that the cofibre of $(1_X\wedge j)u' : L\wedge K\rightarrow\Sigma X$ is $Y$ given by the cofibration\\
{\bf (9.4.7)}\qquad $L\wedge K\stackrel{(1_X\wedge j)u'}{\longrightarrow}\Sigma X\stackrel
{\omega}{\longrightarrow} Y\stackrel{\overline{\mu}_2(1_Y\wedge i'i)}{\longrightarrow}\Sigma L\wedge K$\\
In addition, by the commutativity of the above rectangle we have\\
{\bf (9.4.8)}\qquad $\omega\wedge 1_M = \alpha_{Y\wedge M}m_M(\widetilde{\psi}\wedge 1_M)$.

{\bf Proposition 9.4.9}\quad Under the supposition (I) of the main Thoerem C we have

(1) $Ext_A^{s+1,tq+r}(H^*K,H^*M)$ = 0 for $r = 1,2$,

(2) $Ext_A^{s+1,tq+rq+1}(H^*K,H^*K)$ = 0 for $r = -1,0,1,2$.

{\bf Proof}: (1) By the supposition, $Ext_A^{s+1,tq-q+r}(Z_p,Z_p)$ = 0 for $r = -1,0,1,\\2,3$, then
$(j')_*Ext_A^{s+1,tq+r}(H^*K,Z_p)\subset Ext_A^{s+1,tq-q-r-1}(H^*M,Z_p)$ = 0
for $r = 1,2,3$ and so $Ext_A^{s+1,tq+r}(H^*K,Z_p) = (i')_*Ext_A^{s+1,tq+r}(H^*M,Z_p)$ = 0
for $r = 1,2,3$ (cf. Prop. 9.3.0(1)) and the result follows.

(2) Consider the following exact sequence ($r = -1,0,1,2$)

$0 = Ext_A^{s+1,tq+(r+1)q+2}(H^*K,H^*M)\stackrel{(j')^*}{\longrightarrow}Ext_A^{s+1,
tq+rq+1}(H^*K,H^*K)$

$\qquad\quad\stackrel{(i')^*}{\longrightarrow}Ext_A^{s+1,tq+rq+1}(H^*K,H^*M)$\\
induced by (9.1.2). The right group is zero for $r = 0,1,2$( cf. (1) and Prop. 9.3.2(1)(2)) and also is zero for $r = -1$
which is obtained by the supposition on $Ext_A^{s+1,tq-q+r}(Z_p,Z_p)$ = 0 for $r = -1,0,1,2$. The left group
is zero for $r = -1,0,1$
(cf. (1) and Prop. 9.3.2). The left group also is zero for $r = 2$, this is because
$Ext_A^{s+1,tq+rq+u}(Z_p,Z_p)$ = 0
for $r = 2,3, u = 0,1,2,3$ by the supposition. Then the middle group is zero as desired.
Q.E.D.

{\bf Proposition 9.4.10}\quad Under the supposition (I) of the main Theorem C we have

(1)\quad $Ext_A^{s+1,tq+q+1}(H^*W\wedge K, H^*X\wedge M)$ = 0.

(2)\quad $Ext_A^{s+1,tq+2q+1}(H^*Y,H^*M)\cong Z_p\{((1_Y\wedge j)\alpha_{Y\wedge M})_*(\tilde{\sigma})\}$,

$ Ext_A^{s+1,tq+q}(H^*Y,H^*Y)\cong Z_p\{(\overline{u})^*((1_Y\wedge j)\alpha_{Y\wedge M})_*(\tilde{\sigma})\}$,

(3)\quad $Ext_A^{s+1,tq+3q}(H^*X,H^*M)\cong Z_p\{((1_X\wedge j)\alpha_{X\wedge M})_*(\tilde{\sigma})\}$

{\bf Proof}: \quad (1) Consider the following exact sequence

$0 = Ext_A^{s+1,tq+3q+1}(H^*W\wedge K, H^*M)\stackrel{m_M(\tilde{\psi}\wedge 1_M)^*}
{\longrightarrow} Ext_A^{s+1,tq+q+1}(H^*W\wedge K, H^*X\wedge M)
\stackrel{(u')^*}{\longrightarrow} Ext_A^{s+1,tq+q+1}(H^*W
\wedge K, H^*L\wedge K)$\\
induced by (9.4.2). The right group is zero by Prop. 9.4.9(2) and (9.1.12)(9.1.3). The left group
also is zero by
$Ext_A^{s+1,tq+rq+1}(H^*K,H^*M)$ = 0(for $r = 1,2,3$)(cf. the proof of Prop. 9.4.9(2))
. Then the middle group is zero as desired.

(2) Since $\overline{u}(1_Y\wedge j)\alpha_{Y\wedge M}\in [\Sigma^{q-1}M,M]\cong Z_p\{ij\alpha, \alpha ij\}$,
then $\overline{u}(1_Y\wedge j)\alpha_{Y\wedge M} = \lambda_1 ij\alpha + \lambda_1\alpha ij$
where the scalar $\lambda_1,\lambda_2\in Z_p$ satisfy $\lambda_1j\alpha ij\alpha + \lambda_2 j\alpha^2ij$ = 0.
Consider the following exact sequence

$Ext_A^{s+1,tq+2q}(Z_p,H^*M)\stackrel{(\overline{w})_*}{\longrightarrow}
Ext_A^{s+1,tq+2q+1}(H^*Y,H^*M)$

$\qquad\quad\stackrel{(\overline{u})_*}{\longrightarrow} Ext_A^{s+1,tq+q}(H^*M,H^*M)\stackrel{(j\alpha)_*}{\longrightarrow}$\\
induced by (9.1.5). The left group is zero which can be obtained
by the supposition on $Ext_A^{s+1,tq+2q+k}(Z_p,Z_p)$ = 0 (for $k =
0,1$).  By Prop. 9.3.1(2), the right group has two generators
$(ij)_*\alpha_*(\tilde{\sigma})$ and
$\alpha_*(ij)_*(\tilde{\sigma})$ . Then
$(\overline{u})_*\\Ext_A^{s+1,tq+2q+1}(H^*Y, H^*M)$ has unique
generator $(\overline{u})_*((1_Y\wedge j)\alpha_{Y\wedge
M})_*(\tilde{\sigma})$ so that the first result follows. For the
second result, consider the following exact sequence

$Ext_A^{s+1,tq+q}(Z_p,Z_p)\stackrel{(\overline{w})_*}{\longrightarrow} Ext_A^{s+1,tq+q+1}(H^*Y,Z_p)$

$\qquad\quad \stackrel{(\overline{u})_*}{\longrightarrow} Ext_A^{s+1,tq}(H^*M,Z_p)
\stackrel{(j\alpha)_*}{\longrightarrow}$\\
induced by (9.1.5). By the supposition, the left group has unique generator
 $h_0\sigma = (j\alpha i)_*(\sigma)$
so that im $(\overline{w})_*$ = 0. The right group is zero or has (one or two) generator
$i_*(\sigma')$ such that $(j\alpha)_*i_*(\sigma') = h_0\sigma'\neq 0$.
Then the middle group is zero and so the second result follows.

(3) Since $\widetilde{\psi}(1_X\wedge j)\alpha_{X\wedge M}\in [\Sigma^{q-1}M,M]
\cong Z_p\{ij\alpha, \alpha ij\}$, then $\widetilde{\psi}(1_X\wedge j)\alpha_{X\wedge M} = \lambda_3 ij\alpha
+ \lambda_4 \alpha ij$, where the scalar $\lambda_3,\lambda_4\in Z_p$  satisfy
$\lambda_3(1_Y\wedge j)\alpha_{Y\wedge M} ij\alpha + \lambda_4 (1_Y\wedge j)\alpha_{Y\wedge M}\alpha ij$ = 0.
Then, similar to that in (2), $(\widetilde{\psi})_*\\Ext_A^{s+1,tq+3q}(H^*X,H^*M)$
has unique generator $(\widetilde{\psi})_*((1_X\wedge j)\alpha_{X\wedge M})_*(\tilde{\sigma})$
so that $Ext_A^{s+1,tq+3q}(H^*X,H^*M)$ has unique generator $((1_X\wedge j)\alpha_{X\wedge M})_*(\tilde{\sigma})$
, this is because $Ext_A^{s+1,tq+3q+1}(H^*Y,H^*M)$ = 0 which can be obtained by the supposition (I)(b) on
$Ext_A^{s+1,tq+rq+k}(Z_p,Z_p)$ = 0 for $r = 1,2, k = -1,0,1,2$).  Q.E.D.

{\bf Proposition 9.4.11}\quad Under the supposition (I) of the main Theorem C we have

(1)\quad $Ext_A^{s,tq-2q}(H^*M,H^*X\wedge M)\cong Z_p\{m_M(\widetilde{\psi}\wedge 1_M)^*(\tilde{\sigma})\}$.

(2)\quad $Ext_A^{s,tq+rq+u}(H^* K,H^*M)$ = 0 for $r = -1,1,2,3, u = 0,1,2$

\qquad $Ext_A^{s,tq}(H^*K,H^*K)\cong Z_p\{\sigma_K\}$ satisfying $(i')^*(\sigma_K) = (i')_*(\tilde{\sigma})$,

(3)\quad $Ext_A^{s,tq}(H^*L\wedge K,H^*L\wedge K)\cong Z_p\{\sigma_{L\wedge K}\}$

\qquad satisfying $(j''\wedge 1_K)_*(\sigma_{L\wedge K}) = (j''\wedge 1_K)^*(\sigma_K)$

\qquad $Ext_A^{s,tq+rq+u}(H^*L\wedge K,H^*M)$ = 0 for $r = 1,2,3,u = 0,1,2$,

(4)\quad $Ext_A^{s,tq+rq+u}(H^*W\wedge K,H^*M)$ = 0  for $r = 1,2,3,u = 0,1,2$,

\qquad $Ext_A^{s,tq+q}(H^*W\wedge K,H^*X\wedge M)$ = 0

{\bf Proof}: \quad (1) Consider the following exact sequence

$Ext_A^{s,tq}(H^*M, H^*M\wedge M)\stackrel{(\tilde{\psi}\wedge 1_M)^*}{\longrightarrow}
Ext_A^{s,tq-2q}(H^*M,H^*X\wedge M)$

$\qquad\quad\stackrel{(\tilde{u}w_2\wedge 1_M)^*}{\longrightarrow} Ext_A^{s,tq-2q}
(H^*M, H^*Y\wedge M)$\\
induced by (9.2.12).  By the suppopsition on
$Ext_A^{s,tq-rq+u}(Z_p,Z_p)$ = 0 with ($r = 1,2$, $u = 0,1,2$) and
the degree of the top cell of $Y\wedge M$ is $q + 3$ we know that
the right group is zero.
 Since $(\overline{m}_M)^*Ext_A^{s,tq}(H^*M,H^*M\wedge M)\subset
Ext_A^{s,tq+1}(H^*M,H^*M)$ = 0 ( cf. Prop. 9.3.0(2)), then the left group
has unique generator $(m_M)^*(\tilde{\sigma})$ and so the result follows.

(2) Consider the following exact sequence   ($r = -1,1,2,3, u = 0,1,2$)

$Ext_A^{s,tq+rq+u}(H^*M,H^*M)\stackrel{(i')_*}{\longrightarrow} Ext_A^{s,tq+rq+u}(H^*K,H^*M)$

$\qquad\quad\stackrel{(j')_*}{\longrightarrow} Ext_A^{s,tq+(r-1)q+u-1}
(H^*M,H^*M)\stackrel{\alpha_*}{\longrightarrow}$\\
induced by (9.1.2). The left group is zero for $r = -1,1,2,3,u = 0,1,2$, this is obtained from
the supposition  I(c) on
$Ext_A^{s,tq+rq+k}(Z_p,Z_p)$ = 0 (for $r = -1,1,2,3,k = -1,0,1,2,3$). By the supposition and Prop.
9.3.0(2),the right group is zero except for
$r = 1, u = 0,1$ it has unique generator
$(ij)_*(\tilde{\sigma})$ or $\tilde{\sigma}$ respectively. However,
it satisfies $\alpha_*(ij)_*(\tilde{\sigma})\neq 0$, $\alpha_*(\tilde{\sigma})\neq 0$
then, the middle group is zero as desired. Consider the following exact sequence

$0 = Ext_A^{s,tq+q+1}(H^*K,H^*M)\stackrel{(j')^*}{\longrightarrow} Ext_A^{s,tq}(H^*K,H^*K)$

$\qquad\quad\stackrel{(i')^*}{\longrightarrow} Ext_A^{s,tq}(H^*K,H^*M)\stackrel
{\alpha^*}{\longrightarrow}$\\
induced by (9.1.2). The left group is zero as shown above.  The right group has
unique generator $(i')_*(\tilde{\sigma})$ , this is because $(j')_*Ext_A^{s,tq}\\(H^*K,H^*M)
\subset Ext_A^{s,tq-q-1}(H^*M,H^*M)$ = 0  and $Ext_A^{s,tq}(H^*M,H^*M)
\cong Z_p\{\tilde{\sigma}\}$. Then the middle has unique generatot $\sigma_K$ as desired.

(3) Consider the following exact sequence  ($r = -1,0$)

$Ext_A^{s,tq+(r+1)q}(H^*K,H^*K)\stackrel{(j''\wedge 1_K)^*}{\longrightarrow} Ext_A^{s,tq+rq}(H^*K,H^*L\wedge K)$

$\qquad\quad\stackrel{(i''\wedge 1_K)^*}{\longrightarrow}
Ext_A^{s,tq+rq}(H^*K,H^*K)\stackrel{(\alpha_1\wedge 1_K)^*}{\longrightarrow}$\\
induced by (9.1.3). The left group is zero for  $r = 0$ , this is because by (2)
$(i')^*Ext_A^{s,tq+q}(H^*K,H^*K)\subset
Ext_A^{s,tq+q}(H^*K,H^*M)$ = 0 and  $Ext_A^{s,tq+2q+1}\\(H^*K,H^*M)$ = 0.
Moreover, by (2), the left group has unique generator  $\sigma_K$ for $r = -1$. The right group
is zero  $r = -1$, this is because
by (2) $(i')^*Ext_A^{s,tq-q}(H^*K,H^*K)
\subset Ext_A^{s,tq-q}(H^*K,H^*M)$ = 0  and $Ext_A^{s,tq+1}(H^*K,\\H^*M)$ = 0.
The right group has unique generator $\sigma_K$ for $r = 0$ which satisfies
$(\alpha_1\wedge 1_K)^*(\sigma_K)\neq 0\in Ext_A^{s+1,tq+q}(H^*K,H^*K)$
, this is because $(i')^*(\alpha_1\wedge 1_K)^*(\sigma_K) = (\alpha_1\wedge 1_M)^*(i')^*(\sigma_K)
= (\alpha_1\wedge 1_M)^*(i')_*(\tilde{\sigma}) = (i')_*(\alpha_1\wedge 1_M)_*(\tilde{\sigma})\neq 0\in
Ext_A^{s+1,tq+q}(H^*K,H^*M)$. Then the middle group is zero for $r = 0$ and has unique generator
$(j''\wedge 1_K)^*(\sigma_K)$ for $r = -1$ so that the first result
can be obtained by the following exact sequence

$0 = Ext_A^{s,tq}(H^*K,H^*L\wedge K)\stackrel{(i''\wedge 1_K)_*}{\longrightarrow}
Ext_A^{s,tq}(H^*L\wedge K,H^*L\wedge K)$

$\qquad\quad\stackrel{(j''\wedge 1_K)_*}{\longrightarrow} Ext_A^{s,tq-q}
(H^*K,H^*L\wedge K)\stackrel{(\alpha_1\wedge 1_K)_*}{\longrightarrow}$\\
induced by (9.1.3). For the second result, look at the following exact sequence ($r = 1,2,3, u = 0,1,2$)

$0 = Ext_A^{s,tq+rq+u}(H^*K, H^*M)\stackrel{(i''\wedge 1_K)_*}{\longrightarrow} Ext_A^{s,tq+rq+u}(H^*L\wedge K,H^*M)$

$\qquad\quad\stackrel{(j''\wedge 1_K)_*}{\longrightarrow} Ext_A^{s,tq+(r-1)q+u}
(H^*K,H^*M)\stackrel{(\alpha_1\wedge 1_K)_*}{\longrightarrow}$\\
induced by (9.1.3). By (2), the left group is zero for $r = 1,2,3, u = 0,1,2$ and the right group
also is zero for $r = 2,3, u = 0,1,2$. By Prop. 9.3.0 and the supposition , the right group
also is zero for $r = 1,u = 1,2$.
For $r = 1, u = 0$, The right group  has unique generator $(i')_*(\tilde{\sigma})$ which satisfies
$(\alpha_1\wedge 1_K)_*(i')_*(\tilde{\sigma})\neq 0$.
Then the middle group is zero for  $r = 1,2,3,u = 0,1,2$.

(4) Consider the following exact sequence ($r = 1,2,3,u = 0,1,2$)

$0 = Ext_A^{s,tq+rq+u}(H^*L\wedge K,H^*M)\stackrel{(w\wedge 1_K)_*}{\longrightarrow}
Ext_A^{s,tq+rq+u}(H^*W\wedge K,H^*M)$

$\qquad\quad \stackrel{(j''u\wedge 1_K)_*}{\longrightarrow} Ext_A^{s,tq+(r-2)q+u}(H^*K,H^*M)
\stackrel{(\phi\wedge 1_K)_*}{\longrightarrow}$\\
induced by (9.1.12). By (3), the left group is zero for $r = 1,2,3,u = 0,1,2$. By (2), the right group
is zero for $r = 1,3,u = 0,1,2$ and by Prop. 9.3.0 and the supposition, it also is zero for $r = 2, u = 1$.
  For  $r = 2, u = 0$, the right group has unique generator $(i')_*(\tilde{\sigma})$ which satisfies
$(\phi\wedge 1_K)_*(i')_*(\tilde{\sigma})\neq 0\in Ext_A^{s+1,tq+2q}(H^*L\wedge K,H^*M)$.
Then the middle group is zero for $r = 1,2,3,u = 0,1,2$ as desired.

Since $(\widetilde{u}w_2\overline{w}\wedge 1_M)^*Ext_A^{s,tq+q}(H^*W\wedge K,H^*X\wedge M)\subset
Ext_A^{s,tq+q}(H^*W\wedge K,H^*M)$ = 0, then, by (9.1.5), $(\widetilde{u}w_2\wedge 1_M)^*
Ext_A^{s,tq+q}(H^*W\wedge K,H^*X\wedge M) = (\overline{u}\wedge 1_M)^*Ext_A^{s,tq+2q+1}(H^*W\wedge K,H^*M\wedge M)$ = 0.
and by using (9.2.12) we know that $Ext_A^{s,tq+q}(H^*W\wedge K,H^*X\wedge M) =
(\widetilde{\psi}\wedge 1_M)^*Ext_A^{tq+3q}(H^*W\wedge K,H^*M\wedge M)$ = 0.  Q.E.D.

\quad

The proof of the main Theorem C will be done by some argument processing in the
Adams resolution (cf. 9.2.9) of some spectra related to the sphere spectrum $S$.
Before proving the main Theorem C , we first prove the following Lemmas.

{\vskip 2mm}

{\bf Lemma 9.4.12}\quad  Under the supposition (I)(II) of the main Theorem C we have

(1) \quad  Let $\widetilde{h_0\sigma}\in [\Sigma^{tq+q+1}M,KG_{s+1}\wedge M]$
be a $d_1$-cycle which represents $\alpha_*(\tilde{\sigma})\in Ext_A^{s+1,tq+q+1}(H^*M,H^*M)$,
 then $(\bar c_{s+1}\wedge 1_M)\widetilde{h_0\sigma} = (1_{E_{s+2}}\wedge\alpha )(\kappa\wedge 1_M)$
(up to scalar), where $\kappa\in \pi_{tq+1}E_{s+2}$ such that
$\bar a_{s+1}\cdot\kappa = \bar c_s\cdot \sigma$  and $\sigma\in \pi_{tq}KG_s
\cong Ext_A^{s,tq}(Z_p,Z_p)$.

(2) \quad $(1_{E_{s+2}}\wedge\phi\wedge 1_M)(\kappa\wedge 1_M)$ = 0, $(1_{E_{s+2}}\wedge \alpha_1\wedge 1_M)(\kappa\wedge 1_M)$ = 0.

{\bf Proof}:   (1)  Since $(1_{KG_{s+1}}\wedge i')\widetilde{h_0\sigma}$ is a $d_1$-boundary,
then $(\bar c_{s+1}\wedge 1_K)\\(1_{KG_{s+1}}\wedge i')(\widetilde{h_0\sigma})$ = 0 so that
$(\bar c_{s+1}\wedge 1_M)\widetilde{h_0\sigma} = (1_{E_{s+2}}\wedge \alpha )
f'$ for some $f'\in [\Sigma^{tq+1}M, E_{s+2}\wedge M]$.  It follows that $(\bar a_{s+1}
\wedge 1_M)(1_{E_{s+2}}\wedge\alpha )f'$ = 0  and so $(\bar a_{s+1}\wedge 1_M)f' = (1_{E_{s+1}}\wedge j')
f'_2$ with $f'_2\in [\Sigma^{tq+q+1}M, E_{s+1}\wedge K]$. The $d_1$-cycle
$(\bar b_{s+1}\wedge 1_K)f'_2$ represents an element in $Ext_A^{s+1,tq+q+1}
(H^*K,H^*M)$ and this group is zero by Prop. 9.3.2(2),  then $(\bar b_{s+1}\wedge 1_K)f'_2 =
(\bar b_{s+1}\bar c_s\wedge 1_K)g'_0$ for some $g'_0\in [\Sigma^{tq+q+1}M,
KG_s\wedge K]$.  Consequently we have, $f'_2 = (\bar c_s
\wedge 1_K)g'_0 + (\bar a_{s+1}\wedge 1_K)f'_3$ for some $f'_3\in [\Sigma^{tq+
q+2}M, E_{s+2}\wedge K]$  and so $(\bar a_{s+1}\wedge 1_M)f' = (\bar a_{s+1}
\wedge 1_M)(1_{E_{s+2}}\wedge j')f'_3 + (\bar c_s\wedge 1_M)(1_{KG_s}\wedge j')g'_0 = (\bar a_2
\wedge 1_M)(1_{E_{s+2}}\wedge j')f'_3 + (\bar c_s\wedge 1_M)(\sigma\wedge 1_M) =
(\bar a_{s+1}\wedge 1_M)(1_{E_{s+2}}\wedge j')f'_3 +(\bar a_{s+1}\wedge 1_M)
(\kappa\wedge 1_M)$, where the $d_1$-cycle $(1_{KG_s}\wedge j')g'_0\in [\Sigma^{tq}M, KG_s
\wedge M]$ represents an element in $Ext_A^{s,tq}(H^*M,H^*M)$ and this group has
unique generator $\tilde{\sigma}$ so that it equals to $\sigma\wedge 1_M$ (mod $d_1$-boundary).
Hence we have $f' = (1_{E_{s+1}}\wedge j')f'_3 + (\kappa\wedge 1_M)
+ (\bar c_{s+1}\wedge 1_M)\tilde{g}_1$ for some $\tilde{g}_1\in [\Sigma^{tq+1}M,
KG_{s+1}\wedge M]$ and so $(\bar c_{s+1}\wedge 1_M)\widetilde{h_0\sigma}
= (1_{E_{s+2}}\wedge\alpha )f' = (1_{E_{s+1}}\wedge\alpha )(\kappa\wedge 1_M)$
which shows the result.

(2) Since $Ext_A^{s+1,tq+rq}(Z_p,Z_p)$ is zero for $r = 2$ and has
unique generator $h_0\sigma = (j'')_*(\phi)_*(\sigma)$ for$r = 1$,
then $Ext_A^{s+1,tq+2q}(H^*L,Z_p) \cong Z_p\{(\phi)_*(\sigma)\}$
and $Ext_A^{s+1,tq+2q}(H^*W,Z_p)$ = 0. By this and a similar proof
as given in (1) we know that $(1_{E_{s+2}}\wedge\phi)\kappa =
(\bar c_{s+1}\wedge 1_L) \sigma\phi$ (up to scalar), where
$\sigma\phi\in \pi_{tq+2q}(KG_{s+1}\wedge L)$ is a $d_1$-cycle
which represents $(\phi)_*(\sigma)\in
Ext_A^{s+1,tq+2q}(H^*L,\\Z_p)$. Then, by the supposition (II) of
the main Theorem C we have
 $(1_{E_{s+2}}\wedge\phi\wedge 1_M)(\kappa\wedge 1_M) = (\bar c_{s+1}\wedge 1_{L\wedge M})
 (\sigma\phi\wedge 1_M)$ = 0 so that the result follows. Q.E.D.

{\bf  Lemma 9.4.13}\quad Under the supposition (I) of the main Theorem C we have

(1)\quad $ Ext_A^{s,tq}(H^*X\wedge M,H^*X\wedge M)\cong Z_p\{[\sigma\wedge 1_{X\wedge M}]\}$.

(2) \quad  For any $d_1$-cycle $g_0\in [\Sigma^{tq+q}X, KG_{s+1}\wedge X]$, $g_0 =
\lambda'(h_0\sigma\wedge 1_X)$ (mod $d_1$-boundary)
with $\lambda'\in Z_p$ and $(\psi_{X\wedge M})_*[h_0\sigma\wedge 1_{X\wedge M}]\neq 0\in Ext_A^{s+1,tq+1}(H^*Y\wedge M,H^*X\wedge M)$.

{\bf Proof}\quad  (1) Consider the following exact sequence

$Ext_A^{s,tq+2q}(H^*L\wedge K,H^*M)\stackrel{m_M(\tilde{\psi}\wedge 1_M)^*}
\to\longrightarrow Ext_A^{s,tq}(H^*L\wedge K,H^*X\wedge M)$

$\qquad\quad\stackrel{(u')^*}{\longrightarrow} Ext_A^{s,tq}(H^*L\wedge K,H^*L\wedge K)\stackrel
{((1_L\wedge i')(\phi\wedge 1_M))^*}{\longrightarrow}$\\
induced by (9.4.2).  By Prop. 9.4.11(3), the left group is zero and
the right group has unique generator $\sigma_{L\wedge K}$ which satisfies
$((1_L\wedge i')(\phi\wedge 1_M))^*(\sigma_{L\wedge K})
\neq 0\in Ext_A^{s+1,tq+2q}(H^*L\wedge K,H^*M)$, this is because $(j''\wedge 1_K)_*((1_L\wedge i')(\phi\wedge 1_M))^*(\sigma_{L\wedge K})
= ((1_L\wedge i')(\phi\wedge 1_M))^*(j''\wedge 1_K)_*(\sigma_{L\wedge K}) = ((1_L\wedge i')(\phi\wedge 1_M)^*(j''\wedge 1_K)^*(\sigma_K)
= ((\alpha_1\wedge 1_K)i')^*(\sigma_K) = (\alpha_1\wedge 1_M)^*(i')_*(\tilde{\sigma})
= (i'(\alpha_1\wedge 1_M))_*(\tilde{\sigma})\neq 0\in Ext_A^{s+1,tq+q}(H^*K,H^*M)$.
Then the middle group is zero. Look at the following exact sequence

$Ext_A^{s,tq}(H^*L\wedge K,H^*X\wedge M)\stackrel{(u')_*}{\longrightarrow}
Ext_A^{s,tq}(H^*X\wedge M,H^*X\wedge M)$

$\qquad\quad\stackrel{m_M(\tilde{\psi}\wedge 1_M)_*}{\longrightarrow} Ext_A^{s,tq-2q}(H^*M,H^*X\wedge M)
\stackrel{((1_L\wedge i')(\phi\wedge 1_M))_*}{\longrightarrow}$\\
induced by (9.4.2). As shown above, the left group is zero.   By Prop. 9.4.11(1), the
right group has unique generator $m_M(\widetilde{\psi}\wedge 1_M)^*(\tilde{\sigma})\\ = m_M(\widetilde{\psi}\wedge 1_M)^*[\sigma\wedge 1_M]
= [(\sigma\wedge 1_M)m_M(\widetilde{\psi}\wedge 1_M)] = [(1_{KG_s}\wedge m_M(\widetilde{\psi}\wedge 1_M))(\sigma\wedge 1_{X\wedge M})]
= m_M(\widetilde{\psi}\wedge 1_M)_*[\sigma\wedge 1_{X\wedge M}]$
and it satisfies $((1_L\wedge i')(\phi\wedge 1_M))_*m_M(\widetilde{\psi}\wedge 1_M)^*(\tilde{\sigma})
= ((1_L\wedge i')(\phi\wedge 1_M))_*m_M(\widetilde{\psi}\wedge 1_M)_*[\sigma\wedge 1_{X\wedge M}]$ = 0.
Then the middle group has unique generator $[\sigma\wedge 1_{X\wedge M}]$ as desired.

(2) Note that $(\widetilde{\psi})_*(\widetilde{u}w_2)^*Ext_A^{s+1,tq+q}(H^*X,H^*X)\subset
Ext_A^{s+1,tq-q-1}(H^*M,\\H^*Y)$. Similar to that in Prop. 9.3.0(1), by the supposition
we know that $Ext_A^{s+1,tq}(H^*M,H^*M)$ is zero or has (one or two ) generator
$\tilde{\sigma}'$, then\\ $Ext_A^{s+1,tq-q-1}(H^*M,H^*Y)$ is zero or has (one or two) generator
$(\overline{u})^*(\tilde{\sigma}')$ and it satisfies
 $((1_Y\wedge j)\alpha_{Y\wedge M})_*(\overline{u})^*(\tilde{\sigma}')
= ((1_Y\wedge j)\alpha_{Y\wedge M})_*(\overline{u})_*[\sigma'\wedge 1_Y] = (1_Y\wedge \alpha_1)_*[\sigma'\wedge 1_Y]
= [h_0\sigma'\wedge 1_Y]\neq 0$, then $(\widetilde{\psi})_*(\widetilde{u}w_2)^*Ext_A^{s+1,tq+q}(H^*X,H^*X)$ = 0
and so we have $(\widetilde{u}w_2)^*Ext_A^{s+1,tq+q}(H^*X,H^*X) = (\widetilde{u}w_2)_*Ext_A^{s+1,tq+q}(H^*Y,\\H^*Y)$ = 0
, this is because $Ext_A^{s+1,tq+q}(H^*Y,H^*Y)\cong Z_p\{((1_Y\wedge j)\alpha_{Y\wedge M})_*(\overline{u})^*\\(\tilde{\sigma})\}$ ( cf. Prop. 9.4.10(2)).
 Then $Ext_A^{s+1,tq+q}(H^*X,H^*X) = (\widetilde{\psi})^*Ext_A^{s+1,tq+3q}\\(H^*X,H^*M)$ and it has unique generator
$(\widetilde{\psi})^*((1_X\wedge j)\alpha_{X\wedge M})_*(\tilde{\sigma})
= ((1_X\wedge j)\alpha_{X\wedge M})_*[(\sigma\wedge 1_M)\widetilde{\psi}]
= ((1_X\wedge j)\alpha_{X\wedge M})_*[(1_{KG_{s+1}}\wedge\widetilde{\psi})(\sigma\wedge 1_X)]
= ((1_X\wedge j)\alpha_{X\wedge M})_*m_M(\widetilde{\psi}\wedge 1_M)_*(1_X\wedge i)_*[\sigma\wedge 1_X]
= (1_X\wedge j\alpha i)_*[\sigma\wedge 1_X] = [h_0\sigma\wedge 1_X]$ (cf. Prop. 9.4.10(3))
Then the first result follows.
For the second result , by (9.4.6),  the $d_1$-cycle $(1_{KG_{s+1}}\wedge m_M(\overline{u}\wedge 1_M)\psi_{X\wedge M})(h_0\sigma\wedge 1_{X\wedge M})
= (1_{KG_{s+1}}\wedge m_M(\widetilde{\psi}\wedge 1_M))(h_0\sigma\wedge 1_{X\wedge M})
= (h_0\sigma\wedge 1_M)m_M(\widetilde{\psi}\wedge 1_M)$ and it represents an element
$m_M(\widetilde{\psi}\wedge 1_M)^*[h_0\sigma\wedge 1_M]
= m_M(\widetilde{\psi}\wedge 1_M)^*(\alpha_1\wedge 1_M)_*(\tilde{\sigma})\neq 0$
so that the second result follows. Q.E.D.

\vspace {2mm}

{\bf Proof the main Theorem C} \quad  By Lemma 9.4.12(1), it suffices to prove
$(\bar c_{s+1}\wedge 1_M)\widetilde{h_0\sigma} = (1_{E_{s+1}}\wedge \alpha)(\kappa\wedge 1_M) = 0$.
The proof is divided into the following two steps.

{\bf  Step 1} \qquad  To prove \quad $(\kappa\wedge 1_{X\wedge M})(1_X\wedge\alpha)$ = 0.

By (9.4.3), $(\phi\wedge 1_M)m_M(\widetilde{\psi}\wedge 1_M) = (u''\wedge 1_M)(1_X\wedge\alpha)$,
then by Lemma 9.4.12(2) we have $(1_{E_{s+2}}\wedge u''\wedge 1_M)(1_{E_{s+2}}\wedge 1_X\wedge\alpha)(\kappa\wedge 1_{X\wedge M})
= (1_{E_{s+2}}\wedge \phi\wedge 1_M)(\kappa\wedge 1_M)m_M(\widetilde{\psi}\wedge 1_M)$
= 0.
Moreover, by (9.2.16)  we have $(1_{E_{s+2}}\wedge 1_X\wedge\alpha)(\kappa\wedge 1_{X\wedge M}) =
(1_{E_{s+2}}\wedge \widetilde{u}w_3\wedge 1_M)f$ for some $f\in [\Sigma^{tq+q+1}X\wedge M, E_{s+2}\wedge W\wedge M]\cap (ker d)$ (cf. Cor. 6.4.15).
By composing $(1_{E_{s+2}}\wedge 1_X\wedge i'i\wedge 1_M)$ on the above equation we have
$(1_{E_{s+2}}\wedge \widetilde{u}w_3\wedge 1_{K\wedge M})(1_{E_{s+2}}\wedge 1_W\wedge i'i\wedge 1_M)f
= (1_{E_{s+2}}\wedge (1_X\wedge (i'i\wedge 1_M)\alpha)(\kappa\wedge 1_{X\wedge M})
= (1_{E_{s+2}}\wedge 1_X\wedge\overline{m}_K i'(\alpha_1\wedge 1_M)))(\kappa\wedge 1_{X\wedge M})$ = 0
,where we use the result on $(1_{E_{s+2}}\wedge \alpha_1\wedge 1_M)(\kappa\wedge 1_M)$ = 0
in Lemma 9.4.12(2).
Consequently, by (9.2.16), $(1_{E_{s+2}}\wedge 1_W\wedge i'i\wedge 1_M)f = (1_{E_{s+2}}
\wedge w'(\pi\wedge 1_L)\wedge 1_{K\wedge M})f_2$ = 0 (for some $f_2\in [\Sigma^{tq+1}X\wedge M,
E_{s+2}\wedge L\wedge K\wedge M]$), this is because $ \pi\wedge 1_K $ = 0.
Then  by (9.1.4) we have
$ f = (1_{E_{s+2}}\wedge 1_W\wedge\epsilon\wedge 1_M)f_3 = (1_{E_{s+2}}\wedge 1_W
\wedge\alpha m_M(\overline{u}\wedge 1_M)f_3$  for some $f_3\in [\Sigma^{tq+q+2}X\wedge M, E_{s+2}
\wedge W\wedge Y\wedge M]\cap (ker d)$ ( cf. Cor. 6.4.15)  and so\\
{\bf (9.4.14)}\quad\qquad $(1_{E_{s+2}}\wedge 1_X\wedge\alpha)(\kappa\wedge 1_{X\wedge M})$

$\quad = (1_{E_{s+2}}\wedge\widetilde{u}w_3\wedge 1_M)(1_{E_{s+2}}\wedge 1_W\wedge
\alpha m_M(1_M\wedge\overline{u}))f_3$

$\quad = (1_{E_{s+2}}\wedge\alpha_{X\wedge M}(j''u\wedge 1_M))(1_{E_{s+2}}\wedge 1_W\wedge
m_M(1_M\wedge \overline{u}))f_3$ ( cf. (9.4.3))

By (9.4.14), $(\bar a_{s+1}\wedge 1_{X\wedge M})
(1_{E_{s+2}}\wedge \widetilde{u}w_3\wedge 1_M)(1_{E_{s+2}}\wedge (1_W\wedge\alpha m_M(\overline{u}\wedge 1_M))f_3$ =
$(\bar a_{s+1}\wedge 1_{X\wedge M})(1_{E_{s+2}}\wedge 1_X\wedge\alpha)(\kappa\wedge 1_{X\wedge M}) =
(\bar c_s\wedge 1_{X\wedge M})(1_{KG_s}\wedge 1_X\wedge\alpha)(\sigma\wedge 1_{X\wedge M})$ = 0
,this is because $\alpha$ induces zero homomorphism in $Z_p$-cohomology.
  Then, by (9.2.16)
and $w'(\pi\wedge1_L)\wedge 1_M = (w\wedge 1_M)(1_L\wedge\alpha)$  we have\\
{\bf (9.4.15)}\quad $(\bar a_{s+1}\wedge 1_{W\wedge M})(1_{E_{s+2}}\wedge 1_W\wedge
\alpha m_M(\overline{u}\wedge 1_M))f_3 $

$ = (1_{E_{s+1}}\wedge (1_W\wedge\alpha)(w\wedge 1_M))f_5$\\
with $f_5\in [\Sigma^{tq}X\wedge M, E_{s+1}\wedge L\wedge M]\cap (ker d)$ ( cf. Cor. 6.4.15 ).

By (9.4.15)(9.1.2) , $(\bar a_{s+1}\wedge 1_{W\wedge
M})(1_{E_{s+2}}\wedge 1_W\wedge m_M(\overline{u}\wedge 1_M))f_3 =
(1_{E_{s+1}}\wedge w\wedge 1_M)f_5 + (1_{E_{s+1}}\wedge 1_W\wedge
j')f_6$ for some $f_6\in [\Sigma^{tq+q+1}X\wedge M, E_{s+1}\wedge
W\wedge K]\cap (ker d)$ ( cf. Prop. 6.5.26). Since
$(1_W\wedge\alpha_1)w = w(1_L\wedge\alpha_1) = w\cdot \phi j''$ =
0,  then $w = (1_W\wedge j'')\psi_W$ , where $\psi_W\in
[\Sigma^qL, W\wedge L]$.  So we have $ w\wedge 1_M =  (1_W\wedge
j'')\psi_W\wedge 1_M = (1_W\wedge m_M(\overline{u}\wedge 1_M))
((1_W\wedge\bar h)\psi_W\wedge 1_M)$.  Hence, $ - (\bar
a_{s+1}\wedge 1_{W\wedge Y\wedge M})f_3 = (1_{E_{s+2}}\wedge
(1_W\wedge h)\psi_W\wedge 1_M)f_5 + (1_{E_{s+1}}\wedge 1_W\wedge
(1_Y\wedge i)r)f_6 + (1_{E_{s+1}}\wedge 1_W\wedge (r\wedge
1_M)\overline{m}_K)f_7$  and by Prop. 6.5.26, $f_7 = f_8(1_X\wedge
i') + f_9(1_X\wedge i'ij)$, where $f_8\in [\Sigma^{tq+q}X\wedge K,
E_{s+1}\wedge W\wedge K]\cap (ker d)$ and $f_9\in
[\Sigma^{tq+q+1}X\wedge K, E_{s+1}\wedge W\wedge K]\cap (ker d)$.
Since $d((1_Y\wedge i)r) = ((r\wedge 1_M)d(1_K\wedge i) = (r\wedge
1_M)(1_K\wedge m_M)(T_{K,M}\wedge 1_M)(1_M \wedge 1_K\wedge
i)\overline{m}_K = (r\wedge 1_M)(1_K\wedge m_M(1_M\wedge
i)\overline{m}_K = (r\wedge 1_M)\overline{m}_K$, by applying the
derivation $d$ using Theorem 6.4.8(1) we have $-
(1_{E_{s+1}}\wedge 1_W\wedge (r\wedge 1_M)\overline{m}_K)f_6 -
(1_{E_{s+1}}\wedge 1_W\wedge (r\wedge
1_M)\overline{m}_K)f_9(1_X\wedge i')$
= 0 (Note : $f_6$  has odd degree)  and so\\
{\bf (9.4.16)} \qquad $ - (\bar a_{s+1}\wedge 1_{W\wedge Y\wedge M})f_3 = (1_{E_{s+1}}\wedge
(1_W\wedge\bar h)\psi_W\wedge 1_M)f_5$

$\quad + (1_{E_{s+1}}\wedge 1_W\wedge (1_Y\wedge i)r)f_6
+ (1_{E_{s+1}}\wedge 1_W\wedge (r\wedge 1_M)\overline{m}_K)f_8(1_X\wedge i')$

$\quad  - (1_{E_{s+1}}\wedge 1_W\wedge (r\wedge 1_M)\overline{m}_K)f_6(1_X\wedge ij)$

Note that the $d_1$-cycle
$(\bar b_{s+1}\wedge 1_{W\wedge K})f_6\in [\Sigma^{tq+q+1}X\wedge M, KG_{s+1}\wedge W\wedge K]\cap (ker d)$
represents an element in $Ext_A^{s+1,tq+q+1}(H^*W\wedge K, H^*X\wedge M)$ and by Prop. 9.4.10(1) this group is zero,
then $(\bar b_{s+1}\wedge 1_{W\wedge K})f_6 = (\bar b_{s+1}\bar c_s\wedge 1_{W\wedge K})g$
for some $g\in [\Sigma^{tq+q+1}X\wedge M, KG_s\wedge W\wedge K]\cap (ker d)$ (cf. Prop. 6.5.26)   and so
$f_6 = (\bar c_s\wedge 1_{W\wedge K})g + (\bar a_{s+1}\wedge 1_{W\wedge K})f'$
with $f'\in [\Sigma^{tq+q+2}X\wedge M, E_{s+2}\wedge W\wedge K]\cap (ker d)$ ( cf. Prop. 6.5.26). Then we have\\
{\bf (9.4.17)}\qquad $- (\bar a_{s+1}\wedge 1_{W\wedge Y\wedge M})f_3
= (1_{E_{s+1}}\wedge (1_W\wedge \bar h)\psi_W\wedge 1_M)f_5$

$\qquad  + (\bar a_{s+1}\wedge 1_{W\wedge Y\wedge M})(1_{E_{s+2}}
\wedge 1_W\wedge (1_Y\wedge i)r)f' $

$\qquad + (\bar c_s\wedge 1_{W\wedge Y\wedge M})
(1_{KG_s}\wedge 1_W\wedge (1_Y\wedge i)r)g$

$\qquad  + (\bar a_{s+1}\wedge 1_{W\wedge Y\wedge M})(1_{E_{s+2}}\wedge 1_W\wedge (r\wedge 1_M)\overline{m}_K)f'(1_X\wedge ij)$

$\qquad - (\bar c_s\wedge 1_{W\wedge Y\wedge M})(1_{KG_s}\wedge 1_W\wedge (r\wedge 1_M)\overline{m}_K)g(1_X\wedge ij)$

$\qquad + (1_{E_{s+1}}\wedge 1_W\wedge (r\wedge 1_M)\overline{m}_K)f_8(1_X\wedge i')$

Let  $P$ be the cofibre of $(1_W\wedge\bar h)\psi_W : \Sigma^{q+1}L\rightarrow W\wedge Y$
given by the cofibration\\
{\bf (9.4.18)}\qquad $\Sigma^{q+1}L\stackrel{(1_W\wedge \bar h)\psi_W}{\longrightarrow}
W\wedge Y\stackrel{w_5}{\longrightarrow} P\stackrel{u_5}{\longrightarrow} \Sigma^{q+2}L$\newline
 Then the cofibre of $w_5(1_W\wedge r) : W\wedge K\rightarrow P$  is $\Sigma X$
given by the cofibration\\
{\bf (9.4.19)}\qquad $ W\wedge K\stackrel{w_5(1_W\wedge r)}{\longrightarrow}V
\stackrel{w_6}{\longrightarrow} \Sigma X\stackrel{u_6}{\longrightarrow}\Sigma W\wedge K$\\
 This can be seen by the following homotopy commutative diagram of $3\times 3$-Lemma

 \quad

$\qquad\qquad W\wedge K\quad\stackrel{w_5(1_W\wedge r)}{\longrightarrow}\quad P\quad
\stackrel{u_5}{\longrightarrow}\quad\Sigma^{q+2}L$

$\qquad\qquad\quad\searrow 1_W\wedge r\quad\nearrow w_5\quad\searrow w_6\quad\nearrow u''$

$\qquad\qquad\qquad\quad W\wedge Y\qquad\qquad\qquad \Sigma X$

$\qquad\qquad\quad\nearrow ^{(1_W\wedge\bar h)\psi_W}\searrow 1_W\wedge\epsilon
\nearrow \widetilde{u}w_3\searrow u_6$

$\qquad\qquad\Sigma^{q+1}L\stackrel{w'(\pi\wedge 1_L)}{\longrightarrow}\quad\Sigma W
\quad\stackrel{1_W\wedge i'i}{\longrightarrow}\quad\Sigma W\wedge K$

\quad

Note that $u_6 = \mu_{X\wedge M}(1_X\wedge i)$,  then by composing
$(\bar b_{s+1}\wedge 1_{P})(1_{E_{s+1}}\wedge w_5\wedge j)$ on the
left hand side of (9.4.17) and composing
 $(1_X\wedge i)$ on the right hand side  we have
$(\bar b_{s+1}\wedge 1_P)(1_{E_{s+1}}\wedge w_5(1_W\wedge r))f_8(1_X\wedge i'i)$ = 0
and so $(\bar b_{s+1}\wedge 1_{W\wedge K})f_8(1_X\wedge i'i) = (1_{KG_{s+1}}\wedge u_6)g_0 = (1_{KG_{s+1}}
\wedge \mu_{X\wedge M}(1_X\wedge i))g_0 = (1_{KG_{s+1}}\wedge\mu_{X\wedge M})(g_0\wedge 1_M)(1_X\wedge i)$
with  $d_1$-cycle $g_0\in [\Sigma^{tq+q}X, KG_{s+1}\wedge X]$. Moreover , by Lemma 9.4.13(2),
$g_0 = \lambda_1(h_0\sigma\wedge 1_X)$ (mod $d_1$-boundary), where $\lambda_1\in Z_p$.
On the other hand, by applying the derivation $d$ to  $(\bar b_{s+1}\wedge 1_{W\wedge K})f_8(1_X\wedge i'ij) = (1_{KG_{s+1}}\wedge\mu_{X\wedge M})(g_0\wedge 1_M)(1_X\wedge ij)$
we have

\quad\\
{\bf (9.4.20)}\qquad $(\bar b_{s+1}\wedge 1_{W\wedge K})f_8(1_X\wedge i') =
= (1_{KG_{s+1}}\wedge\mu_{X\wedge M})(g_0\wedge 1_M)$ ,

$\qquad\qquad g_0 = \lambda_1 (h_0\sigma\wedge 1_{X})\in [\Sigma^{tq+q}X, KG_{s+1}\wedge X]$
\qquad (mod $d_1$-boundary)

\quad

Consider the following commutative diagram of exact sequences

\quad

$\Sigma^{q+1}L\wedge M\quad\stackrel{w\wedge 1_M}{\longrightarrow}\quad \Sigma^{q+1}W\wedge M
\quad\stackrel{j''u\wedge 1_M}{\longrightarrow}\quad \Sigma^{3q+1}M\stackrel{\phi\wedge 1_M}{\longrightarrow}
\Sigma^{q+2}L\wedge M$

$\qquad\big\uparrow 1_{L\wedge M}\qquad\qquad\quad\qquad\big\uparrow 1_W\wedge m_M(\overline{u}\wedge 1_M)
\qquad\big\uparrow u_7\qquad\qquad\big\uparrow 1_{L\wedge M}$

$\Sigma^{q+1}L\wedge M\stackrel{(1_W\wedge\bar h)\psi_W\wedge 1_M}{\longrightarrow} W\wedge Y\wedge M
\quad\stackrel{w_5\wedge 1_M}{\longrightarrow}\quad P\wedge M\stackrel{u_5\wedge 1_M}
{\longrightarrow}\Sigma^{q+2}L\wedge M$\\
of the cofibrations (9.1.12)(9.4.18).  Since the left rectangle homotopy commutes
then there exists  $u_7\in [\Sigma^{-3q-1}P\wedge M, M]$
such that all the above rectangle homotopy commute.  That is we have\\
{\bf (9.4.21)}\qquad $u_7(w_5\wedge 1_M) = (j''u\wedge 1_M)(1_W\wedge m_M(\overline{u}\wedge 1_M)),\quad$

$\quad\qquad (\phi\wedge 1_M)u_7 = \pm\quad u_5\wedge 1_M$\\
where $u_7\in [\Sigma^{-3q-1}P\wedge M, M]$.  By the above two equations
, we have the following homotopy commutative diagram of $3\times 3$-Lemma in which
we use the cofibrations (9.2.12)(9.4.18)(9.1.23)

$\qquad P\wedge M\quad\stackrel{u_5\wedge 1_M}{\longrightarrow}\quad\Sigma^{q+2}L\wedge M
\quad\stackrel{w\wedge 1_M}{\longrightarrow}\quad \Sigma^{q+2} W\wedge M$

$\qquad\qquad\searrow u_7\quad\nearrow \phi\wedge 1_M\quad\searrow
^{((1_W\wedge \bar h)\psi_W\wedge 1_M}\nearrow _{1_W\wedge
m_M(\overline{u}\wedge 1_M)}\searrow ^{j''u\wedge 1_M}$\newline
{\bf (9.4.22)}$\qquad\qquad\Sigma^{3q+1}M\quad\qquad\qquad\qquad
\Sigma W\wedge Y\wedge M\qquad\qquad\quad \Sigma^{3q+2}M$

$\qquad\qquad\nearrow j''u\wedge 1_M\searrow ^{(\phi_W\wedge 1_K)i'}
\nearrow _{1_W\wedge (r\wedge 1_M)\overline{m}_K}\searrow w_5\wedge 1_M\quad\nearrow u_7$

$\quad\Sigma^{q+1}W\wedge M\stackrel{\tilde{\lambda}(1_W\wedge \alpha'i')}{\longrightarrow}
\Sigma^2W\wedge K\qquad\quad\longrightarrow\qquad\quad\Sigma P\wedge M$\\
Then there is a cofibration\\
{\bf (9.4.23)}\qquad $\Sigma^{3q-1}M\stackrel{(\phi_W\wedge 1_K)i'}
{\longrightarrow} W\wedge K\stackrel{(w_5\wedge 1_M)(1_W\wedge (r\wedge 1_M)\overline{m}_K)}
{\longrightarrow}$

$\qquad\qquad \Sigma^{-1}P\wedge M\stackrel{u_7}{\longrightarrow}\Sigma^{3q}M$\newline
in which $\phi_W\in [\Sigma^{3q-1}S, W]$ such that $u\cdot \phi_W = \phi\in [\Sigma^{2q-1}S, L]$.
Since $(\phi\wedge 1_K)i'\cdot u_7 = (u\cdot\phi_W\wedge 1_K)i'\cdot u_7$ = 0,
then by (9.4.2) we have\\
{\bf (9.4.24)}\qquad $u_7 = m_M(\widetilde{\psi}\wedge 1_M)u_8$\\
where $u_8\in [\Sigma^{-q-1}P\wedge M,X\wedge M]$.
On the other hand, by (9.4.8), $(\omega\wedge 1_M)u_8(w_5(1_W\wedge r)\wedge 1_M) = \alpha_{Y\wedge M}
m_M(\widetilde{\psi}\wedge 1_M)u_8(w_5(1_W\wedge r)\wedge 1_M) = \alpha_{Y\wedge M} u_7\\(w_5(1_W\wedge r)\wedge 1_M)
= \alpha_{Y\wedge M} j'(j''u\wedge 1_K)(1_W\wedge m_K)$ = 0 (cf. (9.1.27)).
Then, by (9.4.7),  $u_8(w_5(1_W\wedge r)\wedge 1_M) = ((1_X\wedge j)u'\wedge 1_M)\Delta_1$
with $\Delta_1\in [\Sigma^{-q}W\wedge K\wedge M, L\wedge K\wedge M]\cap (ker d)$.
By composing  $\mu_{X\wedge M}(1_X\wedge i)\wedge 1_M$ on the above equation and using
 (9.4.19) we have
$((1_X\wedge j)u'\wedge 1_M)\Delta_1(\mu_{X\wedge M}(1_X\wedge i)\wedge 1_M)$ = 0
and so by (9.4.7)(9.4.6),  $\Delta_1(\mu_{X\wedge M}(1_X\wedge i)\wedge 1_M) = (\overline{\mu}_2(1_Y\wedge i'i)\wedge 1_M)\psi_{X\wedge M}$.
 Then $(j''\wedge 1_{K\wedge M})\Delta_1(\mu_{X\wedge M}(1_X\wedge i)\wedge 1_M)
= ((j''\wedge 1_K)\overline{\mu}_2(1_Y\wedge i'i)\wedge 1_M)\psi_{X\wedge M} =
(i'\overline{u}\wedge 1_M)\psi_{X\wedge M} = (i'i\wedge 1_M)m_M(\overline{u}\wedge 1_M)\psi_{X\wedge M}
+ (i'\wedge 1_M)\overline{m}_M(j\overline{u}\wedge 1_M)\psi_{X\wedge M}$
and so $(j''\wedge 1_{K\wedge M})\Delta_1(\mu_{X\wedge M}(1_X\wedge i)\cdot\widetilde{u}w_2\wedge 1_M)$ = 0.
 Consequently we have  $(j''\wedge 1_{K\wedge M})\Delta_1(\mu_{X\wedge M}(\widetilde{u}w_2\wedge 1_M)\wedge 1_M)
\in (1_Y\wedge j\wedge 1_M)^*[\Sigma^{-2q}Y\wedge M, K\wedge M]$ = 0, this is because
the degree of the top cell of $Y\wedge M$ is $q+3$.  Then $(j''\wedge 1_{K\wedge M})\Delta_1(\mu_{X\wedge M}\wedge 1_M)\in
(\widetilde{\psi}\wedge 1_{M\wedge M})^*[M\wedge M\wedge M, K\wedge M]$ and so
 $(j'(j''\wedge 1_K)\wedge 1_M)\Delta_1(\mu_{X\wedge M}\wedge 1_M)$ = 0
and by (9.4.4) we have $(j'(j''\wedge 1_K)\wedge 1_M)\Delta_1 = \Delta_2 (j'(j''u\wedge 1_K)\wedge 1_M)
= \lambda (j'(j''u\wedge 1_K)\wedge 1_M)$ with $\lambda\in Z_p$, this is because
$\Delta_2\in [M\wedge M,M\wedge M]\cap (ker d)\cong Z_p\{1_{M\wedge M}\}$.
Hence, $(\widetilde{\psi}\wedge 1_M)u_8(w_5(1_W\wedge r)\wedge 1_M) =
(\widetilde{\psi}(1_X\wedge j)u'\wedge 1_M)\Delta_1 = (j'(j''\wedge 1_K)\wedge 1_M)\Delta_1 =
\lambda (j'(j''u\wedge 1_K)\wedge 1_M)$ and by (9.4.21)(9.4.24) we know that $\lambda = 1$ so that\\
{\bf (9.4.25)}\qquad $m_M(\widetilde{\psi}\wedge 1_M)u_8(w_5\wedge 1_M)(1_W\wedge
(1_Y\wedge i)r) = j'(j''u\wedge 1_K)$

$\qquad\quad  = (j\widetilde{\psi}\wedge 1_M)u_8(w_5\wedge 1_M)(1_W\wedge (r\wedge 1_M)\overline{m}_K)$ ,

$\qquad\quad (j\widetilde{\psi}\wedge 1_M)u_8(w_5\wedge 1_M)(1_W\wedge (1_Y\wedge i)r) = ijj'(j''u\wedge 1_K)$\\
where we use $(jj'\wedge 1_M)\overline{m}_K = j'$ in the above equation.
By composing $(1_{E_{s+1}}\wedge u_8(w_5\wedge 1_M))$ (it has odd degree) on (9.4.17) we have\\
{\bf (9.4.26)}\qquad $(\bar a_{s+1}\wedge 1_{X\wedge M})(1_{E_{s+2}}\wedge u_8(w_5\wedge 1_M))f_3$

$ = - (\bar a_{s+1}\wedge 1_{X\wedge M})(1_{E_{s+2}}\wedge u_8(w_5\wedge 1_M)(1_W\wedge (1_Y\wedge i)r)f'$

$ \qquad  - \overline{\lambda}(\bar a_{s+1}\wedge 1_{X\wedge M})(1_{E_{s+2}}\wedge u'(u\wedge 1_K))f'(1_X\wedge ij)$

$\qquad + (\bar c_s\wedge 1_{X\wedge M})(1_{KG_s}\wedge u_8(w_5\wedge 1_M)(1_W\wedge (1_Y\wedge i)r)g$

$\qquad - (\bar c_s\wedge 1_{X\wedge M})(1_{KG_s}\wedge u_8(w_5\wedge 1_M)(1_W\wedge (r\wedge 1_M)\overline{m}_K))g(1_X\wedge ij)$

$\qquad + \overline{\lambda}(1_{E_{s+1}}\wedge u'(u\wedge 1_K))f_8(1_X\wedge i')$\newline
where we use $u_8(w_5\wedge 1_M)(1_W\wedge (r\wedge 1_M)\overline{m}_K) = \overline{\lambda} u'(u\wedge 1_K)$,
for some nonzero $\overline{\lambda}\in Z_p$.  Moreover, by (9.4.20)(9.4.6), $(\bar b_{s+1}\wedge 1_{L\wedge K})(1_{E_{s+1}}\wedge u\wedge 1_K)f_8(1_X\wedge i')
= (1_{KG_{s+1}}\wedge (u\wedge 1_K)\mu_{X\wedge M})(g_0\wedge 1_M) = (1_{KG_{s+1}}\wedge
\overline{\mu}_2(1_Y\wedge i')\psi_{X\wedge M})(g_0\wedge 1_M) = \lambda_1(1_{KG_{s+1}}\wedge
\overline{\mu}_2(1_Y\wedge i')\psi_{X\wedge M})(h_0\sigma\wedge 1_{X\wedge M}) =
\lambda_1 (h_0\sigma\wedge 1_{L\wedge K})\overline{\mu}_2(1_Y\wedge i')\psi_{X\wedge M}$
(mod $d_1$-boundary).  Then $[(\bar b_{s+1}\wedge 1_{L\wedge K})(1_{E_{s+1}}\wedge u\wedge 1_K)
f_8(1_X\wedge i')] = \lambda_1(\phi\wedge 1_K)_*(j''\wedge 1_K)_*[(\sigma\wedge 1_{L\wedge K})
\overline{\mu}_2(1_Y\wedge i')\psi_{X\wedge M}] = \lambda_1(\phi\wedge 1_K)_*(j''\wedge 1_K)_*(\overline{\mu}_2(1_Y\wedge i'))_*
(\psi_{X\wedge M})_*[\sigma\wedge 1_{X\wedge M}]
= \lambda_1(\phi\wedge 1_K)_*(i')_*(m_M(\overline{u}\wedge 1_M))_*(\psi_{X\wedge M})_*[\sigma\wedge 1_{X\wedge M}]
= \lambda_1((1_L\wedge i')(\phi\wedge 1_M))_*(m_M(\widetilde{\psi}\wedge 1_M)_*[\sigma\wedge 1_{X\wedge M}]$ = 0 $\in
Ext_A^{s+1,tq}(H^*L\wedge K,H^*X\wedge M)$.  That is we have $(\bar b_{s+1}\wedge 1_{L\wedge K})(1_{E_{s+1}}\wedge u\wedge 1_K)f_8(1_X\wedge i')
= (\bar b_{s+1}\bar c_s\wedge 1_{L\wedge K})g_3$ with $g_3\in [\Sigma^{tq}X\wedge M, KG_s\wedge L\wedge K]
\cap (ker d)$ ( cf. Prop. 9.5.26) and so $(1_{E_{s+1}}\wedge u\wedge 1_K)f_8(1_X\wedge i')
= (\bar c_s\wedge 1_{L\wedge K})g_3 + (\bar a_{s+1}\wedge 1_{L\wedge K})f'_2$
with $f'_2\in [\Sigma^{tq+1}X\wedge M, E_{s+2}\wedge L\wedge K]\cap (ker d)$ (cf. Prop.  6.5.26). Hence, (9.4.26) becomes\\
{\bf (9.4.27)}\qquad $(\bar a_{s+1}\wedge 1_{X\wedge M})(1_{E_{s+2}}\wedge u_8(w_5\wedge 1_M))f_3$

$ = - (\bar a_{s+1}\wedge 1_{X\wedge M})(1_{E_{s+2}}\wedge u_8(w_5\wedge 1_M)(1_W\wedge (1_Y\wedge i)r))f'$

$\quad - \overline{\lambda}(\bar a_{s+1}\wedge 1_{X\wedge M})(1_{E_{s+2}}\wedge u'(u\wedge 1_K))f'(1_X\wedge ij)$

$\quad + (\bar c_s\wedge 1_{X\wedge M})(1_{KG_s}\wedge u_8(w_5\wedge 1_M)(1_W\wedge (1_Y\wedge i)r))g$

$\quad - (\bar c_s\wedge 1_{X\wedge M})(1_{KG_s}\wedge u_8(w_5\wedge 1_M)(1_W\wedge (r\wedge 1_M)\overline{m}_K))g(1_X\wedge ij)$

$\quad + \overline{\lambda}(\bar c_s\wedge 1_{X\wedge M})(1_{KG_s}\wedge u')g_3 + \overline{\lambda}(\bar a_{s+1}\wedge 1_{X\wedge M})(1_{E_{s+2}}\wedge u')f'_2$

 By (9.4.27), $(1_{KG_s}\wedge u_8(w_5\wedge 1_M)(1_W\wedge (1_Y\wedge i)r)g -
(1_{KG_s}\wedge u_8(w_5\wedge 1_M)(1_W\wedge (r\wedge 1_M)\overline{m}_K)g(1_X\wedge ij)
+ \overline{\lambda}(1_{KG_s}\wedge u')g_3\in [\Sigma^{tq}X\wedge M, KG_s\wedge X\wedge M]$
is a  $d_1$-cycle which represents an element in $Ext_A^{s,tq}(H^*X\wedge M,H^*X\wedge M)\cong Z_p\{[\sigma\wedge 1_{X\wedge M}]\}$
( cf. Lemma 9.4.13).  Then we have\\
{\bf (9.4.28)}\qquad $(1_{KG_s}\wedge u_8(w_5\wedge 1_M)(1_W\wedge (1_Y\wedge i)r)g
+ \overline{\lambda}(1_{KG_s}\wedge u')g_3$

$\quad - (1_{KG_s}\wedge u_8(w_5\wedge 1_M)(1_W\wedge (r\wedge 1_M)\overline{m}_K))g(1_X\wedge ij)$

$\qquad = \bar\lambda_0(\sigma\wedge 1_{X\wedge M})$
\quad  (mod $d_1$-boundary).\\  Now we consider the cases of  $\bar\lambda_0 \neq 1$ or
$\bar\lambda_0 = 1$ separately .

If $\bar\lambda_0\neq 1$, then by (9.4.27) and $\bar c_s\cdot\sigma = \bar a_{s+1}\cdot\kappa$ we have

$\qquad(1_{E_{s+2}}\wedge u_8(w_5\wedge 1_M))f_3 = - (1_{E_{s+2}}\wedge u_8(w_5\wedge 1_M)(1_W\wedge (1_Y\wedge i)r)f'$

$\qquad\quad - \bar \lambda (1_{E_{s+2}}\wedge u'(u\wedge 1_K))f'(1_X\wedge ij) + \bar\lambda
(1_{E_{s+2}}\wedge u')f'_2$

$\qquad\quad + \bar\lambda_0 (\kappa\wedge 1_{X\wedge M}) + (\bar c_{s+1}\wedge 1_{X\wedge M})g_4$\newline
with $g_4\in [\Sigma^{tq+1}X\wedge M, KG_{s+1}\wedge X\wedge M]$ and  by composing
$(1_{E_{s+2}}\wedge 1_X\wedge\alpha) = (1_{E_{s+2}}\wedge\alpha_{X\wedge M}m_M(\widetilde{\psi}\wedge 1_M))$
we obtain that $(1_{E_{s+2}}\wedge 1_X\wedge\alpha)(\kappa\wedge 1_{X\wedge M}) = (1_{E_{s+2}}
\wedge\alpha_{X\wedge M}(j''u\wedge 1_M)(1_W\wedge m_M(\overline{u}\wedge 1_M))f_3
= (1_{E_{s+2}}\wedge\alpha_{X\wedge M}\cdot m_M(\widetilde{\psi}\wedge 1_M)u_8(w_5\wedge 1_M))f_3
= \bar\lambda_0(1_{E_{s+2}}\wedge 1_X\wedge\alpha)(\kappa\wedge 1_{X\wedge M})$
so that the result of the step 1 follows.

If $\bar\lambda_0 = 1$, then by composing $(1_{KG_s}\wedge m_M(\widetilde{\psi}\wedge 1_M))$ on (9.4.28)
and using (9.4.25) we have $(1_{KG_s}\wedge j'(j''u\wedge 1_K))g = (1_{KG_s}
\wedge m_M(\widetilde{\psi}\wedge 1_M)u_8(w_5\wedge 1_M)(1_W\wedge (1_Y\wedge i)r)g
= (\sigma\wedge 1_M)m_M(\widetilde{\psi}\wedge 1_M)$ (mod $d_1$-boundary).
Moreover, by composing $(1_{KG_s}\wedge j\widetilde{\psi}\wedge 1_M)$ on (9.4.28) and using
 (9.4.25) we have

$(1_{KG_s}\wedge j\widetilde{\psi}\wedge 1_M)(\sigma\wedge 1_{X\wedge M})$

$ =  (1_{KG_s}\wedge (j\widetilde{\psi}\wedge 1_M)u_8(w_5\wedge 1_M)(1_W\wedge (1_Y\wedge i)r)g$

$\quad - (1_{KG_s}\wedge (j\widetilde{\psi}\wedge 1_M)u_8(w_5\wedge 1_M)(1_W\wedge (r\wedge 1_M)\overline{m}_K))g(1_X\wedge ij)$

$\quad + \overline{\lambda}(1_{KG_s}\wedge (j\widetilde{\psi}\wedge 1_M)u')g_3$

$ = (1_{KG_s}\wedge ij(j'(j''\wedge 1_K)(u\wedge 1_K))g$

$\quad - (1_{KG_s}\wedge j'(j''u\wedge 1_K))g(1_X\wedge ij) + \bar \lambda (1_{KG_s}\wedge j'(j''\wedge 1_K))g_3$\quad
  by (9.4.25)

$ = (1_{KG_s}\wedge ij)(\sigma\wedge 1_M)m_M(\widetilde{\psi}\wedge 1_M)  - (\sigma\wedge 1_M)m_M(\widetilde{\psi}\wedge 1_M)(1_X\wedge ij)$

$\quad + \overline{\lambda}(1_{KG_s}\wedge j'(j''\wedge 1_K))g_3$

$ = (1_{KG_s}\wedge j\widetilde{\psi}\wedge 1_M)(\sigma\wedge 1_{X\wedge M})
+ \overline{\lambda}(1_{KG_s}\wedge j'(j''\wedge 1_K))g_3$\quad  by (9.4.6)\\
(mod $d_1$-boundary), then $(1_{KG_s}\wedge j'(j''\wedge 1_K))g_3$ = 0 and so
 $g_3 = (1_{KG_s}\wedge\overline{\mu}_2(1_Y\wedge i'))g_5$ (mod $d_1$-boundary)
for some $g_5\in [\Sigma^{tq+q+1}X\wedge M, KG_s\wedge Y\wedge M]$.  So, by (9.4.6)(9.4.20)
$(1_{KG_{s+1}}\wedge\overline{\mu}_2(1_Y\wedge i')\psi_{X\wedge M})(g_0\wedge 1_M) =
(1_{KG_{s+1}}\wedge (u\wedge 1_K)\mu_{X\wedge M})(g_0\wedge 1_M) =
(\bar b_{s+1}\wedge 1_{L\wedge K})(1_{E_{s+1}}\wedge u\wedge 1_K)f_8(1_X\wedge i') = (\bar b_{s+1}\bar c_s\wedge 1_{L\wedge K})g_3
= (\bar b_{s+1}\bar c_s\wedge 1_{L\wedge K})(1_{KG_s}\wedge\overline{\mu}_2(1_Y\wedge i'))g_5$
so that $(1_{KG_{s+1}}\wedge \psi_{X\wedge M})(g_0\wedge 1_M) = (\bar b_{s+1}\bar c_s\wedge 1_{Y\wedge M})g_5$
, this shows $\lambda_1(\psi_{X\wedge M})_*[h_0\sigma\wedge 1_{X\wedge M}] =
(\psi_{X\wedge M})_*[g_0\wedge 1_M] = 0 \in Ext_A^{s+1,tq+1}(H^*Y\wedge M,.H^*X\wedge M)$
and by Lemma 9.4.13(2) we have $\lambda_1$ = 0 .  Then $[g_0\wedge 1_M]$ = 0 and so
$(\bar b_{s+1}\wedge 1_{W\wedge K})f_8(1_X\wedge i') = (\bar b_{s+1}\bar c_s\wedge 1_{W\wedge K})g_6$
for some $g_6\in [\Sigma^{tq+q}X\wedge M, KG_s\wedge W\wedge K]$  and $f_8(1_X\wedge i')
= (\bar c_s\wedge 1_{W\wedge K})g_6 + (\bar a_{s+1}\wedge 1_{W\wedge K})f'_3$
with $f'_3\in [\Sigma^{tq+q+1}X\wedge M,E_{s+2}\wedge W\wedge K]$.   Then, by composing
$(1_{E_{s+1}}\wedge w_5\wedge 1_M)$ on (9.4.17) we have

$- (\bar a_{s+1}\wedge 1_{P\wedge M})(1_{E_{s+2}}\wedge w_5\wedge 1_M)f_3$

$ = (\bar a_{s+1}\wedge 1_{P\wedge M})(1_{E_{s+2}}\wedge (w_5\wedge 1_M)(1_W\wedge (1_Y\wedge i)r)f'$

$\quad + (\bar a_{s+1}\wedge 1_{P\wedge M})(1_{E_{s+2}}\wedge (w_5\wedge 1_M)(1_W\wedge (r\wedge 1_M)\overline{m}_K)f'(1_X\wedge ij)$

$\quad + (\bar a_{s+1}\wedge 1_{P\wedge M})(1_{E_{s+2}}\wedge (w_5\wedge 1_M)(1_W\wedge (r\wedge 1_M)\overline{m}_K)f'_3
+ (\bar c_s\wedge 1_{P\wedge M})g_7$\\
where the $d_1$-cycle $g_7 = (1_{KG_s}\wedge (w_5\wedge 1_M)(1_W\wedge (1_Y\wedge i)r)g
- (1_{KG_s}\wedge (w_5\wedge 1_M)(1_W\wedge (r\wedge 1_M)\overline{m}_K)g(1_X\wedge ij)
+ (1_{KG_s}\wedge (w_5\wedge 1_M)(1_W\wedge (r\wedge 1_M)\overline{m}_K)g_6\in
[\Sigma^{tq+q+1}X\wedge M, KG_s\wedge P\wedge M]$ which represents an element in
$Ext_A^{s,tq+q+1}(H^*P\wedge M,H^*X\wedge M)$.  However, this group is zero
, this can be obtained by the following exact sequence

\quad

$0 = Ext_A^{s,tq+q}(H^*W\wedge K,H^*X\wedge M)\stackrel{((w_5\wedge 1_M)(1_W\wedge (r\wedge 1_M)\overline{m}_K)_*}
{\longrightarrow} $

$\qquad Ext_A^{s,tq+q+1}(H^*P\wedge M,H^*X\wedge M)\stackrel{(u_7)_*}{\longrightarrow} $

$\qquad Ext_A^{s,tq-2q}(H^*M,H^*X\wedge M)\stackrel{((1_W\wedge i')(\phi_W\wedge M))_*}{\longrightarrow}$\\
induced by (9.4.23), where the left group is zero by Prop. 9.4.11(4) and by Prop. 9.4.11(1)
the right group has unique generator
 $m_M(\widetilde{\psi}\wedge 1_M)^*(\tilde{\sigma})$,
which satisfies $((1_W\wedge i')(\phi_W\wedge 1_M))_*m_M(\widetilde{\psi}\wedge 1_M)^*(\tilde{\sigma})\neq 0\in Ext_A^{s+1,tq+q}
(H^*W\wedge K,H^*X\wedge M)$.

Then, $(\bar c_s\wedge 1_{P\wedge M})g_7$ = 0  and so $- (1_{E_{s+2}}\wedge w_5\wedge 1_M)f_3 = (1_{E_{s+2}}\wedge (w_5\wedge 1_M)u_8
(1_W\wedge (1_Y\wedge i)r)f' - (1_{E_{s+2}}\wedge (w_5\wedge 1_M)(1_W\wedge (r\wedge 1_M)\overline{m}_K)f'(1_X\wedge ij)
+ (1_{E_{s+2}}\wedge (w_5\wedge 1_M)(1_W\wedge (r\wedge 1_M)\overline{m}_K))f'_3
+ (\bar c_{s+1}\wedge 1_{P\wedge M})g_8$ for some $g_8\in [\Sigma^{tq+q+2}X\wedge M, KG_{s+1}
\wedge P\wedge M]$.  By composing $(1_{E_{s+2}}\wedge\alpha_{X\wedge M}\cdot u_7)$  we have
$(1_{E_{s+2}}\wedge 1_X\wedge\alpha)(\kappa\wedge 1_{X\wedge M}) = (1_{E_{s+2}}
\wedge\alpha_{X\wedge M}(j''u\wedge 1_M)(1_W\wedge m_M(\overline{u}\wedge 1_M))f_3
= (1_{E_{s+2}}\wedge\alpha_{X\wedge M}\cdot u_7(w_5\wedge 1_M))f_3$ = 0
. This shows the result of step 1.

\quad

{\bf  Step 2}\qquad  To prove \quad $(\bar c_{s+1}\wedge 1_M)\widetilde{h_0\sigma} =
(\kappa\wedge 1_M)\alpha $ = 0.

\quad

 By (9.4.3)(9.4.4), $\mu_{X\wedge M}(1_X\wedge\alpha i) = \mu_{X\wedge M}
\alpha_{X\wedge M}\widetilde{\psi}$ = 0
and so by (9.1.15) $\mu_{X\wedge M} = \mu_{X\wedge K'}(1_X\wedge v)$ ,  where $\mu_{X\wedge K'}
\in [X\wedge K', W\wedge K]$.
 We claim that $X\wedge K'$ splits into $W\wedge K\vee \Sigma^qY$,  that is, there is
a split cofibration $\Sigma^qY\rightarrow X\wedge K'\rightarrow W\wedge K$, this can be seen by
the following homotopy commutative diagram
of $3\times 3$-Lemma and using $(1_Y\wedge j)\alpha_{Y\wedge M} j' =
 r (1_K\wedge\alpha_1)$

\quad

\quad

$\qquad\qquad X\wedge M\quad\stackrel{\mu_{X\wedge M}}{\longrightarrow}\quad W\wedge K\quad
\stackrel{0}{\longrightarrow}\quad\Sigma^{q+1}Y$

$\qquad\qquad\quad\searrow 1_X\wedge v\nearrow \mu_{X\wedge K'}
\searrow ^{j'(j''u\wedge 1_K)}\nearrow _{(1_Y\wedge j)\alpha_{Y\wedge M}}$

$\qquad\qquad\qquad X\wedge K'\qquad\qquad\qquad \Sigma^{3q+1}M$

$\qquad\qquad\quad\nearrow \widetilde{\tau}_{X\wedge K'}\searrow 1_X\wedge
\quad\nearrow\widetilde{\psi}\quad\searrow\alpha_{X\wedge M}$

$\qquad\qquad\Sigma^qY\quad \stackrel{\tilde{u}w_2}{\longrightarrow}\qquad\Sigma^{q+1}X\quad
\stackrel{1_X\wedge\alpha i}{\longrightarrow}\quad\Sigma X\wedge M$\\
Hence, there is a split cofibration $\Sigma^qY\stackrel{\tau_{X\wedge K'}}{\longrightarrow}
 X\wedge K'\stackrel{\mu_{X\wedge K'}}{\longrightarrow} W\wedge K$ and so there are
 $\nu_{X\wedge K'} :  X\wedge K'\rightarrow\Sigma^q Y$ and $\widetilde{\nu}_{X\wedge K'}
: W\wedge K\rightarrow X\wedge K'$ such that $\nu_{X\wedge K'}\cdot \tau_{X\wedge K'} = 1_Y,\quad
\mu_{X\wedge K'}\cdot \widetilde{\nu}_{X\wedge K'} = 1_{W\wedge K},$
$ \widetilde{\tau}_{X\wedge K'}\cdot \nu_{X\wedge K'} +
\widetilde{\nu}_{X\wedge K'}\cdot \mu_{X\wedge K'} = 1_{X\wedge K'}$.

By the result of step 1 we have $(\kappa\wedge 1_{M\wedge X\wedge K'})(\alpha\wedge 1_{X\wedge K'})$ = 0,
then $(\kappa\wedge 1_{M\wedge Y})(\alpha\wedge 1_{Y}) = (1_{E_{s+2}}\wedge 1_M\wedge\nu_{X\wedge K'})
(\kappa\wedge 1_{M\wedge X\wedge K'})(\alpha\wedge 1_{X\wedge K'})(1_M\wedge\tau_{X\wedge K'})$ = 0.
Moreover, by using the splitness in (9.1.32) we have $(\bar c_{s+1}\wedge 1_M)\widetilde{h_0\sigma} =
(\kappa\wedge 1_M)\alpha = (1_{E_{s+2}}\wedge 1_M\wedge\widetilde{\nu})(\kappa\wedge 1_{M\wedge Y\wedge K'})(\alpha\wedge 1_{Y\wedge K'})
(1_M\wedge\widetilde{\tau})$ = 0  which shows the main Theorem C.  Q.E.D.

\vspace{2mm}

{\bf  Remark.}\quad In the proof of the main Theorem C, We only use the supposition (II) for our geometric input
to obtain that
 $(1_{E_{s+2}}\wedge\phi\wedge 1_M)(\kappa\wedge 1_M)m_M(\widetilde{\psi}\wedge 1_M)$ = 0.
 Then , the geometric supposition (II) of the main Theorem C can be weakened
 to be the supposition on
$m_M(\widetilde{\psi}\wedge 1_M)^*(\phi\wedge 1_M)_*(\tilde{\sigma})\in Ext_A^{s+1,tq}
(H^*L\wedge M,H^*X\wedge M)$ is a permanent cycle in the ASS.

\vspace{2mm}

Using some new cofibrations in this section , we also can give an alternative proof of
Theorem 9.3.9( and so the main Theorem B). We first do some preminalaries.

Since $\alpha'\alpha'i' = 0$, then by (9.1.23), there exists
$\alpha''_{Y\wedge M}\in [\Sigma^{q-2}Y\wedge M,K]$  such that
$\alpha''_{Y\wedge M}(r\wedge 1_M)\overline{m}_K = \alpha'$.  By
applying the derivation $d$, $d(\alpha''_{Y\wedge M})(r\wedge
1_M)\overline{m}_K = - d(\alpha') = 0$ and so $d(\alpha''_{Y\wedge
M})\in (m_M(\overline{u}\wedge 1_M))^*[\Sigma^{2q}M,K]$ = 0.
$\alpha_{Y\wedge M}(1_Y\wedge i)r\in [\Sigma^{q-2}K,K] \cong
Z_p\{\alpha''\}$ and so $\alpha''_{Y\wedge M}(1_Y\wedge i)r =
\lambda\alpha''$ for some $\lambda\in Z_p$. Note that
$d((1_Y\wedge i)r) = (r\wedge 1_M)d(1_K\wedge i) = (r\wedge
1_M)\overline{m}_K$, then by applying the derivation $d$, we have
$\alpha' = \alpha''_{Y\wedge M} (r\wedge M)\overline{m}_K =
\lambda d(\alpha'') = - \lambda \alpha'$ and so $\lambda = - 1$.
By (9.1.8), $\bar h i'' = \overline{w}$, $r i' = \overline{w}\cdot
j$( up to sign), then $(r\wedge 1_M)\overline{m}_K i' = - (r
i'\wedge 1_M)\overline{m}_M = \pm (\overline{w}\wedge 1_M) = \pm
(\bar h i''\wedge 1_M)$ and so $\alpha''_{Y\wedge M}(\bar h
i''\wedge 1_M) = \lambda_0\alpha''_{Y\wedge M} (r\wedge
1_M)\overline{m}_K i' = \lambda_0 \alpha'i' = \lambda_0 i'
((\alpha_1)_L i''\wedge 1_M)$ and we have $\alpha''_{Y\wedge M}
(\bar h\wedge 1_M) = \lambda_0i'((\alpha_1)_L\wedge 1_M)$, where
$\lambda_0 = \pm 1$. On the other hand, $i'((\alpha_1)_L\wedge
1_M) (1_L\wedge j')(i''\wedge 1_K) = i'(\alpha_1\wedge 1_M)j' =
i'(ij\alpha - \alpha ij)j'$ = 0, then $i'((\alpha_1)_L\wedge 1_M)
(1_L\wedge j') = \lambda'\alpha'' (j''\wedge 1_K)$ with
$\lambda'\in Z_p$. By composing the map $\widetilde{\Delta}$ in
Theorem 6.5.18 we have $\lambda'\alpha'i'ijj' =
\lambda'\alpha''i'j' = \lambda'\alpha''(j''\wedge
1_K)\widetilde{\Delta} = i'((\alpha_1)_L\wedge 1_M) (1_L\wedge
j')\widetilde{\Delta} = - i'((\alpha_1)_L i''\wedge 1_M)ijj' = -
\alpha'i'ijj'$  so that $\lambda' = -1$.  Concludingly , there is
$\alpha''_{Y\wedge M}\in [\Sigma^{q-2}Y\wedge M, K]$ such that\\
{\bf (9.4.29)}\qquad $\alpha''_{Y\wedge M}(r\wedge 1_M)\overline{m}_K = \alpha',\quad\qquad \alpha''_{Y\wedge M}(1_Y\wedge i)r
= - \alpha'',$

$\qquad\qquad d(\alpha''_{Y\wedge M}) = 0, \quad\qquad \alpha''_{Y\wedge M}(\bar h\wedge 1_M) = \lambda_0 i'((\alpha_1)_L\wedge 1_M),$

$\qquad\qquad i'((\alpha_1)_L\wedge 1_M)(1_L\wedge j') = - \alpha''(j''\wedge 1_K)$\\
where $\lambda_0 = \pm 1$.

Note that  the cofibre of $\alpha''_{Y\wedge M} : \Sigma^{q-2}Y\wedge M\to K$ is $X\wedge M$ given by the cofibration\\
{\bf (9.4.30)}\qquad $\Sigma^{q-2}Y\wedge M\stackrel{\alpha''_{Y\wedge M}}{\longrightarrow} K\stackrel{u'(i''\wedge 1_K)}{\longrightarrow}
X\wedge M\stackrel{\psi_{X\wedge M}}{\longrightarrow}\Sigma^{q-1}Y\wedge M$\\
and the above map $\psi_{X\wedge M}\in [X\wedge M,\Sigma^{q-1}Y\wedge M]$ and
$u'\in [L\wedge K, X\wedge M]$ is just the map in (9.4.2) and
(9.4.6). This can be seen by the equation $m_M(\overline{u}\wedge 1_M)\psi_{X\wedge M}
= m_M(\widetilde{\psi}\wedge 1_M)$ in (9.4.6),(9.4.2) and the following homotopy commutative
diagram of $3\times 3$-Lemma

$\qquad\quad X\wedge M\quad\stackrel{m_M(\tilde{\psi}\wedge 1_M)}{\longrightarrow}\quad\Sigma^{2q}M\qquad\stackrel{\alpha'i'}
{\longrightarrow}\qquad\Sigma^{q+1}K$

$\qquad\qquad\quad\searrow\psi_{X\wedge M}\quad\nearrow ^{m_M(\overline{u}\wedge 1_M)}\quad\searrow ^{(\phi\wedge 1_K)i'}\quad
\nearrow j''\wedge 1_K\searrow ^{(r\wedge 1_M)\overline{m}_K}$\\
{\bf (9.4.31)}\quad\qquad$\Sigma^{q-1}Y\wedge M\qquad\qquad\qquad\quad\Sigma L\wedge K\qquad\qquad\qquad\Sigma^{q}Y\wedge M$

$\qquad\qquad\nearrow ^{(r\wedge 1_M)\overline{m}_K}\quad\searrow ^{\alpha''_{Y\wedge M}}\quad\nearrow i''\wedge 1_K\quad
\searrow u'\qquad\quad\nearrow \psi_{X\wedge M}$

$\qquad\quad\Sigma^qK\qquad\stackrel{\alpha'}{\longrightarrow}\qquad\quad\Sigma K\qquad\stackrel{u'(i''\wedge 1_K)}{\longrightarrow}
\qquad\Sigma X\wedge M$\\
and by this we have the following relation\\
{\bf (9.4.32)}\qquad $\psi_{X\wedge M}u' = (r\wedge 1_M)\overline{m}_K(j''\wedge 1_K)$.

{\bf Proposition 9.4.33}\quad  Let $p\geq 5$ and $V$ be any spectrum,
then for any map $f\in [\Sigma^*K, V\wedge K]$ we have
$(1_V\wedge\alpha')d(f) = d(f)\alpha' = 0$.

{\bf  Proof}:\quad  By (6.5.12), $\alpha\wedge 1_K =
\overline{m}_K'\alpha'm_K'$, where $m_K' = m_KT : M\wedge K\to K$,
$\overline{m}_K' = T\overline{m}_K : \Sigma K\to M\wedge K$.
$d(f)\alpha'm_K' = (1_V\wedge m_K')(T'\wedge 1_K)(1_M\wedge
f)\overline{m}_K'\alpha' m_K' = (1_V\wedge m_K')(T'\wedge
1_K)(1_M\wedge f)(\alpha\wedge 1_K) = (1_V\wedge m_K')(T'\wedge
1_K)(\alpha\wedge 1_V\wedge 1_K) (1_M\wedge f) = (1_V\wedge
m_K'(\alpha\wedge 1_K))(T'\wedge 1_K)(1_M\wedge f)$ = 0,  where
$T' : M\wedge V\to V\wedge M$is the switching map. Q.E.D.

{\bf Proposition 9.4.34}\quad Under the supposition (I) of the main Theorem B we have

(1) $Ext_A^{s,tq-1}(H^*K,H^*K)$ = 0.

(2) $Ext_A^{s,tq}(H^*Y\wedge M,H^*K)$ has unique generator $(1_Y\wedge i)_*r_*[\sigma\wedge 1_K]$.

{\bf  Proof}:\quad (1) Consider the following exact sequence

$\qquad Ext_A^{s,tq+q}(H^*M,H^*M)\stackrel{(i')_*}{\longrightarrow} Ext_A^{s,tq+q}(H^*K,H^*M)$

$\qquad\qquad\stackrel{(j')_*}{\longrightarrow} Ext_A^{s,tq-1}(H^*M,H^*M)\stackrel{\alpha_*}{\longrightarrow}$\\
induced by (9.1.2). By the supposition (I), the right group has
unique generator $j^*i_*(\sigma)$  which satisfies
$\alpha_*j^*i_*(\sigma) = (ij)_*\alpha_*(\tilde{\sigma})\neq 0$.
Then im $(j')_*$ = 0. By the supposition (I), the left group is
zero or has two generators $(ij)_*\alpha_*(\widetilde{\tau}'),
(ij)^*\alpha_*(\widetilde{\tau}')$ (this can be obtained by a
similar proof as given in Prop. 9.3.1(2)), then $
Ext_A^{s,tq+q}(H^*K,H^*M)
 = (i')_*Ext_A^{s,tq+q}(H^*M,H^*M)$ is zero or has unique generator $(i')_*(\alpha_1\wedge 1_M)_*
(\widetilde{\tau}')$. Look at the following exact sequence

$\qquad Ext_A^{s,tq+q}(H^*K,H^*M)\stackrel{(j')^*}{\longrightarrow} Ext_A^{s,tq-1}(H^*K,H^*K)$

$\qquad\qquad \stackrel{(i')^*}{\longrightarrow} Ext_A^{s,tq-1}(H^*K,H^*M)\stackrel{\alpha^*}{\longrightarrow}$\\
induced by (9.1.2). By the supposition (I), the right group has
unique generator $j^*(i'i)_*(\sigma)$ which satisfies
$\alpha^*j^*(i'i)_*(\sigma) = (i')_*(ij)_*\alpha_*(\tilde{\sigma})
\neq 0\in Ext_A^{s+1,tq+q}(H^*K,H^*M)$  so that im $(i')^*$ = 0.
The left group is zero or has unique generator
$(i')_*(\alpha_1\wedge 1_M)_*(\widetilde{\tau}')$
 and so im $(j')^*$ = 0 and the middle group is zero as desired.

(2)  For any $g\in Ext_A^{s,tq}(H^*Y\wedge
M,H^*K)$,$m_M(\overline{u}\wedge 1_M))_*(g)\in
Ext_A^{s,tq-q-1}\\(H^*M,H^*K)\cong Z_p
\{(j')^*(\tilde{\sigma})\}$, this can be obtained from
$Ext_A^{s,tq}(H^*M,H^*M)\\\cong Z_p\{\tilde{\sigma}\}$ in Prop.
9.3.0(2) and $Ext_A^{s,tq-q-1}(H^*M,H^*M)$ = 0, where the last is
obtained by the supposition (I) on $Ext_A^{s,tq-q+u}\\(Z_p,Z_p)$ =
0(for $ u = 0,-1,1$). Then $(m_M(\overline{u}\wedge 1_M))_*(g) =
\lambda'(j')^*[\sigma\wedge 1_M] = \lambda'[(\sigma\wedge 1_M)j']
= \lambda'[(1_{KG_s}\wedge j')(\sigma\wedge 1_K)] = \lambda'
(j')_*[\sigma\wedge 1_K] = \lambda' (m_M(\overline{u}\wedge
1_M))_*(1_Y\wedge i)_*r_*[\sigma\wedge 1_K]$ and so $g = \lambda'
(1_Y\wedge i)_*r_*[\sigma\wedge 1_K]$ ( with $\lambda'\in Z_p$)
modulo $((r\wedge 1_M)\overline{m}_K)_*Ext_A^{s,tq-1} (H^*K,H^*K)$
= 0.  Q.E.D.

\vspace{2mm}

{\bf An alternative proof of Theorem 9.3.9}: \quad By the supposition (II) of the main Theorem B
we have $(1_{E_{s+2}}\wedge\alpha')(\kappa\wedge 1_K) = (\bar c_{s+1}\wedge 1_K)
(h_0\sigma\wedge 1_K)$ = 0, then $(\kappa\wedge 1_K) = (1_{E_{s+2}}\wedge j''\wedge 1_K)f$ and we have $d((1_{E_{s+2}}\wedge j''
\wedge 1_K)f)$ = 0.  That is we have\\
{\bf (9.4.35)}\qquad $\kappa\wedge 1_K = (1_{E_{s+2}}\wedge j''\wedge 1_K)f$\qquad\quad $d(f) = (1_{E_{s+2}}\wedge i''\wedge 1_K)f'$\\
for some $f\in [\Sigma^{tq+q+1}K, E_{s+2}\wedge L\wedge K], f'\in [\Sigma^{tq+q+2}K, E_{s+2}\wedge K]$.

By (9.4.29)(9.4.35),$(1_{E_{s+2}}\wedge\alpha''_{Y\wedge M}(\bar h\wedge 1_M)(1_L\wedge j'))f = \lambda_0(1_{E_{s+2}}
\wedge i'((\alpha_1)_L\wedge 1_M)(1_L\wedge j'))f = - \lambda_0(1_{E_{s+2}}\wedge\alpha''(j''\wedge 1_K))f = - \lambda_0
(1_{E_{s+2}}\wedge\alpha'')(\kappa\wedge 1_K)$, where $\lambda_0 = \pm 1$.  That is we have\\
{\bf (9.4.36)}\qquad $(1_{E_{s+2}}\wedge\alpha''_{Y\wedge M}(\bar h\wedge 1_M)(1_L\wedge j'))f = - \lambda_0(1_{E_{s+2}}
\wedge\alpha'')(\kappa\wedge 1_K)$\qquad\quad  where $\lambda_0 = \pm 1$

It follows that $(\bar a_{s+1}\wedge 1_K)(1_{E_{s+2}}\wedge\alpha''_{Y\wedge M}(\bar h\wedge 1_M)(1_L\wedge j'))f = - \lambda_0
(\bar a_{s+1}\wedge 1_K)(1_{E_{s+2}}\wedge\alpha'')(\kappa\wedge 1_K) = - \lambda_0(\bar c_s\wedge 1_K)(1_{KG_s}\wedge\alpha'')
(\sigma\wedge 1_K)$ = 0 since $\alpha''$ induces zero homomorphism in $Z_p$-cohomology.
  Then,  by (9.4.30),$(\bar a_{s+1}\wedge 1_{Y\wedge M})
(1_{E_{s+2}}\wedge (\bar h\wedge 1_M)(1_L\wedge j'))f = (1_{E_{s+1}}\wedge\psi_{X\wedge M})f_2$ with $f_2\in [\Sigma^{tq+q-1}K,\\
E_{s+1}\wedge X\wedge M]$. Consequently, $(\bar b_{s+1}\wedge 1_{Y\wedge M})(1_{E_{s+1}}\wedge\psi_{X\wedge M})f_2$ = 0 and so by (9.4.30) we have
$(\bar b_{s+1}\wedge 1_{X\wedge M})f_2 = (1_{KG_{s+1}}\wedge u'(i''\wedge 1_K))g$, with $d_1$-cycle $g\in [\Sigma^{tq+q-1}K,KG_{s+1}
\wedge K]$ and this $d_1$-cycle represents an element in $Ext_A^{s+1,tq+q-1}(H^*K,H^*K)\cong Z_p\{(h_0\sigma)''\}$. Then $[g] = \lambda'(h_0\sigma)''
= \lambda'(\alpha'')_*[\sigma\wedge 1_K]$ for some $\lambda'\in Z_p$ so that

$\quad [(\bar b_{s+1}\wedge 1_{X\wedge M})f_2] = (u'(i''\wedge 1_K))_*[g] = \lambda'(u'(i''\wedge 1_K))_*(\alpha'')_*
[\sigma\wedge 1_K]$

$ = \lambda'(u'(i''\wedge 1_K))_*(\alpha''_{Y\wedge M})_*
((1_Y\wedge i)r)_*[\sigma \wedge 1_K]$ = 0\\
 Hence, $(\bar b_{s+1}\wedge 1_{X\wedge M})f_2 = (\bar b_{s+1}\bar c_s\wedge 1_{X\wedge M})g_2$ for some $g_2\in [\Sigma^{tq+q-1}K,
KG_s\wedge X\wedge M]$ and so $f_2 = (\bar c_s\wedge 1_{X\wedge
M})g_2 + (\bar a_{s+1}\wedge 1_{X\wedge M})f_3$ with $f_3\in
[\Sigma^{tq+q-1}K,\\E_{s+2}\wedge X\wedge M]$ and we have

$(\bar a_{s+1}\wedge 1_{Y\wedge M})(1_{E_{s+2}}\wedge (\bar h\wedge 1_M)(1_L\wedge j'))f$

$ = (\bar a_{s+1}\wedge 1_{Y\wedge M})(1_{E_{s+2}}\wedge\psi_{X\wedge M})f_3 + (\bar c_s\wedge 1_{Y\wedge M})(1_{KG_s}
\wedge\psi_{X\wedge M})g_2$

$ = (\bar a_{s+1}\wedge 1_{Y\wedge M})(1_{E_{s+2}}\wedge\psi_{X\wedge M})f_3 + \lambda (\bar c_s\wedge 1_{Y\wedge M})
(1_{KG_s}\wedge (1_Y\wedge i)r)(\sigma\wedge 1_K)$

$ = (\bar a_{s+1}\wedge 1_{Y\wedge M})(1_{E_{s+2}}\wedge\psi_{X\wedge M})f_3 $

$ \qquad + \lambda (\bar a_{s+1}\wedge 1_{Y\wedge M})
(1_{E_{s+2}}\wedge (1_Y\wedge i)r)(\kappa\wedge 1_K)$\\
with $\lambda\in Z_p$, where the $d_1$-cycle $(1_{KG_s}\wedge\psi_{X\wedge M})g_2\in [\Sigma^{tq}K,KG_s\wedge Y\wedge M]$ represents an element
$\lambda ((1_Y\wedge i)r)_*[\sigma\wedge 1_K]\in Ext_A^{s,tq}(H^*Y\wedge M,H^*K)$( cf. Prop. 9.4.34(2)) and so it equals to $\lambda
(1_{KG_s}\wedge (1_Y\wedge i)r)(\sigma\wedge 1_K)$ ( mod $d_1$-boundary). Then we have
$(1_{E_{s+2}}\wedge (\bar h\wedge 1_M)(1_L\wedge j'))f
= (1_{E_{s+2}}\wedge\psi_{X\wedge M})f_3 + \lambda(1_{E_{s+2}}\wedge (1_Y\wedge i)r)(\kappa\wedge 1_K) + (\bar c_{s+1}\wedge
1_{Y\wedge M})g_3$ for some $g_3\in [\Sigma^{tq+1}K, KG_{s+1}\wedge Y\wedge M]$.  By composing $1_{E_{s+2}}\wedge 1_Y\wedge\alpha$ and using(9.4.8) we have $(1_{E_{s+2}}\wedge\omega\wedge 1_M)f_3 = (1_{E_{s+2}}\wedge\alpha_{Y\wedge M}m_M(\widetilde{\psi}\wedge 1_M))f_3
= (1_{E_{s+2}}\wedge (1_Y\wedge\alpha)\psi_{X\wedge M})f_3 = - \lambda (1_{E_{s+2}}\wedge (1_Y\wedge\alpha i)r)(\kappa\wedge 1_K)
= - \lambda (1_{E_{s+2}}\wedge (r\wedge 1_M)\overline{m}_K\alpha')(\kappa\wedge 1_K)$ = 0 and by (9.4.7) we have $f_3 = (1_{E_{s+2}}
\wedge (1_x\wedge j)u'\wedge 1_M)f_4$ with $f_4\in [\Sigma^{tq+q+1}K, E_{s+2}\wedge L\wedge K\wedge M]$. That is we have\\
{\bf (9.4.37)}\qquad $(1_{E_{s+2}}\wedge (\bar h\wedge 1_M)(1_L\wedge j'))f = (1_{E_{s+2}}\wedge\psi_{X\wedge M}((1_X\wedge j)u'
\wedge 1_M))f_4$

$\qquad + \lambda (1_{E_{s+2}}\wedge (1_Y\wedge i)r)(\kappa\wedge 1_K) + (\bar c_{s+1}\wedge 1_{Y\wedge M})g_3$

$ = (1_{E_{s+2}}\wedge\psi_{X\wedge M}u')f_5 + (1_{E_{s+2}}\wedge\psi_{X\wedge M}(1_X\wedge ij)u'(1_L\wedge m_K))f_4$

$\qquad + \lambda (1_{E_{s+2}}\wedge (1_Y\wedge i)r)(\kappa\wedge 1_K) + (\bar c_{s+1}\wedge 1_{Y\wedge M})g_3$\\
where we use $f_4 = (1_{E_{s+2}}\wedge (1_L\wedge\overline{m}_K)(1_{L\wedge K}\wedge j))f_4 + (1_{E_{s+2}}\wedge (1_{L\wedge K}
\wedge i)(1_L\wedge m_K))f_4$ and write $(1_{E_{s+2}}\wedge 1_{L\wedge K}\wedge j)f_4 = f_5$.

By composing $1_{E_{s+2}}\wedge\alpha''_{Y\wedge M}$ on (9.4.37) and using (9.4.36)(9.4.30) we have
$- \lambda_0(1_{E_{s+2}}\\\wedge\alpha'')(\kappa\wedge 1_K)
= (1_{E_{s+2}}\wedge\alpha''_{Y\wedge M}(\bar h\wedge 1_M)(1_L\wedge j'))f = \lambda (1_{E_{s+2}}\wedge\alpha''_{Y\wedge M}
(1_Y\wedge i)r)(\kappa\wedge 1_K) = - \lambda (1_{E_{s+2}}\wedge \alpha'')(\kappa\wedge 1_K)$.  If $\lambda\neq \lambda_0$,
then $(1_{E_{s+2}}\wedge\alpha'')(\kappa\wedge 1_K)$ = 0 and the Theorem follows.  So, we suppose that
$\lambda = \lambda_0$.

By (9.1.8) we have $\overline{u}\bar h = i\cdot j''$ so that
$m_M(\overline{u}\wedge 1_M)(\bar h\wedge 1_M) = j''\wedge 1_M$( up to sign).
Then, what happen is either $m_M(\overline{u}\wedge 1_M)(\bar h\wedge 1_M) =
\lambda_0(j''\wedge 1_M)$ or $m_M(\overline{u}\wedge 1_M)
(\bar h\wedge 1_M) = - \lambda_0(j''\wedge 1_M)$. Now we consider
this two cases separately.

{\bf  Case 1}\qquad\qquad $m_M(\overline{u}\wedge 1_M)(\bar h\wedge 1_M) = \lambda_0(j''\wedge 1_M)$.

In this case, by composing $1_{E_{s+2}}\wedge m_M(\overline{u}\wedge 1_M)$ on (9.4.37) we have\\
{\bf (9.4.38)}\qquad $\lambda_0(1_{E_{s+2}}\wedge j')(\kappa\wedge 1_K) = \lambda_0(1_{E_{s+2}}\wedge j'(j''\wedge 1_K))f$

$ = \lambda_0(1_{E_{s+2}}\wedge (j''\wedge 1_M)(1_L\wedge j'))f = (1_{E_{s+2}}\wedge m_M(\overline{u}\wedge 1_M)(\bar h\wedge 1_M)
(1_L\wedge j'))f$

$ = (1_{E_{s+2}}\wedge m_M(\overline{u}\wedge 1_M)\psi_{X\wedge M}(1_X\wedge ij)u'(1_L\wedge m_K))f_4$

$\qquad + \lambda_0(1_{E_{s+2}}\wedge j')(\kappa\wedge 1_K) + (\bar c_{s+1}\wedge 1_M)(1_{KG_{s+1}}\wedge m_M(\overline{u}\wedge 1_M))g_3$

$ = - (1_{E_{s+2}}\wedge j'(j''\wedge 1_K)(1_L\wedge m_K))f_4$

$\qquad + \lambda_0(1_{E_{s+2}}\wedge j')(\kappa\wedge 1_K) + (\bar c_{s+1}\wedge 1_M)(1_{KG_{s+1}}\wedge m_M(\overline{u}\wedge 1_M))g_3$\\
and so $(1_{E_{s+2}}\wedge j'(j''\wedge 1_K)(1_L\wedge m_K))f_4 = (\bar c_{s+1}\wedge 1_M)(1_{KG_{s+1}}\wedge m_M(\overline{u}\wedge 1_M))g_3$,
where we use $m_M(\overline{u}\wedge 1_M)\psi_{X\wedge M}(1_X\wedge ij)u' = m_M(\widetilde{\psi}\wedge 1_M)(1_X\wedge ij)u'
= - (j\widetilde{\psi}\wedge 1_M)u' = - j'(j''\wedge 1_K)$ which is obtained by (9.4.6) and the right rectangle of the diagram (9.4.1). Moreover,
by applying the derivation $d$ to the equation (9.4.37) we have\\
{\bf (9.4.39)}\qquad $(1_{E_{s+2}}\wedge (\bar h\wedge 1_M)(1_L\wedge j')(i''\wedge 1_K))f' = (1_{E_{s+2}}\wedge\psi_{X\wedge M}u')d(f_5)$

$\quad + (1_{E_{s+2}}\wedge\psi_{X\wedge M}(1_X\wedge ij)u')d((1_{E_{s+2}}\wedge 1_L\wedge m_K)f_4)$

$\quad + (1_{E_{s+2}}\wedge \psi_{X\wedge M}u'(1_L\wedge m_K))f_4 - \lambda_0(1_{E_{s+2}}\wedge (r\wedge 1_M)\overline{m}_K)(\kappa\wedge 1_K)$

$\quad + (\bar c_{s+1}\wedge 1_{Y\wedge M})d(g_3)$\\
By (9.4.32) we have $\psi_{X\wedge M}u' = (r\wedge
1_M)\overline{m}_K(j''\wedge 1_K)$ so that $ (j\overline{u}\wedge
1_M)\psi_{X\wedge M}u'
= j'(j''\wedge 1_K)$.  Then, by composing $1_{E_{s+2}}\wedge\phi\cdot j\overline{u}\wedge 1_M$ on (9.4.39), it becomes\\
{\bf (9.4.40)}\qquad $\lambda_0(1_{E_{s+2}}\wedge (\phi\wedge 1_M)j')(\kappa\wedge 1_K) = (1_{E_{s+2}}\wedge (\phi\wedge 1_M)j'(j''\wedge 1_K))d(f_5)$

$\quad + (1_{E_{s+2}}\wedge (\phi\wedge 1_M)ijj'(j''\wedge 1_K))d((1_{E_{s+2}}\wedge 1_L\wedge m_K)f_4) = 0$\\
here we use $(1_{E_{s+2}}\wedge (\phi\cdot j\overline{u}\wedge 1_M)\psi_{X\wedge M}u'(1_L\wedge m_K))f_4 = (1_{E_{s+2}}\wedge
j'(j''\wedge 1_K)(1_L\wedge m_K))f_4 = (\bar c_{s+1}\wedge 1_{L\wedge M})(1_{KG_{s+1}}\wedge (\phi\wedge 1_M)m_M(\overline{u}
\wedge 1_M))g_3$ = 0 and by $1_L\wedge\alpha_1 = \phi\cdot j''$(up to nonzero scalar) we obtain that
$(1_{E_{s+2}}\wedge (\phi\wedge 1_M)j'(j''\wedge 1_K))d(f_5)
= (1_{E_{s+2}}\wedge (1_L\wedge j'\alpha'))d(f_5)$ = 0 ( cf. Prop. 9.4.33) and so  the first term of the right hand side of (9.4.40) is zero.
the second term of the right hand side of (9.4.40) is zero by the same reason.

It follows from (9.4.40) that $(1_{E_{s+2}}\wedge\phi\wedge 1_K)(\kappa\wedge 1_K) = (1_{E_{s+2}}\wedge (1_L\wedge\mu(i'i\wedge 1_K))(\phi\wedge 1_K)
(\kappa\wedge 1_K)(jj'\wedge 1_K)\nu$ = 0 and so we have $(1_{E_{s+2}}\wedge\alpha'')(\kappa\wedge 1_K) = (1_{E_{s+2}}\wedge
\overline{\Delta}(\phi\wedge 1_K))(\kappa\wedge 1_K)$ = 0 and the Theorem follows.

{\bf Case 2} \qquad\quad $m_M(\overline{u}\wedge 1_M)(\bar h\wedge 1_M) = - \lambda_0(j''\wedge 1_M)$.

In this case, the left hand side of (9.4.38) changes sign, then $ - (1_{E_{s+2}}\wedge j'(j''\wedge 1_K)(1_L\wedge m_K))f_4 + (\bar c_{s+1}\wedge 1_M)
(1_{KG_{s+1}}\wedge m_M(\overline{u}\wedge 1_M))g_3 = -2 \lambda_0(1_{E_{s+2}}\wedge j')(\kappa\wedge 1_K)$ and by composing
$1_{E_{s+2}}\wedge \phi\cdot j\overline{u}\wedge 1_M$ on (9.4.39) we have

$\qquad (\lambda_0 - 2 \lambda_0)(1_{E_{s+2}}\wedge (\phi\wedge 1_M)j')(\kappa\wedge 1_K)$

$ = \lambda_0(1_{E_{s+2}}\wedge (\phi\wedge 1_M)j')(\kappa\wedge 1_K) - (1_{E_{s+2}}\wedge (\phi\wedge 1_M)j'(j''\wedge 1_K)(1_L\wedge m_K))f_4$

$ = (1_{E_{s+2}}\wedge (\phi\wedge 1_M)j'(j''\wedge 1_K))d(f_5)$

$\qquad + (1_{E_{s+2}}\wedge (\phi\wedge 1_M)ijj'(j''\wedge 1_K))d((1_{E_{s+2}}\wedge 1_L\wedge m_K)f_4)$ = 0\\
so that the Theorem follows by the same reason. Q.E.D.

\quad

\begin{center}

{\bf \large \S 5.\quad A sequence of $h_0\sigma$ new families in the stable homotopy groups of spheres}

\end{center}

\vspace{2mm}

In this section, the convergence of a sequence of $h_0\sigma$ and $h_0\sigma'$ new families will
be derived by the main Theorem A in \S 2 and the main Theorem C in \S 4, where
$\sigma$ and $\sigma'$ is a pair of $a_0$-related elements.

\vspace{2mm}

{\bf  Theorem 9.5.1}\quad   Let $p\geq 7, n\geq 2$,  then\\
\centerline{$h_0h_n\in Ext_A^{2,p^nq+q}(Z_p,Z_p),\quad h_0b_{n-1}\in Ext_A^{3,p^nq+q}(Z_p,Z_p)$}\\
are permanent cycles in the ASS and they converge  in the ASS
to homotopy elements of order $p$ in $\pi_{p^nq+q-2}S, \pi_{p^nq+q-3}S$ respectively.

{\bf  Proof} : By [12] Theorem 1.2.14  we have $d_2(h_n) = a_0b_{n-1}\in Ext_A^{3,p^nq+1}\\(Z_p,Z_p),n\geq 1$,
 where $d_2 : Ext_A^{1,p^nq}(Z_p,Z_p)\to Ext_A^{3,p^nq+1}(Z_p,Z_p)$ is a secondary differential
 in the ASS. That is,  $h_n$ and $b_{n-1}$ is a pair of $a_0$-related elements
 so that the main Theorem A can apply to  $(\sigma,\sigma') = (h_n,b_{n-1}), (s,tq) = (1,p^nq)$.  We only need to check
 the supposition (I)(II)(III) in the main Theorem A hold.  By knowledge on the $Z_p$-base
 of $Ext_A^{s,*}(Z_p,Z_p)$ for $(s\leq 3)$ we know that the the supposition (I)(II) of the main
 Theorem A hold for $(\sigma,\sigma') = (h_n,b_{n-1}), (s,tq) = (1,p^nq)$.
  On the other hand, from  some results on $Ext_A^{4,*}(Z_p,Z_p)$ in [17] we know that
 the following hold.

$Ext_A^{4,p^nq+rq+1}(Z_p,Z_p) = 0 ( r = 1,3,4),$

\quad $Ext_A^{4,p^nq+rq}(Z_p,Z_p) = 0 (r = 2,3)$,

\quad $Ext_A^{4,p^nq+2q+1}(Z_p,Z_p)\cong Z_p\{\widetilde{\alpha}_2b_{n-1}\}$,

\quad $Ext_A^{4,p^nq+2}(Z_p,Z_p)\cong Z_p\{a_0^2b_{n-1}\},Ext_A^{4,p^nq+1}(Z_p,Z_p) = 0$.\\
 That is, the supposition )III) of the main Theorem A hold for $(\sigma,\sigma') =
 (h_n,b_{n-1}) , (s,tq) = (1, p^nq)$.  Then, by the main Theorem A we obtain that
 $h_0b_{n-1}\in Ext_A^{3,p^nq+q}(Z_p,Z_p)$,
$i_*(h_0h_n)\in Ext_A^{2,p^nq+q}\\(H^*M,Z_p)$ are permanent cycles in the ASS.
  By Remark 9.2.35, the main Theorem A also obtains that $(1_L\wedge i)_*\phi_*(h_n)
\in Ext_A^{2,p^nq+2q}(H^*L,Z_p)$ is a permanent cycle in the ASS so that  the main Theorem C
can apply to obtain the result of the Theorem, this is because by knowledge on the
$A_p$-base of $Ext_A^{s,*}(Z_p,Z_p)$ for $s= 1,2,3$ we can easy to see that the supposition
(I) of the main Theorem C hold for $(\sigma, \sigma', s,tq) = (h_n,b_{n-1}, 1, p^nq)$.
Q.E.D.

\vspace{2mm}

Now  we apply the main Theorem A and the main Theorem C to $(\sigma,\sigma') = (h_nh_m, h_nb_{m-1}-h_mb_{n-1}),
(s,tq) = (2,p^nq+p^mq)$ to obtain another sequence of $h_0\sigma$ families in the stable
homotopy groups of spheres.
For checking the supposition (I)(II)(III) of the main Theorem A, we first prove the
following Proposition.

\vspace{2mm}

{\bf  Proposition 9.5.2}\quad  Let $p\geq 7, n\geq m + 2\geq 4, tq = p^nq + p^mq$,  then\\
(1) $Ext_A^{4,tq+rq+u}(Z_p,Z_p)$ = 0 for $r = 2,3,4, u = -1,0$ or
$r = 3,4 , u = 1$,

 $ Ext_A^{4,tq+q}(Z_p,Z_p)\cong Z_p\{h_0h_nb_{m-1}, h_0h_mb_{n-1}\},$

 $ Ext_A^{4,tq}(Z_p,Z_p)\cong Z_p\{b_{n-1}b_{m-1}\},$

$ Ext_A^{4,tq+1}(Z_p,Z_p)\cong Z_p\{a_0h_nb_{m-1}, a_0h_mb_{n-1}\}$\\
(2) $Ext_A^{5,tq+rq+1}(Z_p,Z_p) = 0$ for $r = 1,3,4$,

 $Ext_A^{5,tq+rq}(Z_p,Z_p)$ = 0 for $ r = 2,3$,

 $ Ext_A^{5,tq+2q+1}(Z_p,Z_p)\cong Z_p\{\widetilde{\alpha}_2h_nb_{m-1},
\widetilde{\alpha}_2h_mb_{n-1}\}$,

 $ Ext_A^{5,tq+2}(Z_p,Z_p)\cong Z_p\{a_0^2h_nb_{m-1}, a_0^2h_mb_{n-1}\}$,

 $Ext_A^{5,tq+1}(Z_p,Z_p)\cong Z_p\{a_0b_{n-1}b_{m-1}\}$,

 $a_0^2b_{n-1}b_{m-1}\neq 0 \in Ext_A^{6,tq+2}(Z_p,Z_p)$

{\bf  Proof} :  By Theorem  5.5.3, there is a May spectral
sequence (MSS) $\{E_{r}^{s,t,*} , d_{r}\}$ which converges to
$Ext_{A}^{s,t}(Z_{p},Z_{p})$ and whose  $E_{1}$-term is

$\quad E_{1}^{*,*,*} = E(h_{i,j}\mid i > 0, j\geq 0)\otimes P(b_{i,j}\mid i > 0, j\geq 0)
\otimes P(a_{i}\mid i\geq 0)$,\\
where  $E$ denotes the exterior algebra and $P$ denotes a polynomial algebra,
, $h_{i,j}\in $
$E_{1}^{1,2(p^{i}-1)p^{j},2i-1}$, $b_{i,j} \in \\E_{1}^{2,2(p^{i}-1)p^{j+1},p(2i-1)},
a_{i}\in E_{1}^{1,2p^{i}-1, 2i+1}$.  Consider the following second degrees (mod $p^nq$)
of the generators in the $E_1^{*,*,*}$-term, where $0\leq i \leq n , n\geq m + 2\geq 4$

$\qquad \mid h_{s,i}\mid = (p^{s+i-1} + \cdots + p^i)q$\quad (mod $p^{n}q ) ,
\quad 0 \leq i < s+i-1 < n$

$\qquad\qquad\quad = (p^{n-1} + \cdots +  p^i)q $ \quad (mod $p^nq$ ),
\quad $0 \leq i < s+i-1 = n$,

$\qquad \mid b_{s,i-1}\mid = (p^{s+i-1} + \cdots + p^i)q$ \quad (mod $p^nq$),
$1\leq i < s+i-1 < n$,

$\qquad\qquad\quad = (p^{n-1} + \cdots + p^i)q$ \quad (mod $p^nq$), $1\leq i < s+i-1 = n$.

$\qquad \mid a_{i+1}\mid   = (p^{i} + \cdots + 1)q + 1$ \quad (mod $p^{n}q $ ),
\quad $1 \leq i< n$

$\qquad \mid a_{i+1}\mid = (p^{n-1} + \cdots + 1)q + 1$ \quad (mod $p^nq$ ),\quad $i = n$.\\
 For degree $k = tq+rq+u$  such that $0\leq r\leq 4,
-1\leq u\leq 2$  we have $ k \equiv p^mq + rq + u$ (mod $p^nq$ ). Then, for $3\leq w\leq 5$,
$E_1^{w,tq+rq+u,*}$  has no such generators which have one of the above elements as a factor,
this is because such  a generator will have second degree  $(c_{n}p^{n-1} + \cdots + c_{1}p
+ c_{0})q + d$ (mod $p^{n}q$ ), where $c_{i}\neq 0 (1\leq i\leq m-1$ or $m < i < n$), $0\leq c_{l} < p , l = 0,\cdots
, n, 0\leq d\leq 5$. In addition,  the second degree  $\mid b_{1,i-1}\mid
= p^iq$ (mod $p^nq$)( $1\leq i\leq n$),$\mid h_{1,i}\mid = p^iq $
(mod $p^nq$)( $0\leq i\leq n$). Then , exclude the above factor
and the factor which has second degree $\geq tq+pq$,
we know that the only possibility of the factor of the generators in $E_{1}^{w,tq+rq+u,*}$
are  $a_1, a_0, h_{1,0}$ , $h_{1,n}, h_{1.m},$ $b_{1,n-1}, b_{1,m-1}$ .

Then, by degree reasons we have

$E_{1}^{4,tq+rq+1,*}$ = 0 for $ r = 3,4$, \quad $E_1^{4,tq+rq+u,*}$ = 0
for $r = 2,3,4, u = -1,0$

$E_1^{4,tq,*} =  Z_p\{b_{1,n-1}b_{1,m-1}\}, \quad E_1^{4,tq+1,*}\cong Z_p
\{a_0h_{1,n}b_{1,m-1}, a_0h_{1,m}b_{1,n-1}\}$,

$E_1^{4,tq+2,*} = Z_p\{a_0^2h_{1,n}h_{1,m}\}$,

$E_1^{4,tq+2q+1,*} =  Z_p\{h_{1,0}a_1h_{1,n}h_{1,m}\}$,

$E_1^{4,tq+q,*} =
Z_p\{h_{1,0}h_{1,n}b_{1,m-1}, h_{1,0}h_{1,m}b_{1,n-1}\}$,

$ E_1^{3,tq+1,*} =  Z_p\{a_0h_{1,n}h_{1,m}\}$, \quad $E_1^{3,tq.*} = Z_p\{h_{1,n}b_{1,m-1}, h_{1,m}b_{1,n-1}\}$,

$E_1^{3,tq+q,*} =  Z_p\{h_{1,0}h_{1,n}h_{1,m}\}$, \quad $
E_1^{3,tq+2q+1,*} = 0$\newline Note that the differentials in the
MSS is derivative,  that is, $d_{r}(xy) = d_{r}(x) y + (-1)^{s} x
d_{r}(y)$ for $x\in E_{1}^{s,t,*}, y \in E_{1}^{s',t',*}$. In
addition, ,$a_0, h_{1,n}, b_{1,n-1}, h_{1,0}a_1$ are permanent
cycles in the MSS which converge to $ a_0, h_{n}, b_{n-1},
\widetilde{\alpha}_2\in Ext_A^{*,*}\\(Z_p,Z_p)$ respectively.

 Then, the differential
$d_rE_r^{3,tq+sq+u,*}$ = 0 for all $r\geq 1$ and $s = u = 0$ or $s
= 1 , u = 0$ or $s = 0,  u = 1$ or $s = 2, u = 1$ so that
$b_{1,n-1}b_{1,m-1},  a_0h_{1,n}b_{1,m-1},$\\$
a_0h_{1,m}b_{1,n-1},  h_{1,0}h_{1,n}b_{1,m-1},
h_{1,0}h_{1,m}b_{1,n-1}\in E_r^{4,*,*}$ are not $d_r$-boundary in
the MSS and so   $b_{n-1}b_{m-1}, a_0h_nb_{m-1}, a_0h_mb_{n-1},
h_0h_nb_{m-1}, \\h_0h_mb_{n-1}$ are all nontrivial in $Ext_A^{4,*}
(Z_p,Z_p)$. This shows (1).

Similarly, by degree reasons we have

$E_1^{5,tq+q+1,*}\cong Z_p\{a_0h_{1,0}h_{1,n}b_{1,m-1}, a_0h_{1,0}h_{1,m}b_{1,n-1},
a_1b_{1,n-1}b_{1,m-1}\}$

$E_1^{5,tq+rq+1,*}$ = 0  for $ r = 3,4$,\quad $E_1^{5,tq+rq,*} = 0$  for $ r = 2,3$

$E_1^{5,tq+2q+1,*}\cong Z_p\{h_{1,0}a_1h_{1,n}b_{1,m-1}, h_{1,0}a_1h_{1,m}b_{1,n-1}\}$

$E_1^{5,tq+1,*} = Z_p\{a_0b_{1,n-1}b_{1,m-1}\}, \quad E_1^{4,tq+2q+1,*}\cong Z_p\{h_{1,0}a_1h_{1,n}h_{1,m}\}$

$E_1^{5,tq+2,*} = Z_p\{a_0^2h_{1,m}b_{1,n-1}, a_0^2h_{1,n}b_{1,n-1}\}$\\
All the generators of $E_1^{5,tq+q+1,*}$ dy in the MSS, this is because
$ a_0h_{1,0}h_{1,n}b_{1,m-1}\\ = - d_1(a_1h_{1,n}b_{1,m-1})$,
 $ a_0h_{1,0}h_{1,m}b_{1,n-1} = - d_1(a_1h_{1,m}b_{1,n-1})$ and
\\\centerline{$d_1(a_1b_{1,n-1}b_{1,m-1}) = - a_0h_{1,0}b_{1,n-1}b_{1,m-1}\neq 0
\in E_1^{5,tq+q+1,*}$}\\  Then, $Ext_A^{5,tq+q+1}(Z_p,Z_p)$ = 0.  In addition,
similar to that given in (1) we have,
$d_rE_r^{4,tq+u,*} = 0$, $d_rE_r^{4,tq+2q+1,*} = 0$  for all $r\geq 1$, $u = 1,2$
.  Then, the generators in $E_1^{5,*,*}$ converge in the  MSS to
$\widetilde{\alpha}_2h_nb_{m-1},
\widetilde{\alpha}_2h_mb_{n-1}, a_0b_{n-1}b_{m-1}, \\a_0^2h_mb_{n-1}, a_0^2h_nb_{m-1}$
respectively.  For the last result,  note that $d_rE_r^{5,tq+2,*}$ = 0 for all
$r\geq 1$ , then $a_0^2b_{n-1}b_{m-1}\neq 0\in Ext_A^{6,tq+2}(Z_p,Z_p)$.
 This shows  (2).  Q.E.D.

\vspace{2mm}

{\bf  Theorem 9.5.3}\quad  Let $p\geq 7, n\geq m+2\geq 4$,  then\\
\centerline{$h_0h_nh_m\in Ext_A^{3,p^nq+p^mq+q}(Z_p,Z_p),$}\\
\centerline{$h_0(h_nb_{m-1}-h_mb_{n-1})
\in Ext_A^{4,p^nq+p^mq+q}(Z_p,Z_p)$}\\
are permanent cycles  in the ASS and they converge to homotopy elements of order $p$ in
$\pi_{p^nq+p^mq+q-3}S$ and $\pi_{p^nq+p^mq+q-4}S$ respectively.

{\bf Proof} : By [12]p.11 Theorem 1.2.14, there is a nontrivial secondary differential
$d_2(h_n) = a_0b_{n-1} ( n\geq 1$ ) and it follows that
$d_2(h_nh_m) = d_1(h_n)h_m + (-1)^{1+p^nq}h_nd_2(h_m) = a_0h_mb_{n-1} - a_0h_nb_{m-1}$.
 That is, $(h_nh_m , h_mb_{n-1}-h_nb_{m-1})$ is a pair of $a_0$-related elements.
 By applying the main Theorem A to
$(\sigma, \sigma') = (h_nh_m , h_mb_{n-1}-h_nb_{m-1}), (s,tq) = (2,p^nq+p^mq)$ we have
$h_0(h_mb_{n-1}-h_nb_{m-1})\\\in Ext_A^{4,p^nq+p^mq+q}(Z_p,Z_p)$ and $i_*(h_0h_nh_m)
\in Ext_A^{3,p^nq+p^mq+q}(Z_p,Z_p)$ are permanent cycles in the ASS, this is
because by knowledge of $Z_p$-base of $Ext_A^{s,*}(Z_p,\\Z_p)$ for  $s\leq 3$ we know that
the supposition (I)(II)(III) of the main Theorem A hold.
By Remark 9.2.35, the main Theorem A also obtains that $(1_L\wedge i)_*\phi_*(h_nh_m)
\in Ext_A^{3,p^nq+p^mq+2q}(H^*L,Z_p)$ is a permanent cycle in the ASS so that by the
main  Theorem C , the result of the Theorem follows. This is because the supposition (I)
of the main Theorem C hold by the knowledge of the $Z_p$-base of
$Ext_A^{s,*}(Z_p,Z_p)$ for $s = 1,2,3$. Q.E.D.

\vspace{2mm}

From Theorem 9.5.1 and Theorem 9.5.3 , we obtain four families of $h_0\sigma$ new families.
In fact, there are many pairs of $a_0$-related elements so that we can
expect to obtain some other sequence of $h_0\sigma$ new families in the stable homotopy
groups of speheres. We have the following conjectures.

\vspace{2mm}

{\bf  Conjecture 9.5.4}\quad  Let $p\geq 7, n\geq 3$, then there is a secondary
differential $d_2(g_n) = a_0l_n\in Ext_A^{4,p^{n+1}q+2p^nq+1}(Z_p,Z_p), n\geq 3$ (up to nonzero scalar)  and\\
\centerline{$h_0g_n\in Ext_A^{3,p^{n+1}q+2p^nq+q}(Z_p,Z_p)$}\\
\centerline{$h_0l_n\in Ext_A^{4,p^{n+1}q+2p^nq+q}(Z_p,Z_p)$}\\
are permanent cycles in the ASS and they converge to homotopy elements of order $p$ in
$\pi_{p^{n+1}q+2p^nq+q-3}S$ and $\pi_{p^{n+1}q+2p^nq+q-4}S$
respectively, where $g_n\in Ext_A^{2,p^{n+1}q+2p^nq}(Z_p,Z_p), l_n\in Ext_A^{3,p^{n+1}q+2p^nq}(Z_p,Z_p)$.

\vspace{2mm}

{\bf  Conjecture 9.5.5}\quad  Let $p\geq 7, n\geq 3$, then there is a secondary
differential $d_2(k_n)
= a_0l'_n\in Ext_A^{4,2p^{n+1}q+p^nq+1}(Z_p,Z_p)$ ( up to nonzero scalar), $n\geq 3$
and\\
\centerline{$h_0k_n\in Ext_A^{3,2p^{n+1}q+p^nq+q}(Z_p,Z_p)$}\\
\centerline{$h_0l'_n\in Ext_A^{4,2p^{n+1}q+p^nq+q}(Z_p,Z_p)$}\\
are permanent cycles in the ASS and they converge to homotopy
elementws of order $p$ in $\pi_{2p^{n+1}+p^nq+q-3}S$ and
$\pi_{2p^{n+1}q+p^nq+q-4}S$ , where $k_n\in
Ext_A^{2,2p^{n+1}q+p^nq}(Z_p,Z_p), l'_n\in
Ext_A^{3,2p^{n+1}q+p^nq}(Z_p,Z_p)$.

\vspace{2mm}

{\bf Remark 9.5.6}\quad  By [10][25],  there is Thom map $\Phi : Ext_{BP_*BP}^{s,*}(BP_*,BP_*)\\
\rightarrow Ext_A^{s,*}(Z_p,Z_p)$( $s = 2,3$) such that $\Phi(\beta_{p^{n-1}/p^{n-1}-1}) = h_0h_n$,
$\Phi(\beta_{p^{n-1}/p^{n-1}}) = b_{n-1}$,
$\Phi(\beta_{p^{m-1}/p^{m-1}-1}\beta_{p^{n-1}/p^{n-1}} - \beta_{p^{n-1}/p^{n-1}-1}\beta_{p^{m-1}/p^{m-1}})
= h_0(h_mb_{n-1} \\- h_nb_{m-1})$,
$\Phi(\gamma_{p^{n-2}/p^{n-2}-p^{m-1},p^{m-1}-1}) = h_0h_nh_m$.
Then,  the $h_0h_n, h_0b_{n-1}, \\h_0(h_mb_{n-1}-h_nb_{m-1}), h_0h_nh_m$-map
obtained by Theorem 9.5.1 and Theorem 9.5.3 are represented  by
$\beta_{p^{n-1}/p^{n-1}-1}$ + other terms $\in Ext_{BP_*BP}^{2,p^nq+q}(BP_*,BP_*)$,
$\alpha_1\beta_{p^{n-1}/p^{n-1}}$ + other terms $\in Ext_{BP_*BP}^{3,p^nq+q}(BP_*,BP_*)$,
$\beta_{p^{m-1}/p^{m-1}-1}\beta_{p^{n-1}/p^{n-1}} - \beta_{p^{n-1}/p^{n-1}-1}
\cdot\beta_{p^{m-1}/p^{m-1}}$ + other terms $\in Ext_{BP_*BP}^{4,p^nq+p^mq+q}(BP_*,BP_*)$,\\
$\gamma_{p^{n-2}/p^{n-2}-p^{m-1},p^{m-1}-1}$ + other terms $\in Ext_{BP_*BP}^{3,p^nq+p^mq+q}(BP_*,BP_*)$
respectively in the Adams-Novikov spectral sequence.

\quad

\begin{center}

{\bf \large\S 6. \quad A sequence of $h_0\sigma\widetilde{\gamma}_s, g_0\sigma
\widetilde{\gamma}_s$ new families in the stable homotopy groups of spheres}

\end{center}

\vspace{2mm}

In this section, we use the main Theorem B to obtain
$i'_*i_*(g_0h_n),\\ i'_*i_*(g_0b_{n-1}),
i'_*i_*(g_0h_nh_m),i'_*i_*(g_0(h_nb_{m-1}-b_mb_{n-1}))$ et al converge to the corresponding nontrivial
homotopy elements in the homotopy groups of Smith-Toda spectrum $V(1)$.
In base of these results, we obtain a sequence of $g_0\sigma\widetilde{\gamma}_s,
h_0\sigma\widetilde{\gamma}_s$ new families in the stable homotopy groups of spheres.

\vspace{2mm}

{\bf  Theorem 9.6.1}\quad Let $p\geq 5, n\geq 2$,  then\\
\centerline{$i'_*i_*(g_0h_n)\in Ext_A^{3,p^nq+pq+2q}(H^*K,Z_p)$}\\
\centerline{$i'_*i_*(g_0b_{n-1})\in Ext_A^{4,p^nq+pq+2q}(H^*K,Z_p)$}\\
are permanent cycles in the ASS and they converge to the corresponding homotopy element in
$\pi_{p^nq+pq+2q-3}K, \pi_{p^nq+pq+2q-4}K$ respectively.

{\bf  Proof} : We first apply the main Theorem B to $(\sigma, s, tq) = (h_n, 1, p^nq)$.
By Theorem 9.5.1, the supposition (II) of the main Theorem B holds. Moreover,
by knowledge on the $Z_p$-base of $Ext_A^{s,*}(Z_p,Z_p)$ for $s = 1,2,3$ we know that
the supposition (I) of the main Theorem B holds, then the first result of the
Theorem follows by the main Theorem B.

For the second result, we apply the main Theorem B to $(\sigma , s, tq) = (b_{n-1}, 2,p^nq)$.
Similarly by Theorem 9.5.1, the supposition (II) of the main Theorem B holds.
Noreover, by knowledge on the $Z_p$-base of $Ext_A^{s,*}(Z_p,Z_p)$
for $s= 2,3$ and some result in [17] on $Ext_A^{4,*}(Z_p,Z_p)$ we know that
the supposition (I) of the main Theorem B holds. Then , the second result
also follows by the main Theorem B. Q.E.D.

{\bf alternative Proof :} \quad It is known from the proof of Theorem 9.5.1 that
the supposition (I)(II)(III) if the main Theorem A hold for $(\sigma,\sigma') = (h_n,b_{n-1}),(s,tq) = (1,p^nq)$.
Then , applying the main Theorem B' in $\S 3$ to
$(\sigma,\sigma') = (h_n,b_{n-1}),(s,tq) = (1,p^nq)$ we obtain the two results of the
Theorem. Q.E.D.

\vspace{2mm}

{\bf  Theorem 9.6.2}\quad  Let $p\geq 7, n\geq m+2\geq 4$,  then\\
\centerline{$i'_*i_*(g_0h_nh_m)\in Ext_A^{4,p^nq+p^mq+pq+2q}(H^*K,Z_p)$}\\
\centerline{$i'_*i_*(g_0(h_nb_{m-1}-h_mb_{n-1}))\in Ext_A^{5,p^nq+p^mq+pq+2q}(H^*K,Z_p)$}\\
are permanent cycles in the ASS and they converge to nontrivial homotopy elements in
$\pi_{p^nq+p^mq+pq+2q-4}K , \pi_{p^nq+p^mq+pq+2q-5}K$  respectively.

{\bf  Proof} : We first apply the main Theorem B to $(\sigma , s, tq) =
(h_nh_m, 2, p^nq+p^mq)$.  By Theorem
9.5.3, the supposition (II) of the main Theorem B holds. By knowledge on the $Z_p$-base of
$Ext_A^{s,*}(Z_p,Z_p)$ for $s = 2,3$ and some result in [17] on $Ext_A^{4,*}(Z_p,Z_p)$
we know that the supposition (I) of the main Theorem B also holds. Then the first
result follows by the main Theorem B. Moreover, we apply the main Theorem B to
$(\sigma , s, tq) = (h_nb_{m-1}-h_mb_{n-1}, 3, p^nq+p^mq)$.
Similarly by Theorem 9.5.3, the supposition (II) of the main Theorem B holds.
By knowledge on the $Z_p$-base of $Ext_A^{3,*}(Z_p,Z_p)$ for
$s = 3,4$ and the result on $Ext_A^{5,p^nq+p^mq+2q+1}(Z_p,Z_p)\cong Z_p\{
\widetilde{\alpha}_2h_nb_{m-1},\widetilde{\alpha}_2h_mb_{n-1}\}$ in Prop. 9.5.2 we
know that the supposition (I) of the main Theorem B holds. Then, the second result
 follows  immediately by the main Theorem B. Q.E.D.

{\bf  alternative Proof :}\quad It is known from the proof of Theorem 9.5.2
that the supposition (I)(II)(III) of the main Theorem A hold for $(\sigma,\sigma') = (h_nh_m, h_nb_{m-1}-h_mb_{n-1}), (s,tq) = (2,p^nq+p^mq)$.
Then by applying the main Theorem B' in \S 3 to $(\sigma,\sigma') = (h_nh_m, h_nb_{m-1}-h_mb_{n-1}),(s,tq) =
(2,p^nq+p^mq)$, we obtain the two result of the Theorem. Q.E.D.

Using the notation in the cofibration (6.2.7)--(6.2.10), we know that \\
\centerline{$\widetilde{\gamma}_s = ((j_1j_2j_3)_*(\gamma)_*^s(i_3i_2i_1)_*(1)
\in Ext_A^{s,sp^2q+ (s-1)pq+(s-2)q+s-3}(Z_p,Z_p)$}\\
converges to the following third periodicity element in the ASS\\
\centerline{$\gamma_s = j_1j_2j_3\gamma^si_3i_2i_1\in\pi_{sp^2q+(s-1)pq
+(s-2)q-3}S$}\\ where $3\leq s < p$ and $1\in Ext_A^{0,0}(Z_p,Z_p)$.  Now we consider
the products
$g_0\sigma\widetilde{\gamma}_s, h_0\sigma\widetilde{\gamma}_s$,
in $Ext_A^{*,*}(Z_p,Z_p)$ and we will prove that they converge
to the corresponding homotopy element of order $p$ in the stable homotopy groups
of spheres, where $\sigma = h_n , b_{n-1}, h_nh_m,$ or $h_nb_{m-1}-h_mb_{n-1}$.

\vspace{2mm}

{\bf  Theorem 9.6.3}\quad  Let $p\geq 7, n\geq 3,3\leq s < p$, then the products\\
\centerline{$g_0h_n\widetilde{\gamma}_s\neq 0\in Ext_A^{s+3,p^nq+sp^2q+spq+sq+s-3}(Z_p,Z_p)$}\\
\centerline{$g_0b_{n-1}\widetilde{\gamma}_s\neq 0\in Ext_A^{s+4, p^nq+sp^2q+spq+sq+s-3}(Z_p,Z_p)$}\\
are permanent cycles in the ASS and they converge to the corresponding homotopy
elements of order $p$ in the stable homotopy groups of spheres.

{\bf  Proof}: By Theorem  9.6.1, there is a nontrivial $f\in\pi_{p^nq+pq+2q-3}K$ such that
it is represented by $i'_*i_*(g_0h_n)\in Ext_A^{3,p^nq+pq+2q}(H^*K,Z_p)$ in the ASS.
Let $\tilde{f} =
j_1j_2j_3\gamma^si_3 f$ be the following composition ($tq = p^nq+pq+2q-3$)

$\qquad\tilde{f} : \Sigma^{tq}S\stackrel{f}{\longrightarrow}V(1)\stackrel{i_3}
{\longrightarrow}V(2)\stackrel{\gamma^s}{\longrightarrow}\Sigma^{-s(p^2q+pq+q)}V(2)$

$\qquad\qquad\quad\stackrel{j_1j_2j_3}
{\longrightarrow}\Sigma^{-s(p^2+p+1)q+(p+2)q+q+3}S$\\
Since $f$ is represented by $(i_2)_*(i_1)_*(g_0h_n)\in Ext_A^{3,p^nq+pq+2q}(H^*K,Z_p)$
in the ASS, then the above $\tilde{f}$ is represented by
$$c = (j_1j_2j_3)_*(\gamma_*)^s(i_3i_2i_1)_*(g_0h_n)\in Ext_A^{s+3,p^nq+s(p^2+p+1)q+s-3}(Z_p,Z_p)$$
By knowledge of Yoneda products we know that the above element $c$ is just the products
$g_0h_n\widetilde{\gamma}_s
\in Ext_A^{s+3,p^nq+s(p^2+p+1)q+s-3}(Z_p,Z_p)$.  Then, to obtain the first result,
it suffices to prove the product
$g_0h_n\widetilde{\gamma}_s$ is nonzero in the Ext group and it is not a $d_r$-boundary in the ASS,
that is, we still need to prove $Ext_A^{s+3-r, p^nq+s(p^2+p+1)q+s-2-r}(Z_p,Z_p)$ is zero for $r\geq 2$.
We may prove this two facts by an argument in the  May spectral sequence.
By degree reasons, $h_n, g_0, \widetilde{\gamma}_s$ is represented by
$h_{1,n}, h_{2,0}h_{1,0}, h_{2,1}h_{1,2}h_{3,0}a_3^{s-3}
\in E_1^{*,*,*}$ in the MSS respectively. Then, the products $g_0h_n\widetilde{\gamma}_s$ is represented by\\
\centerline{$h_{1,n}h_{2,0}h_{1,0}h_{2,1}h_{1,2}h_{3,0}a_3^{s-3}\in E_1^{s+3,
p^nq+s(p^2+p+1)q+s-3,*}$}\\
in the MSS and so we can do some computation in the degree to prove
$E_1^{s+2,p^nq+s(p^2+p+1)q+s-3,*} = 0$ and $E_1^{s+3-r, p^nq+s(p^2+p+1)q+s-2-r,*} = 0 (r\geq 2)$
so that the first result follows. We leave this computation to the reader. The proof and computation
for the second result is similar. Q.E.D.

\vspace{2mm}

By using  Theorem 9.6.2, Theorem 9.5.1 and Theorem 9.5.3, similar to that given in
the proof of Theorem 9.6.3, we can obtain the following Theorem 9.6.4--9.6.6.

\vspace{2mm}

{\bf  Theorem 9.6.4}\quad  Let $p\geq 7, n\geq m+2\geq 5, 3\leq s < p$,  then\\
\centerline{$g_0h_nh_m\widetilde{\gamma}_s\neq 0\in Ext_A^{s+4, p^nq+p^mq+s(p^2+p+1)q+s-3}(Z_p,Z_p)$}\\
\centerline{$g_0(h_nb_{m-1}-h_mb_{n-1})\widetilde{\gamma}_s\neq 0\in Ext_A
^{s+5, p^nq+p^mq+s(p^2+p+1)q+s-3}(Z_p,Z_p)$}\\
are permanent cycles in the ASS and they converge to the corresponding homotopy elements of order $p$
in the stable homotopy groups of spheres.

\vspace{2mm}

{\bf  Theorem 9.6.5}\quad  Let $p\geq 7, n\geq 3, 3\leq s < p$,  then the products\\
\centerline{$h_0h_n\widetilde{\gamma}_s\neq 0\in Ext_A^{s+2, p^nq+sp^2q+(s-1)(p+1)q+s-3}(Z_p,Z_p)$}\\
\centerline{$h_0b_{n-1}\widetilde{\gamma}_s\neq 0\in Ext_A^{s+3, p^nq+sp^2q+(s-1)(p+1)q+s-3}(Z_p,Z_p)$}\\
are permanent cycles in the ASS and they converge to the corresponding elements of order $p$
in the stable homotopy groups of spheres.

\vspace{2mm}

{\bf  Theorem 9.6.6}\quad  Let $p\geq 7, n\geq m+2\geq 5, 3\leq s < p$,  then the products\\
\centerline{$h_0h_nh_m\widetilde{\gamma}_s\neq 0\in Ext_A^{s+3,p^nq+p^mq+sp^2q+(s-1)(p+1)q+s-3}(Z_p,Z_p)$}\\
\centerline{$h_0(h_nb_{m-1}-h_mb_{n-1})\widetilde{\gamma}_s\neq 0\in Ext_A
^{s+4,p^nq+p^mq+sp^2q+(s-1)(p+1)q+s-3}(Z_p,Z_p)$}\\
are permanent cycles in the ASS and they converge to the corresponding homotopy elements of order $p$
in the stable homotopy groups of spheres.

\vspace{2mm}

{\bf  Remark 9.6.7}\quad  The new families obtained in Theorem 9.6.5 and Theorem 9.6.6

are the composition products of $h_0h_n$-element ,
$h_0b_{n-1}$-element in Theorem 9.5.1,
$h_0h_nh_m$-element,$h_0(h_nb_{m-1}-h_mb_{n-1})$- element  in
Theorem 9.5.3 and $\gamma_s = j_1j_2j_3\gamma^si_3i_2i_1\in
\pi_{sp^2q+(s-1)pq+(s-2)q-3}S$. However, the new families obtained
in  Theorem 9.6.3 and Theorem 9.6.4 are indecomposable elements in
the stable homotopy groups of spheres, that is, they are not
compositions of some other elements of lower degrees in the stable
homotopy groups of spheres. This is because $g_0\in
Ext_A^{2,pq+2q}(Z_p,Z_p)$ dies in the ASS, that is, it support a
nontrivial differential in the Adams spectral sequence : $d_2(g_0)
= b_0\widetilde{\alpha}_2$ (up to nonzero scalar) $\in
Ext_A^{4,pq+2q+1}(Z_p,Z_p)$ which can be easily proved as follows.
Since $\widetilde{\alpha}_2,  b_0$ converge in the ASS to
$\alpha_2 = j\alpha^2 i, \beta_1 = jj'\beta i'i\in \pi_*S$,then
the composition products of $\beta_1\alpha_2 \in \pi_{pq+2q-3}S$
must be represented by  $b_0\widetilde{\alpha}_2\in
Ext_A^{4,pq+2q+1}(Z_p,Z_p)$ in the ASS. However, it is easily seen
that $\beta_1\alpha_2 = jj'\beta\i'ij\alpha^2 i = 0$ and
$b_0\widetilde{\alpha}_2\neq 0 \in Ext_A^{4,pq+2q+1}(Z_p,Z_p)$,
then $b_0\widetilde{\alpha}_2$ must be a $d_r$-boundary. By degree
reason, the only possibility is $b_0\widetilde{\alpha}_2 =
d_2(g_0)$(up to nonzero scalar).

\vspace{2mm}

{\bf  Conjecture 9.6.8}\quad  By the conjecture 9.5.4--9.5.5, we can conjecture that, for
$p\geq 7, n\geq 3, 3\leq s < p$,
the products $h_0g_n\widetilde{\gamma}_s, h_0l_n\widetilde{\gamma}_s, g_0g_n\widetilde{\gamma}_s,
g_0l_n\widetilde{\gamma}_s, \\h_0k_n\widetilde{\gamma}_s, h_0l'_n\widetilde{\gamma}_s,
g_0k_n\widetilde{\gamma}_s, g_0l'_n\widetilde{\gamma}_s$ are permanent cycles in the ASS
and they converge  to the corresponding homotopy elements of order $p$ in the stable homotopy
groups of spheres. In addition, all results or conjectures in this section also hold
when we replace the products with
$\widetilde{\gamma}_s$ to be the products with $\widetilde{\beta}_s (2\leq s < p)$.
That is, we can obtain a sequence of $h_0\sigma\widetilde{\beta}_s, g_0\sigma
\widetilde{\beta}_s$-elements, where $\widetilde{\beta}_s = (j_1j_2)_*\beta_*^s(i_1i_1)_*(1)
\in Ext_A^{s,spq+(s-1)q+s-2}(Z_p,Z_p), 2\leq s < p$.

\quad

\begin{center}

{\bf \large \S 7. \quad  Third periodicity families in the stable homotopy groups of spheres}

\end{center}

\vspace{2mm}

In this section, we will first prove the convergence of $h_n$-elements in the homotopy groups
of Smith-Toda spectrum $V(1)$ and in base of this we obtain the convergence of
third periodicity $\gamma_{p^n/s}$ families $(1\leq s\leq p^n-1)$ in the Adams-Novikov spectral sequence.

\vspace{2mm}

{\bf  Theorem 9.7.1} ([9]  Theorem II)\quad  Let $p\geq 5, n\geq 0$, \\
\centerline{$h_n\in Ext_{BP_*BP}^{1,p^nq}(BP_*, BP_*K)$}\\
be the element represented by $[t_1^{p^n}]$ in the cobar complex. Then
this $h_n$ is a permanent cycle in the Adams-Novikov spectral sequence and it
converges to a nontrivial homotopy element in
$\pi_{p^nq-1}K$.

The proof of the above $h_n$-Theorem will be the main content of this section.
 By Theorem 8.1.6(b)(ii),  there is a relation\\
{\bf (9.7.2)}\qquad $h_n = c_2(p^{n-2}) + v_2^{p^{n-2}}h_{n-2}\in Ext_{BP_*BP}^{1,p^nq}
(BP_*BP_*K)$\\  By [10].p.502  Cor. 7.8, the image of $v_2^sc_2(p^{n-2}) (p^{n-2} > s \geq 1)$
under the boundary homomorphism (or connecting homomorphism)
$$j'_*:Ext_{BP_*BP}^{1,*}(BP_*,BP_*K)\to Ext_{BP_*BP}^{2,*}(BP_*, BP_*M)$$
and $$j_* : Ext_{BP_*BP}^{2,*}(BP_*,BP_*M)\to Ext_{BP_*BP}^{3,*}(BP_*,BP_*)$$
is just the third periodicity family $\gamma_{p^{n-2}/p^{n-2}-s}\neq 0\in
Ext_{BP_*BP}^{3,*}(BP_*,\\ BP_*)$. Then, by  Theorem 9.7.1, the relation (9.7.2) and Theorem 7.3.2,
we can obtain the following convergence Theorem of third periodicity families in the
stable homotopy groups of spheres immediately.

\vspace{2mm}

{\bf  Theorem 9.7.3} ([9] Theorem I) \quad  Let $p\geq 5, n\geq 1$ and
$1\leq s\leq p^n-1$, then the following third periodicity family\\
\centerline{$\gamma_{p^n/s} \in Ext_{BP_*BP}^{3,*}(BP_*,BP_*)$}\\
is a permanent cycle in tha ASS and it converge to an element of order $p$ in
$\pi_*S$ which has degree $p^{n+2}q+(p^n-s)(p+1)q-q-3$.

To prove Theorem 9.7.1, we first prove the following weaker  Theorem.

\vspace{2mm}

{\bf  Theorem 9.7.4} ([9] Theorem 4.1)  \quad  Let $p\geq 5, n\geq 0$,
$h_n\in Ext_A^{1,p^nq}(Z_p,\\Z_p)$
be the element represented by $\xi^{p^n}$ in the cobar complex,  then
$i'_*i_*(h_n)\in Ext_A^{1,p^nq}(H^*K,Z_p)$ is a permanent cycle in the ASS
and it converge to a nontrivial homotopy element in
$\pi_{p^nq-1}K$.

\vspace{2mm}

The proof of Theorem 9.7.4 will be the main content of the rest of this section.
The proof need some preminilaries on low dimensional Ext groups
and an argument processing in the Adams resolution of some spectra related to $S$.
We first prove some results on Ext groups.

\vspace{2mm}

{\bf  Theorem 9.7.5}\quad  Let $p\geq 3, n\geq 2,$  , then

 (1)\quad $Ext_A^{s,p^nq+r}(H^*K,Z_p) = 0$ for $s = 2,3, r = 1,2$,

\qquad $Ext_A^{s,p^nq+1}(H^*K,H^*M) = 0$,

\qquad $Ext_A^{3,p^nq+q}(H^*K,H^*K)\cong Z_p\{(h_0b_{n-1})'\}$.

(2) \quad $Ext_{A}^{s-1,p^{n}q+q+s-3}(H^{*}Y,Z_{p})$ = 0 for s = 1,2,3.

 (3) \quad $Ext_{A}^{1,p^{n}q}(H^{*}K, H^{*}Y)$ = 0,\\
where $Y$ is the spectrum in the cofibration (9.1.4).

{\bf  Proof} : (1) Consider the following exact sequence ($s = 1,2,3, r = 1,2)$

$Ext_A^{s,p^nq+r}(Z_p,Z_p)\stackrel{i_*}{\rightarrow}Ext_A^{s,p^nq+r}(H^*M,Z_p)
\stackrel{j_*}{\rightarrow}Ext_A^{s,p^nq+r-1}(Z_p,Z_p)\stackrel{p_*}{\rightarrow}$\\
induced by (9.1.1).
By knowledge of $Z_p$-base of $Ext_A^{s,*}(Z_p,Z_p)$ for $s =  1,2,3$ we know that
the right group is zero except for $(s,r) = (1,1),(2,1),(2,2), (3,2)$ it has unique generator
$h_n, b_{n-1}, a_0h_n, a_0b_{n-1}$.
However,$p_*(h_n) =  a_0h_n\neq 0, p_*(b_{n-1}) = a_0b_{n-1}\neq 0, p_*(a_0h_n) =
a_0^2h_n\neq 0, p_*(a_0b_{n-1}) = a_0^2b_{n-1}\neq 0$,  then the above $p_*$ is monic and so
im $j_*$ = 0.  In addition, the left group is zero except for $(s,r) = (2,1), (3,1), (3,2)$
it has unique generator $a_0h_n = p_*(h_n), a_0b_{n-1} = p_*(b_{n-1}), a_0^2h_n = p_*(a_0h_n)$ respectively.
 Then we have im $i_*$ = 0 and obtain that $Ext_A^{s,p^nq+r}(H^*M,Z_p) = 0$ for
 $s = 1,2,3, r = 1,2$.

Look at the following exact sequence ($s = 2,3, r = 1,2)$

$0 = Ext_A^{s,p^nq+r}(H^*M,Z_p)\stackrel{i'_*}{\longrightarrow}Ext_A^{s,p^nq+r}
(H^*K,Z_p)$

$\qquad\qquad\quad \stackrel{j'_*}{\longrightarrow}Ext_A^{s,p^nq-q+r-1}(H^*M,Z_p)$\\
induced by (9.1.2). The left group is zero as shown above.  The right group also is zero,
this is because $Ext_A^{s,p^nq-q+r-1}(Z_p,Z_p)$ = 0
for $s = 2,3, r = 1,2,3$(cf. Chap. 5). Then, the middle group is zero for $s = 2,3,  r = 1,2$
and so $Ext_A^{s,p^nq+1}(H^*K,H^*M)$ = 0 $(s = 2,3$).

For the last result, consider the following exact sequence

$Ext_A^{3,p^nq+2q+1}(H^*K,H^*M)\stackrel{(j')^*}{\longrightarrow}Ext_A^{3,p^nq+q}
(H^*K,H^*K)$

$\qquad\qquad\quad\stackrel{(i')^*}{\longrightarrow}Ext_A^{3,p^nq+q}(H^*K,H^*M)
\stackrel{\alpha^*}{\longrightarrow}$\\
induced by (9.1.2). The left group is zero by Prop. 9.3.2(1) and the right group
has unique generator $i'_*(\alpha_1\wedge 1_M)_*(\tilde{b}_{n-1}))$ (cf. Prop. 9.3.1)
such that $\alpha^*(i')_*(\alpha_1\wedge 1_M)_*(\tilde{b}_{n-1})$ = 0.
Then, the middle group has unique generator $(h_0b_{n-1})'$ such that
$(i')^*(h_0b_{n-1})' = i'_*(\alpha_1\wedge 1_M)(\tilde{b}_{n-1})$. Q.E.D.

(2) The result is obvious for  $s = 1$.   For $s = 2,3$, consider the following
exact sequence

$\stackrel{(i'i)_*}{\longrightarrow} Ext_{A}^{s-1,p^{n}q+q+s-3}(H^{*}K,Z_{p})\stackrel{r_{*}}{\longrightarrow}
Ext_{A}^{s-1,p^{n}q+q+s-3}(H^{*}Y,Z_{p})$

$\qquad\qquad\stackrel{\epsilon_{*}}{\longrightarrow} Ext_{A}^{s,p^{n}q+q+s-3}(Z_{p},
Z_{p})\stackrel{(i'i)_*}{\longrightarrow} $\\
induced by (9.1.4). The left group is zero for  $s = 2$, this is because
 $Ext_{A}^{1,t}(Z_{p},Z_{p})$ = 0 for  $ t = -1,-2$ (mod $q$). The left group
  has unique genertor $(i'i)_{*}(h_{0}h_{n})$ for $s = 3$ so that  im $r_{*}$ = 0.
  The right group is zero for $s = 2$ and has unique generator $h_{0}b_{n-1}$ for $s = 3$ which satisfies
$(i'i)_{*}(h_{0}b_{n-1})\neq 0$, then im $\epsilon_{*}$
= 0 and so the middle group is zero for $s = 1,2,3$.

(3) Observe the following exact sequence

$ 0 = Ext_{A}^{0,p^{n}q}(H^{*}K,Z_{p})\stackrel{\epsilon^{*}}{\longrightarrow}
Ext_{A}^{1,p^{n}q}(H^{*}K,H^{*}Y)$

$\qquad\qquad\stackrel{r^{*}}{\longrightarrow} Ext_{A}^{1,p^{n}q}(H^{*}K,H^{*}K)
\stackrel{(i'i)^{*}}{\longrightarrow} $\\
induced by (9.1.4). The left group clearly is zero and the right group
has unique generator
 $(h_{n})'$ (cf. Prop. 9.3.6) which satisfies $(i'i)^{*}(h_{n})' = (i'i)_{*}(h_{n})
\neq 0 \in Ext_{A}^{1,p^{n}q}(H^{*}K,Z_{p})$ , then the middle group is zero as desired.
Q.E.D.

\vspace{2mm}

{\bf  Prop. 9.7.6} \quad   Let $p\geq 3, n\geq 2,$  then

 (1) $\quad Ext_{A}^{2,p^{n}q}(Z_{p},H^{*}M)$ = 0,
$Ext_{A}^{3,p^{n}q+1}(Z_{p}, H^{*}M)$ = 0.

 (2) $\quad Ext_{A}^{2,p^{n}q}(Z_{p}, H^{*}K)$ = 0 ,  $Ext_{A}^{3,p^{n}q+1}(Z_{p},
H^{*}K)$ = 0.

 (3) $\quad Ext_A^{2,p^nq+q-u}(Z_p,H^*K) = 0$ for $u = 0,1$ ,

\qquad\quad $Ext_A^{3,p^nq+q}(Z_p,H^*K)$ = 0.

{\bf  Proof}:  (1) Consider the following exact sequences

$$ Ext_{A}^{2,p^{n}q+1}(Z_{p},Z_{p})\stackrel{j^{*}}{\rightarrow}
Ext_{A}^{2,p^{n}q}(Z_{p},H^{*}M)\stackrel{i^{*}}{\rightarrow} Ext_{A}^{2,
p^{n}q}(Z_{p},Z_{p})\stackrel{p^*}{\rightarrow} $$
$$ Ext_{A}^{3,p^{n}q+2}(Z_{p},Z_{p})\stackrel{j^{*}}{\rightarrow}
Ext_{A}^{3,p^{n}q+1}(Z_{p},H^{*}M)\stackrel{i^{*}}{\rightarrow} Ext_{A}
^{3,p^{n}q+1}(Z_{p},Z_{p})\stackrel{p^*}{\rightarrow} $$
induced by (9.1.1). The upper left group has unique generator
 $a_{0}h_{n}$ which satisfies $j^{*}(a_{0}h_{n}) = j^{*}p^{*}(h_{n})$ = 0
 and the upper right group has unique generator
 $b_{n-1}$ satisfying $p^{*}(b_{n-1}) =
a_{0}b_{n-1} \neq 0 \in Ext_{A}^{3,p^{n}q+1}(Z_{p},Z_{p})$ (cf. Theorem 5.4.1),
then we have $Ext_{A}^{2,p^{n}q}(Z_{p},H^{*}M)$ = 0. The lower left group
has unique generator $a_{0}^{2}h_{n}$ satisfying
$j^{*}(a_{0}^{2}h_{n}) = j^{*}p^{*}$ $(a_{0}h_{n})$ = 0  and the lower right group
has unique generator $a_{0}b_{n-1}$ such that $p^{*}(a_{0}b_{n-1})
= a_{0}^{2}b_{n-1} \neq 0 \in Ext_{A}^{4,p^{n}q+2}(Z_{p},Z_{p})$ (cf. Prop. 9.5.2(2))
, then $Ext_{A}^{3,p^{n}q+1}(Z_{p},H^{*}M)$ = 0.

(2) Consider the following exact sequences

$0 = Ext_{A}^{2,p^{n}q+q+1}(Z_{p},H^{*}M)\stackrel{(j')^{*}}{\longrightarrow}
Ext_{A}^{2,p^{n}q}(Z_{p},H^{*}K)$

$\qquad\qquad\stackrel{(i')^{*}}{\longrightarrow}Ext_{A}
^{2,p^{n}q}(Z_{p},H^{*}M) = 0$

$0 = Ext_{A}^{3,p^{n}q+q+2}(Z_{p},H^{*}M)\stackrel{(j')^{*}}{\longrightarrow}
Ext_{A}^{3,p^{n}q+1}(Z_{p},H^{*}K)$

$\qquad\qquad\stackrel{(i')^{*}}{\longrightarrow}Ext_{A}^{3,p^{n}q+1}(Z_{p},H^{*}M) = 0$\\
induced by (9.1.2). Both two right groups are zero by (1) and both two left groups are also zero
, this is because
 $Ext_{A}^{2,p^{n}q+q+r}(Z_{p},Z_{p})$ = 0 for $r = 1,2$ (cf. Chapter 5)  and
$Ext_{A}^{3,p^{n}q+q+r}(Z_{p},Z_{p})$ = 0 for $ r = 2,3 $ (cf. Theorem 5.4.1), then
the result follows.

(3) Consider the following exact sequence

$\stackrel{\alpha^*}{\longrightarrow}Ext_{A}^{2,p^{n}q+2q}(Z_{p},H^{*}M)\stackrel{(j')^{*}}{\longrightarrow}
Ext_{A}^{2,p^{n}q+q-1}(Z_{p},H^{*}K)$

$\qquad\qquad\stackrel{(i')^{*}}{\longrightarrow}Ext_{A}^{2,p^{n}q+q-1}(Z_{p},H^{*}M)
\stackrel{\alpha^{*}}{\longrightarrow}$\\
induced by (9.1.2). The left group is zero, this is because
  $Ext_{A}^{2,p^{n}q+2q+r}(Z_{p}, Z_{p})$ = 0 for $r = 0,1 $ (cf. Chapter 5). The right
  group has unique generator
$j^{*}(h_{0}h_{n})$ since $Ext_{A}^{2,p^{n}q+q-1}(Z_{p},Z_{p})$ =
0 and $Ext_{A}^{2,p^{n}q+q}(Z_{p},Z_{p})\cong
Z_{p}\{h_{0}h_{n}\}$.   In addition, we claim that
$\alpha^{*}j^{*}(h_{0}h_{n})$ = $\frac{1}{2}\cdot
j^{*}(\widetilde{\alpha}_{2} h_{n}) \neq 0 \in
Ext_{A}^{3,p^{n}q+2q}(Z_{p},\\H^{*}M)$ . To prove this,  it
suffices to prove $\alpha^*j^*(h_0) =
\frac{1}{2}j^*(\widetilde{\alpha}_2 )\in Ext_A^{2,2q}(Z_p,
\\H^*M)$. Since $i^*\alpha^*j^*(h_0) = \alpha_1^*(h_0) = h_0^2 =
0$, then$\alpha^*j^*(h_0) = \lambda j^*(\widetilde{\alpha}_2)$ for
some scalar $\lambda\in Z_p$. Since both sides of the equation
detect the corresponding homotopy elements, then the relation
$\alpha_1j\alpha = \frac{1}{2} \alpha_2 j$ implies $\lambda =
\frac{1}{2}$.  This shows the above claim and so the above
 $\alpha^{*}$  is monic, im $(i')^{*}$ = 0
and we have $Ext_{A}^{2,p^{n}q+q-1}(Z_{p},H^{*}K)$ = 0.
The proof of the case for $u = 0$ is similar.

For the second result, consider the following exact sequence

$Ext_{A}^{3,p^{n}q+2q+1}(Z_{p},H^{*}M)\stackrel{(j')^{*}}{\longrightarrow}
Ext_{A}^{3,p^{n}q+q}(Z_{p},H^{*}K)$

$\qquad\qquad\stackrel{(i')^{*}}{\longrightarrow}Ext_{A}^{3,p^{n}q+q}(Z_{p},H^{*}M)
\stackrel{\alpha^{*}}{\longrightarrow}$\\
induced by (9.1.2). The left group has unique generator
$j_{*}\alpha_{*}\alpha_{*} (\tilde{h}_{n}) =
j_{*}\alpha_{*}\alpha^{*}(\tilde{h}_{n})$, this is because
$Ext_{A} ^{3,p^{n}q+2q+1}(Z_{p},Z_{p})$ has unique generator
$\widetilde{\alpha}_2h_n\\ = j_*\alpha_*\alpha_*i_*(h_n) =
i^*j_*\alpha_*\alpha_*(\tilde{h}_n)$ and $Ext_{A}^{3,p^{n}q+2q+2}
(Z_{p},Z_{p})$ = 0 (cf. Theorem 5.4.1), then im $(j')^{*}$ = 0.
The right group has unique generator
 $j_{*}\alpha_{*}(\tilde{b}_{n-1})$ since $Ext_{A}^{3,p^{n}q+q}(Z_{p},Z_{p})$
has unique generator $h_{0}b_{n-1} = j_{*}\alpha_{*}i_{*}(b_{n-1}) = j_{*}
\alpha_{*}i^{*}(\tilde{b}_{n-1})$ and $Ext_{A}^{3,p^{n}q+q+1}(Z_{p},Z_{p})$
= 0 (cf. Theorem 5.4.1).   In addition, $\alpha^{*} j_{*}\alpha_{*}\cdot (\tilde{b}_{n-1}) = j_{*}\alpha_{*}
\alpha_{*}(\tilde{b}_{n-1}) \neq 0 \in Ext_{A}^{4,p^{n}q+2q+1}(Z_{p},H^{*}M)$,
this is because $i^{*}j_{*}\alpha_{*}\alpha_{*}\cdot (\tilde{b}_{n-1}) = \widetilde{\alpha}_{2}b_{n-1} \neq 0$
(cf. Prop. 9.5.2(2)). Then the above  $\alpha^{*}$ is monic, im $(i')^{*}$ = 0 and so
the middle group is zero as desired. Q.E.D.

\vspace{2mm}

{\bf  Proposition 9.7.7} \quad  Let $p\geq 3, n\geq 2$, then

 (1)\quad $Ext_{A}^{2,p^{n}q}(Z_{p},H^{*}X)$ = 0, $Ext_{A}^{3,p^{n}q+1}(Z_{p},
H^{*}X)$ = 0.

 (2)\quad $Ext_{A}^{3,p^{n}q+q}(H^{*}X,H^{*}K) \cong Z_{p}\{w_{*}
(h_{0}b_{n-1})'\}.$

 (3)\quad $Ext_{A}^{1,p^{n}q+q-1}(H^{*}X,Z_{p})\cong Z_{p}\{\tau_{*}(h_{n}) \}$, \\
where $X$ is the spectrum in the cofibration (9.3.7) , $\tau : \Sigma^{q-1}S\rightarrow X$ is a map satisfying
 $u \tau = i'i : S\rightarrow K$ which is obtained by $\alpha'' i'i = 0$  and (9.3.7).

{\bf Proof} : (1) Consider the following exact sequences

$0 = Ext_{A}^{2,p^{n}q+q-1}(Z_{p},H^{*}K)\stackrel{u^{*}}{\longrightarrow}
Ext_{A}^{2,p^{n}q}(Z_{p},H^{*}X)$

$\qquad\qquad\quad\stackrel{w^{*}}{\longrightarrow}Ext_{A}^{2,p^{n}q}(Z_{p},H^{*}K) = 0$

$0 = Ext_{A}^{3,p^{n}q+q}(Z_{p},H^{*}K)\stackrel{u^{*}}{\longrightarrow}
Ext_{A}^{3,p^{n}q+1}(Z_{p},H^{*}X)$

$\qquad\qquad\quad\stackrel{w^{*}}{\longrightarrow}Ext_{A}^{3,p^{n}q+1}(Z_{p},H^{*}K) = 0$\\
induced by (9.3.7). By Prop. 9.7.6(2)(3), Both sides four groups are zero so that
the result follows.

(2) We first claim that $Ext_{A}^{s,p^{n}q+1}(H^{*}K,H^{*}K)$ = 0( $s = 2,3$),
then the result follows by the following exact sequence

$\stackrel{(\alpha'')_*}{\longrightarrow}Ext_{A}^{3,p^{n}q+q}(H^{*}K,H^{*}K)
\stackrel{w_*}{\longrightarrow}Ext_A^{3,p^nq+q}(H^*X,H^*K)$

$\qquad\qquad\quad\stackrel{u_*}{\longrightarrow}Ext_A^{3,p^nq+1}(H^*K,H^*K) = 0$\\
induced by (9.3.7), where the left group has unique generator
 $(h_{0}b_{n-1})'$(cf. Prop. 9.7.5(1)). To prove the above claim,
consider the following exact sequence

$\stackrel{\alpha^*}{\longrightarrow}Ext_{A}^{s,p^{n}q+q+2}(H^{*}K,H^{*}M)
\stackrel{(j')^{*}}{\longrightarrow}Ext_{A}^{s,p^{n}q+1}(H^{*}K,H^{*}K)$

$\qquad\qquad\quad\stackrel{(i')^{*}}{\longrightarrow}Ext_{A}^{s,p^{n}q+1}
(H^{*}K,H^{*}M)\stackrel{\alpha^*}{\longrightarrow}$\\
induced by (9.1.2)> The right group is zero for  $s = 2,3$ (cf. Prop. 9.7.5(1)) and
the left group is zero by Prop. 9.3.2(2). This shows the above claim.

(3)  Since $\alpha'' i'i = 0$, then, by (9.3.7),  there is $\tau \in \pi_{q-1}X$ such that
$u \tau = i'i : S\rightarrow K$. Consider the following exact sequence

$ 0 = Ext_{A}^{1,p^{n}q+q-1}(H^{*}K,Z_{p})\stackrel{w_{*}}{\longrightarrow}
Ext_{A}^{1,p^{n}q+q-1}(H^{*}X,Z_{p})$

$\qquad\qquad\quad\stackrel{u_{*}}{\longrightarrow}Ext_{A}^{1,p^{n}q}(H^{*}K,
Z_{p})\stackrel{\alpha''_{*}}{\longrightarrow} $\\
induced by (9.3.7). The left group is zero since  $Ext_{A}^{1,t}(Z_{p},Z_{p})$
= 0 for $t \equiv -1, -2$ (mod $q$). The right group has unique generator
$(i'i)_{*}(h_{n})$ which satisfies $\alpha''_{*}(i'i)_{*}(h_{n})$ = 0, then
the middle group has unique generator
 $\tau_{*}(h_{n})$ such that $u_{*}\tau_{*}(h_{n}) = (i'i)_{*}(h_{n})$. Q.E.D.

Since $u \tau\cdot p = i'i\cdot p = 0$, then by (9.3.7) we have $\tau \cdot p = wi'i \alpha_{1}$
(uo to nonzero scalar), this is because $\pi_{q-1}K\cong Z_{p}\{i'i(\alpha_{1})\}.$
Then, by  $Ext_{A}^{2,p^{n}q+1}(H^{*}K,Z_{p})$ = 0(cf. Prop. 9.7.5(1)) and the
Ext exact sequence induced by (9.3.7) we have\\
{\bf (9.7.8)}\quad $ \tau_{*}(a_{0}b_{n-1}) = \tau_{*}p_{*}(b_{n-1}) = w_{*}(i'i)_{*}
(\alpha_{1})_{*}(b_{n-1})$

$\qquad = w_{*}(i'i)_{*}(h_{0}b_{n-1}) \neq 0 \in Ext_{A}^{3,p^{n}q+q}
(H^{*}X,Z_{p})$.

\vspace{2mm}

{\bf  Proposition 9.7.9} \quad  Let $p\geq 3, n\geq 2$, then

 (1)\quad $Ext_{A}^{1,p^{n}q+q-1}(H^{*}K,H^{*}K)$ = 0, $Ext_{A}^{1,p^{n}q}(H^{*}K,
H^{*}X)$ = 0.

 (2)\quad $Ext_{A}^{1,p^{n}q-q+1}(H^{*}K,H^{*}X)\cong Z_{p}\{u^{*}(h_{n})'\}$.

 (3)\quad $Ext_{A}^{2,p^{n}q+q}(H^{*}X, Z_{p}) \cong Z_{p}\{w_{*}(i'i)_{*}
(h_{0}h_{n}) \}$.

 (4)\quad $Ext_{A}^{2,p^{n}q+1}(H^{*}X, H^{*}K)$ = 0.

{\bf Proof}: (1) Consider the following exact sequence

$\stackrel{\alpha^*}{\longrightarrow}Ext_{A}^{1,p^{n}q+2q}(H^{*}K,H^{*}M)\stackrel{(j')^{*}}{\longrightarrow}
Ext_{A}^{1,p^{n}q+q-1}(H^{*}K,H^{*}K)$

$\qquad\qquad\quad\stackrel{(i')^{*}}{\longrightarrow}Ext_{A}^{1,p^{n}q+q-1}(H^{*}K,H^{*}M)$\\
induced by (9.1.2). The right group is zero since
$Ext_{A}^{1,p^{n}q-2}(H^{*}M,H^{*}M)$ = 0,
$Ext_{A}^{1,p^{n}q+q-1}(H^{*}M,H^{*}M)$ = 0 which is obtained by
$Ext_{A}^{1,t}(Z_{p},Z_{p})$ = 0 for $ t \equiv -1,-2$ (mod $q$)
and  $Ext_{A}^{1,p^{n}q+q+t}(Z_{p},Z_{p})$ = 0 for $t = -1, 0, 1
$. The left group also is zero since $Ext_{A}
^{1,p^{n}q+q-1}(H^{*}M, H^{*}M)$ = 0  and
$Ext_{A}^{1,p^{n}q+2q}(H^{*}M,\\H^{*}M)$ = 0 which is obtained by
the same reason as above. Then the middle group is zero.

The second result follows by the following exact sequence

$0 = Ext_{A}^{1,p^{n}q+q-1}(H^{*}K,H^{*}K)\stackrel{u^{*}}{\longrightarrow}
Ext_{A}^{1,p^{n}q}(H^{*}K,H^{*}X)$

$\qquad\qquad\quad\stackrel{w^{*}}{\longrightarrow}Ext_{A}^{1,p^{n}q}(H^{*}K,H^{*}K)\stackrel
{(\alpha'')^{*}}{\longrightarrow} $\\
induced by (9.3.7), where the right group has unique generator
 $(h_{n})'$
which satisfies $(\alpha'')^{*}(h_{n})' = (h_{0}h_{n})'' \neq 0
\in Ext_{A}^{2,p^{n}q+q-1} (H^{*}K,H^{*}K)$ ( cf. (9.3.8)).

(2) Consider the following exact sequence

$\stackrel{\alpha^*}{\longrightarrow}Ext_{A}^{1,p^{n}q+2}(H^{*}K,H^{*}M)\stackrel{(j')^{*}}{\longrightarrow}
Ext_{A}^{1,p^{n}q-q+1}(H^{*}K,H^{*}K)$

$\qquad\qquad\quad\stackrel{(i')^{*}}{\longrightarrow}Ext_{A}^{1,p^{n}q-q+1}(H^{*}K,
H^{*}M)\stackrel{\alpha^*}{\longrightarrow}$\\
induced by (9.1.2). The left group is zero since
$Ext_{A}^{1,p^{n}q-q+1}(H^{*}M, H^{*}M)$ = 0  and
$Ext_{A}^{1,p^{n}q+2}(H^{*}M,H^{*}M)$ = 0 which is obtained by
$Ext_{A}^{1,p^{n}q-q+t}(Z_{p},Z_{p})$ = 0 for $t = 0,1,2 $ and
$Ext_{A} ^{1,p^{n}q+t} (Z_{p},Z_{p})$ = 0  ($t = 1,2 $).  The
right group also is zero since
$Ext_{A}^{1,p^{n}q-2q}(H^{*}M,H^{*}M)$ = 0 and $Ext_{A}^{1,
p^{n}q-q+1}(H^{*}M,H^{*}M)$ = 0 .  Then we have
$Ext_A^{1,p^{n}q-q+1}(H^{*}K,H^{*}K)$ = 0.

The desired result can be obtained by the following exact sequence

$\stackrel{(\alpha'')^*}{\longrightarrow}Ext_{A}^{1,p^{n}q}(H^{*}K,H^{*}K)\stackrel{u^{*}}{\longrightarrow}
Ext_{A}^{1,p^{n}q-q+1}(H^{*}K,H^{*}X)$

$\qquad\qquad\quad\stackrel{w^{*}}{\longrightarrow}Ext_{A}^{1,p^{n}q-q+1}(H^{*}K,H^{*}K) = 0$\\
induced by (9.3.7), where the left group has unique generator
 $(h_{n})'$ (cf. Prop. 9.3.6).

(3) Consider the following exact sequence

$Ext_A^{1,p^nq+1}(H^*K,Z_p)\stackrel{\alpha''_{*}}{\longrightarrow}Ext_{A}^{2,p^{n}q+q}(H^{*}K, Z_{p})
\stackrel{w_{*}}{\longrightarrow}$

$\qquad\qquad Ext_{A}^{2,p^{n}q+q}(H^{*}X,Z_{p})\stackrel{u_{*}}{\longrightarrow}Ext_{A}^{2,p^{n}q+1}(H^{*}K,Z_{p}) = 0$\\
induced by (9.3.7). The right group is zero by Prop. 9.7.5(1) and the left group
has unique generator $(i'i)_{*}(h_{0}h_{n})$ since
 $Ext_{A}^{2,p^{n}q+q}(Z_{p},Z_{p}) \cong Z_{p}\{ h_{0}h_{n} \}$ and
 $Ext_{A}^{2,t}(Z_{p},Z_{p})$ = 0  for $ t\equiv -1, -2 $ (mod $q$).  In addition,
  $Ext_{A}^{1,p^{n}q+1}(H^{*}K,\\Z_{p})$
= 0, this is because $Ext_{A}^{1,p^{n}q+1}(H^{*}M,Z_{p})$ = 0 ( cf. the proof of Prop. 9.7.5(1)) and
$Ext_{A}^{1,p^{n}q-q}(H^{*}M,Z_{p})$ = 0 ( cf Chapter 5),  then $w_*$ is monic so that
the result follows.

(4) Consider the following exact sequence

$ 0 = Ext_{A}^{2,p^{n}q+1}(H^{*}K,H^{*}K)\stackrel{w_{*}}{\longrightarrow}
Ext_{A}^{2,p^{n}q+1}(H^{*}X,H^{*}K)$

$\qquad\qquad\quad\stackrel{u_{*}}{\longrightarrow}Ext_{A}^{2,p^{n}q-q+2}(H^{*}K,H^{*}K) = 0$\\
induced by (9.3.7). The left group is zero as pointed out in the proof of Prop.
9.7.7(2). The right group also is zero since
 $Ext_{A}^{2,p^{n}q+3}(H^{*}K,H^{*}M)$
= 0  and $Ext_{A}^{2,p^{n}q-q+2}(H^{*}K,H^{*}M)$ = 0 which is obtained by
$Ext_{A}^{2,p^{n}q+t}\\(Z_{p},Z_{p})$ = 0  for $t = 2,3$ and $Ext_{A}^{2,p^{n}q-rq+t}
(Z_{p},Z_{p})$ = 0 for $r = 1,2$ , $t = 1,2,3$. Then the middle group is zero as desired. Q.E.D.

\vspace{2mm}

Now we proceed to prove the main Theorem 9.7.4 in this section.
The proof will be done by some argument processing in the Adams resolution
(9.2.9) . We first prove the following Proposition and Lemmas.

\vspace{2mm}

{\bf Proposition 9.7.10}\quad  Let $p\geq 5, n\geq 2$,
$(h_0h_n)''\in [\Sigma^{p^nq+q-1}K, KG_2\wedge K]$ be $d_1$-cycle
which represents the element $(h_0h_n)'' = (\alpha'')^*(h_n)'\in
\\Ext_A^{2,p^nq+q-1}(H^*K,H^*K)$(cf. (9.3.8)), then there exist
$\eta''_{n,2}\in [\Sigma^{p^nq+q-1}K, \\E_2\wedge K]$ and
$(\eta''_{n,2})_Y \in [\Sigma^{p^nq+q-1}Y, E_2\wedge K]$ such that
$(\bar b_2\wedge 1_K)\eta''_{n,2} = (\bar b_2\wedge
1_K)(\eta''_{n,2})_Y\cdot r = (h_0h_n)''\in [\Sigma^{p^nq+q-1}K,
KG_2\wedge K]$ where $r : K\to Y$ is the map in (9.1.4).

{\bf  Proof} : Applying Theorem 9.3.9 to $(\sigma, s, tq) = (h_n, 1,p^nq)$, or applying the mian
Theorem B' and its proof  to $(\sigma,\sigma') = (h_n,b_{n-1}),
(s,tq) = (1,p^nq)$, we have $(\bar c_2\wedge 1_K)(h_0h_n)'' = 0$,.
Them there exists $\eta''_{n,2}\in [\Sigma^{p^nq+q-1}K, KG_2\wedge K]$ such that
$(\bar b_2\wedge 1_K)\eta''_{n,2} = (h_0h_n)''\in [\Sigma^{p^nq+q-1}K, KG_2
\wedge K]$. For the second result, note that
$(h_0h_n)'' i'i\in [\Sigma^{p^nq+q-1}S, KG_2\wedge K] = 0$, this is because
$\pi_{p^nq+tq-u}KG_2\cong Ext_A^{2,p^nq+tq-u}(Z_p,Z_p)$ = 0  for $t = 0,1, u = 1,2,3$.
Then there is $(h_0h_n)''_Y\in [\Sigma^{p^nq+q-1}Y, KG_2\wedge K]$ such that
$(h_0h_n)'' = (h_0h_n)''_Y\cdot r$,  where $r : K\to Y$ is the map in (9.1.4). Then
by Theorem 9.3.9 we have
$(\bar c_2\wedge 1_K)(h_0h_n)''_Y\cdot r = 0$ and so , by the cofibration (9.1.4),
$(\bar c_2\wedge 1_K)(h_0h_n)''_Y = f'_0\epsilon = (1_{E_2}\wedge\epsilon\wedge 1_K)(1_Y\wedge f'_0) = 0$,
for some $f'_0\in [\Sigma^{p^nq+q}S, E_3\wedge K]$, where we use $\epsilon\wedge 1_K = \mu (i'i\wedge 1_K)
(\epsilon\wedge 1_K) = 0$ which is obtained by (9.1.26) and the cofibration (9.1.4). Hence, there is
$(\eta''_{n,2})_Y\in [\Sigma^{p^nq+q-1}Y, E_2\wedge K]$ such that $(\bar b_2\wedge 1_K)
(\eta''_{n,2})_Y = (h_0h_n)''_Y\in [\Sigma^{p^nq+q-1}Y, KG_2\wedge K]$. Q.E.D.

\vspace{2mm}

{\bf  Lemma 9.7.11}\quad Let $p\geq 5, n\geq 2$  and $(\eta''_{n})_{Y} = (\bar a_{0}\bar a_{1}\wedge 1_{K})(\eta''_{n,2})_{Y}
\in [\Sigma^{p^{n}q+q-3}Y, K]$ be the map obtained in Prop. 9.7.10, then
$w (\eta''_{n})_{Y}\cdot r = \lambda' w (\zeta_{n-1}\wedge 1_{K}) + (\bar a_{0}\bar a_{1}
\bar a_{2}\bar a_{3}\wedge 1_{X})f''_{1}$  for some $f''_{1} \in
[\Sigma^{p^{n}q+q+1}K,E_{4}\wedge X]$  and nonzero $\lambda' \in Z_{p}$, where
$\zeta_{n-1} \in \pi_{p^{n}q+q-3}S$ is the element obtained in Theorem 9.5.1
which is represented by $h_{0}b_{n-1} \in Ext_{A}^{3,p^{n}q+q}(Z_{p},Z_{p})$ in the ASS.

{\bf Proof} :  By Prop. 9.7.10  and (9.3.8)(9.3.7), $(\bar b_{2}\wedge 1_{X})(1_{E_{2}}\wedge w )(\eta''_{n,2})_{Y}\cdot r$
= $(1_{KG_{2}}\wedge w)(h_{0}h_{n})'' = (\bar b_{2}\bar c_{1}\wedge 1_{X})g''$ with
$g'' \in [\Sigma^{p^{n}q+q-1}K, KG_{1}\wedge X]$,  then\\
{\bf (9.7.12)} \qquad $ (1_{E_{2}}\wedge w )(\eta''_{n,2})_{Y}\cdot r = (\bar c_{1}\wedge 1_{X})g''
+ (\bar a_{2}\wedge 1_{X})f''_{0}$\\
for some $f''_{0} \in [\Sigma^{p^{n}q+q}K, E_{3}\wedge X]$.  The $d_{1}$-cycle
$(\bar b_{3}\wedge 1_{X})f''_{0} \in [\Sigma^{p^{n}q+q}K, KG_{3}\wedge X]$
represents  an element in $Ext_{A}^{3,p^{n}q+q}(H^{*}X,H^{*}K)$ and this group
has unique generator
 $w_{*}[h_{0}b_{n-1}\wedge 1_{K}]$ ( cf. Prop. 9.7.7(2)), then

$(\bar b_{3}\wedge 1_{X})f''_{0} = \lambda' (1_{KG_{3}}\wedge w)(h_{0}b_{n-1}
\wedge 1_{K}) + (\bar b_{3}\bar c_{2}\wedge 1_{X})\tilde{g}_{0}$

$\qquad =  \lambda' (\bar b_{3}\wedge 1_{X})(1_{E_{3}}
\wedge w)(\zeta_{n-1,3}\wedge 1_{K}) + (\bar b_{3}\bar c_{2}\wedge 1_{X})
\tilde{g}_{0}$\\
with $\lambda' \in Z_p$  and $\tilde{g}_0\in [\Sigma^{p^nq+q}K, KG_2\wedge X]$, where we use
$(\bar b_3\wedge 1_K)(\zeta_{n-1,3}\wedge 1_K) = h_0b_{n-1}\wedge 1_K$ (cf. Theorem 9.5.1)
.  Then $f''_{0} = \lambda' (1_{E_{3}}\wedge w)(\zeta_{n-1,3}\wedge 1_{K}) +
(\bar c_{2}\wedge 1_{X})\tilde{g}_{0} + (\bar a_{3}\wedge 1_{X})f''_{1}$
for some $f''_{1} \in [\Sigma^{p^{n}q+q+1}K, E_{4}\wedge X]$  and so we have
$(\bar a_{2}\wedge 1_{X})f''_{0} = \lambda' (\bar a_{2}\wedge 1_{X})(1_{E_{3}}\wedge w)
(\zeta_{n-1,3}\wedge 1_{K}) + (\bar a_{2}\bar a_{3}\wedge 1_{X})f''_{1}$
and  (9.7.12)  becomes\\
{\bf (9.7.13)} \qquad $ (1_{E_{2}}\wedge w)(\eta''_{n,2})_{Y}\cdot r = (\bar c_{1}\wedge 1_{X})g''$

$\qquad\qquad + \lambda' (\bar a_{2}\wedge 1_{X})(1_{E_{3}}\wedge w)(\zeta_{n-1,3}\wedge 1_{K}) +
(\bar a_{2}\bar a_{3}\wedge 1_{X})f''_{1}$\\
with $g'' \in [\Sigma^{p^{n}q+q-1}K , KG_{1}\wedge X]$, $f''_{1} \in [\Sigma^{p^{n}q+q+1}K,
E_{4}\wedge X]$ and $\lambda'\in Z_{p}$.

To prove the Lemma, it suffices to prove the scalar $\lambda'$ in (9.7.13) is nonzero.
Suppose in contrast that $\lambda'$ = 0, then by (9.7.13)(9.1.4)  we have\\
{\bf (9.7.14)}\qquad $ (\bar a_{2}\bar a_{3}\wedge 1_{X}) f''_{1} i'i = - (\bar c_{1}\wedge 1_{X}) g'' i'i$\\
This will yield a contradiction as shown below.

Note that the $d_{1}$-cycle $g'' i'i \in \pi_{p^{n}q+q-1}KG_{1}\wedge X$ represents
an element in $Ext_{A}^{1.p^{n}q+q-1}(H^{*}X , Z_{p}) \cong Z_{p}\{ \tau_{*}(h_{n})\}$
( cf. Prop. 9.7.7(3)).  Then $g'' i'i = \lambda_{0}(1_{KG_{1}}\wedge\tau )(h_{n})$
, where $h_{n} \in \pi_{p^{n}q}KG_{1} \cong Ext_{A}^{1,p^{n}q}(Z_{p},Z_{p})$
and $\lambda_{0} \in Z_{p}$.   Consequently, (9.7.14)  becomes\\
{\bf (9.7.15)}\qquad $ (\bar a_{2}\bar a_{3}\wedge 1_{X})f''_{1} i'i = -\lambda_{0}
(\bar c_{1}\wedge 1_{X})(1_{KG_{1}}\wedge\tau )(h_{n})$

The equation  (9.7.15) means  the secondary differential $-
\lambda_0 d_2(\tau_*(h_n))$ = 0.
 However, by  [12] p.11  Theorem 1.2.14, $d_{2}(h_{n}) = a_{0}b_{n-1} \neq 0 \in Ext_{A}
^{3,p^{n}q+1}(Z_{p},Z_{p})$,  where $d_{2} :
Ext_{A}^{1,p^{n}q}(Z_{p}, Z_{p})\rightarrow
Ext_{A}^{3,p^{n}q+1}(Z_{p},Z_{p})$ is the secondary differential
in the ASS. This implies that $d_{2}(\tau_{*}(h_{n})) =
\tau_{*}d_{2}(h_{n}) = \tau_{*} (a_{0}b_{n-1}) =
w_{*}(i'i)_{*}(h_{0}b_{n-1})\neq 0 \in Ext_{A}^
{3,p^{n}q+q}(H^{*}X, Z_{p})$ ( cf. (9.7.8)).  This shows that
$\lambda_{0}$ = 0 and so by (9.7.15)  we have $(\bar a_{2}\bar
a_{3}\wedge 1_{X})f''_{1} i'i$ = 0.

It follows that $(\bar a_{3}\wedge 1_{X})f''_{1}i'i = (\bar c_{2}\wedge 1_{X})
g''_{2}$ = 0 , where the $d_{1}$-cycle $g''_{2} \in \pi_{p^{n}q+q}KG_{2}\wedge X$
represents an element in $Ext_{A}^{2,p^{n}q+q}(H^{*}X, Z_{p}) \cong Z_{p}\{w_{*}
(i'i)_{*}\\(h_{0}h_{n})\}$ (cf. Prop. 9.7.9(3)) and the generator of this group
is a permanent cycle in the ASS (cf. Theorem 9.5.1) so that we have
$(\bar c_2\wedge 1_X)g''_2 = 0$.
 Then $f''_{1}  i'i =
(\bar c_{3}\wedge 1_{X})g''_{3} = (\bar c_3\wedge 1_X)g''_4i'i$ for some
$g''_{3} \in \pi_{p^{n}q+q+1}KG_{3}\wedge X$  and $g''_4\in [\Sigma^
{p^nq+q+1}K,KG_3\wedge X]$
, this is because $g''_3\cdot\epsilon$ = 0 which is obtained by the fact that
 $\epsilon : Y\rightarrow
\Sigma S$ induces zero homomorphism in $Z_{p}$-cohomology.
Consequently we have $f''_{1} = (\bar c_{3}\wedge 1_{X})g''_{4} + f''_{2} r$
with $f''_{2} \in [\Sigma^{p^{n}q+q+1}Y, E_{4}\wedge X]$  and $(\bar a_{2}
\bar a_{3}\wedge 1_{X})f''_{1} = (\bar a_{2}\bar a_{3}\wedge 1_{X})f''_{2} r$.
Hence, if $\lambda'$ = 0, (9.7.13)  becomes
$$ (1_{E_{2}}\wedge w)(\eta''_{n,2})_{Y}\cdot r = (\bar a_{2}\bar a_{3}\wedge 1_{X})
f''_{2} r +  (\bar c_{1}\wedge 1_{X})g''_{5} r$$
where $g''_{5} \in [\Sigma^{p^{n}q+q-1}Y, KG_{2}\wedge X]$  such that $ g''_{5} r
= g''$.  Moreover,  by the above equation we have\\
 {\bf (9.7.16)}\qquad $ (1_{E_{2}}\wedge w)(\eta''_{n,2})_{Y} = (\bar a_{2}\bar a_{3}
\wedge 1_{X})f''_{2} + (\bar c_{1}\wedge 1_{X})g''_{5} + f''_{3}\epsilon $\\
with $f''_{3} \in \pi_{p^{n}q+q}E_{2}\wedge X$.  Since $\epsilon : Y\rightarrow \Sigma S$
induces zero homomorphism in  $Z_{p}$-cohomology  , then the right hand side of
(9.7.16) has  filtration $\geq 3$.  However, $(\eta''_{n})_{Y} =  (\bar a_{0}\bar a_{1}\wedge 1_{K})(\eta''_{n,2})_{Y}$
has  filtration 2, this is because it is represented by   $(h_{0}h_{n})''_{Y} \in Ext_{A}
^{2,p^{n}q+q-1}(H^{*}K,H^{*}Y)$ in the ASS.  Moreover,  by the following exact sequence

$ 0 = Ext_{A}^{1,p^{n}q}(H^{*}K,H^{*}Y)\stackrel{\alpha''_{*}}{\longrightarrow}
Ext_{A}^{2,p^{n}q+q-1}(H^{*}K,H^{*}Y)$

$\qquad\qquad\quad\stackrel{w_{*}}{\longrightarrow}Ext_{A}^{2,p^{n}q+q-1}
(H^{*}X,H^{*}Y)\stackrel{(\alpha'')_*}{\longrightarrow}$\\
induced by (9.3.7) we know that  $w_{*}(h_{0}h_{n})''_{Y} \neq 0$, where the left group is zero
by Prop. 9.7.5(3).
 That is to say, $(1_{E_2}\wedge w) (\eta''_{n})_{Y}$
has  filtration 2 which is represented by
 $w_{*}(h_{0}h_{n})''_{Y}$ in the ASS. This shows that the equation (9.7.16) is a contradiction
 so that the scalar $\lambda'$ must be nonzero. Q.E.D.

\vspace{2mm}

{\bf  Lemma 9.7.17}\quad  Let $w : K\rightarrow X$ be the map in the cofibration (9.3.7) and
$W$ is the cofibre of $wi'i : S\rightarrow X$ given by the cofibration $S\stackrel
{wi'i}{\longrightarrow}X\stackrel{w_1}{\longrightarrow}W\stackrel{u_1}
{\longrightarrow}\Sigma S$ , then

\quad

(1)\quad $Ext_A^{s-1,p^nq+q+s-3}(H^*W, Z_p)$ = 0 for  $s = 1,3$ and has unique generator
$(w_1)_*\tau_*(h_n)  = (\tau)^*(w_1)_*[h_n\wedge 1_X]$ for $s = 2$.

(2)\quad $Ext_A^{1,p^nq}(H^*W, H^*X)\cong Z_p\{(w_1)_*(h_n)'_X = (w_1)_*
[h_n\wedge 1_X]\}$.

(3)\quad $(w_1)_*[a_0b_{n-1}\wedge 1_X]\neq 0 \in Ext_A^{3,p^nq+1}(H^*W, H^*X)$.

{\bf  Proof} :  (1) Note that  $W$ also is the cofibre of $r\alpha'' : \Sigma^{q-2}K
\rightarrow Y$, this can be seen by the following homotopy commutative diagram
of $3\times 3$-Lemma.

\quad

\quad

$\qquad\qquad\quad S\qquad\stackrel{wi'i}{\longrightarrow}\qquad X\qquad\stackrel
{u}{\longrightarrow}\quad\Sigma^{q-1}K$

$\qquad\qquad\qquad \searrow i'i\qquad\nearrow w\quad\searrow w_{1}\qquad\nearrow u_{2}$

$\qquad\qquad\qquad\qquad K \quad\qquad\qquad\qquad  W$

$\qquad\qquad\qquad \nearrow \alpha''\qquad\searrow r\quad\nearrow w_{2}\qquad\searrow u_{1}$

$\qquad\qquad \Sigma^{q-2}K\quad \stackrel{r\alpha''}{\longrightarrow}\qquad Y\qquad
\stackrel{\epsilon}{\longrightarrow}\qquad\Sigma S$\\
That is , we have a cofibration\\
$\qquad \Sigma^{q-2}K\stackrel{r\alpha''}{\longrightarrow}Y\stackrel{w_{2}}
{\longrightarrow}W\stackrel{u_{2}}{\longrightarrow}\Sigma^{q-1}K$
and it induces the following exact sequence

$  Ext_{A}^{s-1,p^{n}q+q+s-3}(H^{*}Y,Z_{p})\stackrel{(w_{2})_{*}}
{\longrightarrow}Ext_{A}^{s-1,p^{n}q+q+s-3}(H^{*}W,Z_{p})$

$\qquad\qquad\quad\stackrel{(u_{2})_{*}}{\longrightarrow}
Ext_{A}^{s-1,p^{n}q+s-2}(H^{*}K,Z_{p})\stackrel
{(r\alpha'')_*}{\longrightarrow}$\\
The left group is zero for $s = 1,2,3$ ( cf. Prop. 9.7.5(2)). The right group
also is zero for $s = 1, 3$  (cf. Prop. 9.7.5(1)) and has unique generator
$(i'i)_{*}(h_{n}) = u_{*}\tau_{*}(h_{n}) = (u_{2})_{*}(w_{1})_{*}(\tau)_{*}
(h_{n})$ for $s = 2$. Then the result follows.

(2) Consider the following exact sequence

$ Ext_{A}^{1,p^{n}q}(H^{*}X,H^{*}X)\stackrel{(w_{1})_{*}}{\longrightarrow}
Ext_{A}^{1,p^{n}q}(H^{*}W,H^{*}X)$

$\qquad\qquad\quad\stackrel{(u_{1})_{*}}{\longrightarrow}Ext_{A}^{2,p^{n}q}
(Z_{p},H^{*}X) = 0$\\
The right group is zero by Prop. 9.7.7(1) and by Prop. 9.7.9(2)(1) we know that the left
group has unique generator
$(h_{n})'_{X} = [h_{n}\wedge 1_{X}]$ which satisfies $u_{*}(h_{n})'_{X} = u^{*}(h_{n})'
\in Ext_{A}^{1,p^{n}q-q+1}(H^{*}K,H^{*}X)$
.  Then the middle group has unique generator $(w_{1})_{*}(h_{n})'_{X}$.

(3)  Since $\alpha'' : \Sigma^{q-2}K\to K$ is not an $M$-module map, then , as the
cofibre of $\alpha''$, the spectrum $X$ is not an $M$-moduld spectrum, that is,
the map $p\wedge 1_X\neq 0\in [X,X]$.  So $[a_0\wedge 1_X]
= (p\wedge 1_X)_*[\tau\wedge 1_X]\neq 0\in Ext_A^{1,1}(H^*X,X^*X)$ ( where $\tau$ is the unit in $\pi_0KG_0$),
 and so $[a_0b_{n-1}\wedge 1_X] = (p\wedge 1_X)_*[b_{n-1}\wedge 1_X] = [a_0\wedge 1_X][b_{n-1}\wedge 1_X]
\neq 0\in Ext_A^{3,p^nq+1}(H^*X,X^*X)$ which can be obtained by
by knowledge of Yoneda products and $a_0b_{n-1}\neq 0\in Ext_A^{3,p^nq+1}(Z_p,Z_p)$.
Note to the following exact sequence

$0 = Ext_{A}^{3,p^{n}q+1}(Z_{p},H^{*}X)\stackrel{(wi'i)_{*}}{\longrightarrow}
Ext_{A}^{3,p^{n}q+1}(H^{*}X,H^{*}X)$

$\qquad\qquad\quad\stackrel{(w_{1})_{*}}{\longrightarrow}Ext_{A}^{3,p^{n}q+1}
(H^{*}W,H^{*}X)$\\
where the left group is zero by Prop. 9.7.7(1), then $(w_1)_*$ is monic and so
the reslt follows. Q.E.D.

\vspace{2mm}

{\bf  Remark}\quad The result on $[a_0b_{n-1}\wedge 1_X]\neq 0$ in
Lemma 9.7.17(3) also can be proved by some computation in Ext
groups as follows.   Suppose in contrast that $(p\wedge
1_X)_*[b_{n-1}\wedge 1_X] = [a_0b_{n-1}\wedge 1_X] = 0$, then by
(9.1.1), $[b_{n-1}\wedge 1_X] = (j\wedge 1_X)_*(x_1)$ with $x_1\in
Ext_A^{2,p^nq+1}(H^*M\wedge X,H^*X)$. Recall that $X$ is the
spectrum in (9.3.7), then we have $w^*(1_M\wedge u)_*(x_1)\in
Ext_A^{2,p^nq-q+1}(H^*M\wedge K,H^*K)$ = 0 which can be obtained
by $Ext_A^{2,p^nq-q+1}\\(H^*M\wedge K,H^*M)$ = 0 and
$Ext_A^{2,p^nq+2}(H^*M\wedge K,H^*M)$ = 0. By the Ext exact
sequence induced by (9.3.7) we have $(1_M\wedge u)_*(x_1)\in
u^*Ext_A^{2,p^nq+1}(H^*M\wedge K,H^*K)$. However, this group has
unique generator $(\overline{m}_K)_*(b_{n-1})'$,  then $(1_M\wedge
u)_*(x_1) = \lambda u^* (\overline{m}_K)_*(b_{n-1})'$ for some
$\lambda\in Z_p$.  By applying $(1_M\wedge\alpha'')_*$
 we have $\lambda u^*(1_M\wedge\alpha'')_*(\overline{m}_K)_*(b_{n-1})'$ = 0 and so
 $\lambda (1_M\wedge\alpha'')_*(\overline{m}_K)_*(b_{n-1})'\in (\alpha'')^* Ext_A^{2,p^nq+1}(H^*M\wedge K,H^*K)$.
 Then $\lambda (\alpha_1\wedge 1_K)_*(b_{n-1})' = \lambda (m_K)_*(1_M\wedge \alpha'')_*(\overline{m}_K)_*(b_{n-1})'$ = 0
 which shows that $\lambda = 0$.  So $x_1 = (1_M\wedge w)_*(x_2)$ for some
 $x_2\in Ext_A^{2,p^nq+1}(H^*M\wedge K, H^*X)$.
Similarly we can prove that $w^*(x_2)$ = 0.  Then $x_1\in u^*(1_M\wedge w)_*Ext_A^{2,p^nq+q}(H^*M\wedge K,H^*K)$.
 However, $Ext_A^{2,p^nq+q}(H^*M\wedge K,H^*K)$ has two generators $(i\wedge 1_K)_*(h_0h_n)'$ ,
  $(\overline{m}_K)_*(h_0h_n)''$ and $u^*(h_0h_n)'' = u^*(\alpha'')^*(h_n)'$ = 0, then we have
$x_1 = u^*(1_M\wedge w)_*(i\wedge 1_K)_*(h_0h_n)'$( up to scalar) and so $[b_{n-1}\wedge 1_X]
= (j\wedge 1_X)_*(x_1)$ = 0. This is a contradiction and shows that $[a_0b_{n-1}\wedge 1_X]\neq 0
\in Ext_A^{3,p^nq+1}(H^*X,X^*X)$.

\vspace{2mm}

{\bf Proof of Theorem 9.7.4 :}\quad The result for  $n = 0,1$ is wellknown, then we assume that
 $n\geq 2$. By Lemma 9.7.11  and (9.1.4) we have\\
{\bf (9.7.18)}\qquad $ \lambda' wi'i\zeta_{n-1} = \lambda' w(\zeta_{n-1}\wedge 1_{K})i'i =
- (\bar a_{0}\bar a_{1}\bar a_{2}\bar a_{3}\wedge 1_{X})f''_{1} i'i.$

Moreover, by the cofibration in Lemma 9.7.17 we have $(\bar a_{0}\bar a_{1}\bar a_{2}\bar a_{3}\wedge 1_{W})
(1_{E_{4}}\wedge w_{1})f''_{1} i'i $ = 0 and so $(\bar a_{0}\bar a_{1}\bar a_{2}\bar a_{3}
\wedge 1_{W})(1_{E_{4}}\wedge w_{1})f''_{1} u = k \tau'$ for some $k \in [\Sigma^{p^{n}q-2}X' , W]$
, where $i'i = u \tau : \Sigma^{q-1}S \stackrel{\tau}{\rightarrow}X\stackrel{u}
{\rightarrow}\Sigma^{q-1}K$ which is obtained by $\alpha'' i'i$ = 0 and (9.3.7) and
$\tau' : X\rightarrow X'$
is the map in the following cofibration\\
{\bf (9.7.19)}\qquad $\Sigma^{q-1}S \stackrel{\tau}{\rightarrow}X\stackrel{\tau'}
{\rightarrow}X'\stackrel{\tau''}{\rightarrow}\Sigma^{q}S.$

We claim that $k \in [\Sigma^{p^{n}q-2}X' , W]$  has  filtration $\geq 4$, this can be proved as follows.

By Lemma 9.7.17(1) and (9.7.18) we have $(\tau'')^{*}Ext_{A}^{s-1,p^{n}q+q+s-3}(H^{*}W,\\
Z_{p}) = 0\subset Ext_{A}^{s,p^{n}q+s-2}(H^{*}W,H^{*}X')$ and so
$(\tau')^{*} : Ext_{A}^{s,p^{n}q+s-2}(H^{*}W,
\\H^{*}X')\rightarrow Ext_{A}^{s,p^{n}q+s-2}$ $(H^{*}W, H^{*}X)$
for $ s = 1,2,3$ is monic. Then,  the fact that $k \tau'$ has
filtration$\geq 4$ implies that $k \in [\Sigma^{p^{n}q-2}X' , W]$
also has filtration $\geq 4$. This shows the above claim and so $
k = (\bar a_{0}\bar a_{1} \bar a_{2}\bar a_{3}\wedge 1_{W})k_{3}$
for some $k_{3} \in [\Sigma^{p^{n}q+2}X' , E_{4}\wedge W]$ and
$(\bar a_{0}\bar a_{1}\bar a_{2}\bar a_{3}\wedge
1_{W})(1_{E_{4}}\wedge w_{1})f''_{1}u = (\bar a_{0}\bar a_{1}\bar
a_{2}\bar a_{3}\wedge 1_{W})k_{3}\tau'$.   It follows that
\\{\bf (9.7.20)}\qquad $ (\bar a_{2}\bar a_{3}\wedge 1_{W})(1_{E_{4}}\wedge w_{1})f''_{1}u =
(\bar a_{2}\bar a_{3}\wedge 1_{W})k_{3}\tau' + (\bar c_{1}\wedge 1_{W})\bar g$

$\quad = (\bar a_2\bar a_3\wedge 1_W)k_3\tau' + \lambda_1(\bar c_1
\wedge 1_W)(1_{KG_1}\wedge w_1)(h_n\wedge 1_X)$\\
where the  $d_{1}$-cycle $\bar g = \lambda_1 (1_{KG_1}\wedge w_1)(h_n
\wedge 1_X)\in [\Sigma^{p^{n}q}X, KG_{1}\wedge W]$ with $\lambda_{1}
\in Z_{p}$ which is obtained by Lemma 9.7.17(2).

The equation (9.7.20) means that the  differential
$d_{2}(\lambda_{1}(w_{1})_{*} [h_{n}\wedge 1_{X}])$ = 0. However,
$d_{2}((w_{1})_{*}[h_{n}\wedge 1_{X}]) =
(w_{1})_{*}[a_{0}b_{n-1}\wedge 1_{X}] \neq 0$ ( cf. Lemma
9.7.17(3)). Then the scalar  $\lambda_{1}$ = 0 and we have
 $\bar g$ = 0, $(\bar a_{2}\bar a_{3}\wedge 1_{W})(1_{E_{4}}\wedge w_{1})f''_{1} u
= (\bar a_{2}\bar a_{3}\wedge 1_{W})k_{3}\tau'$  and $(\bar a_{2}\bar a_{3}\wedge 1_{K})
(1_{E_{4}}\wedge u)f''_{1} i'i = (\bar a_{2}\bar a_{3}\wedge 1_{K})(1_{E_{4}}\wedge
u_{2}w_{1})f''_{1} u\tau$ = 0.  Consequently we have $(\bar a_{3}\wedge 1_{K})(1_{E_{4}}
\wedge u)f''_{1}i'i = (\bar c_{2}\wedge 1_{K})\bar g_{2} = 0$ , this is because
the   $d_{1}$-cycle
$\bar g_{2} \in \pi_{p^{n}q+1}KG_{2}\wedge K$ represents an element in
$Ext_{A}^{2,p^{n}q+1}(H^{*}K,Z_{p})$ = 0 (cf. 9.7.5(1)) .  Then, $(1_{E_{4}}\wedge u)f''_{1}i'i
= (\bar c_{3}\wedge 1_{K})\bar g_{3}$ for some $\bar g_{3} \in [\Sigma^{p^{n}q+2}S, KG_{3}\wedge K]$.
Since $(1_{KG_{3}}\wedge \alpha'')\bar g_{3}$ = 0, then $\bar g_{3} = (1_{KG_{3}}\wedge u)
\bar g_{4}$ for some $\bar g_{4} \in [\Sigma^{p^{n}q+q+1}S, KG_{3}\wedge X]$ and so we have
 $(1_{E_{4}}\wedge u)f''_{1} i'i = (\bar c_{3}\wedge 1_{K})(1_{KG_{3}}
\wedge u)\bar g_{4}$ , $f''_{1} i'i = (\bar c_{3}\wedge 1_{X})\bar g_{4} +
(1_{E_{4}}\wedge w)\bar f_{2}$ with $\bar f_{2} \in [\Sigma^{p^{n}q+q+1}S,
E_{4}\wedge K]$.

Hence, by (9.7.18) we have $ \lambda' wi'i\zeta_{n-1} = - (\bar a_{0}\bar a_{1}\bar a_{2}\bar a_{3}\wedge 1_{X})
f''_{1}i'i = - (\bar a_{0}\bar a_{1}\bar a_{2}\\\bar a_{3}\wedge 1_{X})(1_{E_{4}}\wedge w)\bar f_{2}$
and by (9.3.7), $\lambda' i'i \zeta_{n-1} = - (\bar a_{0}\bar a_{1}\bar a_{2}\bar a_{3}
\wedge 1_{K})\bar f_{2} + \alpha'' \omega_{n}$ with $\omega_{n} \in \pi_{p^{n}q-1}K$.
Since $\lambda'i'i\zeta_{n-1}$ is a map of filtration 3 which is represented by
$\lambda' (i'i)_{*}(h_{0}b_{n-1}) \in Ext_{A}^{3,p^{n}q+q}(H^{*}K,Z_{p})$ in the ASS,
then $\alpha''\omega_{n}$
has  filtration 3  and so $\omega_{n} \in \pi_{p^{n}q-1}K$  has filtration $\leq 2$.
 However, by Prop. 9.7.5(1) we have $Ext_{A}^{2,p^{n}q+1}(H^{*}K,Z_{p})$ = 0 , then
$\omega_{n} \in \pi_{p^{n}q-1}K$  must be represented by the unique generator
$(i'i)_{*}(h_{n}) \in Ext_{A}^{1,p^{n}q}(H^{*}K,Z_{p})$ (up to nonzero scalar).
This shows the Theorem. Q.E.D.

\vspace{2mm}

{\bf  Remark}\quad The element $\omega_{n} \in \pi_{p^{n}q-1}K$
obtained in Theorem 9.7.4 can be extended to $(\omega_{n})' \in
[\Sigma^{p^{n}q-1}K,K]$ such that $(\omega_{n})'i'i = \omega_{n}$.
Then, $(\omega_{n})'$ is represented by $(h_{n})' \in
Ext_{A}^{1,p^{n}q}(H^{*}K,H^{*}K)$ in the ASS and
$\alpha''(\omega_{n})'$, $(\omega_{n})'\alpha'' \in
[\Sigma^{p^{n}q+q-3}K,K]$  is represented by $\alpha''_{*}(h_{n})'
= (\alpha'')^{*}(h_{n})' = (h_{0}h_{n})'' \\\in
Ext_{A}^{2,p^{n}q+q-1}(H^{*}K,H^{*}K)$. By Theorem 9.7.4 and Lemma
9.7.11  we have $\alpha''(\omega_{n})' = (\omega_{n})'\alpha'' +
\lambda' \zeta_{n-1}\wedge 1_{K}$ (modulo higher filtration).

By [10] p.511, there is a map $\phi_{*}\phi : BP_{*}BP \rightarrow
A_{*}$ such that $t_{n} \mapsto$  the conjugate of $\xi_{n}$,
where $A_{*} = E[\tau_{0},\tau_{1},\tau_{2},...]\otimes P[\xi_{1},
\xi_{2}, ... ]$ is the dual of the Steenrod algebra $A$.  Then
$\phi_{*}\phi$  induces the Thom  map $\Phi :
Ext_{BP_{*}BP}^{1,p^{n}q}(BP_{*},\\BP_{*}K) \rightarrow
Ext_{A}^{1,p^{n}q} (H^{*}K,Z_{p})$  such that the image of $h_{n}
\in Ext_{BP_{*}BP}^{1,p^{n}q}(BP_{*},\\BP_{*}K)$ is $\Phi (h_{n})
= (i'i)_{*}(h_{n}) \in Ext_{A}^{1,p^{n}q}(H^{*}K,Z_{p})$. Then,
the element $\omega_{n} \in \pi_{p^{n}q-1}K$ obtained in
Theorem 9.7.4 is represented by $h_{n}$  +  other terms $\in Ext_{BP_{*}BP}^{1,p^{n}q}(BP_{*},\\
BP_{*}K)$  in the Adams-Novikov spectral sequence. To know what the elements in the other terms,
we first prove the following Lemma.

\vspace{2mm}

{\bf  Lemma 9.7.21}\quad  By degree reason, $Ext_{BP_{*}BP}^{1,p^{n}q}(BP_{*},BP_{*}K)$
is generated (additively) by the following $v_{2}$-torsion elements $c_{2}(p^{n-2})$ and
 $v_{2}$-torsion free elements $h_{n}, v_{2}^{p^{n-2}(p-1)}h_{n-2} , v_{2}^{a_{i}p^{i}}
h_{i}$ , where $ i\geq 0, a_{i} = (p^{2k}-1)/(p+1)$ , $n - i = 2k \geq 4$.   In addition,
there is a relation $h_{n} = c_{2}(p^{n-2}) + v_{2}^{p^{n-2}(p-1)}h_{n-2} \in
Ext_{BP_{*}BP}^{1,p^{n}q}(BP_{*},BP_{*}K)$.

{\bf  Proof} :   By [19]  Theorem 1.1  and 1.5, $Ext^{1,*}(BP_{*},BP_{*}K)$
 is a  $Z_{p}[v_{2}]$-module which is generated by $v_{2}$-torsion elements $c_{2}(ap^{s})$
and $v_{2}$-torsion free elements $w_{2}, h_{i}$, where  $a \neq 0$ (mod $p$), $s \geq 0$
and $i \geq 0$.  Moreover, the internal degree $\mid h_{i}\mid = p^{i}q$, $\mid c_{2}(ap^{s})
\mid = ap^{s}(p^{2}+p+1)q - q(ap^{s})(p+1)q$  and $\mid w_{2}\mid = (p+1)^{2}q$.

Since $\mid v_{2}^{b}w_{2}\mid \equiv 0$ (mod $(p+1)q)$,  then $\mid v_{2}^{b}w_{2}
\mid \neq p^{n}q$.  If $\mid v_{2}^{b}h_{i}\mid = p^{n}q$,  then $b(p+1)q + p^{i}q
= p^{n}q, b(p+1) = p^{i}(p^{n-i}-1)$  and so $(p^{n-i}-1)$  must be divisible by $p+1$.
 Hence $b = 0$  and $i = n$  or $ b = a_{i}p^{i}$ with $a_{i} = (p^{2k}-1)/(p+1)$
and $n-i = 2k \geq 2$.  Then $h_{n}, v_{2}^{p^{n-2}(p-1)}h_{n-2}$  and $v_{2}^{a_{i}p^{i}}h_{i}
(0 \leq i < n-2$)  are the only torsion free elements of  $Ext^{1,p^{n}q}(BP_{*},BP_{*}K)$.

 If $\mid v_{2}^{b}c_{2}(ap^{s})\mid = p^{n}q$,  then $p^{n}q = ap^{s}(p^{2}+p+1)q
- (q(ap^{s})-b)(p+1)q$, $ap^{s}(p^{2}+p+1) = p^{n} + (q(ap^{s}) - b)(p+1)$  and so
the right hand side must be divisible by $p^{2}+p+1$.  So we have\\
{\bf (9.7.22)}\qquad $ap^{s} = p^{n-2} + \frac{(q(ap^{s})-b-p^{n-2})(p+1)}{p^{2}+p+1}. $\\
We claim that  $s \leq n-2$ which will be proved below , then $q(ap^{s})-b$  must be divisible by
$p^{s}$.  However,  by [19] p.132, $q(ap^{s}) = p^{s}$ for $a = 1$  and $q(ap^{s}) = p^{s}$
+  other terms  for $a \geq 2$. Then, the only possibility is  $q(ap^{s})
= p^{s}, a = 1, b = 0$  and $s = n-2$. That is to say, the only $v_2$-torsion elements in
$Ext^{1,p^{n}q}(BP_{*},BP_{*}K)$ is  $c_{2}(p^{n-2})$.

 Now we prove the above claim.   Suppose in contrast that $s \geq n-1$ , then, by (9.7.22) we have
 $\frac{(q(ap^{s})-b-p^{n-2})(p+1)}{p^{2}+p+1} = ap^{s} - p^{n-2} \geq p^{s} - p^{n-2} ,$
$$2p^{s} > q(ap^{s})-b-p^{n-2} \geq \frac{(p^{s}-p^{n-2})(p^{2}+p+1)}{p+1} >
p^{s+1} - p^{n-1}$$
and this is a contradiction which shows the above claim. Q.E.D.

\vspace{2mm}

{\bf Proof of Theorem 9.7.1} \quad  For the Thom $\Phi :
Ext_{BP_*BP}^{1,p^nq} (BP_{*},BP_{*}K)\\\rightarrow
Ext_{A}^{1,p^{n}q}(H^{*}K,Z_{p})$ we have $\Phi (h_{n}) =
(i'i)_{*}(h_{n})$. By this we know that the element $\omega_{n}
\in \pi_{p^{n}q-1}K$ obtained in Theorem 9.7.4 is represented by
$h_{n}$ + ( other terms) $\in
Ext_{BP_{*}BP}^{1,p^{n}q}(BP_{*},BP_{*}K)$ in the Adams-Novikov
spectral sequence.  By Lemma 9.7.21, the other terms are the
linear combination of $v_{2}^{p^{n-2}(p-1)}h_{n-2}$  and
$v_{2}^{a_{i}p^{i}}h_{i}$, where  $ i\geq 0, n-i = 2k \geq 4$  and
$a_{i} = (p^{2k}-1)/(p+1)$. Let $\beta \in [\Sigma^{(p+1)q}K,K]$
be the known $v_{2}$-map,  then $i'i\alpha_{1} \in \pi_{q-1}K$ and
$i'j'\beta i'i \in \pi_{pq-1}K$ is represented by $h_{0}, h_{1}
\in Ext_{BP_{*}BP}^{1,*}(BP_{*},BP_{*}K)$ respectively.
 That is, $h_{0}, h_{1} \in Ext_{BP_{*}BP}^{1,*}(BP_{*},BP_{*}K)$ are permanent cycles in the
 Adams-Novikov spectral sequence.  Suppose inductively that $h_i\in Ext_{BP_*BP}^{1,p^iq}(BP_*,BP_*K)$ for
$i\leq n-1 (n\geq 2)$ are permanent cycles in the Adams-Novikov spectral sequence.  Since
$\omega_n\in \pi_{p^nq-1}K, \omega_{n+1}\in \pi_{p^{n+1}q-1}K$ are represented by the linear combination of\\
\centerline{$h_n + v_2^{p^{n-2}(p-1)}h_{n-2}$ and $v_2^{a_ip^i}h_i\in
Ext_{BP_*BP}^{1,p^nq}(BP_*,BP_*K),$}\\
\centerline{$h_{n+1} + v_2^{p^{n-1}(p-1)}h_{n-1}$ and $v_2^{a_ip^i}h_i \in Ext_{BP_*BP}^{1,p^{n+1}q}
(BP_*,BP_*K)$}\\
 then $h_n,h_{n+1}\in Ext_{BP_*BP}^{1,*}(BP_*,BP_*K)$ also are permanent cycles.
 This completes the induction and the result of the Theorem follows. Q.E.D.
 \vspace{2mm}

{\bf  Conjecture 9.7.22}\quad  Theorem 9.7.4 can be generalzed to
be the following general result.  Let $p\geq 5, s\leq 4$ ,
$Ext_A^{s,tq}(Z_p,Z_p)\cong Z_p\{x\}, Ext_A^{s+1,tq+q}\\(Z_p,Z_p)
\cong Z_p\{h_0x\}, Ext_A^{s+2,tq+2q+1}(Z_p,Z_p)\cong
Z_p\{\widetilde{\alpha}_2x\}$ and some supposition on vanishes of
some Ext groups. \quad If the secondary differential $d_2(x) =
a_0x'\in \\Ext_A^{s+2,tq+1}(Z_p,Z_p)$ with $x'\in
Ext_A^{s+1,tq}(Z_p,Z_p)$, that is, $x$ and $x'$ is a pair of
$a_0$-related elements,
 then there exists $\omega\in \pi_{tq-s}K$
 such that $i'i\xi = \alpha''\cdot \omega$ (mod $F^{s+2}\pi_*K)$ and $\omega\in\pi_{tq-s}K$
 is represented by $(i'i)_*(x)\in Ext_A^{s,tq}(H^*K,Z_p)$ in the ASS, where
 $\xi\in \pi_{tq+q-s-2}S$ is the homotopy element which is represented by
 $h_0x'\in Ext_A^{s+2,tq+q}(Z_p,Z_p)$ in the ASS and $F^{s+2}\pi_*K$ denotes
  the group consisting of all elements in $\pi_*K$
filtration $\geq s+2$.

\quad

\begin{center}

{\bf \large \S 8.\quad Second periodicity families in the stable homotopy groups of spheres}

\end{center}

\vspace{2mm}

 By Theorem 8.1.2 in chapter 8, $Ext_{BP_*BP}^{1,*}(BP_*,BP_*)$ is generated by
 $\alpha_{tp^n/n+1} (n\geq 0, p$
 not divisible by $t\geq 1$) and It was proved by Novikov that all these first periodicity families
 converge to the im $J\subset \pi_*S$.
In this section, using the $h_0h_{n+1}$-element obtained in Theorem 9.5.1 and
the elements $\beta_{p/r}, 1\leq r\leq p-1$
 and $\beta_{tp/r}, t\geq 2, 1\leq r\leq p$ as our geometric input,
 we prove the following Theorem on the convergence of second periodicity families
 $\beta_{tp^n/r}$ in the Adams-Novikov spectral sequence.
 \vspace{2mm}

{\bf  Theorem 9.8.1}\quad  Let $p\geq 5, n\geq 1, 1\leq s\leq p^{n}-1$  if $t\geq 1$ is not divisible by
$p$  or $1\leq s\leq p^{n}$ if $t\geq 2$ is not divisible by $p$ , then The elements\\
\centerline{$\beta_{tp^{n}/s} \in Ext_{BP_{*}BP}^{2,tp^{n}(p+1)q-sq}
(BP_{*},BP_{*})$}\\
in Theorem 8.1.3 are permanent cycles in the Adams-Novikov spectral sequence and they
converge to the corresponding homotopy elements of order $p$ in
$\pi_{tp^{n}(p+1)q-sq-2}S$.

We will prove Theorem 9.8.1 in case $t\geq 1$ or $t\geq 2$ separately. The proof will be done
by some arguments processing in the cannical Adams-Novikov resolution.
We first do some preminilaries as follows.

Let $M$  be the Moore spectrum whose $BP_*$-homology are $BP_*(M) = BP_*/(p)$. Let
$\alpha : \Sigma^{q}M \rightarrow M$ be the Adams map which induces $BP_*$-homomorphisms
are $v_{1} : BP_{*}/(p)\rightarrow
BP_{*}/(p)$.   Let $K_{r}$  be the cofibre of
$\alpha^{r} : \Sigma^{rq}M \rightarrow M$ given by the cofibration\\
{\bf (9.8.2)}$\qquad\qquad \Sigma^{rq}M \stackrel{\alpha^{r}}{\longrightarrow}M\stackrel
{i'_{r}}{\longrightarrow}K_{r}\stackrel{j'_{r}}{\longrightarrow}\Sigma^{rq+1}M$\\
 The cofibration (9.8.2) induces a short exact sequence of $BP_*$-homology\\
\centerline{$0\rightarrow BP_{*}/(p)\stackrel{v_{1}^{r}}{\longrightarrow}
BP_{*}/(p)\longrightarrow BP_{*}/(p,v_{1}^{r})\rightarrow 0$}\\
Recall from \S 5 in chapter 6, $K_{r}$ is a  $M$-module spectrum and we have the following
derivations\\
{\bf (9.8.3)}\qquad $d(i'_{r})$ = 0,\quad $d(j'_{r})$ = 0,\quad $d(\alpha) = 0,\quad
d(ij) = - 1_M$.

\quad

Moreover, the cofibre of $i'_{s}j'_{r} : K_{r}\rightarrow \Sigma^{rq+1}K_{s}$ is $\Sigma K_{r+s}$
given by the cofibration\\
{\bf (9.8.4)} $\qquad \Sigma^{rq}K_{s}\stackrel{\psi_{s,s+r}}{\longrightarrow}
K_{s+r}\stackrel{\rho_{s+r,r}}{\longrightarrow} K_{r}\stackrel
{i'_{s}j'_{r}}{\longrightarrow}\Sigma^{rq+1}K_{s}$\\
This can be seen by the following homotopy commutative diagram of $3\times 3$-Lemma

$\qquad\qquad K_r\quad\stackrel{i'_sj'_r}{\longrightarrow}\quad\Sigma^{rq+1}K_s\quad
\stackrel{j'_s}{\longrightarrow}\quad \Sigma^{(r+s)q+2}M$

$\qquad\qquad\quad\searrow j'_r\quad\nearrow i'_s\quad\searrow \psi_{s,r+s}\quad\nearrow j'_{r+s}\quad\searrow \alpha^s$\\
{\bf (9.8.5)}\qquad\qquad\quad $\Sigma^{rq+1}M\qquad\qquad\Sigma K_{r+s}\qquad\qquad\qquad\Sigma^{rq+2}M$

$\qquad\qquad\quad\nearrow\alpha^s\quad\searrow\alpha^r\quad\nearrow i'_{r+s}\quad
\searrow\rho_{r+s,r}\quad\nearrow j'_r$

$\qquad\quad\Sigma^{(r+s)q+1}M\stackrel{\alpha^{r+s}}{\longrightarrow}\Sigma M\quad
\quad\stackrel{i'_r}{\longrightarrow}\qquad\Sigma K_r$\\
 Moreover, the cofibration (9.8.4) induces a short exact sequence of $BP_*$-homology
$$0\rightarrow
BP_{*}/(p,v_{1}^{s})\stackrel{v_{1}^{r}}{\longrightarrow}
BP_{*}/(p,v_{1}^{s+r})\longrightarrow BP_{*}/(p,v_{1}^{r})\rightarrow 0$$
and by the homotopy commutative diagram (9.8.5) , we have the following relations\\
{\bf (9.8.6)} $\qquad \psi_{s,s+r}i'_{s} = i'_{s+r}\alpha^{r},\qquad j'_{r}
\rho_{s+r,r} = \alpha^{s}j'_{s+r}$

$\qquad\quad j'_{s+r}\psi_{s,s+r} = j'_{s}, \qquad \rho_{s+r,r}i'_{s+r}
= i'_{r}$.

\vspace{2mm}

{\bf  Proposition 9.8.7}\quad  Let $p\geq 5$  and $f \in [\Sigma^{t}K_{r} ,
S]$ be any map,  then $f = jj'_{r}\overline{f}$ for some $\overline{f}
\in [\Sigma^{t+rq+2}K_{r},K_{r}]$.

{\bf  Proof} :  By Theorem 6.5.16(A) in chapter 6,  there is $\nu_{r} : \Sigma^{rq+2}K_{r}
\rightarrow K_{r}\wedge K_{r}$ such that $(jj'_{r}\wedge
1_{K_{r}})\nu_{r} = 1_{K_{r}}$.  Let $K'_{r}$ be the cofibre of
$jj'_{r} : \Sigma^{-1}K_{r}\rightarrow \Sigma^{rq+1}S$ given by the cofibration
 $\Sigma^{-1}K_{r}\stackrel{jj'_{r}}{\longrightarrow}
\Sigma^{rq+1}S\stackrel{z_{r}}{\longrightarrow}K'_{r}\rightarrow
K_{r}$, then $z_{r}\wedge 1_{K_{r}} = (z_{r} jj'_{r}\wedge 1_{K_{r}})
\nu_{r}$ = 0 $\in [\Sigma^{rq+1}K_{r} , K'_{r}\wedge K_{r}]$.  Consequently,
$z_{r}f = (1_{K'_{r}}\wedge f)(z_{r}\wedge 1_{K_{r}})$ = 0  and so
$f = jj'_{r}\overline{f}$ for some $\overline{f} \in [\Sigma^{t+rq+1}
K_{r}, K_{r}]$.  Q.E.D.

Let

{\bf (9.8.8)}$\quad\qquad\cdots\stackrel{\tilde a_{2}}{\longrightarrow}
\quad\Sigma^{-2}\widetilde{E}_{2}\stackrel{\tilde a_{1}}{\longrightarrow}\quad
\Sigma^{-1}\widetilde{E}_{1}\quad\stackrel{\tilde a_{0}}{\longrightarrow} \qquad\widetilde{E}_{0} = S$

$\quad\qquad\qquad\qquad\qquad\qquad\big\downarrow\tilde b_{2}\qquad\qquad
\quad\big\downarrow \tilde b_{1}\qquad\qquad\qquad \big\downarrow \tilde b_{0}$

$\qquad\qquad\qquad\qquad\Sigma^{-2}BP\wedge
\widetilde{E}_{2}\qquad\Sigma^{-1}BP\wedge \widetilde{E}_{1}\qquad BP\wedge \widetilde{E}_{0} = BP$\\
be the canonnical Adams-Novikov resolution of the sphere spectrum $S$, where
$\widetilde{E}_{s}\stackrel{\tilde b_{s}}{\longrightarrow}BP\wedge \widetilde{E}_{s}\stackrel{\tilde
c_{s}}{\longrightarrow}\widetilde{E}_{s+1}\stackrel{\tilde a_{s}}{\longrightarrow}
\Sigma \widetilde{E}_{s}$ are cofibrations for all  $s\geq 0$ such that $\widetilde{E}_{0} = S,
\tilde b_{s} = \tau\wedge 1_{\widetilde{E}_{s}}$( $s > $0) and $\tilde b_{0} = \tau
: S\rightarrow BP$ is the injection of the bottom cell.   Then $\pi_{t}
BP\wedge \widetilde{E}_{s}$ is the  $E_{1}^{s,t}$-term of the Adams-Novikov spectral sequence,
$(\tilde b_{s+1}\tilde
c_{s})_{*} : \pi_{t}BP\wedge \widetilde{E}_{s}\rightarrow \pi_{t}BP\wedge \widetilde{E}_{s+1}$
are the $d_{1}^{s,t}$-differential and $$E_{2}^{s,t} =
Ext_{BP_{*}BP}^{s,t}(BP_{*},BP_{*}) \Longrightarrow (\pi_{t-s}S)_{p}$$

\vspace{2mm}

{\bf Proposition 9.8.9}\quad Let $p\geq  3, r\geq 1, s\geq 0$ and
$\widetilde{E}_{s}$  be the spectrum in the Adams-Novikov resolution (9.8.8),
$\wedge^{s}BP = BP\wedge\cdots\wedge BP$ be the smash products of $s$ copies of $BP$,
then $(BP\wedge \widetilde{E}_{s})^{*}, (BP\wedge \widetilde{E}_{s})^{*}(M),
(BP\wedge \widetilde{E}_{s})^{*}(K_{r})$ are the direct summand of $(\wedge^{s+1}BP)^{*}, (\wedge
^{s+1}BP)^{*}(M),(\wedge ^{s+1}BP)^{*}(K_{r})$ and we have
$[\Sigma^{t}M, BP\wedge \widetilde{E}_{s}] = (BP\wedge \widetilde{E}_{s})^{-t}(M)$ = 0
for $t\neq -1 $ (mod q), $[\Sigma^{t}K_{r}, BP\wedge \widetilde{E}_{s}] = (BP
\wedge \widetilde{E}_{s})^{-t}(K_{r})$ = 0 for $t \neq -2$ (mod q).

{\bf Proof} :  We first consider the  $BP^{*}$-cohomology.  It is known that
$\pi_{t}BP = BP_{t} = BP^{-t}$,  then $BP^{*} = Z_{(p)}[v_{1},v_{2},
\cdots ]$ , where $\mid v_{i}\mid = - 2(p^{i}-1)$  and $I_{n} = (p, v_{1},
\cdots, v_{n-1}), (p,v_{1}^{r})$  is the invariant ideal of  $BP^{*}$.  Clearly,
there are two exact sequences on $BP^*$-cohomology as follows
$$\qquad 0\rightarrow BP^{*}\stackrel{p}{\longrightarrow}BP^{*}\stackrel
{\rho_{0}}{\longrightarrow}BP^{*}/(p)\rightarrow 0$$
$$\qquad 0\rightarrow \Sigma^{-rq}BP^{*}/(p)\stackrel{v_{1}^{r}}{\longrightarrow}
BP^{*}/(p)\stackrel{\rho_{1}}{\longrightarrow}BP^{*}/(p,v_{1}^{r})
\rightarrow 0$$
where $\rho_{0}, \rho_{1}$ are the projections.

Note that $(\wedge^{s}BP)^{*} = \pi_{*}(\wedge^{s}BP) =
BP_{*}(\wedge^{s-1}BP) = BP_{*}BP\otimes \cdots\otimes BP_{*}BP$ with
$s-1$ copies of $BP_{*}BP$ and $s\geq 2$.  Then we have the follwing short exact sequences
($s\geq 1$)\\
{\bf (9.8.10)} $\qquad 0\rightarrow (\wedge ^{s}BP)^{*}\stackrel{p}{\longrightarrow}
(\wedge^{s}BP)^{*}\longrightarrow (\wedge^{s}BP)^{*}/(p)\rightarrow 0$

$\quad 0\rightarrow \Sigma^{-rq}(\wedge^{s}BP)^{*}/(p)\stackrel
{v_{1}^{r}}{\longrightarrow}(\wedge^{s}BP)^{*}/(p)\longrightarrow
(\wedge^{s}BP)^{*}/(p,v_{1}^{r})\rightarrow 0$

 For any $f \in [\Sigma^{t}M, \wedge^{s}BP] = (\wedge^{s}BP)^{-t}(M)$, if $t\neq 0 $(mod q)
,  then by the sparseness of $(\wedge^sBP)^* = BP_*BP\otimes\cdots\otimes BP_*BP$
( that is, $BP_rBP$ = 0 for $ r \neq 0$ (mod $q$)) we have $fi \in (\wedge^{s}BP)^{-t}$ = 0
;  if $t = 0 $ (mod q), then $fi$ is an element of order $p$ in
$Z_{(p)}$-module $(\wedge^{s}BP)^{-t}$ so that we have $fi$ = 0.  This shows that
 $i^{*} = 0 : (\wedge^{s}BP)^{*}(M)\rightarrow
(\wedge^{s}BP)^{*}$. Similarly we have  $(i'_{r})^{*} = 0 : (\wedge^{s}BP)^{*}(K_{r})
\rightarrow (\wedge^{s}BP)^{*}(M)$ ( $r\geq 1$). Then, the cofibration
(9.1.1) (9.8.2) induces respectively the following short exact sequences for all $(r\geq 1)$
$$\qquad 0\rightarrow (\wedge^{s}BP)^{*}\stackrel{p}{\longrightarrow}
(\wedge^{s}BP)^{*}\stackrel{j^{*}}{\longrightarrow}(\wedge^{s}BP)^{*}
(M)\rightarrow 0 $$
$$\qquad 0 \rightarrow (\wedge^{s}BP)^{*}(M)\stackrel{(\alpha^{r})^{*}}
{\longrightarrow}(\wedge^{s}BP)^{*}(M)\stackrel{(j'_{r})^{*}}{\longrightarrow}
(\wedge^{s}BP)^{*}(K_{r})\rightarrow 0$$
where the degrees $\mid j^{*}\mid = -1$, $\mid (j'_{r})^{*}\mid =
-(rq+1)$.  By comparison to the above two short exact sequences with (9.8.10) we have\\
{\bf (9.8.11)} $\qquad (\wedge^{s}BP)^{*}(M)\cong\Sigma (\wedge^{s}BP)^{*}/(p),$

$\qquad\quad (\wedge^{s}BP)^{*}(K_{r})\cong\Sigma^{rq+2}(\wedge^{s}BP)^{*}/
(p,v_{1}^{r})$.

Let $\mu : BP\wedge BP\rightarrow BP$ be the multiplication of the ring spectrum $BP$ and
 $\tau : S\rightarrow BP$ be the injection of the bottom cell,
,  then we have $\mu (1_{BP}\wedge\tau ) = 1_{BP} = \mu (\tau \wedge
1_{BP})$ so that the cofibration $\widetilde{E}_{s-1}\stackrel{\tilde
b_{s-1}}{\longrightarrow}BP\wedge \widetilde{E}_{s-1}\stackrel{\tilde c_{s-1}}
{\longrightarrow}\widetilde{E}_{s}\stackrel{\tilde a_{s-1}}{\longrightarrow}
\Sigma\widetilde{E}_{s-1}$ induces a split short exact sequence
$$ BP\wedge \widetilde{E}_{s-1}\stackrel{1_{BP}\wedge\tilde b_{s-1}}{\longrightarrow}
BP\wedge BP\wedge \widetilde{E}_{s-1}\stackrel{1_{BP}\wedge \tilde c_{s-1}}
{\longrightarrow}BP\wedge \widetilde{E}_{s}$$
this is because  $(\mu \wedge 1_{\widetilde{E}_{s-1}})(1_{BP}\wedge\tilde b_{s-1}) =
(\mu\wedge 1_{\widetilde{E}_{s-1}})(1_{BP}\wedge \tau\wedge 1_{\widetilde{E}_{s-1}}) =
1_{BP\wedge \widetilde{E}_{s-1}}$.
That is to say, $BP\wedge \widetilde{E}_{s}$ is the direct summand of  $BP\wedge BP\wedge
\widetilde{E}_{s-1}$ and by induction we have $BP\wedge \widetilde{E}_{s}$ is
the direct summand of
 $\wedge^{s+1}BP$.   Hence, $(BP\wedge \widetilde{E}_{s})^{*}, (BP
 \wedge \widetilde{E}_{s})^{*}(M), (BP\wedge \widetilde{E}_{s})^{*}(K_{r})$
  are the direct summand of $(\wedge^{s+1}BP)^{*}, (\wedge^{s+1}BP)^{*}(M), (\wedge^{s+1}BP)^{*}
(K_{r})$ respectively and the last result can be obtained by (9.8.11).  Q.E.D.

\vspace{2mm}

{\bf  Proposition 9.8.12} \quad  Let $p \geq 3, n\geq 1$, then

$\qquad\qquad Ext_{BP_{*}BP}^{0,p^{n}(p+1)q}(BP_{*},BP_{*}/(p,v_{1}
^{p^{n}-1}))$\\
is generated additively by the generators $v_{2}^{p^{n}},
v_{1}^{p^{n}-p^{n-2r}}\tilde{c}_{1}(t_{r}p^{n-2r})$ ( $r\geq 1$),
where $t_{r} = (p^{2r+1}+1)/(p+1)$  and $\tilde{c}_{1}(ap^{s})$  is the generator
in Theorem 8.1.7 in chapter 8 which has degree $sp^s(p+1)q$.

{\bf Proof} :  By Theorem 8.1.7, the desired generators are of the form
 $v_{1}^{b}\tilde{c}_{1}(ap^{s})$ with degrees $bq + ap^{s}
(p+1)q = p^{n}(p+1)q$, $ a\geq 1$ is not divisible by $p$, $0\leq
b < p^{n}-1$  and $b \geq $ max $\{0, p^{n}-1-q_{1}(ap^{s})\}$,
where $q_{1}(ap^{s}) = p^{s}$ if $ a = 1$, $q_{1}(ap^{s}) = p^{s}
+ p^{s-1} -1$  if
 $ a\geq 2$ is not divisible by $p$.

 If $b = 0$,
then the generator is $v_{2}^{p^{n}}$.  Since $b < p^{n}-1$,  then $s < n$ and
$b \equiv 0$ ( mod $p^{s}$) and so $ b\geq p^{n}-1-q_{1}(ap^{s})\geq
p^{n}-p^{n-1}$ if $b\geq 1$.  Let $b = (p-1)p^{n-1} + c_{n-2}p^{n-2}
+ \cdots + c_{s}p^{s}$ be the p-adic expasion of $b$ such that  $0\leq
c_{i}\leq p-1$.  By $b\geq (p^{n}-1)-(p^{s}+p^{s-1}-1)$ or
$b\geq (p^{n}-1)-p^{s}$  we have $c_{n-2}p^{n-2} + \cdots + c_{s}p^{s}
\geq p^{n-1}-p^{s}-p^{s-1}$  or $p^{n-1}-p^{s}-1$.  Consequently we have $c_{n-2}
= \cdots = c_{s} = p-1$.  On the other hand, $ b$  is divisible by $p+1$,  then
$(p-1) - c_{n-2} + \cdots + (-1)^{n-1-s}c_{s}$  = 0  so that $n-1-s$  must be odd.
 Let $n-1-s = 2r-1$,  then we have $s = n-2r, b = p^{n}-p^{n-2r}$ as desired and
 $a = (p^{2r+1}+1)/(p+1)$.  Q.E.D.

\vspace{2mm}

{\bf  Proposition 9.8.13}\quad  Let $p\geq 3, n\geq 1$,  then

(1) $Ext_{BP_{*}BP}^{2,p^{n+1}q+q}(BP_{*},BP_{*})$ is generated additively by
the generators \\$\beta_{p^{n}/p^{n}-1}, \beta_{t_{r}p^{n-2r}/
p^{n-2r}-1}$  for all $r\geq 1$), where $t_{r} = (p^{2r+1}+1)/(p+1)$.

(2) $Ext_{BP_{*}BP}^{1,p^{n+1}q+q}(BP_{*},BP_{*}M)$ is generated additively by
the generators $\beta'_{p^{n}/p^{n}-1}, \beta'_{t_{r}p^{n-2r}/p^{n-2r}-1} $ for all
$r\geq 1$), where $t_{r} = (p^{2r+1}+1)/(p+1)$ and $\beta'_{tp^{n}/s}$ is the generator in
$Ext_{BP_{*}BP}^{1,*}(BP_{*},BP_{*}M)$ such that  $j_{*}(\beta'_{t
p^{n}/s}) = \beta_{tp^{n}/s} \in Ext_{BP_{*}BP}^{2,*}(BP_{*},BP_{*})$.

{\bf Proof} :  By Theorem 8.1.3 in chapter 8,  $Ext_{BP_{*}BP}^{2,*}(BP_{*},
BP_{*})$ is generated additively by the generators $\beta_{ap^{s}/b,c+1}
\in Ext_{BP_{*}BP}^{2,ap^{s}(p+1)q-bq}(BP_{*},BP_{*})$ , where $s\geq
0$,  $ a\geq 1$ is not divisible by $p$, $b\geq 1, c\geq 0$ and subject to

(i)  \quad $b\leq s$  if $a$ = 1.

(ii) \quad $  p^{c}\mid b \leq p^{s-c}+p^{s-c-1} -1$

(iii)\quad $p^{s-c-1} + p^{s-c-2} -1 < b $  if $p^{c+1} \mid b$,\\
 and $\beta_{ap^{s}/b,1} = \beta_{ap^{s}/b}$. Then, for
$\beta_{ap^{s}/b,c+1}\in Ext_{BP_{*}BP}^{2,p^{n+1}q+q}(BP_{*},BP_{*})$
 we have $ap^{s}(p+1)q-bq = p^{n+1}q + q = p^{n}(p+1)q - (p^{n}-1)q$
so that $ap^{s}(p+1)q + (p^{n}-1-b)q = p^{n}(p+1)q$. Similar to that in the proof of Prop.
9.8.12 we have $a = 1, s = n , b =
p^n-1$  or $a =(p^{2r+1}+1)/(p+1), b = p^{n-2r}-1, s = p^{n-2r}$( $r\geq
1$) and consequently  $c = 0$.  This shows  (1) and the proof of  (2) is similar. Q.E.D.

\vspace{2mm}

After finishing the proof of the above Proposition, we proceed to prove
Theorem 9.8.1 in case $t\geq 1$.
The proof will be done by some argument processing in the Adams-Novikov resolution of some spectra
and using  the $h_0h_{n+1}$-map in Theorem 9.5.1 as our geometric input.
We first prove the following Lemma.

\vspace{2mm}

{\bf  Lemma 9.8.14}\quad If $g''$ is the element in $\pi_{p^{n}(p+1)q}BP$ such that
 $\tilde b_{1}\tilde c_{0}g'' j\\j'_{p^{n}-1} = 0 \in
[\Sigma^{p^{n+1}q+q-2}K_{p^{n}-1} , BP\wedge E_{1}]$, then there exists
$\bar g = px_{1} + v_{1}^{p^{n}-1}x_{2} \in \pi_{p^{n}(p+1)q}BP$
such that $\tilde b_{1}\tilde
c_{0}(g'' - \bar g ) = 0 \in \pi_{p^n(p+1)q} BP\wedge \widetilde{E}_{1}$, where
$x_1, x_2$  is some elements in $ \pi_*BP$.

{\bf Proof} :   Let $\mu : BP\wedge BP\rightarrow BP$  be the multiplication of
the ring spectrum $BP$, then  $\mu (\tilde b_{0}\wedge 1_{BP}) = 1_{BP} = \mu (1_{BP}
\wedge \tilde b_{0})$,  where $\tilde b_{0} = \tau : S\rightarrow BP$  is the injection
of the bottom cell as stated above.
Then we have the following split cofibration\\
$\quad BP\stackrel{1_{BP}\wedge\tilde b_{0}}{\longrightarrow}BP\wedge BP
\stackrel{1_{BP}\wedge \tilde c_{0}}{\longrightarrow}BP\wedge \widetilde{E}_{1}
\stackrel{1_{BP}\wedge \tilde a_{0} = 0}{\longrightarrow} \Sigma
BP$

$BP\wedge \widetilde{E}_{1}\stackrel{1_{BP}\wedge\tilde b_{0}\wedge 1_{\widetilde{E}_{1}}}
{\longrightarrow}BP\wedge BP\wedge\widetilde{E}_{1}\stackrel{1_{BP}\wedge \tilde
c_{1}}{\longrightarrow}BP\wedge \widetilde{E}_{2}\stackrel{1_{BP}\wedge \tilde a_{0}
= 0 }{\longrightarrow}\Sigma BP\wedge \widetilde{E}_{1}$\\
and there is $\mu' : BP\wedge \widetilde{E}_{1}\rightarrow BP\wedge BP$ such that
$(1_{BP}\wedge \tilde b_{0})\mu + \mu' (1_{BP}\wedge \tilde c_{0}) =
1_{BP\wedge BP}$.

By $\tilde b_{1}\tilde c_{0}g'' jj'_{p^{n}-1}$ = 0 we have $\tilde
c_{0} g'' jj'_{p^{n}-1} = \tilde a_{1} g'$ with $g' \in$ \\$
[\Sigma^{p^{n+1}q+q-1}K_{p^{n}-1}, \widetilde{E}_{2}]$.  Then
$(1_{BP}\wedge\tilde c_{0}g''jj'_{p^{n}-1}) = (1_{BP}\wedge \tilde
a_{1})(1_{BP}\wedge g')$ = 0 so that$1_{BP}\wedge g''jj'_{p^{n}-1}
= [(1_{BP}\wedge\tilde b_{0}) \mu + \mu' (1_{BP}\wedge \tilde
c_{0})](1_{BP}\wedge g''jj'_{p^{n}-1})$ =
$(1_{BP}\wedge\tilde b_{0})\mu (1_{BP}\wedge g''jj'_{p^{n}-1})$  and we have\\
{\bf (9.8.15)} $\qquad (\tilde b_{0}\wedge 1_{BP})g''jj'_{p^{n}-1} = (1_{BP}
\wedge g''jj'_{p^{n}-1})(\tilde b_{0}\wedge 1_{K_{p^{n}-1}})$

$\qquad\qquad\quad = (1_{BP}\wedge \tilde b_{0})\mu (1_{BP}\wedge g''
jj'_{p^{n}-1})(\tilde b_{0}\wedge 1_{K_{p^{n}-1}})$

$\qquad\qquad\quad = (1_{BP}\wedge\tilde b_{0})g''jj'_{p^{n}-1}$.

Note that $(\tilde b_{0}\wedge 1_{BP})_{*}, (1_{BP}\wedge\tilde b_{0})_{*}
: BP_{*}\rightarrow BP_{*}BP$ are the right and left unit
 $\eta_{R}, \eta_{L} : BP_{*}\rightarrow BP_{*}BP$ respectively, then by (9.8.15) we have
$\eta_{R}(g'') = \eta_{L}(g'')$ mod $(p,v_{1}^{p^{n}-1})$. This means that
 $(p, v_{1}^{p^{n}-1}, g'')$ is a $BP_{*}$ invariant ideal , or equivalently,
$g'' \in Ext_{BP_{*}BP}^{0.p^{n}(p+1)q}(BP_{*},
BP_{*}/(p,v_{1}^{p^{n}-1}))$. Then by Prop.  9.8.12 we have\\
{\bf (9.8.16)}$ \qquad g'' = \lambda v_{2}^{p^{n}} + \Sigma
\lambda_{r}v_{1}^{p^{n}-p^{n-2r}}\tilde{c}_{1}(t_{r}p^{n-2r}) + px_{1}
+ v_{1}^{p^{n}-1}x_{2} \in BP_{*}$\\
where $1\leq \lambda, \lambda_{r} \leq p-1, t_{r} = (p^{2r+1}+1)/(p+1)$
and $x_{1}, x_{2}$ are some elements in $  BP_{*}$.

Let $\bar g = px_{1} + v_{1}^{p^{n}-1}x_{2}$, then $(\tilde b_{0}\wedge
1_{BP})(g'' - \bar g) = \eta_{R}(g'' - \bar g) = \eta_{L}(g'' - \bar
g) = (1_{BP}\wedge\tilde b_{0})(g'' - \bar g)$  so that $\tilde b_{1}\tilde
c_{0}(g'' - \bar g) = (1_{BP}\wedge \tilde c_{0})(\tilde b_{0}\wedge
1_{BP})(g'' - \bar g)$ = 0.  Q.E.D.

\vspace{2mm}

{\bf Proof of Theorem  9.8.1 in case $t\geq 1$}\quad By Theorem 9.5.1, there is
 $\tilde{\eta}_{n+1}\in\pi_{p^{n+1}q+q-1}M$ such that $\eta_{n+1} =
j\tilde{\eta}_{n+1} \in \pi_{p^{n+1}q+q-2}S$ is represented in the ASS by
 $h_{0}h_{n+1} \in Ext_{A}^{2,p^{n+1}q+q}(Z_{p},Z_{p})$.
 By Theorem 8.1.5 in chapter 8, $\Phi (\beta_{p^{n}/p^{n}-1}) = h_{0}
 h_{n+1}$, where $\Phi : Ext_{BP_{*}BP}^{2,*}(BP_{*},BP_{*})\rightarrow
 Ext_{A}^{2,*}(Z_{p}, Z_{p})$ is the Thom map.
 Then $\eta_{n+1} = j\tilde{\eta}_{n+1}$ is represented by $\beta_{p^{n}/p^{n}-1}
 + x \in Ext_{BP_{*}BP}^{2,p^{n+1}q+q}(BP_{*},BP_{*})$  in the Adams-Novikov
 spectral sequence,  where $ x = \Sigma_{r\geq 1}
\lambda_{r}\beta_{t_{r}p^{n-2r}/p^{n-2r}-1}$ with $\lambda_{r}
\in Z_{(p)}$ ( cf. Prop. 9.8.13).  Moreover, $\tilde{\eta}_{n+1}\in
\pi_{p^{n+1}q+q-1}M$ is represented by  $\beta'_{p^{n}/p^{n}-1} +
x' + i_{*}(y)$ in the Adams spectral sequence, where $y \in Ext_{BP_{*}BP}^{1,*}(BP_{*},BP_{*})$,
 and $x' = \Sigma_{r\geq 1}\lambda_{r}\beta'_{t_{r}p^{n-2r}/p^{n-2r}-1}$ ,
$\beta'_{tp^{n}/s}$ are the elements in  $Ext_{BP_{*}BP}^{1,*}(BP_{*},BP_{*}(M))$ such that
$j_{*}\beta'_{tp^{n}/s} = \beta_{tp^{n}/s} \in Ext_{BP_{*}BP}^{2,*}
(BP_{*},BP_{*})$ .   It is known that all the generators in
$Ext_{BP_{*}BP}^{1,*}(BP_{*},BP_{*})$ are permanent cycles in the Adams-Novikov spectral sequence
,  then there exists $\tilde{f} \in \pi_{p^{n+1}q
+q-1}M$ such that it is represented by $\beta'_{p^{n}/P^{n}-1} + x'$.
In addition,  $\tilde{f}$  can be extended by $f \in [\Sigma^{p^{n+1}q+q-1}M,M]
\cap ker d$ such that $\tilde{f} = f i $. Recall from   (9.8.8)

$\qquad \cdots \stackrel{\tilde a_{2}\wedge 1_{M}}{\longrightarrow}
\Sigma^{-2}\widetilde{E}_{2}\wedge M\stackrel{\tilde a_{1}\wedge 1_{M}}{\longrightarrow}
\Sigma^{-1}\widetilde{E}_{1}\wedge M\stackrel{\tilde a_{0}\wedge 1_{M}}{\longrightarrow}
\widetilde{E}_{0}\wedge M = M$

$\qquad\quad\qquad\qquad\qquad \Big\downarrow \tilde b_{2}\wedge
1_{M}\qquad\qquad\Big\downarrow
\tilde b_{1}\wedge 1_{M}\quad\qquad\Big\downarrow\tilde b_{0}\wedge 1_{M}$

$\qquad\qquad\Sigma^{-2}BP\wedge \widetilde{E}_{2}\wedge M \qquad\Sigma^{-1}BP\wedge
\widetilde{E}_{1}\wedge M\qquad BP\wedge M$\\
is the Adams-Novikov resolution of the Moore spectrum $M$.   Then $f i $  can be lifted to
$f_{1}i \in\pi_{p^{n+1}q+q}(\widetilde{E}_{1}\wedge M)$  with
$f_{1} \in [\Sigma^{p^{n+1}q+q}M, E_{1}\wedge M]\cap ker d$  such that $\tilde{a}_0\wedge 1_M)f_1i = f i$
and the $d_{1}$-cycle
$(\tilde b_{1}\wedge 1_{M})f_{1}i \in \pi_{p^{n+1}q+q}BP\wedge \widetilde{E}_{1}\wedge M$ represents
$\beta'_{p^{n}/p^{n}-1} + x' \in Ext_{BP_{*}BP}^{1,p^{n+1}q+q}(BP_{*},
BP_{*}(M))$.  By applying $d$ to the equation
 $(\tilde a_{0}\wedge 1_{M})f_{1}ij = fij$ we have
 $(\tilde a_{0}\wedge 1_{M})f_{1} = f$.

Since $(\tilde a_{0}\wedge 1_{M})f_{1}i = fi \in \pi_{p^{n+1}q+q-1}M$ is represented by
 $\beta'_{p^{n}/p^{n}-1} + x' \in
Ext_{BP_{*}BP}^{1,p^{n+1}q+q}(BP_{*},BP_{*}(M))$ in the Adams-Novikov spectral sequence and
 $v_{1}^{p^{n}-1}(\beta'_{p^{n}/p^{n}-1} + x')$ = 0,  then $(\tilde a_{0}\wedge 1_{M})f_{1}
\alpha^{p^{n}-1}i = f\alpha^{p^{n}-1}i = \alpha^{p^{n}-1}fi$  has $BP$-
filtration $> 1$ so that $(\tilde b_{1}\wedge 1_{M})f_{1}\alpha^{p^{n}-1}i$
is a  $d_{1}$-boundary and it equals to  $(\tilde b_{1}\tilde c_{0}\wedge
1_{M})gi$ for some $g \in [\Sigma^{p^{n}(p+1)q}M, BP\wedge M]$.  Hence,
$(\tilde b_{1}\wedge 1_{M})f_{1}\alpha^{p^{n}-1} = (\tilde b_{1}\tilde c_{0}
\wedge 1_{M})g$ , this is because $\pi_{p^{n}(p+1)q+1}BP\wedge \widetilde{E}_{1}\wedge M$ = 0
which is obtained by the sparseness fo $BP_{*}(\widetilde{E}_{1}\wedge M)$. Consequently we have

 $ f_{1}\alpha^{p^{n}-1} = (\tilde c_{0}\wedge 1_{M}) g +
(\tilde a_{1}\wedge 1_{M})f_{2}$
\quad with $f_{2}\in [\Sigma^{p^{n}(p+1)q+1}M,\widetilde{E}_{2}\wedge M]$\\
and

 $(1_{\widetilde{E}_{1}}\wedge j)f_{1}\alpha^{p^{n}-1} = \tilde c_{0}
g''j + \tilde a_{1}(1_{\widetilde{E}_{2}}\wedge j)f_{2}$\quad  with
$g''\in \pi_{p^{n}(p+1)q}BP$, \\
where $g'' j = (1_{BP}\wedge j)g$, this is because $(1_{BP}\wedge
j)gi \in \pi_{p^{n}(p+1)q-1}BP$ = 0. In addition, by $\tilde
b_{2}(1_{\widetilde{E}_{2}} \wedge j)f_{2}\in
[\Sigma^{p^{n}(p+1)q}M,BP\wedge \widetilde{E}_{2}]$= 0  and
$[\Sigma^{p^{n}(p+1)q+1}M, BP\\\wedge \widetilde{E}_{3}]$ = 0 (
cf. Prop. 9.8.9), then $(1_{\widetilde{E}_{2}}\wedge j)f_{2} =
\tilde a_{2}\tilde a_{3}f_{3}$
for some $f_{3}\in [\Sigma^{p^{n}(p+1)q+2}M, E_{4}]$  and we have\\
{\bf (9.8.17)} $\qquad (1_{E_{1}}\wedge j)f_{1}\alpha^{p^{n}-1} = \bar c_0 g'' j
+ \tilde a_{1}\tilde a_{2}\tilde a_{3}f_{3}$.

 By (9.8.17) we have $\tilde b_{1}\tilde c_{0}g'' jj'_{p^{n}-1}$
= 0, then by Lemma 9.8.14, there is $\bar g = px_{1} + v_{1}^{p^{n}-1}x_{2}
\in \pi_{p^{n}(p+1)q}BP$ with $x_{1}, x_{2} \in \pi_{*}BP$ such that
 $\tilde b_{1}\tilde c_{0}(g'' - \bar g)$ = 0.   Consequently, $\tilde c_{0}(g''
-\bar g) = \tilde a_{1}f_{4}$ for some $f_{4} \in \pi_{p^{n}(p+1)q+1}
\widetilde{E}_{2}$  and  $f_{4} = \tilde a_{2}\tilde a_{3}f_{5}$ with
$f_{5} \in \pi_{p^{n}(p+1)q+3}\widetilde{E}_{4}$ which is obtained by the sparseness of
 $\pi_{*}BP\wedge \widetilde{E}_{s}$.  So, (9.8.17)  becomes\\
{\bf (9.8.18)} $\qquad (1_{\widetilde{E}_{1}}\wedge j)f_{1}\alpha^{p^{n}-1} = \tilde c_{0} \bar gj
+ \tilde a_{1}\tilde a_{2}\tilde a_{3} f_{5} j + \tilde a_{1}\tilde a_{2}\tilde
a_{3}f_{3}$.\\
Note that $\bar g = px_{1} + v_{1}^{p^{n}-1}x_{2}$, then
, $\bar g jj'_{p^{n}-1}$ = 0 $\in BP^{*}(K_{p^{n}-1})\cong
\Sigma^{-(p^{n}-1)q-2}\\BP^{*}/(p,v_{1}^{p^{n}-1})$ so that
$\bar g j = \tilde{g}\alpha^{p^{n}-1}$ for some $\tilde{g}\in
[\Sigma^{p^{n+1}q+q-1}M,BP]$.  Consequently,  by (9.8.4),  the equation (9.8.18)  becomes\\
{\bf (9.8.19)} $\qquad (1_{\widetilde{E}_{1}}\wedge j)f_{1}j'_{1}\rho_{p^{n},1} =
(1_{\widetilde{E}_{1}}\wedge j)f_{1}\alpha^{p^{n}-1}j'_{p^{n}}$

$\qquad\qquad = \tilde c_{0}\tilde{g}\alpha^{p^{n}-1}j'_{p^{n}} +
\tilde a_{1}\tilde a_{2}\tilde a_{3}(f_{5} j + f_{3})j'_{p^{n}}$

$\qquad\qquad = \tilde c_{0}\tilde{g}j'_{1}\rho_{p^{n},1} + \tilde
a_{1}\tilde a_{2}\tilde a_{3}(f_{5}j + f_{3})j'_{p^{n}}$\\
Moreover, by (9.8.4)(9.8.6)  we have $\tilde a_{1}\tilde
a_{2}\tilde a_{3}(f_{5}j + f_{3})j'_{p^{n}-1} = \tilde a_{1}\tilde
a_{2}\tilde a_{3}(f_{5}j +f_{3}) j'_{p^{n}}\psi_{p^{n}-1,p^{n}}$ =
0 and so $(f_{5}j + f_{3})j'_{p^{n}-1}$ = 0 , this is because
$[\Sigma^{p^{n+1}q+q+r}\\K_{p^{n}-1}, BP\wedge
\widetilde{E}_{2+r}]$ = 0 for $ r = -1, 0, 1$( cf. Prop. 9.8.9).
This shows that $(f_{5}j + f_{3}) = f_6\alpha^{p^{n}-1}$ with $f_6
\in [\Sigma^{p^{n+1}q+q+2}M,
\widetilde{E}_{4}]$.  Hence,  the equation (9.8.19) becomes\\
{\bf (9.8.20)} $\qquad (1_{\widetilde{E}_{1}}\wedge j)f_{1}j'_{1}\rho_{p^{n},1} =
\tilde c_{0}\tilde{g}j'_{1}\rho_{p^{n},1} + \tilde a_{1}\tilde a_{2}\tilde a_{2}
f_{6}\alpha^{p^{n}-1}j'_{p^{n}}$

$\qquad\qquad\qquad = \tilde c_{0}\tilde{g}j'_{1}\rho_{p^{n},1} + \tilde
a_{1}\tilde a_{2}\tilde a_{3}f_{6}j'_{1}\rho_{p^{n},1}$\\
and by  (9.8.6) we have \\
{\bf (9.8.21)} $\qquad (1_{\widetilde{E}_{1}}\wedge j)f_{1}j'_{1} = \tilde c_{0}
\tilde{g}j'_{1} + \tilde a_{1}\tilde a_{2}\tilde a_{3}f_{6}j'_{1} + \epsilon
i'_{p^{n}-1}j'_{1}$\\
for some $\epsilon \in [\Sigma^{p^{n+1}q+q-1}K_{p^{n}-1}, \widetilde{E}_{1}]$.
By composing $\tilde a_{0}$ to (9.8.21) we have $j f j'_{1} = \tilde a_{0}\tilde a_{1}
\tilde a_{2}\tilde a_{3}f_{6}j'_{1} + \tilde a_{0}\epsilon i'_{p^{n}-1}j'_{1}
= \tilde a_{0}\tilde a_{1}\tilde a_{2}\tilde a_{3}f_{6}j'_{1} +
jj'_{p^{n}-1}\overline{\epsilon }i'_{p^{n}-1}j'_{1}$ , this is because $\tilde a_{0}
\epsilon = jj'_{p^{n}-1}\overline{\epsilon }$ for some
$\overline{\epsilon }\in [\Sigma^{p^{n}(p+1)q}K_{p^{n}-1}, K_{p^{n}-1}]$
(cf. Prop. 9.8.2).  Consequently we have \\
{\bf (9.8.22)} $\qquad j f i = \tilde a_{0}\tilde a_{1}\tilde a_{2}\tilde a_{3}f_{6}i
+ jj'_{p^{n}-1}\overline{\epsilon }i'_{p^{n}-1}i + f_{7}\alpha
i$\\ with $f_{7} \in [\Sigma^{p^{n+1}q-2}M, S]$.

We claim that the map  $f_{7}\alpha i$ in (9.8.22) has  filtration $\geq 3$
so that by (9.8.22) we obtain that $jj'_{p^{n}-1}\tilde{\epsilon}i'_{p^{n}-1}i \in
\pi_{p^{n+1}q+q-2}S$ is represented by $h_{0}h_{n+1} \in$
$Ext_{A}^{2,p^{n+1}q+q}(Z_{p},Z_{p})$  in the Adams spectral sequence.
This claim will be proved in the last.

Then, $jj'_{p^{n}-1}\tilde{\epsilon}i'_{p^{n}-1}i$ is represented by
$\lambda_{0}\beta_{p^{n}/p^{n}-1} + \Sigma_{r\geq 1}\lambda_{r}\beta_{t_{r}
p^{n-2r}/p^{n-2r}-1} \\\in Ext_{BP_{*}BP}^{2,*}(BP_{*},BP_{*})$ in the Adams-Novikov spectral sequence
so that by Prop. 9.8.12 we know that $\tilde{\epsilon}i'_{p^{n}-1}i \in
\pi_{p^{n}(p+1)q}K_{p^{n}-1}$ is represented by  $v_{2}^{p^{n}}
+ \Sigma_{r\geq 1}\lambda_{r}v_{1}^{p^{n}-p^{n-2r}}\tilde{c}_{1}(t_{r}
p^{n-2r})\in Ext_{BP_{*}BP}^{0,p^{n}(p+1)q}(BP_{*},
BP_{*}K_{p^{n}-1})$, where $\lambda_{r} \in Z_{p}$  and $t_{r} =
(p^{2r+2}+1)/(p+1)$.

By [22] Theorem  C,D, it is known that  $v_{2}^{tp}\in
Ext_{BP_{*}BP}^{0,*}(BP_{*},BP_{*}K_{r})$ , for $t\geq 1, 1\leq
r\leq p-1$, is a permanent cyce in the Adams-Novikov spectral
sequence. Suppose inductively that $v_{2}^{tp^{s}}\in
Ext_{BP_{*}BP}^{0,*}(BP_{*},BP_{*}K_{r})$ (for $t\geq 1, 1\leq
r\leq p^{s}-1$ and $s\leq n-1$) is a permanent cycle in the
Adams-Novikov spectral sequence, then we know that
$v_{1}^{p^{n}-p^{n-2r}} \tilde{c}_{1}(t_{r}p^{n-2r}) \in
Ext_{BP_{*}BP}^{0,p^{n}(p+1)q}\\(BP_{*}, BP_{*}K_{p^{n}-1})$ also
is a permanent cycle for all $r\geq 1$. Moreover, by the
representation of the above $\tilde{\epsilon}i'_{p^n-1}i$ we
obtain that $v_{2}^{p^{n}} \in
Ext_{BP_{*}BP}^{0,p^{n}(p+1)q}(BP_{*},\\BP_{*}K_{p^{n}-1})$ is a
permanent cycle.  Hence, by (9.8.6), there exists $k \in
[\Sigma^{p^{n}(p+1)q}K_{p^{n}-1}, K_{p^{n}-1}]$ such that the
induced $BP_{*}$-homomorphism $k_{*} = v_{2}^{p^{n}}$.  In
addition, the map $\rho_{p^{n}-1,r} : K_{p^{n}-1}\rightarrow
K_{r}$ in (9.8.4) for all  $r\leq p^{n}-1$ is a projection,  then
$\rho_{p^{n}-1,r}k^{t}i'_{p^{n}-1} i \in \pi_{tp^{n} (p+1)q}K_{r}$
is represented by $v_{2}^{tp^{n}}\in Ext_{BP_{*}BP}^{0,*}
(BP_{*},BP_{*}K_{r})$ in the Adams-Novikov spectral sequence. This
completes the induction and
$jj'_{r}\rho_{p^{n}-1,r}k^{t}i'_{p^{n}-1}i \in \pi_{*}S$ is just
the
 $\beta_{tp^{n}/r}$-element of the Theorem.

Now our remaining work is to prove the above claim. We turn to an
argument in the ASS and let $A$ be the mod $p$ Steenrod algebra.
By $Ext_{A}^{1, p^{n+1}q-1}(Z_{p},H^{*}M)\cong
Z_{p}\{j^{*}(h_{n+1})\}$  and the result on
$\beta_{p^{n}/p^{n}}\in
Ext_{BP_{*}BP}^{2,p^{n+1}q}\\(BP_{*},BP_{*})$ support a nontrivial
differential in the Adams-Novikov spectral sequence in [12] p.106
Thoerem 5.4.8(i),
we know that\\
{\bf (9.8.23)} $ \qquad j^{*}(h_{n+1}) \in Ext_{A}^{1,p^{n+1}q-1}(Z_{p},H^{*}M)$
dies in the ASS \\
Then,  the map $f_{7}\in [\Sigma^{p^{n+1}q-2}M,S]$ in (9.8.22)
has filtration $\geq 2$ in the ASS and so $f_{7}\alpha i$ has filtration $\geq 3$.
Moreover, by (9.8.21) we know that $jj'_{p^{n}-1}\overline{\epsilon }i'_{p^{n}-1}i$ and
 $j f i \in \pi_{p^{n+1}q+q-2}S$ must have the same filtration  so that
 it is repesented by $h_{0}h_{n+1} \in Ext_{A}^{2,p^{n+1}q+q}(Z_{p},Z_{p})$
in the ASS. This shows the above claim and the Theorem is proved. Q.E.D.
\vspace{2mm}

{\bf  Remark 9.8.24} \quad We give a detail proof of the result in (9.8.23) as follows.
It will be done by some argument processing in the Adams resolution (9.2.9).
Suppose in contrast that the map $j^*h_{n+1}\in Ext_A^
{1,p^{n+1}q-1}(Z_p,H^*M)$ is a permanent cycle in the ASS, then we have
$\bar{c}_{1}h_{n+1}\cdot j  $ = 0 , where $h_{n+1}\in \pi_{p^{n+1}q}
KG_{1}\cong Ext_{A}^{1,p^{n+1}q}(Z_{p},Z_{p})$.  Consequently $\bar{c}_{1}
h_{n+1} = \bar f\cdot p $ for some $\bar f\in \pi_{p^{n+1}q}
E_{2}$.   On the other hand, $\bar{b}_{2}\bar f \in \pi_{p^{n+1}q}KG_{2}
\cong Ext_{A}^{2,p^{n+1}q}(Z_{p},Z_{p})\cong Z_{p}\{b_{n}\}$ so that we have
$\bar{b}_{2}\bar{f} = \lambda \cdot b_{n}$ with $\lambda \in
Z_{p}$.  However, the scalar  $\lambda$ must be zero, this is because
$b_{n}$ support a nontrivial differential $d_{2p-1}(b_{n}) = d_{2p-1}\Phi
(\beta_{p^{n}/p^{n}}) = \Phi d_{2p-1}(\beta_{p^{n}/p^{n}}) = \Phi
(\alpha_{1}\beta_{p^{n-1}/p^{n-1}}^{p}) = h_{0}b_{n-1}^{p} \neq 0$
(cf. [12] p.206  Theorem 5.4.8(i) ).  Hence $\bar{f} = \bar{a}_{2}
\bar{f}_{1}$ for some $\bar{f}_{1} \in \pi_{p^{n+1}q+1}E_{3}$
and we have $\bar{c}_{1}h_{n+1} = \bar{a}_{2}\bar{f}_{1}\cdot p
= \bar{a}_{2}\bar{a}_{3}\bar{f}_{2}$ with $\bar{f}_{2}
\in \pi_{p^{n+1}q+2}E_{4}$. This means that the secondary differential
$d_{2}(h_{n+1})$ = 0 which contradicts with the following known nontrivial
differential $d_{2}
(h_{n+1}) = a_{0}b_{n} \neq 0 \in Ext_{A}^{3,p^{n+1}q+1}(Z_{p},Z_{p})$
( cf. [12] p.11 Theorem 1.2.14).  So we have $\tilde{c}_{1}h_{n+1}\cdot j
\neq 0$ and so (9.8.23) holds.

\vspace{2mm}

Now we proceed to prove Theorem 9.8.1 in case $t\geq 2$. We first prove
the following Lemmas and Propositions.

\vspace{2mm}

{\bf  Lemma 9.8.25}\quad Let $p\geq 3$ and $v_{1}x \in Ext_{BP_{*}BP}^{0,
tp^{n}(p+1)q}(BP_{*},BP_{*}K_{p^{n}})$,  then $v_{1}x = \Sigma_{r = 1}^{[n/2]}
\lambda_{r}v_{1}^{p^{n}-p^{n-2r}}\tilde{c}_{1}(a_{r}p^{n-2r})$, where
$\lambda_{r}\in Z_{p}, a_{r} = (tp^{2r+1} + tp^{2r}-  p^{2r}+1)/(p+1)$
and $\tilde{c}_{1}(ap^{s})$ is the generator in Theorem 8.1.7 in chapter 6 which has degree
$ap^s(p+1)q$.

{\bf Proof} :  By Theorem 8.1.7 in chapter 8, $v_{1}x $ is a
linear combination of the following generators
$v_{1}^{b}\tilde{c}_{1}(ap^{s})$, where $a\geq 1$ is not divisible
by $p$, $1\leq b < p^{n}, b\geq max\{0, p^{n} - q_{1}(ap^{s}\}$
and $q_{1}(ap^{s}) = p^{s} $ if  $a = 1, q_{1}(ap^{s}) = p^{s} +
p^{s-1} -1$ if $a\geq 2$.

By degree reasons we have $bq + ap^{s}(p+1)q = tp^{n}(p+1)q$,  then $s
< n, b\geq p^{n} - (p^{s}+p^{s-1}-1) > 0$  and so $ b\geq p^{n} - p^{n-1}
$ if $s\leq n-2$.  If $ s = n-1$, then $b$ is divisible by $p^{n-1}(p+1)$
so that $b\geq p^{n} + p^{n-1}$.  So, in any case we have $b\geq p^{n-1}
(p-1)$ and the remaining steps is similar to that given in the proof of
Prop. 9.8.12. Q.E.D.

\vspace{2mm}

{\bf  Proposition 9.8.26}\quad Let $r > s$  and $\rho_{r,s} : K_{r}\rightarrow
K_{s}$ be the map in (9.8.4),  then $d(\rho_{r,s}) = i'_{s}\xi
j'_{r}$ with $\xi \in [\Sigma^{rq+1}M,M]\cap ker d$.

{\bf Proof} :   By (9.8.6)(9.8.3) we have $j'_{s}d(\rho_{r,s}) = d(j'_{s}\rho_{r,s}) =
d(\alpha^{r-s}j'_{r})$ = 0 ,  then $d(\rho_{r,s}) = i'_{s}
\overline{\xi}$ for some $\overline{\xi} \in [\Sigma K_{r},M]$ and
$\overline{\xi} = \xi j'_{r}$ with $\xi \in [\Sigma^{rq+1}M,M]$
, this is because $\overline{\xi}i'_{r} \in [\Sigma M,M]$ = 0.  By Theorem 6.4.14 in chapter 6,
we may assume  $\xi = \xi_{1} + \xi_{2}ij$ with $\xi_{1}, \xi_{2}
\in ker d\cap [\Sigma^{*}M,M]$.   Then  $d(\rho_{r,s}) = i'_{s}\xi_{1}
j'_{r} + i'_{s}\xi_{2}ijj'_{r}$ and by applying the derivation $d$ using (9.8.3) we have
 $i'_{s}\xi_{2} j'_{r}$ = 0.  Consequently we have $i'_{s}\xi_{2}
= \xi_{3}\alpha^{r} = 0$,  this is because $\xi_{3} \in [\Sigma^{2}M,K_{s}]$ = 0.
Then $d(\rho_{r,s}) = i'_{s}\xi_{1}j'_{r}$ with $\xi_{1} \in ker d\cap
[\Sigma^{rq+1}M,M]$.  Q.E.D.

\vspace{2mm}

{\bf  Proof of Thoerem 9.8.1 in case $t\geq 2$}:\quad By Theorem 9.8.1 in case
$t \geq 1$,  there exists $ f' \in [\Sigma^{p^{n}(p+1)q}K_{a},K_{a}]$ such that the
induced $BP_{*}$-homomorphism $f'_{*} = v_{2}^{p^{n}}$, where we briefly write
$p^{n}-1$ as  $a$.  By Theorem 6.5.22 in chapter 6, we may assume $f\in
Mod_{*}$ , this is because the components of $f'$ in $Der_{*}$ and
$Mod_{*}\delta_{0}$ induce zero $BP_*$-homomorphism.

 Write $\delta' = i'_{s}j'_{s} \in [\Sigma^{-sq-1}K_{s},K_{s}]$. By Theorem 6.5.23 in chapter 6,
$\delta' f' - f'\delta' \in Mod_{*}$ and this group is a commutative subring of
$[\Sigma^*K_s, K_s]$. Then we have $f'(\delta' f' -
f'\delta') = (\delta' f' - f'\delta' f')f'$ or equivalently, $(f')^{2}\delta' -
\delta' (f')^{2} = 2((f')^{2}\delta' - f'\delta' f')$. By induction we have
$(f')^{s}\delta' - \delta' (f')^{s} = s((f')^{s}\delta' - (f')^{s-1}\delta' f')$
, $s\geq 1$.  That is\\
{\bf (9.8.27)} $\qquad s\cdot (f')^{s-1}\delta' f' = \delta' (f')^{s} + (s-1)(f')^{s}
\delta'$,\qquad $s\geq 1$

Let $\rho_{a,1} : K_{a}\rightarrow K_{1}$ be the projection in
(9.8.4), then by Theorem in chapter 6, $\rho_{a,1}(f')^{s}i'_{a}i
\in \pi_{*}K_{1}$ can be extended to $k_{s} \in Mod_{*}\subset
[\Sigma^{sp^n(p+1)q}K_1,\\ K_1]$ such that
 $\rho_{a,1}(f')^{s}i'_{a}i = k_{s}i'_{1}i$
and $(k_{s})_{*} = v_{2}^{sp^{n}}$.  Since $j'_{1}k_{s}i'_{1}i =
\alpha^{a-1}j'_{a}(f')^{s}i'_{a}i$  and $(i'_{1}j'_{1}k_{s} - k_{s}i'_{1}
j'_{1})\in Mod_{*}$,  then $(i'_{1}j'_{1}k_{s}-k_{s}i'_{1}j'_{1})i'_{1}i$
= 0  and so $i'_{1}j'_{1}k_{s} = k_{s}i'_{1}j'_{1}$.  By applying the derivation $d$
to the equation $\rho_{a,1}(f')^{s}i'_{a}\delta = k_{s}i'_{1}\delta$
(where we write $\delta = ij$)  we have\\
{\bf (9.8.28)} $\qquad \rho_{a,1}(f')^{s}i'_{a} = k_{s}i'_{1} - d(\rho_{a,1})(f')^{s}
i'_{a}\delta = k_{s}i'_{1} - i'_{1}\xi j'_{a}(f')^{s}i'_{a}\delta$

$\qquad = k_{s}i'_{1} - i'_{1}j'_{a}(f')^{s}i'_{a}\xi\delta$, \quad
$s\geq 1$,\\
where $\xi \in [\Sigma^{aq+1}M,M]\cap ker d$( cf. Prop. 9.8.26).  Let
$t\geq 2$ is not divisible by $p$, then by (9.8.27)(9.8.28) we have $i'_{1}j'_{a}(f')^{t}i'_{a} =
\rho_{a,1}i'_{a}j'_{a}(f')^{t}i'_{a} = t\cdot
\rho_{a,1}(f')^{t-1}i'_{a}j'_{a}f'i'_{a}$ = $t\cdot k^{t-1}i'_{1}j'_{a}
f'i'_{a} - t\cdot i'_{1}j'_{a}(f')^{t-1}i'_{a}\xi\delta j'_{a}f'i'_{a}$  and\\
{\bf (9.8.29)} $\qquad i'_{1}j'_{a}\phi = t\cdot k^{t-1}i'_{1}j'_{a}f'i'_{a}$,\\
where we write $\phi = (f')^{t}i'_{a} + t\cdot (f')^{t-1}i'_{a}\xi\delta j'_{a}f'i'_{a}$.

Let $X$ be the cofibre of $i'_{1}j'_{a}f'i'_{a} : \Sigma^{p^{n+1}q+q-1}M
\rightarrow K_{1}$ given by the cofibration in the upper row of the following homnotopy
commutative diagram ($ m = (t-1)p^{n}(p+1)q+(p^{n}-1)q$)

$\quad \Sigma^{-1}X \quad\stackrel{u}{\longrightarrow} \quad\Sigma^{p^{n+1}q+q-1}M
\stackrel{i'_{1}j'_{a}f'i'_{a}}{\longrightarrow} \quad K_{1}\qquad\stackrel{w}
{\longrightarrow} \qquad X$\\
{\bf (9.8.30)}$\quad\Big\downarrow \bar\phi\quad\qquad\qquad\qquad\Big\downarrow
\phi\qquad\qquad\qquad\quad\Big\downarrow t\cdot k^{t-1}\qquad\quad\Big
\downarrow \bar \phi$

$\quad\Sigma^{-m-1}K_{a+1}\stackrel{\rho_{a+1,a}}{\longrightarrow}\Sigma
^{-m-1}K_{a}\stackrel{i'_{1}j'_{a}}{\longrightarrow}\quad\Sigma^{-m+aq}K_{1}
\stackrel{\psi_{1,a+1}}{\longrightarrow} \Sigma^{-m}K_{a+1}$\\
Note that the above middle rectangle is homotopy commutative by (9.8.29),
then there exists  $\bar \phi$ such that all the above rectangles
commute up to homotopy.

By $wi'_{1}j'_{a}f'i'_{a}$ = 0  we have $wi'_{1}j'_{a}f' = y
j'_{a}$ with $y\in [\Sigma^{p^{n}(p+1)q}M,X]$ so that $uyj'_{a}$ =
0 and $uy = \lambda\cdot\alpha^{a}$ for some $\lambda \in [M,M]\cong
Z_{p}\{1_{M}\}$,  that is we have\\
{\bf (9.8.31)}$\qquad wi'_{1}j'_{a}f' = yj'_{a}$, \quad $uy =
\lambda\cdot\alpha^{a}$ for some $\lambda \in Z_{p}$.\\
On the other hand, $j'_{a}(f')^{s}\psi_{1,a}i'_{1} =
j'_{a}(f')^{s}i'_{a}\alpha^{a-1} =
\alpha^{a-1}j'_{a}(f')^{s}i'_{a} =
j'_{1}\rho_{a,1}\\(f')^{s}i'_{a} = j'_{1}k_{s}i'_{1}$, then
$j'_{a}f'\psi_{1,a} = j'_{1}k_{1} + \eta j'_{1}$ with $\eta\in
[\Sigma ^{*}M,M]$ and so $yj'_{1} = wi'_{1}j'_{a}f'\psi_{1,a} =
wi'_{1}j'_{1}k_{1} + wi'_{1}\eta j'_{1} = w k_{1}i'_{1}j'_{1} +
w i'_{1}\eta j'_{1}$  and we have\\
{\bf (9.8.32)} $\qquad y = w k_{1}i'_{1} + w i'_{1}\eta + z \alpha$ \quad
with $z\in [\Sigma^{*}M,X]$,

$\qquad\quad \bar \phi y = t\cdot\psi_{1,a+1}k_{t-1}k_{1}i'_{1} +
t\cdot\psi_{1,a+1}k_{t-1}i'_{1}\eta + \bar \phi z \alpha$\\
which is obtained by (9.8.31).

We claim that\\
{\bf (9.8.33)} \quad $\bar\phi z\alpha i\in \pi_{tp^{n}(p+1)q+aq}K_{a+1}$  has
$BP$ filtration $> 0$\\
This will be proved in the last. Then $\bar \phi y i = t\cdot\psi_{1,a+1}k_{t-1}k_{1}i'_{1}i$ (modulo
higher  filtration) is represented by $t\cdot v_{1}^{a}v_{2}^{tp^{n}}
\in Ext_{BP_{*}BP}^{0,*}(BP_{*},BP_{*}K_{a+1})$ in the Adams-Novikov spectral sequence.

 Hence , by (9.8.31)(9.8.30)(9.8.28)(9.8.27) we have $\bar\phi y j'_{a} =\bar\phi
wi'_{1}j'_{a}f' = t\cdot\psi_{1,a+1}k_{t-1}i'_{1}j'_{a}f' = t\cdot
\psi_{1,a+1}\rho_{a,1}(f')^{t-1}i'_{a}j'_{a}f' = \psi_{1,a+1}\rho_{a,1}
(i'_{a}j'_{a}(f')^{t} + (t-1)(f')^{t}i'_{a}j'_{a})$ = $ (t-1)\psi_{1,a+1}
\rho_{a,1}(f')^{t}i'_{a}j'_{a}$ and so $\bar \phi y = (t-1)\psi_{1,a+1}
\rho_{a,1}(f')^{t}i'_{a} + \tilde{f}\alpha^{a}$  with $\tilde{f} \in
[\Sigma^{tp^{n}(p+1)q}M,K_{a+1}]$.

By the claim (9.8.33), $\bar\phi yi$ is represented by $t\cdot v_{1}^{a}v_{2}^{tp^{n}}$
in the Adams-Novikov spectral sequence,  then $\tilde{f}\alpha^{a}i$ is represented by
$v_{1}^{a}v_{2}^{tp^{n}}$ and so  $\tilde{f}i \in \pi_{tp^{n}(p+1)q}K_{a+1}$
is represented by $v_{2}^{tp^{n}} + v_{1}x \in Ext_{BP_{*}BP}^{0,tp^{n}(p+1)q}
(BP_{*},BP_{*}K_{a+1})$,  where $v_{1}x = \Sigma_{r=1}^{[n/2]}
\lambda_{r}v_{1}^{p^{n}-p^{n-2r}}\tilde{c}_{1}(a_{r}p^{n-2r})$ which is obtained by Lemma 9.8.25.

By [20], if $t\geq 2$ is not divisible by $p$ and $ 1\leq r\leq p$,
$v_{2}^{tp}\in Ext_{BP_{*}BP}^{0,*}(BP_{*},\\BP_{*}K_{r})$ is a permanent cycle in the
Adams-Novikov spectral sequence.  Suppose inductively that $v_{2}^{tp^{s}}\in Ext_{BP_{*}BP}^{0,*}(BP_{*},
BP_{*}K_{r})$ are permanent cycles for all  $ t\geq 2$ is not divisible by $p$,
$1\leq r\leq p^{s}$ and $ s\leq n-1$.  Then, it is easily seen that
$v_{1}^{p^{n}-p^{n-2r}}\tilde{c}_{1}(a_{r}p^{n-2r})$ is realizable
in $[\Sigma^{tp^{n}(p+1)q}K_{a+1}, K_{a+1}]$ so that the above  by the induction hypothesis
we know that $v_{1}x$ also is a permanent cycle.                                      .
 So, $v_{2}^{tp^{n}}\in Ext_{BP_{*}BP}^{0,*}(BP_{*},BP_{*}K_{a+1})$
is a permanent cycle in the Adams-Novikov spectral sequence and there exists
  $h\in\pi_{tp^{n}(p+1)q}K_{a+1}$ such that the induced $BP_{*}$-homomorphism
$h_{*} = v_{2}^{tp^{n}}$. Hence, for $1\leq r\leq a+1 = p^{n}$,
$jj'_{r}\rho_{a+1,r}h \in \pi_{tp^{n}(p+1)q-rq-2}S$ is just the
$\beta_{tp^{n}/s}$-element of the Theorem.

Now our remaining work is to prove the claim (9.8.33). Recall  as
known above that $j'_{a}f'i'_{a}i\in \pi_{*}M$ is represented by
$\beta'_{p^{n}/p^{n}-1} \in
Ext_{BP_{*}BP}^{1,*}(BP_{*},\\BP_{*}M)$ in the Adams-Novikov
spectral sequence and $\beta'_{p^{n}/p^{n}-1} =
v_{1}\beta'_{p^{n}/p^{n}}$,  then
$(i'_{1})_{*}(\beta'_{p^{n}/p^{n}-1}) = 0 \in
Ext_{BP_{*}BP}^{1,*}(BP_{*},BP_{*}K_{1})$  and so
$i'_{1}j'_{a}f'i'_{a}i\in \pi_{*}K_{1}$  has  $BP$-filtration
$\geq q+1$. Then, in the Adams-Novikov resolution of the spectrum
$K_{1}$ , $i'_{1}j'_{a}f'i'_{a}i$  can be lifted to
$\kappa\in\pi_{*}\widetilde{E}_{q+1}\wedge K_{1}$  such that
$(\tilde a_{0}\wedge 1_{K_{1}})\cdots(\tilde a_{q}\wedge
1_{K_{1}})\kappa = i'_{1}j'_{a}f'i'_{a}i$. Since $K_{1}$ is an
$M$-module spectrum , then $\kappa = \kappa' \cdot i$ with
$\kappa' \in [\Sigma^{*}M,\widetilde{E}_{q+1}\wedge K_{1}]$.
Consequently we have $i'_{1}j'_{a}f'i'_{a} = (\tilde a_{0}\wedge
1_{K_{1}})\cdots (\tilde a_{q}\wedge 1_{K_{1}})\kappa' + \sigma j$
with $\sigma \in [\Sigma^{p^{n+1}q+q}S,K_{1}]$.  Note that
$(\tilde b_{0}\wedge 1_{K_{1}}) \sigma \in
\pi_{p^{n+1}q+q}BP\wedge K_{1}\cong Hom_{BP_{*}BP}^{p^{n+1}q+q}
(BP_{*},BP_{*}(BP\wedge K_{1}))$ is a $d_1$-cycle in the
Adams-Novikov resolution of $K_1$ and it represents an element in
$Ext_{BP_{*}BP}^{0, p^{n+1}q+q}(BP_{*},BP_{*}K_{1})$. However,
this group is zero by degree reason , this  is because
$Ext_{BP_{*}BP}^{0,*}(BP_{*},BP_{*}K_{1})\cong Z_{p}[v_{2}]$).
Then we have $(\tilde b_{0}\wedge 1_{K_{1}})\sigma $ = 0 so that
$\sigma$ can be lifted to $\sigma' \in
\pi_{*}\widetilde{E}_{q+1}\wedge K_{1}$ such that $(\tilde
a_{0}\wedge 1_{K_{1}})\cdots (\tilde a_{q}\wedge 1_{K_{1}})\sigma'
= \sigma$.   So we have $i'_{1}j'_{a}f'i'_{a} = (\tilde
a_{0}\wedge 1_{K_{1}})\cdots (\tilde a_{q}\wedge
1_{K_{1}})(\kappa' +\sigma' j)$. By this we know that the
following short exact sequence induced by the cofibration in the
top row of (9.8.30) is a split exact sequence of
$BP_*BP$-comodule:

$\quad 0 \rightarrow BP_{*}K_{1}\stackrel{u_{*}}{\longrightarrow}
BP_{*}X\stackrel{w_{*}}{\longrightarrow}BP_{*}M\rightarrow 0$\\
where $\mid w_{*}\mid = -(p^{n+1}+1)q$\\
Moreover, this splitness also hold in the following $Ext_{BP_{*}BP}^{0,*}$-stage :

$\quad 0 \rightarrow Ext^{0}K_{1}\stackrel{u_{*}}{\longrightarrow}
Ext^{0}X\stackrel{w_{*}}{\longrightarrow}Ext^{0}M\rightarrow 0$\\
That is to say,  there is an invariant $BP_{*}$-homomorphism $u' :
Ext^{0}X \rightarrow Ext^{0}K_{1}$  and $ w' : Ext^{0}M\rightarrow
Ext^{0}X$ such that $u'u_{*} = 1_{Ext^{0}K_{1}}, w_{*}w' =
1_{Ext^{0}M}$  and $u_{*}u' + w'w_{*} = 1_{Ext^{0}X}$, where we
briefly write $Ext_{BP_{*}BP}^{0,*}(BP_{*},\\BP_{*}X) $  as
$Ext^{0}X$.

To prove the claim (9.8.33), suppose in contrast that $\bar\phi z\alpha i
\in \pi_{*}K_{a+1}$  has $BP$-filtration 0, then, by (9.8.32), it is represented by
$\lambda v_{1}^{a}v_{2}^{tp^{n}}\in Ext^{0}K_{a+1}$  in the Adams-Novikov spectral sequence, where
$\lambda\neq 0 \in Z_{p}$.  Then $z i \in \pi_{tp^{n}(p+1)q+(a-1)q}
X$  must have $BP$ filtration 0 and it is represented by some
$x\in Ext^{0,*}X$ and $(\bar\phi )_{*}(v_{1}x) = \lambda\cdot v_{1}^{a}
v_{2}^{tp^{n}}$.  However,

$\qquad x = u_{*}u'(x) + w'w_{*}(x) = w'w_{*}(x)$\\
, this is because by degree reason we have $u'(x)\in Ext^{0,tp^{n}(p+1)q+(a-1)q}K_{1}$ = 0,
Then

$ x = \lambda' w'(v_{1}^{r})$,\quad for some $\lambda'\in Z_{p}, rq = tp^{n}
(p+1)q+(a-1)q-p^{n+1}q-q$ \\
since $w_{*}(x)\in Ext^{0,r}M\cong Z_{p}\{v_{1}^{r}\}$,
Then\\\centerline{$\lambda v_1^av_2^{tp^n} = (\bar\phi)_*(v_1x) = \lambda' (\bar\phi )_{*}w'(v_{1}^{r+1})$}\\
Moreover, since $w'(v_{1}^{r+1})$ is belong to $Ext^{0,*}X$ the summand
which is isomorphic to  $Ext^{0,*}M$ and
$Ext^{0,*}M$ is a trivial $Z_{p}[v_{2}]$-module, then we have

$\qquad  0 = v_{2}^{p^{n}}\cdot (\bar\phi )_{*}w'(v_{1}^{r+1}) =
\lambda\cdot v_{1}^{a}v_{2}^{(t+1)p^{n}}\in Ext^{0,*}K_{a+1}$\\
This is a contradiction and then shows the claim(9.8.33). Q.E.D.

\vspace{2mm}

After finishing the proof of Thoerem 9.8.1 on second periodicity elements in
the stable homotopy groups of spheres, we state the following Theorem on further result
on second periodicity families in the stable homotopy groups of spheres without proof.
The proof is done in base on the result of Theorem 9.8.1 and using
some properties of the spectrum $M(p^r,v_1^{ap^s})$ which is the geometric realization of $BP_*/(p^r,v_1^{ap^s})$.
The details of the proof can be seen in [23] \S 3.

\vspace{2mm}

{\bf  Theorem 9.8.34}\quad Let $p\geq 5$.  $j = cp^i\leq p^{n-i}-1$ if $t\geq 1$
($cp^i\leq p^{n-i}$ if $t\geq 2$), then the element \\
\centerline{$\beta_{tp^n/j,i+1}\in Ext_{BP_*BP}^{2,*}(BP_*,BP_*)$}\\
is a permanent cycle in the Adams-Novikov spectral sequence and it converges
to the corresponding homotopy element of order $p^{i+1}$ in
$\pi_*S$.

\quad

\begin{center}

{\bf \large REFERENCES}

\end{center}

\vspace{2mm}

[1] Aikawa T. , 3-Dimensional cohomology of the mod p Steenrod algebra
{\it Math. Scanfd. } 47(1980), 91--115.

[2] Cohen R. , Odd primary families in stable homotopy theory. {\it Memoirs of Amer.
Math. Soc.} No. 242(1981).

[3] Cohen R. and Goerss P. , Secondary cohomology operations that detect homotopy classes.
{\it Topology} 22(1984), 177--194.

[4] Hoffman P. , Relations in the stable homotopy of Moore spaces. {\it Proc.
London Math. Soc.} 18(1968), 621--634.

[5] Jinkun Lin and Qibing Zheng , A new family of filtration seven in the stable homotopy of spheres.
{\it Hiroshima Math. J.} 28(1998) 183--205.

[6] Jinkun Lin , A new family of filtration three in the stable homotopy of spheres.
{\it Hiroshima Math. J.} 31(2001) 477--492.

[7] Jinkun Lin , Some new families in the stable homotopy of spheres revisited.
{\it Acta Math. Sinica} 18(2002) 95--106.

[8] Jinkun Lin , Two new families in the stable homotopy groups of sphere
and Moore spectrum.  {\it Chin. Ann. of Math.} 27B(2006) 311-328.

[9] Jinkun Lin, Third periodicity families in the stable homotopy of spheres.
{\it JP Journal of Geometry and Topology} 3(2003), 179-219.

[10] Miller H. R. , Ravenel D.C. and Wilson W.S. . Periodic phenomena in the
Adams-Novikov spectra sequence. {\it Ann. of Math. } 106(1977) 469-516.

[11] Oka S. , Multilpicative structure of finite spectra and stable homotopy
of spheres. {\it Lecture Notes in Math. v.1051} Springer-verlag (1984).

[12] Ravenel D.C. , Complex cobordism and stable homotopy groups of
spheres .  {\it Academic Press Inc. (1986)}

[13] Thomas E. and Zahler R. , Generalized higher order cohomology operations
and stable homotopy groups of spheres. {\it Advances in Math.} 20(1976) 287--328.

[14] Toda H. , Algebra of stable homotopy of $Z_p$-spaces and applications.
{\it J. Math. Kyoto Univ.} 11(1971),197--251.

[15] Toda H. , On spectra realizing exterior part of the Steenrod algebra.
{\it Topology} 10(1971), 53--65.

[16]  Wang X. and  Zheng Q. , The convergence of $\widetilde{\alpha}^{(n)}h_0h_k$
{\it Science in China} 41(1998), 622--636.

[17] Wang X., On the 4-dimensional cohomology of the Steenrod algebra.
{\it Beijing Mathematics} 1(1995), 80-99.

[18] Zhou X., Higher cohomology operations that detect homology class.
{\it Lecture Notes in Math. v.1370} Springer-verlag 1980, 416--436.

[19] Miller H.R. and Wilson W.S. , On Novikov's $Ext^1$ modulo an
invariant prime ideal. {\it Topology} 15(1976), 131-141.

[20] Oka S., Realizing some cyclic $BP_{*}$ modules
and applications to stable homotopy of spheres.
{\it Hiroshima Math. J. } 7(1977), 427-447.

[21]  Oka S., Small ring spectra and p-rank of the
stable homotopy of spheres.  {\it Contemp. Math.} 49(1983), 267-308.

[22] Oka S., A new familiy in the stable homotopy of spheres. {\it
Hiroshima Math. J.} 5(1975),87-114.

[23] Jinkun Lin , Detection of second periodicity families in stable
homotopy of spheres. {\it American J. of Math. } 112(1990), 595-210.

[24] Jinkun Lin, A pull back Theorem in the Adams spectral
sequence,   Acta Math. Sinica. v.34(2008) no.3,471-490

[25] Hirofumi Nakai, The Chrometic $E_1$-term $H^0M_1^2$ for p > 3.
{\it New York Journal of Math. } 6(2000), 21-54.

\end{document}